\documentclass[12pt]{article}
\usepackage{amsthm}
\usepackage{amsmath}
\usepackage{natbib}
\usepackage[colorlinks,citecolor=blue,urlcolor=blue,filecolor=blue,backref=page]{hyperref}
\usepackage{graphicx}
\usepackage{multicol}
\usepackage{amssymb}
\usepackage{rotating}
\usepackage{array}
 \usepackage{color,mathrsfs,enumerate,adjustbox,setspace}
\usepackage[margin=.98in]{geometry}

\newcommand{\blind}{1}

\newcommand{\sone}{\mathbb{S}}
\newcommand{\real}{\mathbb{R}}

\newcommand{\ltwo}{\mathbb{L}^2}
\newcommand{\G}{\mathcal{G}}
\newcommand{\D}{\mathcal{D}}

\newcommand{\Po}{\mathcal{P}}
\newcommand{\Tn}{\mathcal{T}_n}
\newcommand{\Dn}{\Delta_{n-1}}
\newcommand{\Dnn}{\Delta_{n-1}}
\newcommand{\nt}{\lfloor nt \rfloor}
\newcommand{\id}{\gamma_{\textbf{\textit{id}}}}

\newtheorem{proposition}{Proposition}
\newtheorem{theorem}{Theorem}
\newtheorem{corollary}{Corollary}
\theoremstyle{definition}

\newtheorem{algorithm}{Algorithm}

\begin{document}

\def\spacingset#1{\renewcommand{\baselinestretch}%
{#1}\small\normalsize} \spacingset{1}


\if1\blind
{
  \title{\bf Distribution on Warp Maps for Alignment of Open and Closed Curves}
  \author{Karthik Bharath$^1$ and Sebastian Kurtek$^2$\\
    $^1$\small{School of Mathematical Sciences, University of Nottingham, UK}\\
    $^2$\small{Department of Statistics, The Ohio State University, USA}}
    \date{}
  \maketitle
} \fi

\if0\blind
{
  \bigskip
  \bigskip
  \bigskip
  \begin{center}
    {\LARGE\bf Distribution on Warp Maps for Alignment of Open and Closed Curves}
\end{center}
  \medskip
} \fi

\bigskip
\begin{abstract}
Alignment of curve data is an integral part of their statistical analysis, and can be achieved using model- or optimization-based approaches. The parameter space is usually the set of monotone, continuous warp maps of a domain. Infinite-dimensional nature of the parameter space encourages sampling based approaches, which require a distribution on the set of warp maps. Moreover, the distribution should also enable sampling in the presence of important landmark information on the curves which constrain the warp maps. For alignment of closed and open curves in $\mathbb{R}^d, d=1,2,3$, possibly with landmark information, we provide a constructive, point-process based definition of a distribution on the set of warp maps of $[0,1]$ and the unit circle $\sone$ that is (1) simple to sample from, and (2) possesses the desiderata for decomposition of the alignment problem with landmark constraints into multiple unconstrained ones. For warp maps on $[0,1]$, the distribution is related to the Dirichlet process. We demonstrate its utility by using it as a prior distribution on warp maps in a Bayesian model for alignment of two univariate curves, and as a proposal distribution in a stochastic algorithm that optimizes a suitable alignment functional for higher-dimensional curves. Several examples from simulated and real datasets are provided.
\end{abstract}

\noindent%
\vfill
\newpage
\spacingset{1.45} 

\section{Introduction}
\label{introduction}
In contrast to standard multivariate analysis, the concept of phase variation is a unique feature of functional data. For functional data obtained as parametric curves representing geometric objects in high-resolution images, establishing correspondence between points on the curves is an important task. Failure to isolate and quantify variation due to phase or lack of correspondence between points can be detrimental when computing descriptive summaries, and for subsequent inferential tasks. This process of isolation and quantification is referred to as registration or alignment.

There exist several approaches to alignment. One popular approach is continuous monotone alignment, which refers to the alignment of two curves $g_i: D \to \mathbb{R}^d,\ d\geq 1,\ i=1,2$, where $D$ is a compact domain, by estimating a homeomorphic self map $\gamma:D \to D$, known as a warp map, that best matches $g_2 \circ \gamma$ to $g_1$ (or vice versa). A variational formulation quantifies the matching through a cost or energy functional, and alignment is defined as the determination of an optimal warp map $\gamma^*$ from a class $W$ that minimizes the cost functional. Alternatively, a statistical model-based formulation of the alignment task seeks to estimate the warp map $\gamma$ based on a likelihood function defined using some discrepancy measure between $g_1$ and $g_2 \circ \gamma$, conditional on $\gamma$.
%

\subsection{Motivation and related work}
\label{sec:motcont}
The focus of this work is on the construction of a probability distribution on $W$ and development of a simple sampling scheme, which would enable quantification of uncertainty on the optimal alignment of open and closed curves using \emph{any} of the available continuous monotone methods. The need for a distribution enabling stochastic approaches to alignment arises in two contexts: (1) since the parameter space of warp maps is a non-linear function space, deterministic algorithms for solving the variational problem can get stuck in local minima, and (2) in a Bayesian model a prior distribution on $W$ is required. The construction and the desirable features of the distribution are motivated by the following common modeling scenarios that arise in the analysis of biomedical and biological datasets.

Boundaries of data objects such as organs in high-resolution medical images are often represented as parametrized open curves (connected curves that begin and end at different points) or simple closed curves (connected curves that do not cross themselves, and begin and end at the same point). In certain cases, explicit physiological information can constrain the warp maps. 
\begin{figure}[!t]
\begin{center}
	\begin{tabular}{cc}
		\includegraphics[width=2in]{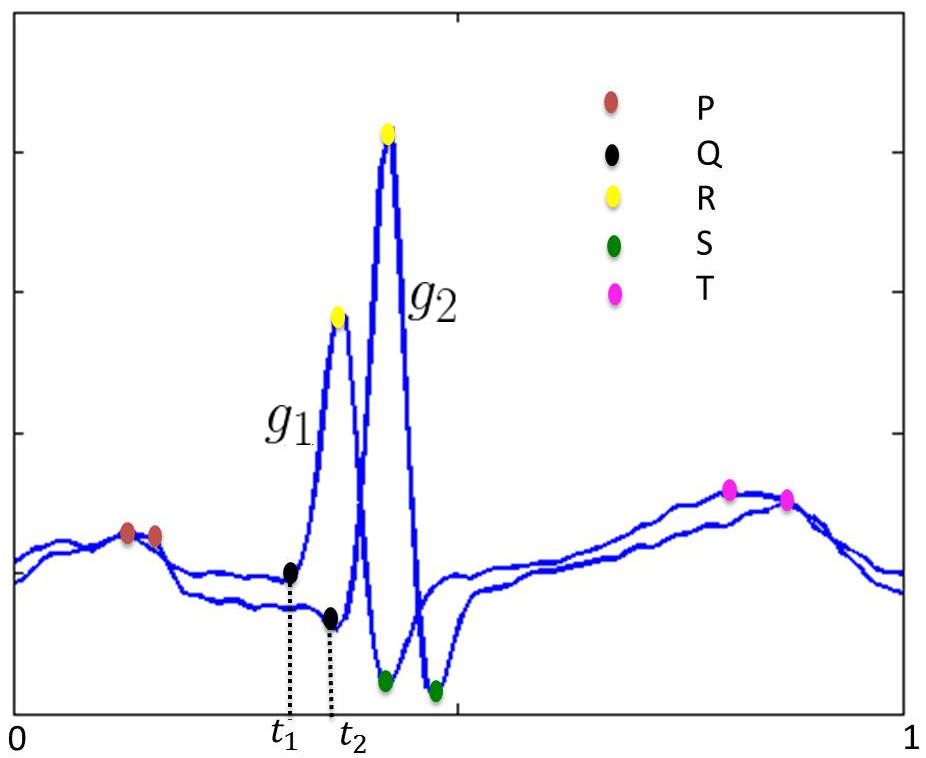}&\includegraphics[width=2in]{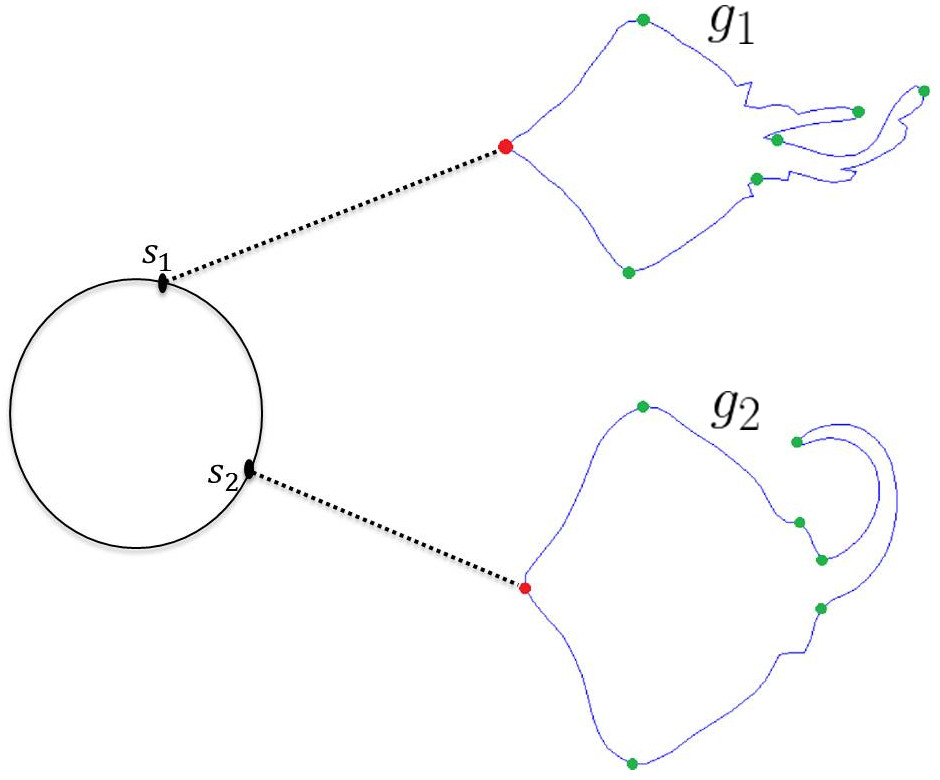}
	\end{tabular}
	\caption{\small Left: Alignment of two ECG cycles with corresponding PQRST complex landmarks. Right: Alignment of the shapes of two stingrays, represented as simple closed curves with landmark features (red and green points).}
	\label{fig1}
\end{center}
\end{figure}
In Figure \ref{fig1}, the left panel shows two open curves (univariate functions) $g_i:[0,1]\to \mathbb{R},\ i=1,2$ representing two electrocardiogram (ECG) cycles with corresponding PQRST complexes marked as fixed landmarks (maxima and minima marked on each cycle). For example, the point Q on $g_1$ at $t_1$ is to be registered to the corresponding Q point on $g_2$ attained at $t_2$ (so do the feature points P, R, S and T). This introduces constraints on any warp map $\gamma: [0,1]\to[0,1]$: $\gamma(t_1)=t_2$ to ensure the matching of the Q point across the two functions (additional constraints are used to ensure the matching of the P, R, S and T points as well); between the landmarks, $\gamma$ is unconstrained. The example in the right panel shows the outlines of two stingrays, represented as embeddings of the circle of unit circumference $\sone$ in $\mathbb{R}^2$, or equivalently as planar closed curves $g_i:\sone \to \mathbb{R}^2,\ i=1,2$. The snout of the stingrays (marked in red) as well as other landmark points (marked in green) are to be matched. This again imposes constraints on a warp map $\gamma: \sone \to \sone$ requiring exact matching of landmarks with unconstrained matching in-between.

There are two complementary requirements of the alignment procedure in the above scenarios: (1) uncertainty around the `best' alignment needs to be captured, and (2) constraints due to landmark points need to be automatically incorporated into the alignment procedure. Operationally, this demands a probability distribution on $W$ that satisfies the following.
\begin{enumerate}[(i)]
	\item An alignment task with $m$ constraints can be decomposed into $m+1$ unconstrained ones by performing \emph{independent} alignment of the curves across $m+1$ subsets of the domain $D$. This requires the initial constrained alignment problem to satisfy the desirable \emph{subset invariance} condition: if $\gamma^*:D\to D$ is the optimal warp map and if $\gamma^*_B:B\to B,\ B \subset D$ is the optimal map when alignment is performed only on a subset $B$ induced by landmarks, then $\gamma^*=\gamma^*_B$ on $B \subset D$ \citep{TY}. In the presence of landmark-induced constraints, a probability distribution $\pi$ on $W$ is said to satisfy subset invariance if its restriction $\pi_B$ to warp maps of $B \subset D$ is a suitably re-scaled version of $\pi$, and is independent of $\pi_{D\backslash B}$.
	\item For efficient exploration of $W$, the distribution should be flexible enough to be centred at any warp map of choice. For example, in the Figure \ref{fig1} examples, as a first step, one can align the two curves \emph{only} at the constraint points with a piecewise linear (PL) warp map. Alignment of the remaining regions can be carried out by sampling in the neighborhood of the PL warp map within a stochastic algorithm, or by employing a prior distribution centred at the PL warp map in a Bayesian model.
\end{enumerate}

Literature on alignment has mostly focused on functional data defined on a closed interval $[a,b]$: for frequentist approaches see \cite{KG, GG, KG2, ZSBM} and \cite{TM}; for Bayesian approaches see \cite{TI,Ian,CSS,YL} and \cite{SK}. Alignment methods for closed curves with $D=\sone$ are conspicuous by their absence within statistics literature; a notable exception is the differential-geometric approach of \cite{AS}. For a good account of variational strategies for alignment of closed and open curves, we refer the reader to \cite{AK}.

Landmark-constrained alignment, under a geometric, non-stochastic framework, was only recently studied \citep{SKBM,bauer}. In those methods, it is not possible to capture or model uncertainty in the optimal alignment. As far as we know, stochastic approaches to curve alignment in the presence of landmarks have not been considered before.
\subsection{Contributions}

In a Bayesian model for alignment, \cite{Ian} define a finite-dimensional Dirichlet distribution on discretized warp maps as a prior distribution; this results in a finite-dimensional specification of the prior instead of a functional one. We first prove that their construction results in a limiting degenerate distribution on $W$ as dimensionality increases (Theorem \ref{nowarp}). We further establish degeneracy of the limiting distribution for a class $\mathcal{C}$ of finite-dimensional distributions that generalize the Dirichlet (Corollary \ref{nowarp2}).

Employing point process methods, we propose a modification of their finite-dimensional specification that results in a constructive definition of a non-degenerate distribution on the set of warp maps of $[0,1]$ and $\sone$ that satisfies requirements (i) and (ii) (Theorem \ref{th1}). We show that the distribution is a canonical one in the sense that it remains the limiting distribution for all finite-dimensional specifications from a class $\mathcal{C}$. The distribution is related to the Dirichlet process \citep{TF} on the set of probability measures on $[0,1]$ (Remark 2). Our approach provides an explicit link between the self-similarity- and Markov-type properties of the Dirichlet process and requirements (i) and (ii) (Proposition \ref{properties}). The distribution can be parametrized by a concentration parameter $\theta>0$ whose value determines how close the probability mass is distributed around the chosen average warp map. The warp maps sampled from the distribution are discontinuous with probability one. \emph{We show that this is unavoidable if one insists on subset invariance in requirement (i)}. The constructive definition, along with the concentration parameter $\theta$, identifies the finite-dimensional projections (coordinates) of the distribution with a Dirichlet distributed random vector, with parameters that depend intimately on the discretization and the choice of the average warp map. This leads to a simple algorithm for sampling PL warp maps.

Finally, we propose a novel stochastic algorithm based on the proposed distribution for alignment of open and closed curves in $\mathbb{R}^k, k=1,2,3$, possibly with landmark constraints. We elucidate on the importance of requirements (i) and (ii) on the distribution through several simulation and real-data examples. In addition, we employ the proposed distribution in a Bayesian model for alignment, similar to the one used by \cite{Ian}. 

The rest of this paper is organized as follows. Section \ref{fixed_part} examines the algorithm of \cite{Ian} and details its shortcomings. Section \ref{random_part} modifies the preceding construction, and proposes a theoretically-justified approach for warp maps of $[0,1]$; the properties of the constructed distribution are studied in Section \ref{properties_section}. Section \ref{circle} extends the construction mechanism to a corresponding distribution on warp maps of $\sone$. Section \ref{examples} presents sample warp maps under different settings, and results from different analyses of real open and closed curve data, possibly with landmark constraints, under a Bayesian model and using a novel stochastic algorithm (Section \ref{simulated_annealing}). Section \ref{discussion} discusses extensions of the proposed methods. The Supplementary Material contains proofs of all results, an alternative construction of a distribution on warp maps of $\sone$, and detailed descriptions of the datasets used in this work.
\section{Construction using fixed partitions and issues}
\label{fixed_part}
We first consider the construction of a distribution on the set of warp maps of a closed subinterval of the real line, which, without loss of generality, can be assumed to be $[0,1]$. The possibility of landmarks on the observed curve data implies that the class of smooth warp maps is inappropriate. Instead, consider the class given by $W_I:=\{\gamma:[0,1] \to [0,1],\ \text{increasing, continuous},\ \gamma(0)=0,\ \gamma(1)=1\}$. In a recent paper on Bayesian alignment of curves, \cite{Ian} proposed a simple method to obtain a continuous random warp map in $W_I$. The following is a summary of their algorithm.
\begin{algorithm}{\textbf{Fixed partition-based sampling of warp maps.}}
	\begin{enumerate}[1.]
		\item Choose a deterministic set of ordered points $0=:t_0<t_1<t_2<\ldots<t_{n-1}<t_{n}:=1$ that induces a partition $\Tn$ of $[0,1]$.
		\item Sample an $n$-dimensional Dirichlet distributed random vector with all parameters set to the same value $\alpha>0$.
		\item Construct a warp map on $[0,1]$ by linear interpolation of the increments.
	\end{enumerate}\label{alg1}
\end{algorithm}
The resulting warp maps are continuous and are elements of $W_I$. The idea behind this approach is based on the fact that for $\gamma \in W_I$, its increments $p_i:=\gamma(t_{i})-\gamma(t_{i-1})$ are positive and satisfy $\sum_{i=1}^n p_i=1$. Hence, $p_i$ can be identified with coordinates of the simplex $\Dnn:=\{x=(x_1,\ldots,x_{n}) \in \mathbb{R}^{n}:\ x_i\geq 0,\ \sum_{i=1}^{n}x_i=1 \}$. The parameter space of warp maps generated in this fashion is essentially finite-dimensional, since a warp map is fully determined by its values at the ordered set of points $t_i$. \cite{Ian} state (without proof) that as $n \to \infty$, this results in a Dirichlet process.

For a fixed partition $\Tn$, Algorithm \ref{alg1} recommends simulating a Dirichlet random vector with all $n$ parameters equal to $\alpha>0$, or equivalently $\alpha(1,1,\ldots,1)$.  The parameter vector is hence independent of the size $n$ of partition. The natural way to incorporate information from the partition $\Tn$ is to assume that the parameter vector is of the form $\alpha^*(\frac{1}{n},\ldots,\frac{1}{n})$ where $\alpha^*=n\alpha$. This can be generalized to an arbitrary deterministic partition $\mathcal{T}=\{B_1,\ldots,B_k,\ k\geq 1\}$ of $[0,1]$ by simulating a $|\mathcal{T}|$-dimensional Dirichlet random vector with parameters $\alpha(\mu(B_1),\ldots,\mu(B_k))$, where $\alpha>0$, $\mu$ is a finite measure on $[0,1]$, and $|T|$ denotes the cardinality of the set $T$. An equi-spaced or uniform partition arises by setting $\mathcal{T}=\mathcal{T}_n$ with $t_i-t_{i-1}=1/n, i=1,\ldots,n$, and taking $\mu$ as the Lebesgue measure.

From a practical perspective, Algorithm \ref{alg1} is appealing since the partition $\Tn$ can be chosen in various ways; the increments can be sampled from distributions different from the Dirichlet on $\Dnn$, and linear interpolation results in continuous maps. The two pertinent questions are: (1) What is the corresponding distribution on $W_I$ as $n \to \infty$? and (2) can the sampled Dirichlet random vector be identified with the finite-dimensional distributions of a stochastic process? The lack of partition information in the parameter vector has a significant implication on the answers to questions (1) and (2).

The algorithm constructs a warp map by linear interpolation of the increments obtained from a fixed partition $\mathcal{T}_n$. The natural setting for its examination is the space $C([0,1])$ of real-valued, continuous functions on $[0,1]$, with the linearly interpolated process based on the partial sum of the increments:
\begin{equation}
\label{LIN}
Y_n(t):=\sum_{i=1}^{\lfloor nt \rfloor}p_i+(nt-\nt)p_{\nt+1}, \quad t \in[0,1].
\end{equation}
Clearly, $Y_n(0)=0$, $Y_n(1)=1$, and $Y_n$ is continuous, increasing in $(0,1)$ and an element of $W_I$. The case when $\alpha=1$ and $(p_1,\ldots,p_{n})$ is uniform on $\Dnn$ (for a given partition of $\Tn)$ is particularly instructive as it captures the key features of the algorithm.
\begin{theorem}
	\label{nowarp}
	Suppose $(p_1,\ldots,p_{n})$ based on a fixed equi-spaced partition $\Tn$ is uniformly distributed on $\Dnn$, and $\id:[0,1]\to [0,1]$ with $\id(t)=t$. In $C([0,1])$, equipped with the uniform topology, $Y_n$ converges in probability to the identity warp map $\id$. The process $\sqrt{n}(Y_n(t)-\id(t))$ converges in distribution to a standard Brownian Bridge process.
\end{theorem}
Theorem \ref{nowarp} states that for large $n$, when $\alpha=1$, the sampling algorithm proposed by \cite{Ian} results in a distribution on $W_I$ that is \emph{degenerate at the identity warp map}. The conclusion also remains true with $\id$ replaced by another deterministic warp map for a fixed non-equi-spaced partition $\mathcal{T}_n$ (see Supplementary Material). Moreover, the fluctuations away from the identity warp can be captured by the behavior of a standard Brownian Bridge. This result suggests that the resulting distribution on the class $W_I$, in the limit, is governed only by the value of the Dirichlet scalar concentration parameter $\alpha$, and concentrates on the identity warp map. See Section \ref{simulations} for numerical illustrations of such degenerate behavior.

In fact, such an uninteresting distribution on $W_I$ resulting from Algorithm \ref{alg1} is not restricted to the case where the increments are uniformly distributed on $\Dn$; \emph{this phenomenon is applicable to a rather large class of distributions on $\Dn$ based on spacings of random variables}. Suppose $x_1,\ldots,x_{n}$ are independent from an absolutely continuous distribution function $F$ on $[0,1]$, with density $f$ and quantile function $Q$ on $(0,1)$. Extend the definition of $Q$ to $[0,1]$ by setting $0=:Q(0)=\lim_{t\to 0} Q(t)$ and $1=:Q(1)=\lim_{t\to 1} Q(t)$. Set $x_{0:n}:=0$, $x_{n:n}:=1$, let $0<x_{1:n}<x_{2:n} <\ldots <x_{n-1:n}<1$ a.s. denote the corresponding order statistics, and let $p_i=x_{i:n}-x_{i-1:n}$ for $i=1,\ldots,n$ be the spacings. Since $\sum_{i=1}^np_i=1$, $(p_1,\ldots,p_n)$ is a random vector on $\Dnn$. When $f$ is the uniform density on $[0,1]$, $(p_1,\ldots,p_n)$ is a Dirichlet distributed random vector with each parameter equalling one. This class of distributions is hence a natural extension of the one used in Algorithm \ref{alg1}.
\begin{corollary}
	\label{nowarp2}
	Suppose $(p_1,\ldots,p_{n}) \in \Dnn$ based on any fixed equi-spaced partition $\Tn$ is the vector of spacings of i.i.d. random variables with a twice differentiable distribution function $F$ and quantile function $Q$. If $f$ is the corresponding probability density, assume that $\inf_{0\leq x \leq 1}f(Q(x))>0$ and $\sup_{0 \leq x \leq 1}|f'(Q(x))|<\infty$. Then, $Y_n$ converges in probability to $Q$ in $C([0,1])$ equipped with the uniform topology.
\end{corollary}
The conditions on $F$ in Corollary \ref{nowarp2} are not too restrictive and are satisfied by several densities with support on $[0,1]$. For example, one can easily check that if $F$ is the distribution function of a non-central Beta \citep{JLH} with both shape parameters equal to one, and a non-centrality parameter $\kappa>0$, then $F$ is twice differentiable with bounded second derivative, and $\inf_{0 \leq x \leq 1}f(x)$ is $e^{-\kappa/2}(\kappa/2+1)>0$, which is attained in the limit at zero. As with Theorem \ref{nowarp} the conclusion of Corollary \ref{nowarp2} holds for a non-equi-spaced partition as well (See Supplementary Material).

\section{Improved construction via random partitions}
\label{random_part}
While the simplicity of Algorithm \ref{alg1} is attractive, the degeneracy of the resulting distribution on $W_I$ (as $n \to \infty$) is disconcerting. In this section, we offer a simple modification of the previous approach that salvages the situation. The sampling method in Algorithm \ref{alg1} depends on the choice of the partition $\Tn$ of $[0,1]$. We demonstrate that choosing a random partition $\mathcal{T}_n(H)$ based on order statistics of an i.i.d. sample from distribution $H$ on $[0,1]$, in conjunction with a point process representation, results in a limit process with sample paths that lie in $W_I$ centred at a desired warp map. The motivation for using partitions induced by order statistics is the fact that conditional on $n$, $t_{i:n}$ have the same distributions as the order statistics of an i.i.d. sample from $h(t)/\int_0^1 h(u)du$, where $h$ is the intensity of a non-homogeneous Poisson process on $[0,1]$. Next, we provide a new algorithm for sampling warp maps in $W_I$, and study the theoretical properties of the associated distribution.
\begin{algorithm}{\textbf{Random partition-based sampling of warp maps on $[0,1]$.}}
	\begin{enumerate}
		\item Choose order statistics $0=:t_{0:n}<t_{1:n}<t_{2:n}<\ldots<t_{n-1:n}<t_{n:n}:=1$ of a random sample from distribution $H$ on $[0,1]$.
		\item Sample an $n$-dimensional Dirichlet distributed random vector $(p_1,\ldots,p_n)$ with parameters set to $\left(t_{1:n}-t_{0:n},\ldots,t_{n:n}-t_{n-1:n}\right)$.
		\item Construct a warp map on $[0,1]$ by linear interpolation of the increments.
	\end{enumerate}\label{alg2}
\end{algorithm}\noindent Algorithm \ref{alg2} is easy to implement, and extends Algorithm \ref{alg1} to one based on random partitions (and a subsequent change in Dirichlet parameters).

\subsection{Theoretical support for Algorithm \ref{alg2}}
As $n \to \infty$, the limiting distribution associated with Algorithm \ref{alg2} can be identified in two ways. The first approach is to start with finite-dimensional distributions at chosen time points and posit the existence of a process with the chosen finite-dimensional projections based on Kolomogorov's consistency theorem (see \cite{RS}). \emph{With this approach, it is then difficult to centre the distribution at a desired warp map}. The approach we adopt in this paper is to constructively define a distribution that arises as a limit (as $n \to \infty$) based on a point process formulation using transformed increments.

We first review the Gamma subordinator process. A process $\G(t),\ t \in [0,1]$ is a Gamma subordinator taking values in $\mathbb{R}_+$ if, for $0\leq s<t\leq 1$, $\G(t)-\G(s)$ is Gamma distributed with shape parameter $t-s$ and scale parameter equal to one. It is a Levy process with Levy measure $\lambda(dy)=y^{-1}e^{-y}dy$ and sample paths that are discrete with probability one, and thereby allows for a point process representation: $\G(t):=\sum_{\nu_x\leq t}\nu_y,\ t \in [0,1]$, where $\nu=(\nu_x,\nu_y)\in [0,1]\times\mathbb{R}_+$ is a Poisson point process with intensity measure $dx \times \lambda(dy)$. For a distribution function (not necessarily one corresponding to a probability measure) $H$ on $[0,1]$, such that $\lim _{x \to 1}H(x)=c<\infty$, consider the time-changed Gamma process $\G{(H(t))}$ whose increments $\G(H(t))-\G(H(s))$ are Gamma distributed with shape parameter $H(t)-H(s)$ and scale parameter equal to one. Then, the normalized Gamma process $t \mapsto \G(H(t))/\G(c)$ is the Dirichlet process $\mathcal{D}(H(t))$ with base measure or parameter $H$, taking values in $[0,1]$. Its sample paths are hence random functions mapping $[0,1]$ to itself. We denote the laws of $\G$ and $\D$ by $\mathbb{G}$ and $\mathbb{D}$, and $\G\circ H$ and $\D \circ H$ by $\mathbb{G}\circ H$ and $\mathbb{D}\circ H$, respectively.

Using the representation of a pure jump Levy process by \cite{FK}, we examine the existence of a limit process with sample paths in $W_I$ with `finite-dimensional' Dirichlet distributions with appropriate parameters. The following result is formulated for the general class of distributions on $\Dn$ induced by spacings of i.i.d. random variables on $[0,1]$ with density $f$ and distribution function $F$, of which the Dirichlet with all parameters set to the same value is a special case. Consider a random partition $\mathcal{T}_n(H)$ of $[0,1]$ based on order statistics $0=:t_{0:n}<t_{1:n}<\ldots<t_{n-1:n}<t_{n:n}:=1$ from an i.i.d sample $\{t_i,\ i=1,\ldots,n-1\}$ with absolutely continuous distribution function $H$ on $[0,1]$. Independent of $t_i$, consider $(p_1,\ldots,p_n) \in \Dnn$ obtained as spacings of an i.i.d sequence $x_i$ from a density $f$ chosen as described in Corollary \ref{nowarp2}. Let $F$ and $Q$ be the corresponding distribution and quantile functions of $x_i$, respectively. Set $v_1=p_1$ and $v_i=p_1+\ldots+p_i,\ i=2,\ldots,n$, and consider the transformed random variables $z_{i,n}=nf(Q(\zeta_{i,n}))v_i$ where $0\leq \zeta_{i,n} \leq 1$ is a deterministic sequence such that $\max_{1\leq i \leq n}|\frac{i}{n}-\zeta_{i,n}|=O(1/n)$. Define $\lambda^{-1}(x):=\inf\{t \in \mathbb{R}:\lambda(t)\geq x\}$.
\begin{theorem}\label{th1}
	Let $\Po_n:=\sum_{i=1}^{n} \delta_{\{t_i,\lambda^{-1}(z_{i,n})\}}$ be a sequence of point processes, where $\delta_{\{x,y\}}$ is the point measure at $(x,y) \in [0,1]\times \mathbb{R}_+$. Assume $f$ to be continuous and positive on $(0,1)$ with $\lim_{x \downarrow 0}f(x)>0$ and finite.
	\begin{enumerate}[(1)]
		\item Suppose $p_i$ for each $i=1,\ldots,n$ possesses a unimodal density. The sequence $\Po_n$ converges in the vague topology to the Poisson point process $\Po$ with intensity measure $H(dt) \times \lambda(dy)$, where $\lambda(dy)=y^{-1}e^{-y}dy$.
		\item If $\Po_n$ converges in the vague topology to $\Po$, then the sequence of processes $G_n(t):=\sum_{i}\lambda^{-1}(z_{i,n})\mathbb{I}_{t_i\leq t},\ t \in [0,1]$ converges weakly to the time-changed Gamma process $\G\circ H$ in the Skorohod $J_1$ topology. However, the linearly interpolated version of $G_n$ converges to $\G\circ H$ in the Skorohod $M_1$ topology.
	\end{enumerate}
\end{theorem}
Evidently, the limit process $\Po$ is the Gamma process $\mathcal{G}\circ H$, and $\Po_n$, when normalized, converges weakly to a limit, which we refer to as the \emph{Dirichlet process} $\mathcal{D} \circ H$ with $dH$ as the base measure, where $dH$ is the Lebesgue-Stieltjes measure corresponding to the distribution function $H$. The probability measure $\mathbb{D}\circ H$ on the class $W_I$ is constructed \emph{only} using the increments $p_i=\gamma(t_{i:n})-\gamma(t_{i-1:n})$ of warp maps on a random partition $\Tn(H)$ based on $t_i$ i.i.d. from $H$. \emph{Thus, Theorem \ref{th1} ensures that the finite-dimensional Dirichlet distributions can be identified with the finite-dimensional distributions of $\mathcal{D}\circ H$}, and form a consistent family in the Kolmogorov sense \citep{RS}. It can be interpreted as follows: conditioned on a partition $\Tn(H)=\{0=:t_{0:n}<t_{1:n}<\ldots<t_{n-1:n}<t_{n:n}:=1\}$,
\begin{align}
\label{dirichlet}
&\mathbb{D}\circ H(\gamma({t_{1:n}})\in dx_1, \ldots, \gamma({t_{n-1:n}})\in dx_{n-1})= \nonumber\\
&=\frac{\Gamma(1)}{\prod_{i=1}^{n} \Gamma(t_{i:n}-t_{i-1:n})}
\prod_{i=1}^{n} (x_{i}-x_{i-1})^{(t_{i:n}-t_{i-1:n})}dx_1\ldots dx_{n-1},
\end{align}
where $t_{n:n}=x_n=1$ and $t_{0:n}=x_0=0$. A few remarks are in order at this stage.

\noindent\textbf{Remark 1.} The Ferguson-Klass representation is based on the transformation of a homogeneous Poisson random measure on $\mathbb{R}_+ \times \mathbb{R}_+$  under the inverse of the tail of a Levy measure $g$ on $\mathbb{R}_+$. Theorem \ref{th1} uses this representation with Levy measure $\lambda$ to (1) retain the simplicity of Algorithm \ref{alg1}, (2) ensure that the distribution on $W_I$ satisfies subset invariance, and (3) ensure that the distribution on $W_I$ can be centred at a desired warp map. Proposition \ref{properties} in the next section clarifies (2) in view of Lukac's characterization of the Gamma distribution; in other words, the subset invariance requirement fixes the Levy measure to be the Gamma tail measure $\lambda(dx)=x^{-1}e^{-x}dx$. For a \emph{fixed} partition, the measure $\mathbb{D}\circ H$ with $H(t)=t$ was constructed by \cite{RS} (Proposition 3.4) using Kolmogorov's extension theorem; as mentioned earlier, such an approach cannot be used to obtain a distribution that can be centred at any desired warp map.

\noindent\textbf{Remark 2.} Let $H$ be the uniform distribution function. The definition of the probability measure $\mathbb{D}\circ H=\mathbb{D}$ on $W_I$ does not identify the increments of a $\gamma$ with the set of probability measures on a finite set. This is in contrast to Ferguson's Dirichlet process, say $\bar{\D}$, which is constructed on the set $\mathbb{P}([0,1])$ of probability measures on $[0,1]$, topologized by weak convergence and indexed by Borel sets of $\mathbb{P}([0,1])$. Suppose $W_I$ is identified with distribution functions on $[0,1]$. Consider the homeomorphism $h:W_I \to \mathbb{P}([0,1])$ that assigns to each $\gamma \in W_I$ its Stieltjes measure $d\gamma$. Equip $W_I$ now with the image of the weak topology on $\mathbb{P}([0,1])$ under the map $h^{-1}: \mathbb{P}([0,1]) \to W_I$. Then, the law of $\mathcal{\bar{D}}$ can be viewed as the push-forward $h_*\mathbb{D}$ of $\mathbb{D}$. On the other hand, if $\bar{h}: \mathbb{P}([0,1])\to W_I$ such that $\bar{h}$ assigns to each $\mu \in \mathbb{P}([0,1])$ a function $\gamma(t):=\sup\{u \in [0,1]: \mu[0,u]\leq t\}$, then $W_I$ is identified with quantile functions on $[0,1]$. Interestingly, the pullback of $\mathbb{D}$ on $W_I$ is topologically very different from the law of Ferguson's process (see p.1131 of \cite {RS}). Importantly, when considering warp maps of $\sone$, the relationships to distribution or quantile functions of measures on $\sone$ are unavailable.

\noindent\textbf{Remark 3.} Part (2) of Theorem \ref{th1} is striking: linear interpolation of the sample paths of $G_n$ does not affect the limit process. This has an important implication for Algorithm \ref{alg2} using random partitions, while not giving up continuity of obtained warp maps. The weaker Skorohod's $M_1$ topology is used since the linear interpolation of $G_n$ implies that we seek convergence to a limit jump process with unmatched jumps in the converging sequence of processes; this cannot be achieved under the usual $J_1$ topology.

\noindent\textbf{Remark 4.} Part (1) of Theorem \ref{th1} states that the limit process is unchanged if the vector of increments of the warp map $\gamma$ is assumed to have a distribution on the simplex $\Dnn$ obtained as the spacings of i.i.d. random variables from an arbitrary density on $[0,1]$. This, in a certain sense, makes the Dirichlet process a natural choice as a probability measure on $W_I$. The assumption of unimodality of the densities of $p_i$ is not critical and can be relaxed at some technical cost.

\section{Properties of the proposed distribution}
\label{properties_section}

Next, we study the theoretical properties of the constructed distribution $\mathbb{D}\circ H$ on warp maps of $[0,1]$. We also consider the case when landmark constraints need to be enforced during the registration process.

\subsection{Automatic regularization}

The limit Gamma process $\G \circ H$ (and hence $\D \circ H$) in Theorem \ref{th1} is a pure jump process with a.s. discrete paths, which leads to a pure jump warp map $\gamma$. Nonetheless, for a fixed $u \in [0,1]$, the function $u \mapsto  \gamma(u)$ is continuous at $u$ $\mathbb{D} \circ H$-almost surely since $\G \circ H$ is a Levy process, and hence stochastically continuous: $|\G(H(t+u)) -\G(H(u))|\overset{P}\to 0$ as $t \to 0$ since $H$ is absolutely continuous. We gather from Part (2) of Theorem \ref{th1} that linear interpolation does not affect the discrete nature of the limit Gamma process. For the alignment problem, the number and frequency of large jumps are pertinent since Algorithm \ref{alg2} is based on choosing a fine partition (large $n$); the occurrence and likelihood of warp maps that contain regions of `large warping' are particularly important. Moreover, this also sheds light on how likely we are to sample warp maps from $\mathbb{D}\circ H$ that deviate considerably from their average $H$.


Based on a vector of increments $(p_1,\ldots,p_n)$ uniformly distributed on $\Dnn$, we consider the `large jumps' process defined as the real-valued partial sum process $Y_{n}(t):=\sum_{i=1}^{\nt}[\xi_i-E(\xi_i\mathbb{I}_{\xi_i\leq 1})],\ 0 \leq t \leq 1$ taking values in $D([0,1])$, where $\xi_{i,n}:=np_i-\log n$ is a triangular array sequence. The definition of $\xi_{i:n}$ is motivated by the fact that the largest jump amongst the $p_i$ is of size $O_P(\log n/n)$ when $(p_1,\ldots,p_n)$ is uniform on $\Dn$ \citep{LD}. Jumps of smaller order can also be studied using intermediate spacings amongst the $p_i$ under a different normalizing transformation (see \cite{NBZ} for details). 
\begin{theorem}
	\label{th2}
	Let $Y(t),\ t\in[0,1]$ be a real-valued Levy jump process with Levy measure $\nu(dy)=e^{-y}dy$. The sequence $Y_n$ converges weakly to $Y$ in $D([0,1])$ equipped with the Skorohod $J_1$ topology.
\end{theorem}
\noindent Theorem \ref{th2} is a functional limit theorem for the process $Y_n$, and describes the probabilistic behavior of Algorithm \ref{alg2} in generation of warp maps which contain regions of high warping. The centering term in Theorem \ref{th2} cannot be dispensed with, although other equivalent terms can be chosen. The arrivals of the large increments are governed by a finite activity Levy process since $\int_0^\infty e^{-u}du<\infty$, which implies that large jumps occur infrequently; $Y$ is hence a compound Poisson process. This ensures that under Algorithm \ref{alg2}, once an $H$ that generates the random partition and sets the average warp map is chosen, \emph{we are not likely to sample warp maps from $\mathbb{D}\circ H$ that contain far more large jumps relative to those in the average warp.} Thus, Algorithm \ref{alg2} offers automatic regularization toward the average warp map.

\subsection{Landmark constraints and a global concentration parameter}
\label{parametric}

In practice, in the presence of landmarks, one decomposes the unconstrained registration problem into multiple sub-problems corresponding to intervals formed due to the landmark constraints. For example, in the case of $m=2$ landmarks, suppose that the landmark locations on the domain of two open curves $g_1$ and $g_2$ are at $t_i,\ t_{k}\in [0,1]$, and at $t_{j},\ t_l\in [0,1]$, respectively, with $1\leq i<k\leq n-1$ and $1\leq j<l\leq n-1$. This induces three intervals of interest on each curve: $\{[0,t_i],[t_i,t_k],[t_k,1]\}$ for $g_1$ and $\{[0,t_j],[t_j,t_l],[t_l,1]\}$ for $g_2$. The points $t_i$ and $t_k$ can now be exactly matched to $t_j$ and $t_l$ using a PL warp map, resulting in two new points $t_i^*$ and $t_k^*$. This leads to three classes of warp maps: $W_1=\{\gamma:[0,t_i^*] \to [0,t_i^*]\}$, $W_2=\{\gamma:[t_i^*,t_k^*] \to [t_i^*,t_k^*]\}$, and $W_3=\{\gamma:[t_k^*,1] \to [t_k^*,1]\}$. The original registration problem on $[0,1]$ has thus been decomposed into three similar ones: (1) match $g_1|_{[0,t_i^*]}$ and $g_2|_{[0,t_i^*]}$, (2) match $g_1|_{[t_i^*,t_k^*]}$ and $g_2|_{[t_i^*,t_k^*]}$, and (3) match $g_1|_{[t_k^*,1]}$ and $g_2|_{[t_k^*,1]}$, where $g|_{[a,b]}$ is the restriction of $g$ to the subinterval $[a,b]$. Note that $g_1|_{[0,t_i^*]}(t_i^*)=g_2|_{[0,t_i^*]}(t_i^*)=g_1|_{[t_i^*,t_k^*]}(t_i^*)=g_2|_{[t_i^*,t_k^*]}(t_i^*)$, and $g_1|_{[t_i^*,t_k^*]}(t_k^*)=g_2|_{[t_i^*,t_k^*]}(t_k^*)=g_1|_{[t_k^*,1]}(t_k^*)=g_2|_{[t_k^*,1]}(t_k^*)$; such relationships are true for all warp maps in the classes $W_1,\ W_2$ and $W_3$.

The above decomposition requires a property of warp maps called subset invariance. The corresponding requirement on a probability measure on $W_I$ is the following. Consider the map $S({[a,b]}):W_I \to W_I$,
\begin{equation}
\label{subset_invariance}
S({[a,b]})(\gamma(u))=\frac{\gamma((1-u)a+ub)-\gamma(a)}{\gamma(b)-\gamma(a)}, \quad u\in [0,1],\ 0\leq a<b\leq 1.
\end{equation}
A parametrized probability measure $\mathbb{P}_\theta,\ \theta \in \Theta$ on $W_I$ is said to satisfy subset invariance if its push-forward under the map $S([a,b])$, $S([a,b])_\#\mathbb{P}_\theta$ is $\mathbb{P}_{\theta(b-a)}$. To facilitate subset invariance, and to ensure that the three sub-problems borrow strength, we introduce a concentration parameter $\theta>0$, which allows the distribution to interpolate between the indicator (step) function and the average warp map. Under the notation employed in the preceding section, consider the parametrized Levy measure $\lambda_\theta(y)=\theta\int_y^\infty e^{-t}t^{-1}dt$. The limit process $\Po$ in Theorem \ref{th1} is then a Gamma process  with Levy measure $\lambda_\theta$, such that, for  any Borel set $A$ and any $0 \leq s<t \leq 1$,
$$P(\G(\theta t)-\G(\theta s)\in A) =\int_A \frac{1}{\Gamma(\theta(t-s))}y^{\theta(t-s)-1}e^{-y}dy.$$
Then $\mathcal{D}_{\theta}(t)=\G(\theta t)/\G(\theta),\ t \in [0,1]$ is the corresponding Dirichlet process. 
\begin{proposition}
	\label{properties}
	The distribution $\mathbb{D}_\theta \circ H$ satisfies the following properties \citep{RS}.
	\begin{enumerate}[(1)]
		\itemsep .5em
		\item \underline{\emph{Concentration around the mean}}: For the partition based on $H$, $E_\theta(\gamma_t)=H(t)$ and Var$_\theta(\gamma_t)=\frac{1}{1+\theta}H(t)(1-H(t))$,\ $t \in [0,1],\ \theta>0$.
		\item \underline{\emph{Subset invariance}}: For the map $S([a,b])$ in Equation (\ref{subset_invariance}), the push-forward measure $S([a,b])_\#\mathbb{D}_\theta \circ H$ equals $\mathbb{D}_{\theta(H(b)-H(a))} \circ H.$
		\item \underline{\emph{Markov-type property}}: For $\theta>0$, the push-forward measures $S({[a,b]})_\#\mathbb{D}_\theta \circ H$ and $S({[0,1]\backslash[a,b]})_\#\mathbb{D}_\theta \circ H)$ depend on each other only at the endpoints $a$ and $b$.
		\item As $\theta \to 0$, $\mathbb{D}_\theta \circ H$ converges
		to a uniform distribution on the subset $\{\gamma:[0,1]\to[0,1]:\gamma(t)=\mathbb{I}_{[0,H(t)]} \}$ of $W_I$.
		\item As $\theta \to \infty$, $\mathbb{D}_\theta \circ H$ converges to the point mass distribution $\delta_{H}$ at $H$.	
		\item Let $\bar{\gamma}:[0,1]\to [0,1]$ be continuous and increasing. For every $\theta>0$,  $\bar{\gamma}^{-1}_*\mathbb{D}_\theta \circ H$ is absolutely continuous with respect to  $\mathbb{D}_\theta \circ H$.
	\end{enumerate}
\end{proposition}

The proof of $(1)$ follows from direct computation using Equation (\ref{dirichlet}). Proofs of $(2)$ and $(3)$ are easily obtained from the representation of the Dirichlet process as a normalized Gamma process, $\mathcal{D}_\theta(t)=\G(\theta t)/\G(\theta)$, and the fact that $\mathcal{D}_\theta(t)$ is independent of $\G(\theta)$ for every $t\in [0,1]$ and $\theta>0$. The independence and self-similarity properties are hence based on Lukac's characterization of the Gamma distribution: if $X_1,\ldots,X_n$ are independent Gamma random variables, then $Y=\sum_{i=1}^n X_i$ and the vector $(X_1/Y,\ldots,X_n/Y)$ are independent. For proofs of $(3)$ and $(4)$, we refer the reader to the proof of Proposition 3.14 by \cite{RS}. Properties $(1),\ (4)$ and $(5)$ ensure that the distribution on $W_I$ is centred at $H$. The parameter $\theta$ behaves like a concentration parameter: varying $\theta$ moves mass away from or toward $H$, and offers a rich class of probability models for $W_I$. Property $(6)$ is a crucial distributional property: $W_I$ is closed under composition, and hence any distribution on $W_I$ should be quasi-invariant with respect to composition. Its proof can be found in Theorem 4.3 of \cite{RS}.

With the introduction of $\theta>0$, the distribution $\mathbb{D}_\theta\circ H$ has finite-dimensional projections interpreted in the following manner: conditioned on a partition $\Tn(H)=\{0=:t_{0:n}<t_{1:n}<\ldots<t_{n-1:n}<t_{n:n}:=1\}$,
\begin{align}
\label{dirichlet}
&\mathbb{D}_\theta\circ H(\gamma({t_{1:n}})\in dx_1, \ldots, \gamma({t_{n-1:n}})\in dx_{n-1})= \nonumber\\
&=\frac{\Gamma(\theta)}{\prod_{i=1}^{n} \Gamma(\theta(t_{i:n}-t_{i-1:n}))}
\prod_{i=1}^{n} (x_{i}-x_{i-1})^{\theta(t_{i:n}-t_{i-1:n})}dx_1\ldots dx_{n-1}.
\end{align}

\noindent\textbf{Remark 5.} Part (2) of Proposition \ref{properties} identifies $\mathbb{D}_\theta\circ H$ with the subset invariance property. \emph{In other words, if subset invariance is a requisite property for a distribution on warp maps, with or without landmarks, it is not possible to construct one that does not concentrate on discontinuous warp maps.} This stands in contrast to the usual smoothness assumptions associated with warp maps \citep{RS2,CSS}, but is rarely an issue in practice. Intuitively, decomposition of alignment based on landmarks can be linked with an independent increments property of the stochastic process. Furthermore, under the popular square-root velocity transform that we employ in the applications for alignment (see Section \ref{applications} for details), \cite{LRK} proved that the optimal warp map is PL when at least one of the two curves to be aligned is PL (see Theorem 6). The lack of smoothness of the sample paths of $\D \circ H$ is thus not unrealistic.

\section{Extension to distributions on warp maps of $\sone$}
\label{circle}

Our aim in this section is to construct a probability measure on the set $W_\sone:=\{\gamma:\sone \to \sone: \text{continuous, orientation preserving}\}$ of warp maps of $\sone$. From a practical perspective, we wish to develop an easy-to-implement sampling method, similar to Algorithm \ref{alg2}, for generating random warp maps of $\sone$ based on a suitable discretization. Our approach `unwraps' $\sone$ at a specific point $c$ and proceeds to identify $W_\sone$ with the product space $W_I \times \sone$ through the identification of $\sone$ with $[0,1]$. This amounts to using the probability measure $\mathbb{D}_\theta\circ H$ on $W_I$ along with one on $\sone$, based on viewing the circle of unit circumference $\sone$ as the quotient group $\sone=\mathbb{R}/2\pi\mathbb{Z}$ with the addition operation inherited from $\mathbb{R}$. We thus move from $[0,1]$ to the circle with unit length $\sone$ by identifying the endpoints of the interval.

Through the identification of $\sone$ with $[0,1]$, every continuous mapping $\beta: \mathbb{R} \to \mathbb{R}$ induces a continuous mapping of $\sone$ onto itself such that $\beta(t+j)=\beta(t)+j$ for all $t \in \mathbb{R}$, where $j$ is an integer ($\beta$ is unique up to addition of an integer and $\beta(t)-t$ is periodic with period $j$). If $\beta$ is monotone increasing and $j=+1$, we say that the induced map on $\sone$ is orientation-preserving (based on a choice of clockwise or anti-clockwise orientation). Specifically, consider the class $W_\mathbb{R}:=\{\beta:\mathbb{R}\to \mathbb{R}: \beta(t+1)=\beta(t)+1, \text{ continuous and non-decreasing}\}.$ Each member $\beta$ of $W_\mathbb{R}$ induces a warp map $\tilde{\beta}:\sone \to \sone$ with $\tilde{\beta}(e^{2\pi it})=e^{2\pi i\beta(t)}$, where $\beta$ is referred to as the lift of $\tilde{\beta}$. This $\beta$ satisfies $\beta(t+1)=\beta(t)+1$ for all $t \in [0,1]$, and consequently we have, for $t \in [0,1]$, $\beta(t)=\gamma(t)+c $, where $\gamma$ is a warp map of $[0,1]$ and $c \in (0,1]$ (through the identification of $[0,1]$ with $\mathbb{R}/2\pi \mathbb{Z}$). This procedure can be viewed as one that produces a warp map of $\sone$ by `unwrapping' $\sone$ at a chosen point $c$ and generating a warp map of $[0,1]$.

The random version of this corresponds to choosing a $c$ according to a non-atomic probability measure $\mu$ on $[0,1]$, independent of $\mathbb{D}_\theta \circ H$, resulting in a product probability measure $\mu \times \mathbb{D}_\theta \circ H$ on $(0,1] \times W_I$.
The procedure outlined above induces a bijection between the set $W_\sone$ and $(0,1] \times W_I$, and results in the following algorithm.
\begin{algorithm}{\textbf{Random partition-based sampling of warp maps on $\sone$.}}
	\begin{enumerate}
		\item Choose $c$ from $\mu$ on $[0,1]$.
		\item Sample $\gamma$ from $\mathbb{D}_\theta \circ H$ using Algorithm \ref{alg2}.
		\item Set $\gamma_s(t):=(\gamma(t)+c)\mod 1$.
	\end{enumerate}\label{alg3}
\end{algorithm}

\begin{figure}[!t]
\begin{center}
	\includegraphics[scale=0.3]{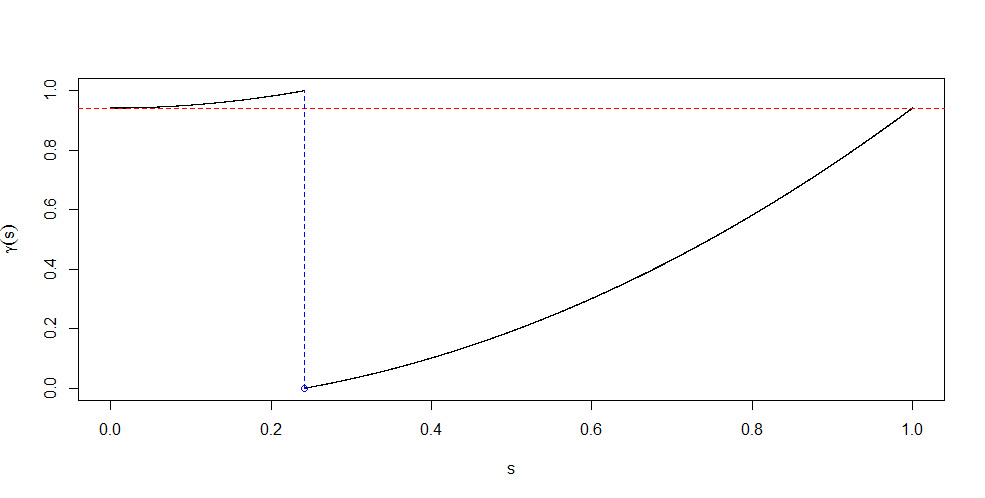}
	\caption{A warp map $\gamma_s$ based on $\gamma(t)=t^2$ and $c=0.94$ with $t_c=0.24$.} 
	\label{warpcircle}
\end{center}
\end{figure}

This method was used in the work of \cite{GMW} while constructing random homeomorphisms of $\sone$. The map $\gamma_s$ is a warp map of $[0,1]$ with a single point of discontinuity $t_c \in [0,1]$ at which $\gamma(t_c)+c=1$, thereby ensuring that $\gamma_s(t_c)=0$. The point $t_c$ is unique to $\gamma_s$.  Figure \ref{warpcircle} offers an illustration of this approach for $\gamma(t)=t^2$ with $c=0.94$, leading to $t_c=0.24$.
Proposition \ref{warp_circ} formalizes this for a random $\gamma_s$ generated in this fashion.
\begin{proposition}
	\label{warp_circ}
	Conditional on $c$ from $\mu$, for each $\gamma_s$ from $\mathbb{D}_\theta\circ H$, the following hold with probability one:
	\begin{enumerate}[(1)]
		\item $\gamma_s(0)=c$.
		\item A unique $t_c$ exists in the interior of $[0,1]$ such that $\lim_{t \uparrow t_c}\gamma_s(t)=0$ and $\lim_{t \downarrow t_c}\gamma_s(t)=1$.
	\end{enumerate}
\end{proposition}
Thus, starting with probability measures $\mathbb{D}_\theta \circ H$ on $W_I$ and $\mu$ on $(0,1]$, the sampling scheme induces the product probability measure $ \mu \times \mathbb{D}_\theta \circ H$ on the set of warp maps of $\sone$. The corresponding measure, independent of the unwrapping point $c$, can be obtained by integrating the product probability measure with respect to $\mu$. Note that the bijection ensures that (trivially) the resulting measure on warp maps of $\sone$ is necessarily absolutely continuous with respect to the probability measure on $W_I$ used in the construction. The Supplementary Material contains an `intrinsic' construction that circumvents the need to unwrap $\sone$. But, the construction provided above is easier to implement in practice, and is thus the only one used in subsequent applications.

\section{Numerical illustrations}
\label{examples}
The degeneracy phenomenon in Algorithm \ref{alg1} can arise due to two reasons: (1) for large values of the concentration parameter $\alpha$ regardless of the size $n$ and manner of construction of the deterministic partition $\Tn$; (2) moderate values for $\alpha$ but large $n$ regardless of the manner of construction of the deterministic partition $\Tn$. High values of the concentration parameter $\theta$ in Algorithm \ref{alg2} also results in a distribution $\mathbb{D}_\theta \circ H$ highly concentrated around $H$. We first illustrate this behavior in a simulation exercise. We then study the utility of Algorithms \ref{alg2} and \ref{alg3} in alignment problems in the context of real data analysis based on finite-dimensional distributions of the process $\mathcal{D}_\theta\circ H$ with law $\mathbb{D}_\theta\circ H$ in view of Theorem \ref{th1}. Note that the theorems merely offer theoretical support for the two algorithms; in particular, there is no need to transform the increments with the inverse Levy measure.

\subsection{Simulation examples}
\label{simulations}
\begin{figure}[!htb]
\begin{center}
\begin{tabular}{|cccc|}
\hline
$n=20$ & $n=100$ & $n=300$ & $n=500$\\
\includegraphics[trim=6cm 2cm 6cm 2cm, clip=true,width=1.3in]{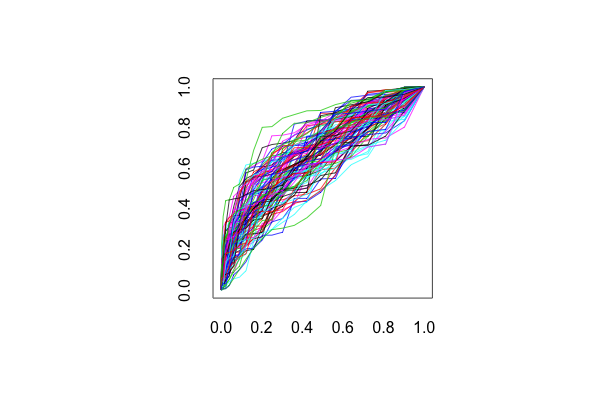}&
\includegraphics[trim=6cm 2cm 6cm 2cm, clip=true,width=1.3in]{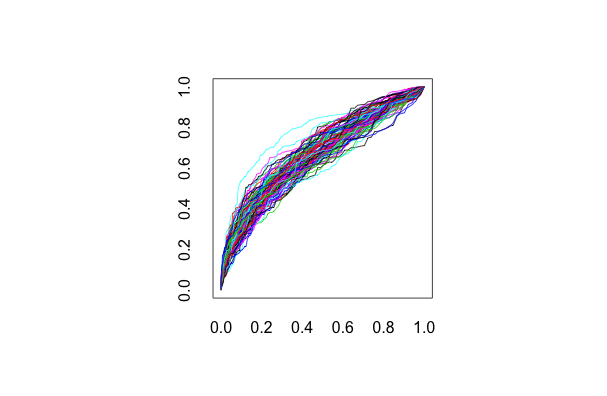}&
\includegraphics[trim=6cm 2cm 6cm 2cm, clip=true,width=1.3in]{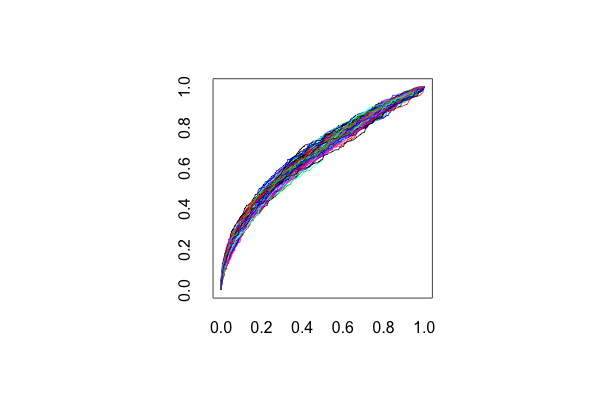}&
\includegraphics[trim=6cm 2cm 6cm 2cm, clip=true,width=1.3in]{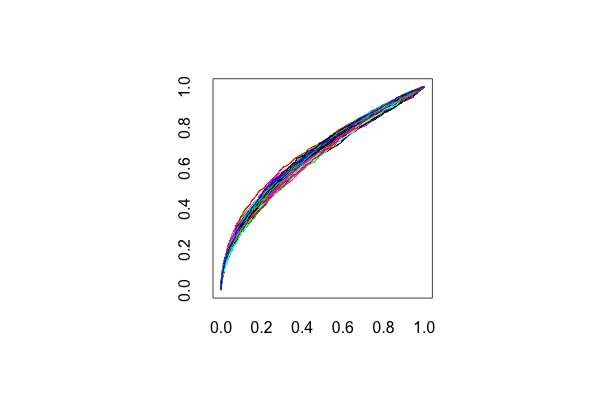}\\
\hline
\includegraphics[trim=6cm 2cm 6cm 2cm, clip=true,width=1.3in]{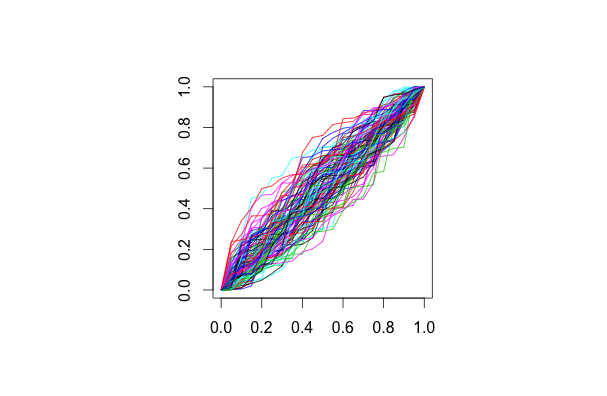}&
\includegraphics[trim=6cm 2cm 6cm 2cm, clip=true,width=1.3in]{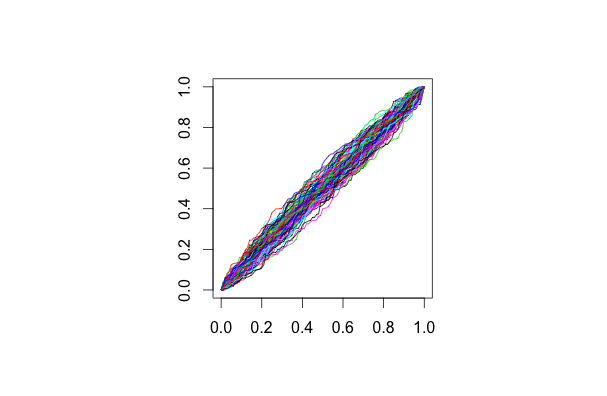}&
\includegraphics[trim=6cm 2cm 6cm 2cm, clip=true,width=1.3in]{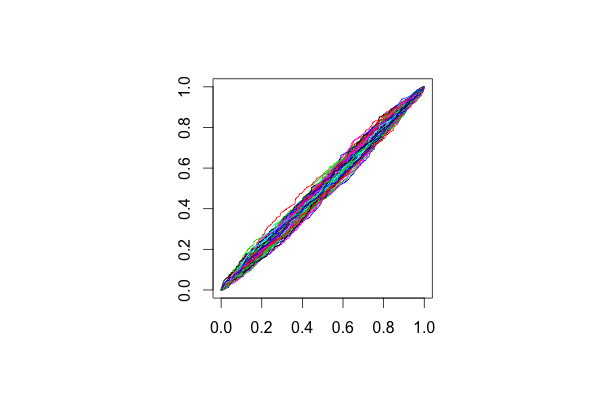}&
\includegraphics[trim=6cm 2cm 6cm 2cm, clip=true,width=1.3in]{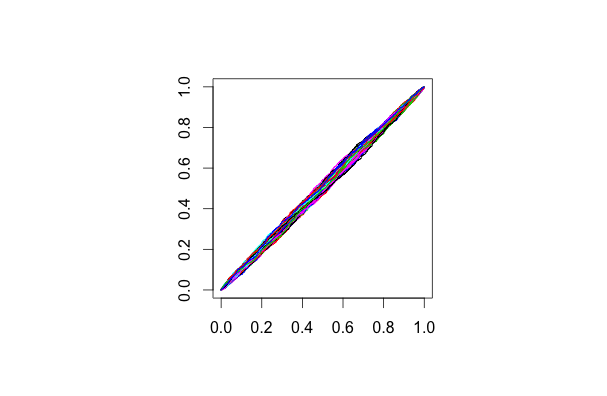}\\
\hline
\end{tabular}
\caption{\small{Sample warp maps from Algorithm \ref{alg1} based on different deterministic partitions $\Tn$ and sizes $n$ with $\alpha=1.2$. Top: $\Tn$ is a non-equi-spaced partition based on a Beta(2,1). Bottom: Equi-spaced partition based on a U[0,1].}}
\label{degenwarp}
\end{center}
\end{figure}

\begin{figure}[!t]
	\begin{center}
		\begin{adjustbox}{max width=\textwidth}
			\begin{tabular}{|c|ccc||ccc|}
				\hline
				&\multicolumn{3}{c||}{(a) Identity average warp}&\multicolumn{3}{|c|}{(b) Non-identity average warp}\\
				\hline
				&\small{$\theta=0.1$}&\small{$\theta=10$}&\small{$\theta=100$}&\small{$\theta=0.1$}&\small{$\theta=10$}&\small{$\theta=100$}\\
				\hline
				\small{$n=5$}&\includegraphics[width=.8in]{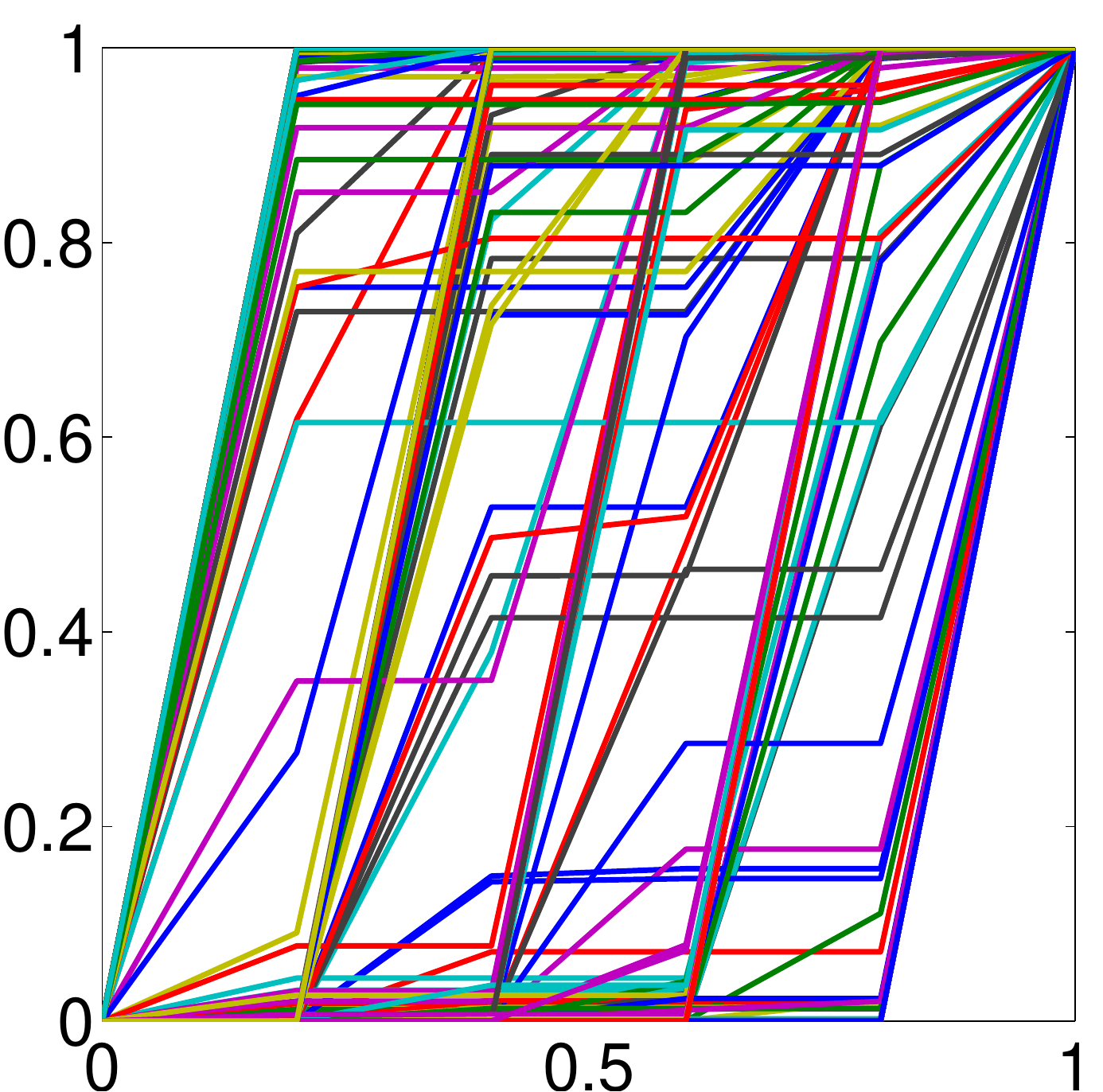}&\includegraphics[width=.8in]{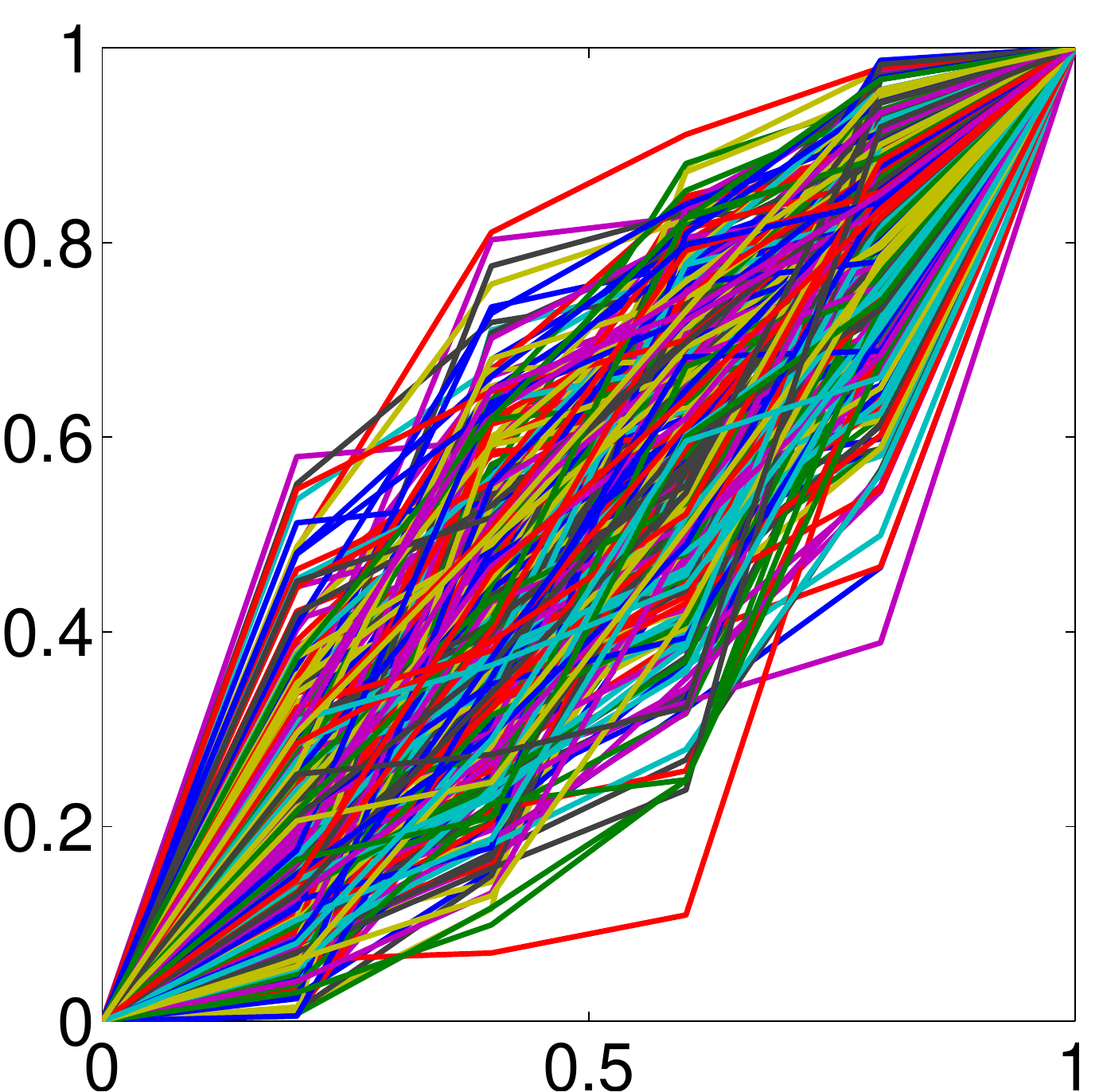}&\includegraphics[width=.8in]{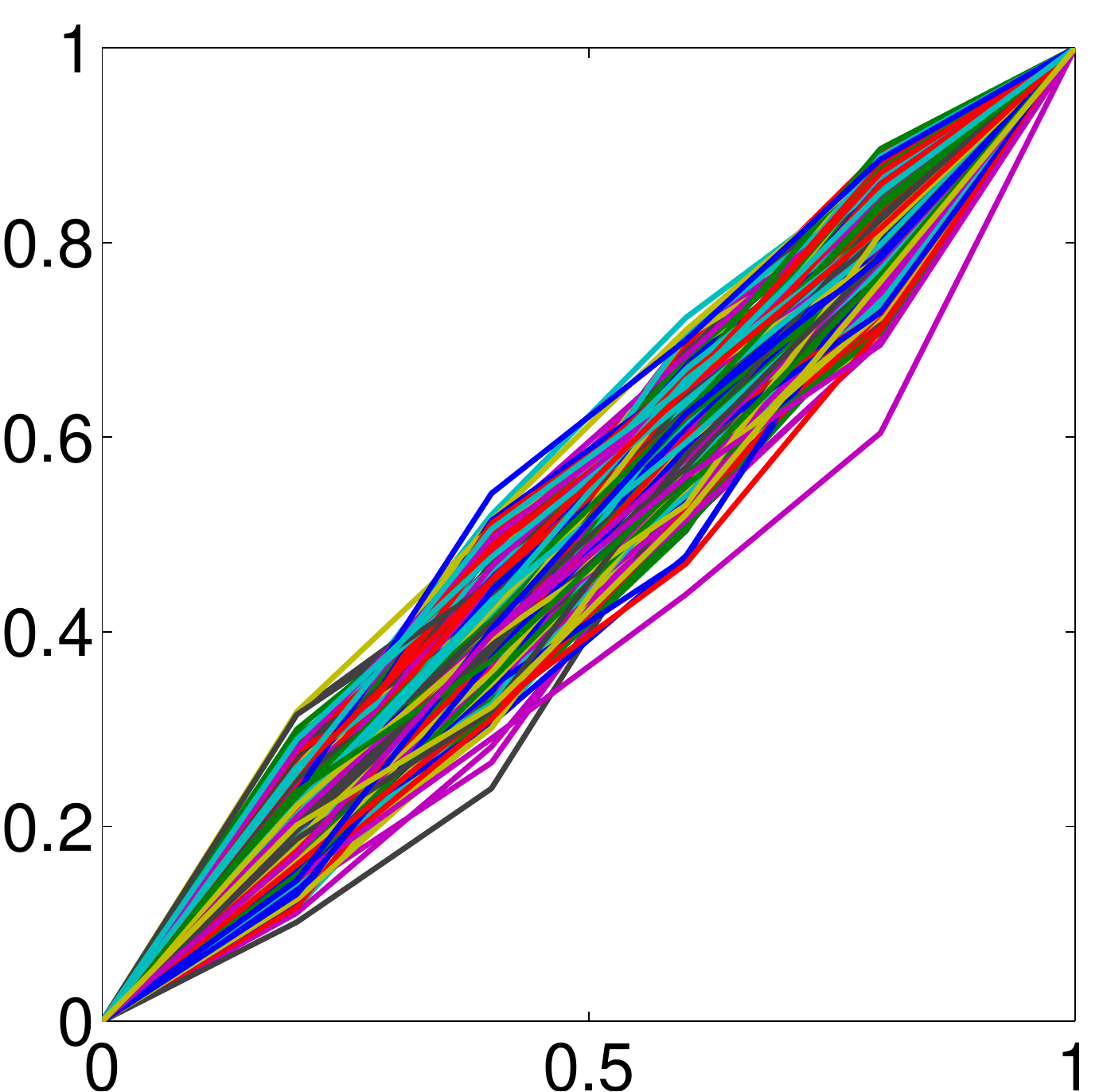}&
				\includegraphics[width=.8in]{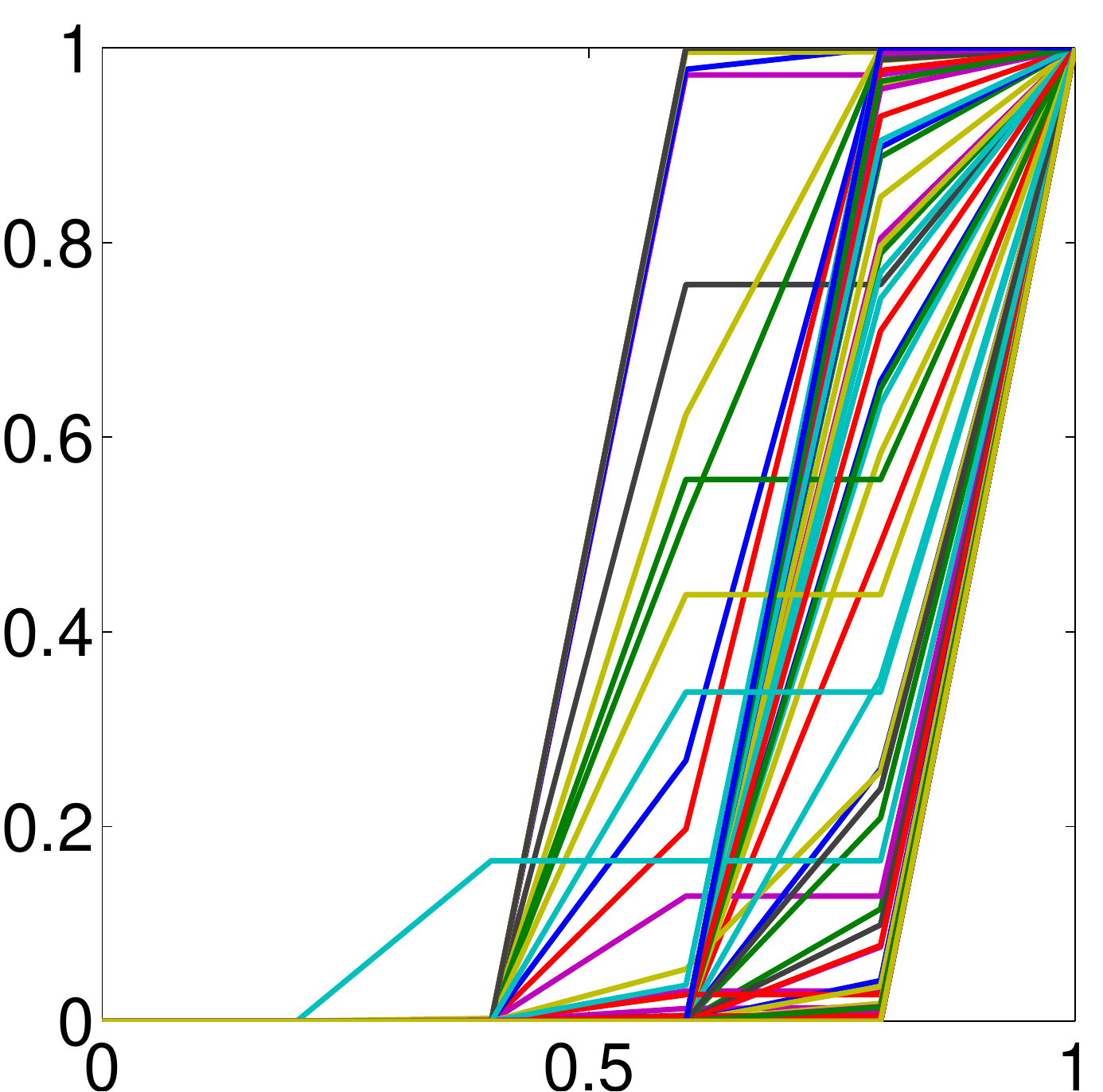}&\includegraphics[width=.8in]{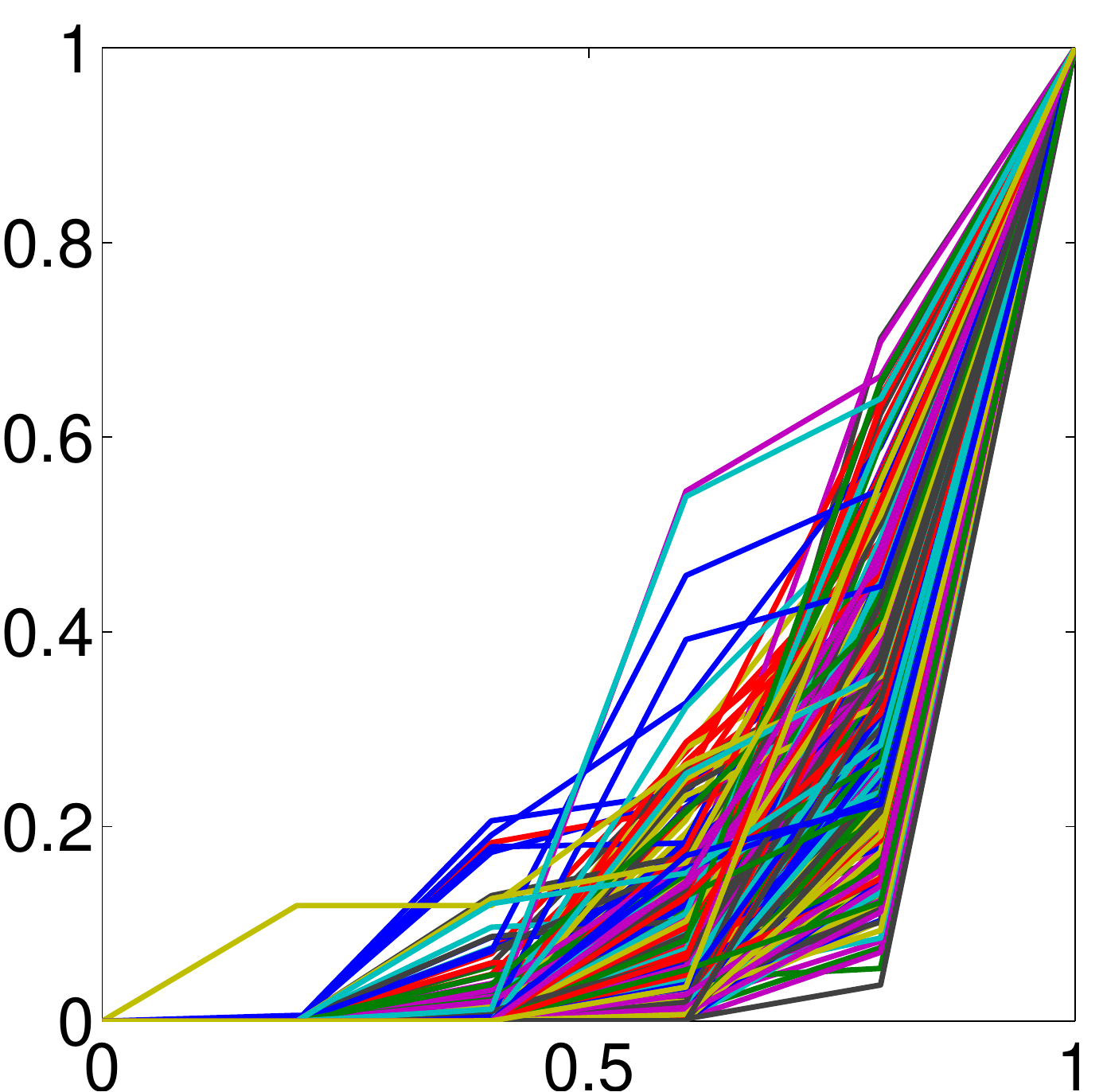}&\includegraphics[width=.8in]{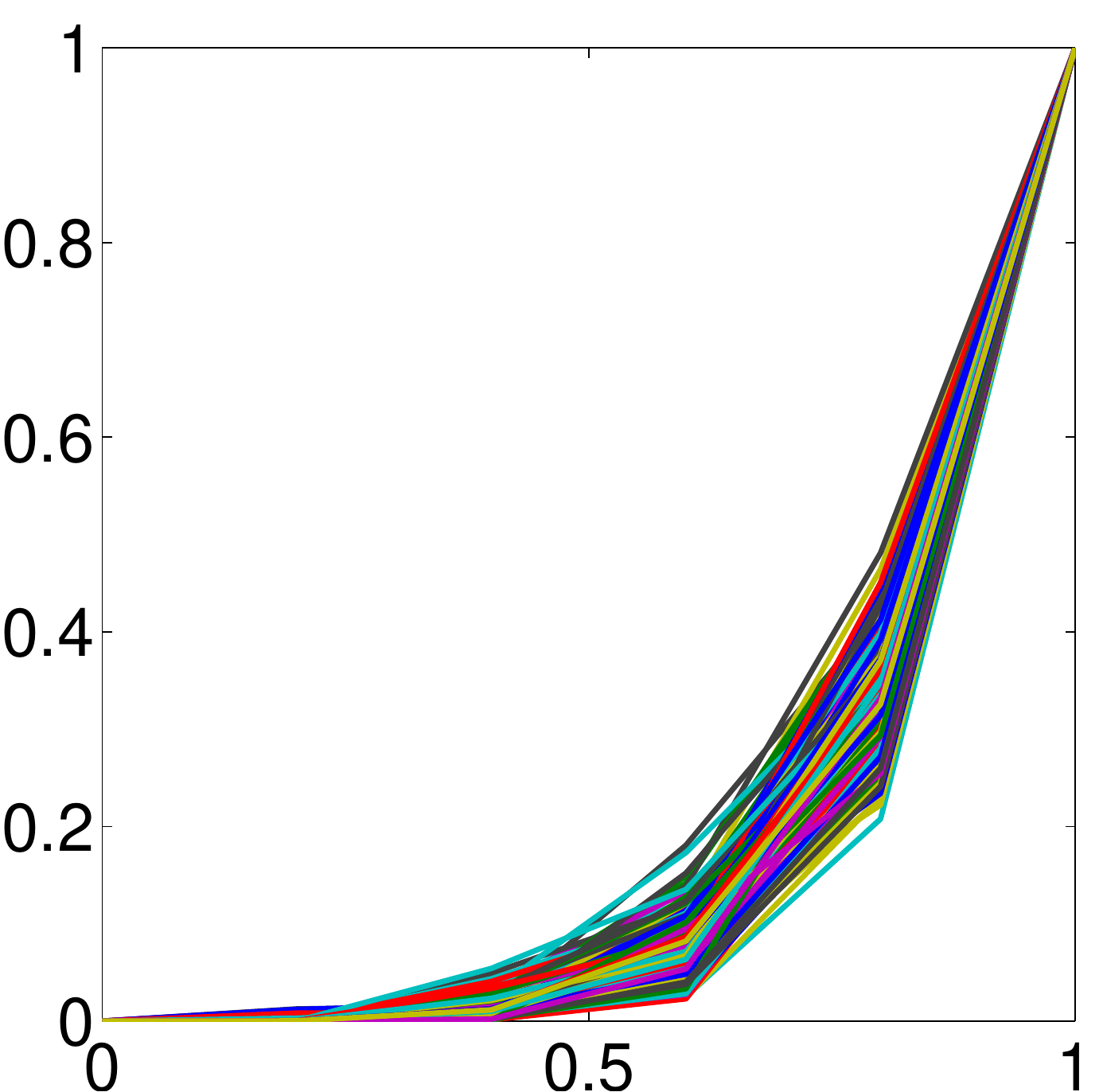}\\
				\small{$n=20$}&\includegraphics[width=.8in]{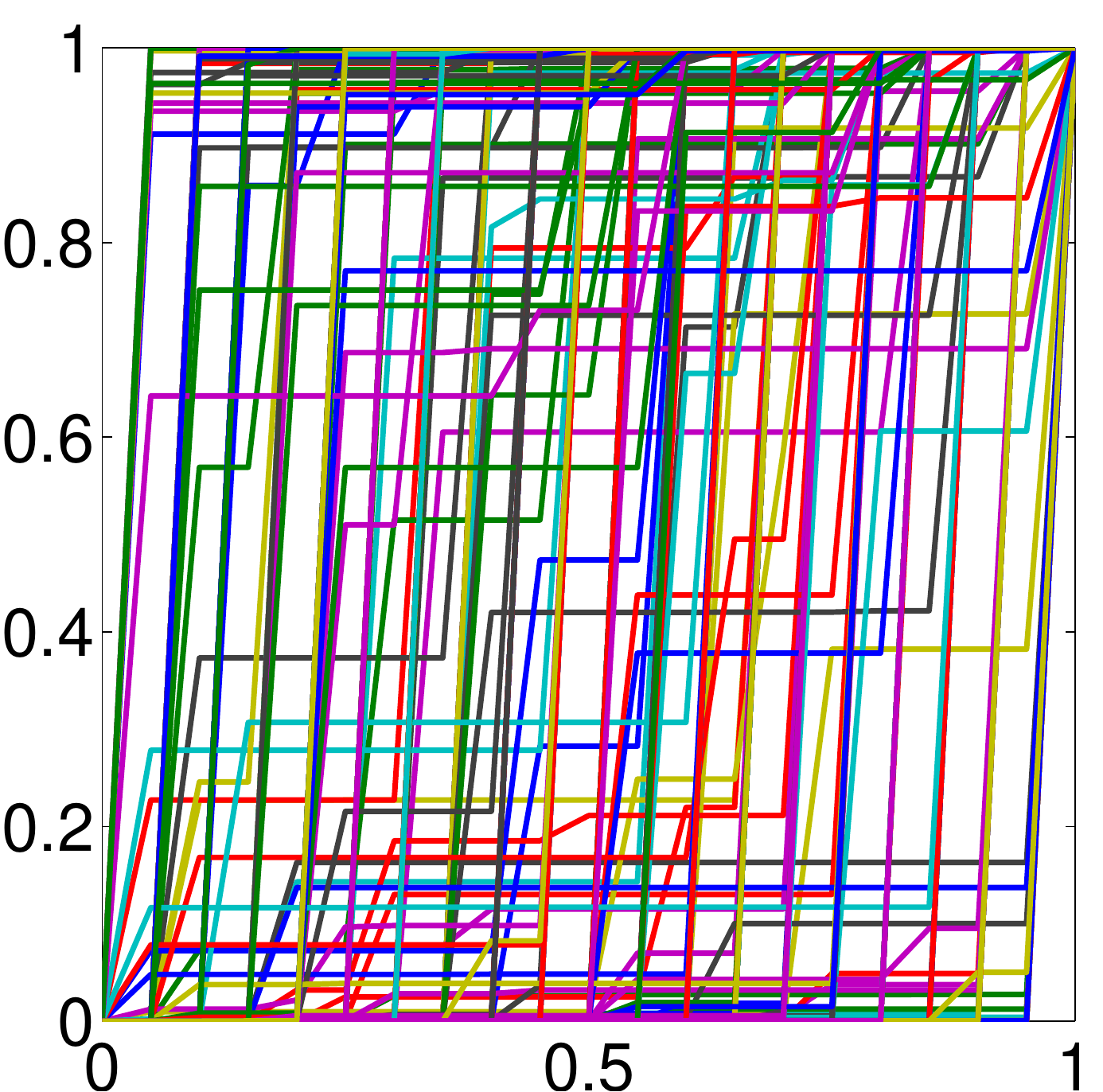}&\includegraphics[width=.8in]{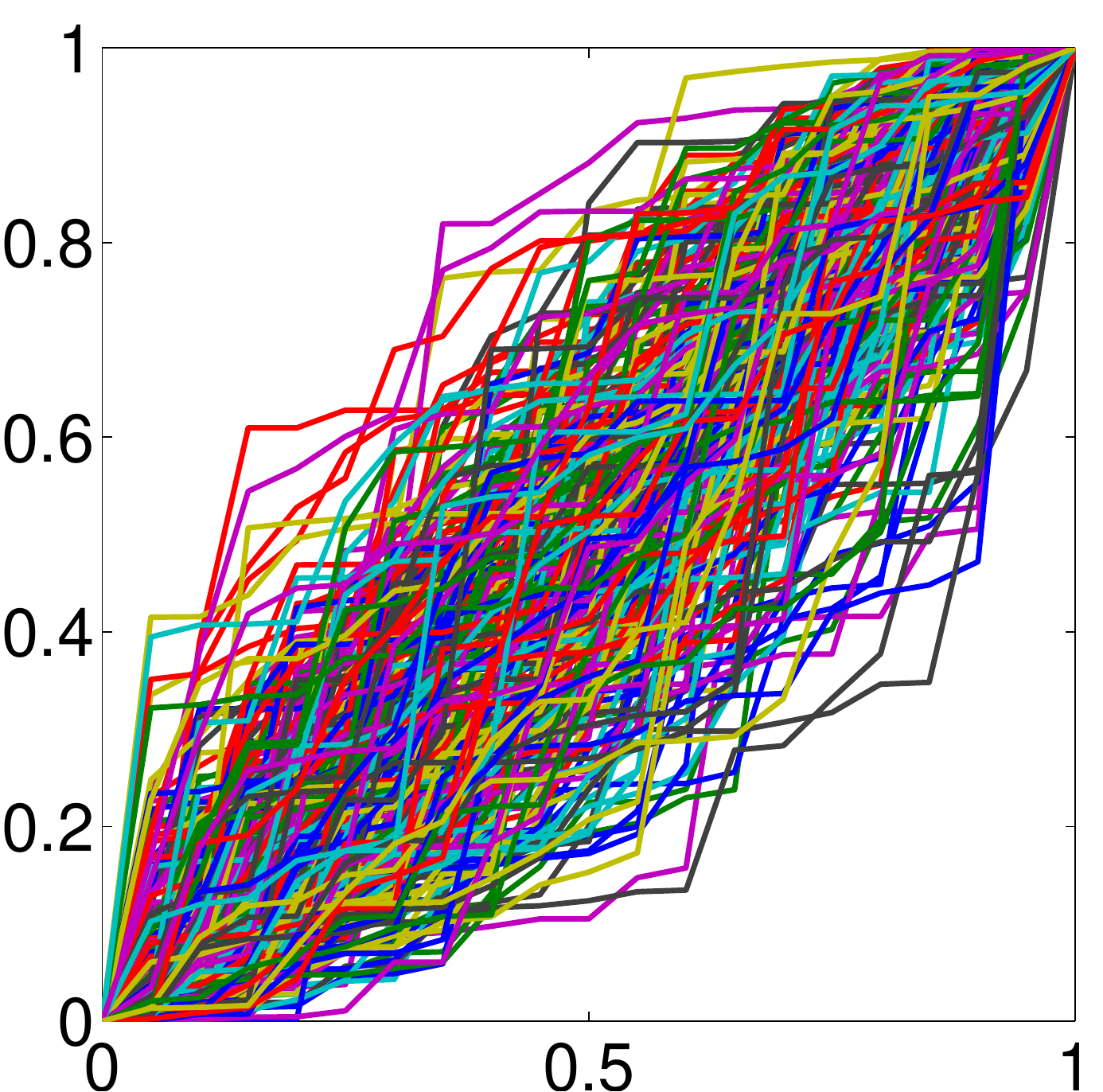}&\includegraphics[width=.8in]{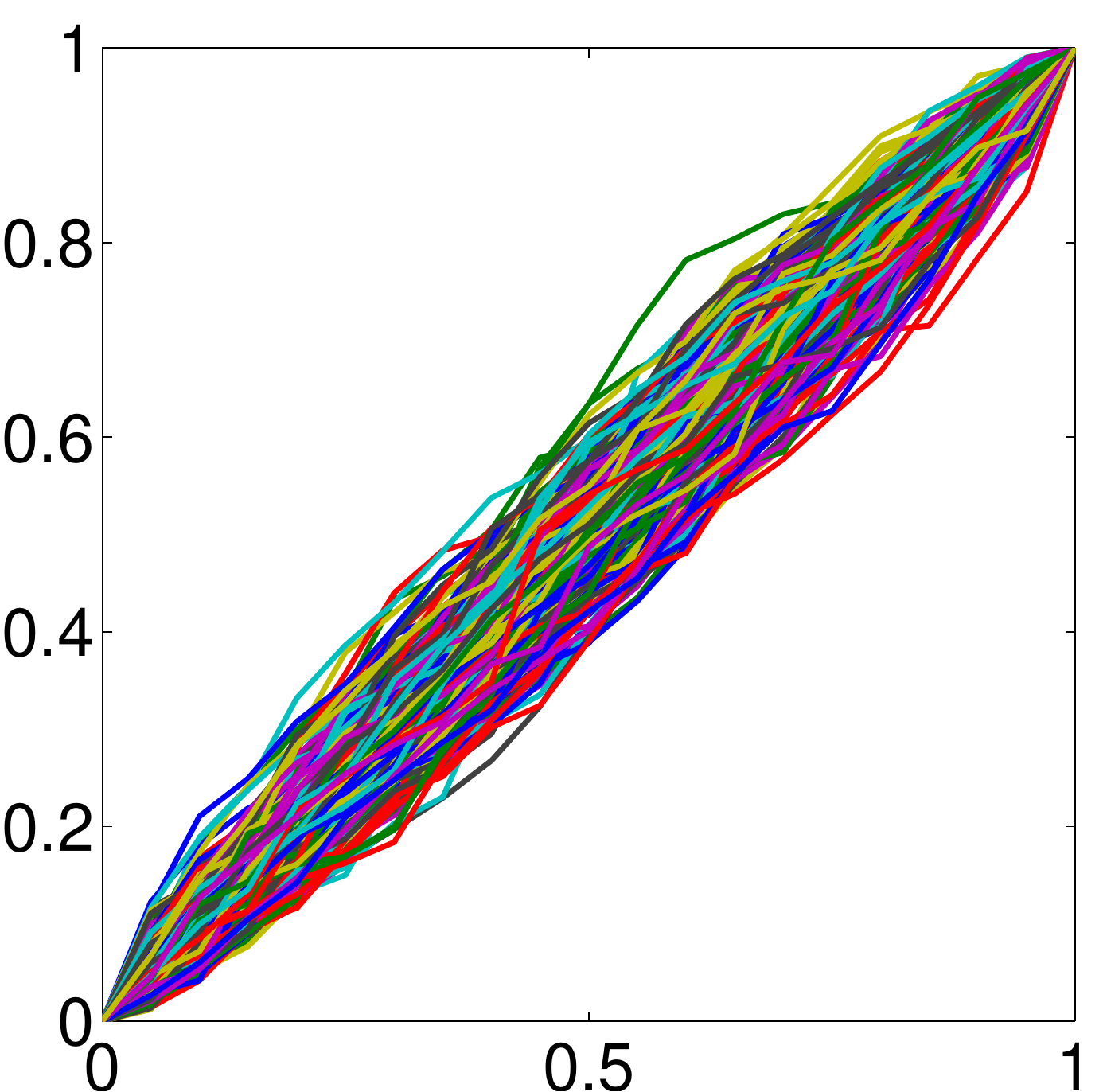}&
				\includegraphics[width=.8in]{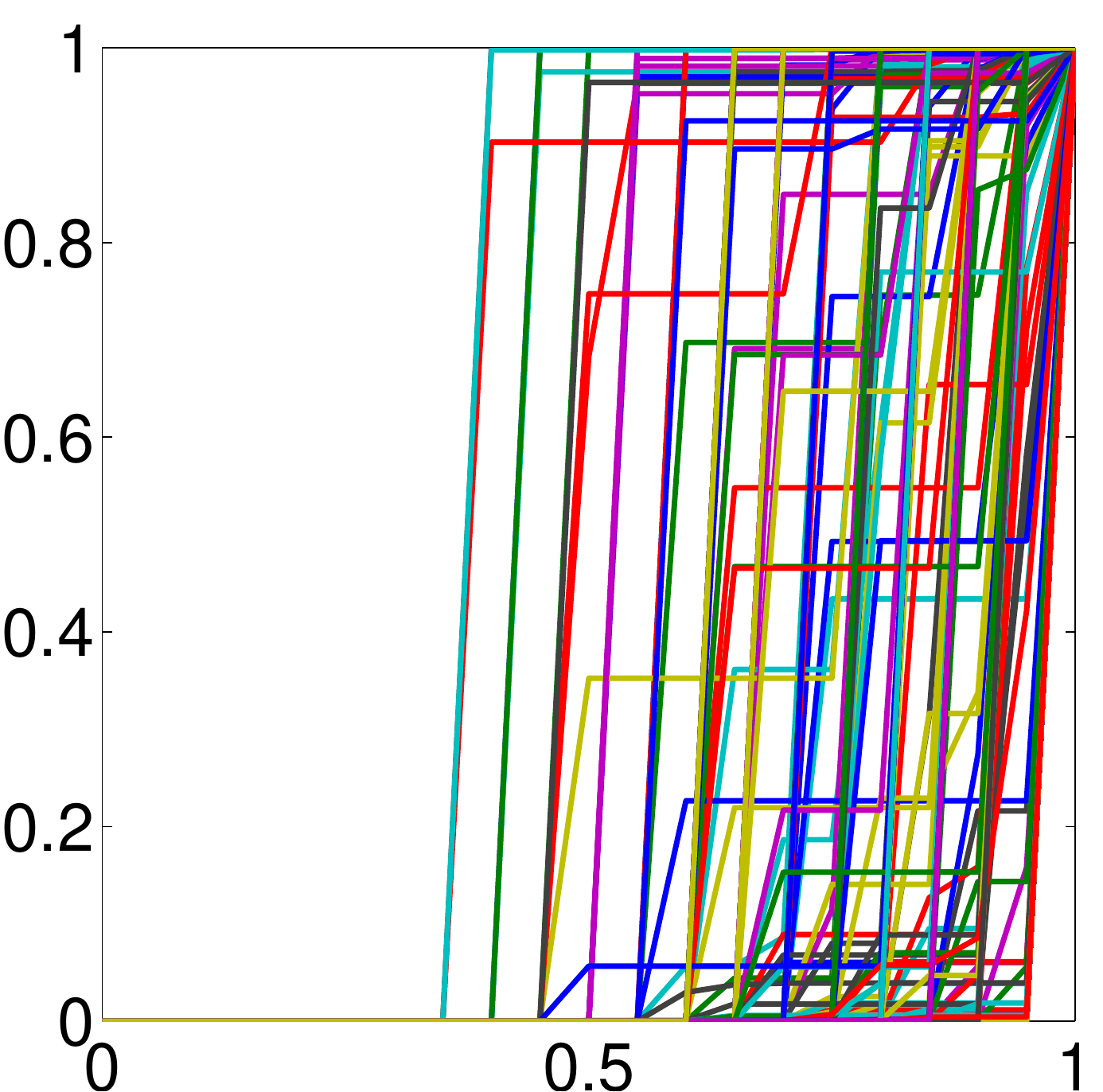}&\includegraphics[width=.8in]{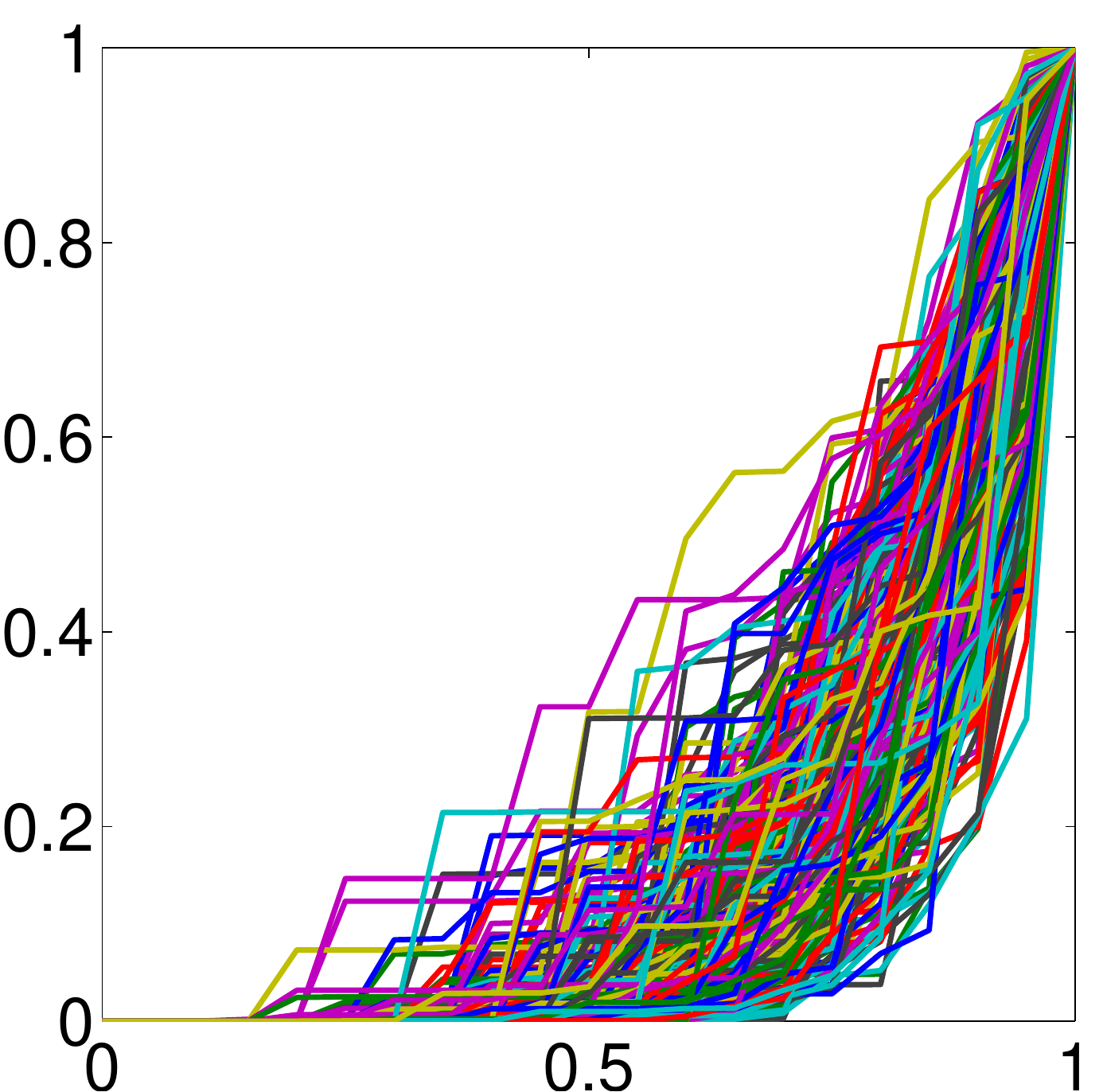}&\includegraphics[width=.8in]{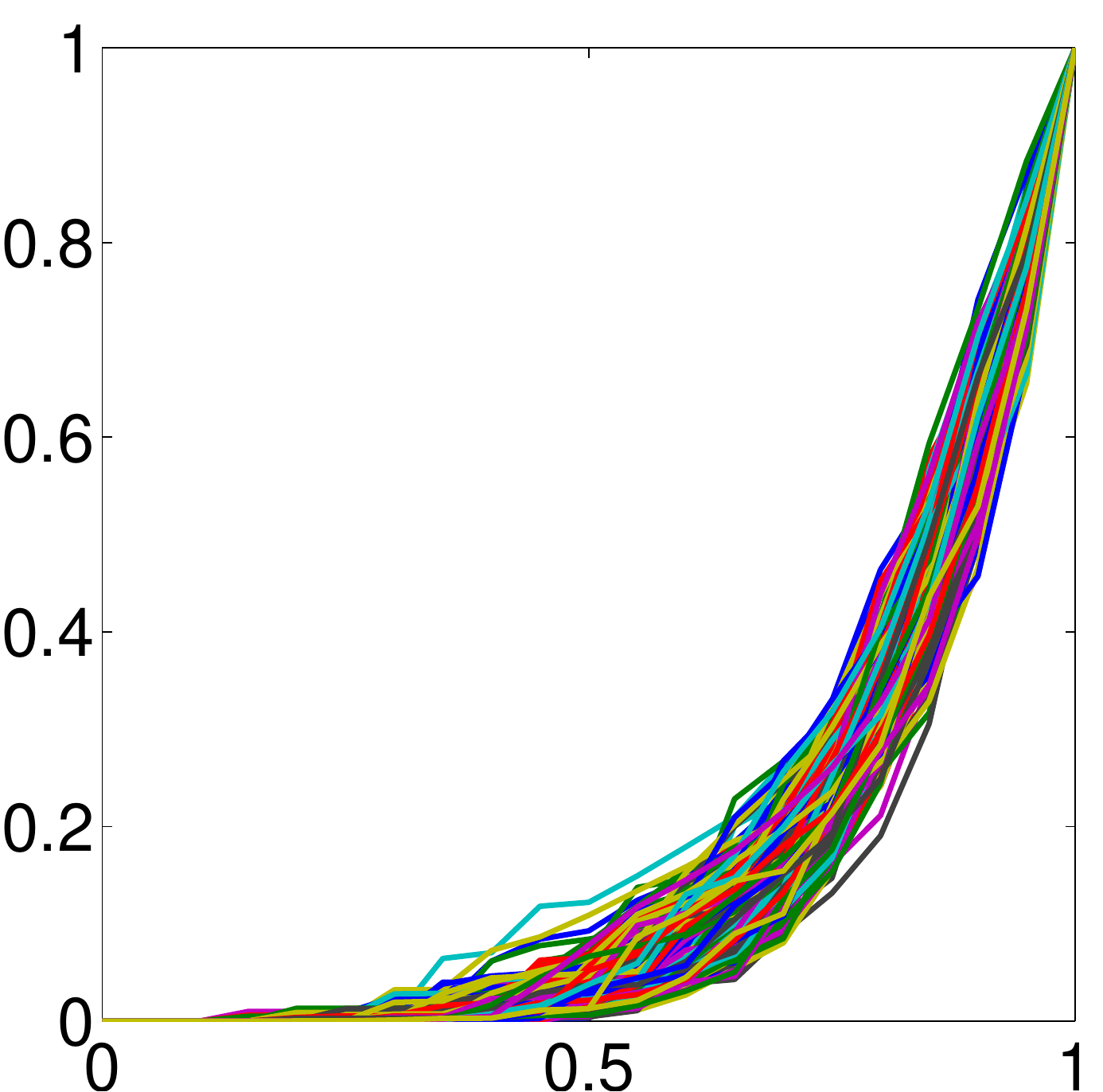}\\
				\small{$n=80$}&\includegraphics[width=.8in]{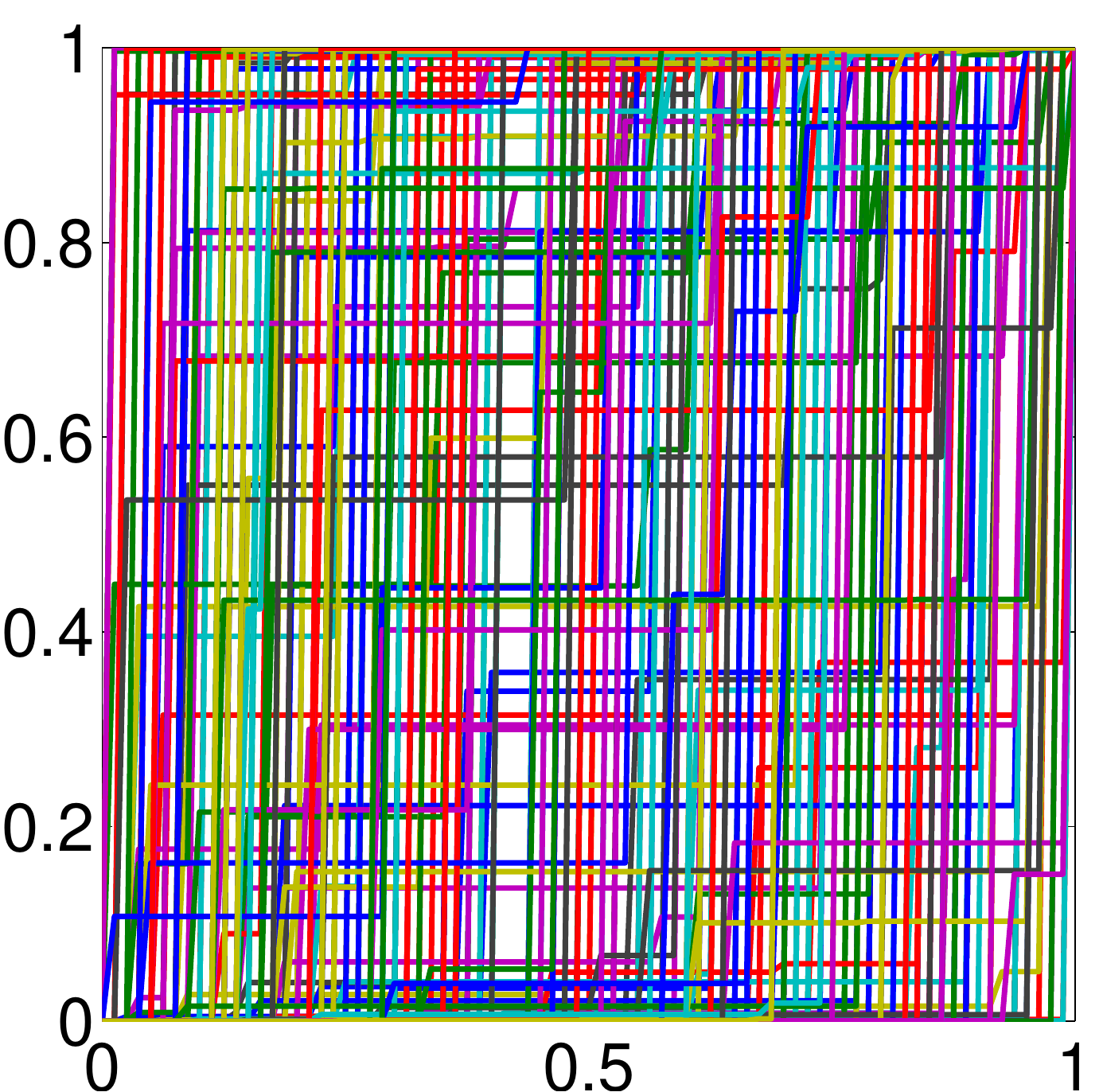}&\includegraphics[width=.8in]{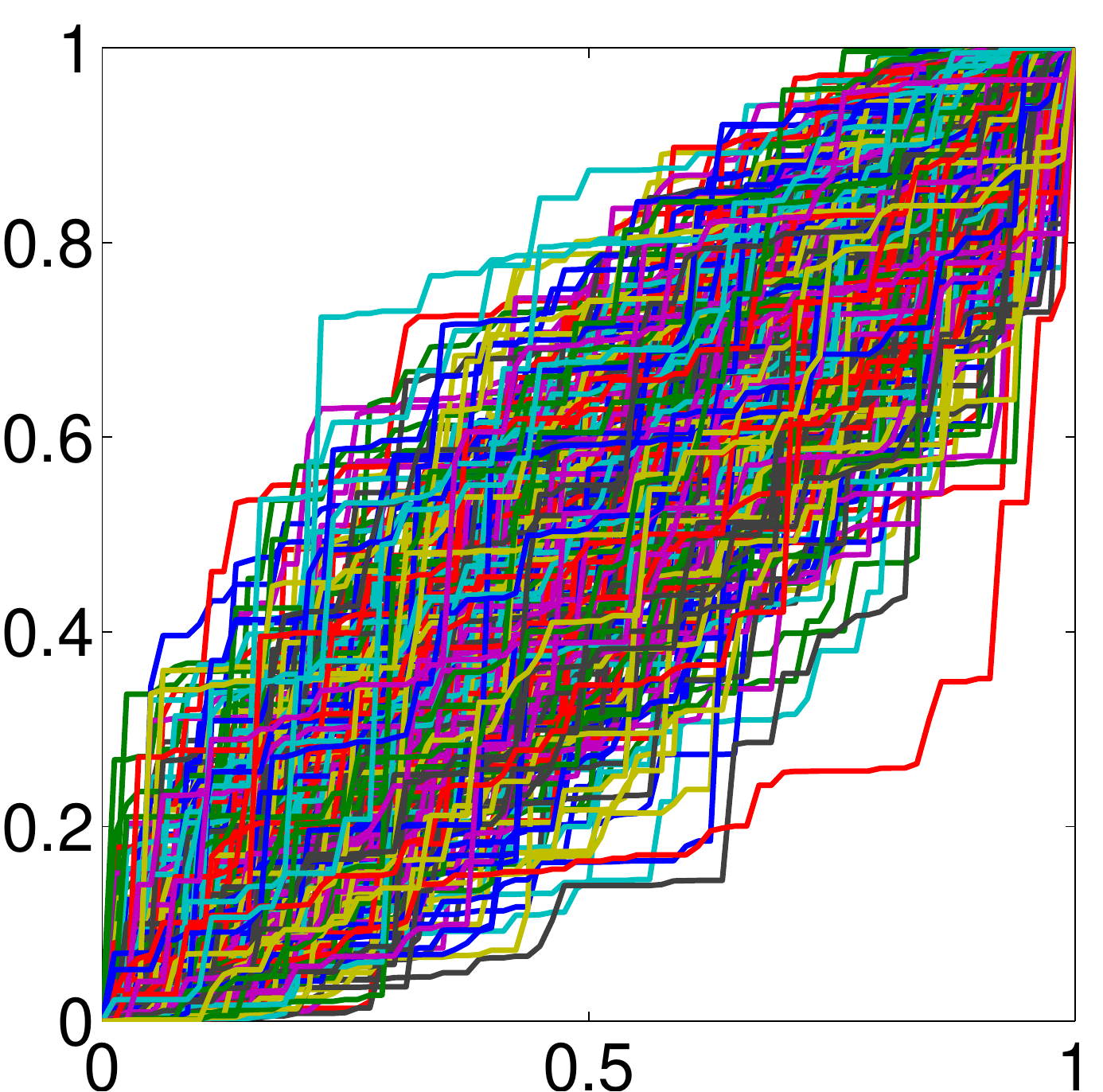}&\includegraphics[width=.8in]{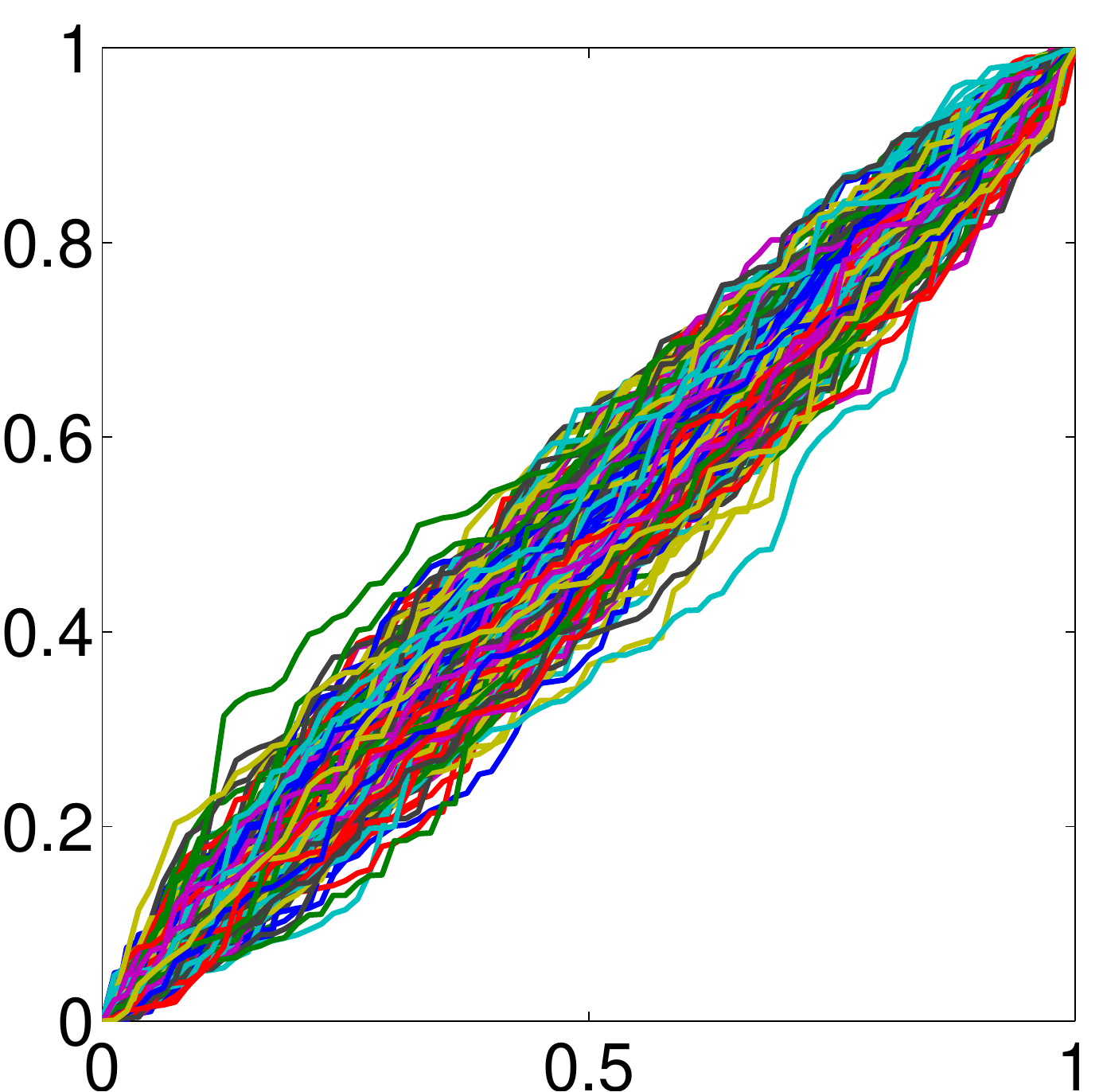}&
				\includegraphics[width=.8in]{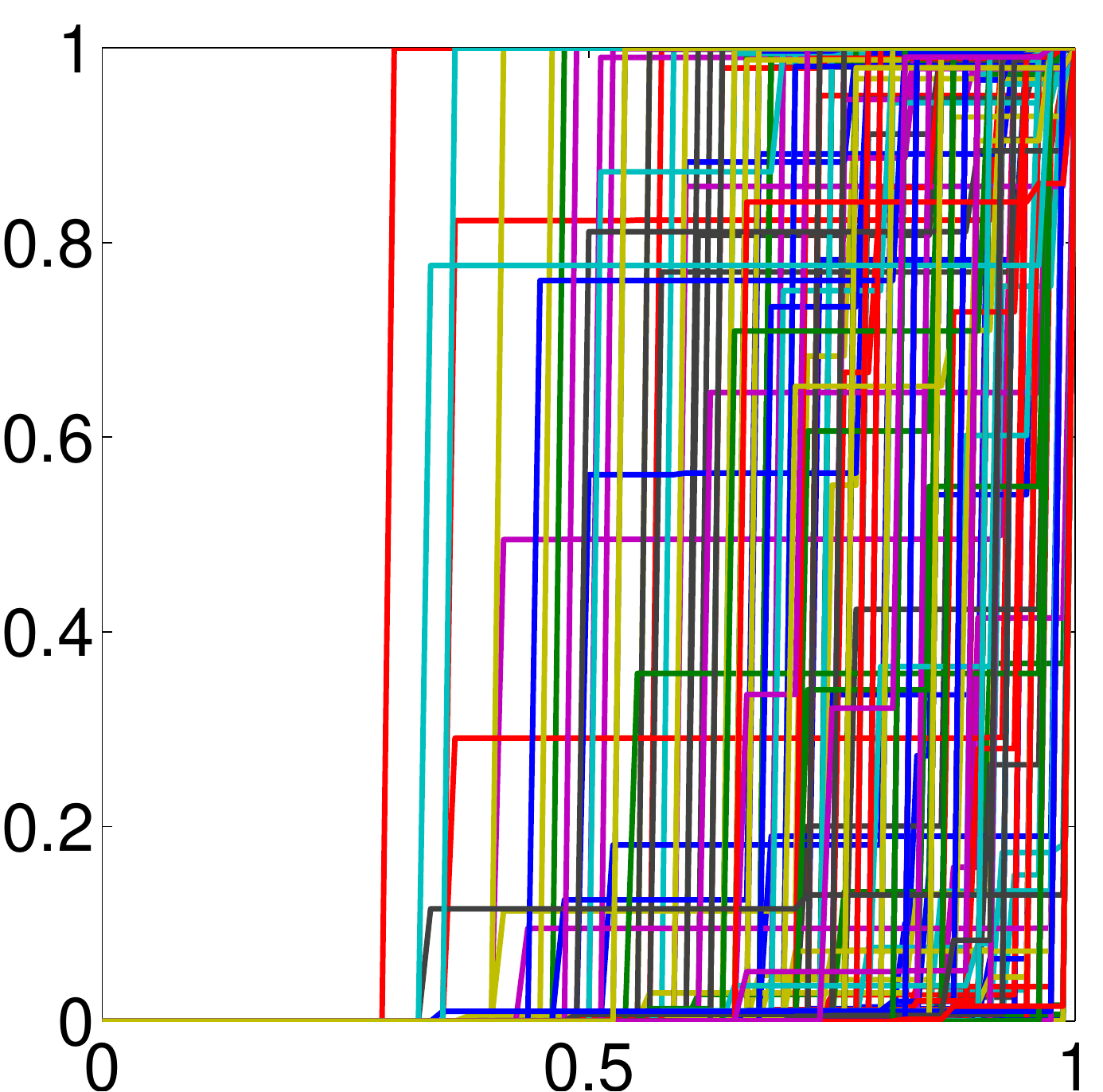}&\includegraphics[width=.8in]{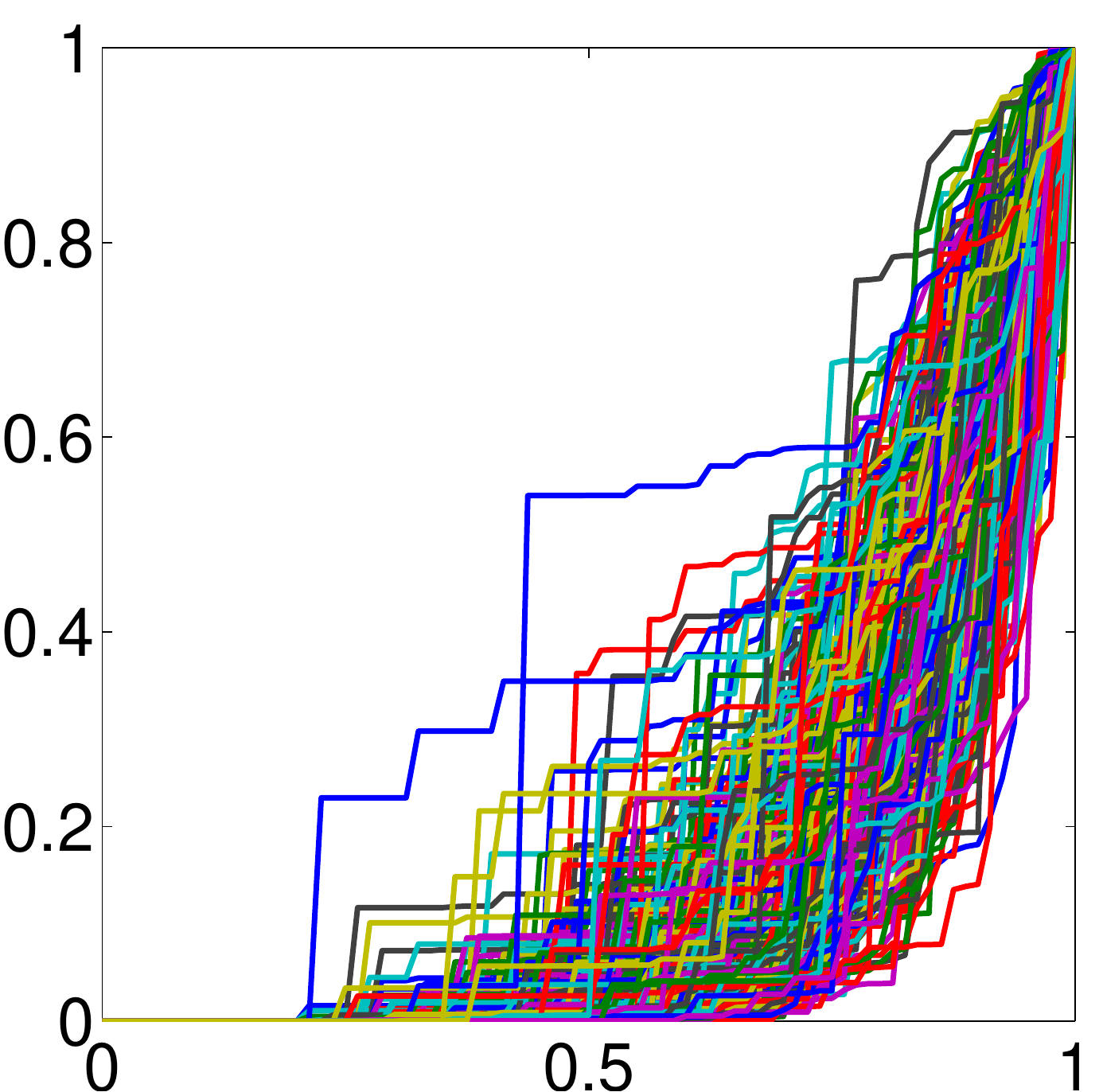}&\includegraphics[width=.8in]{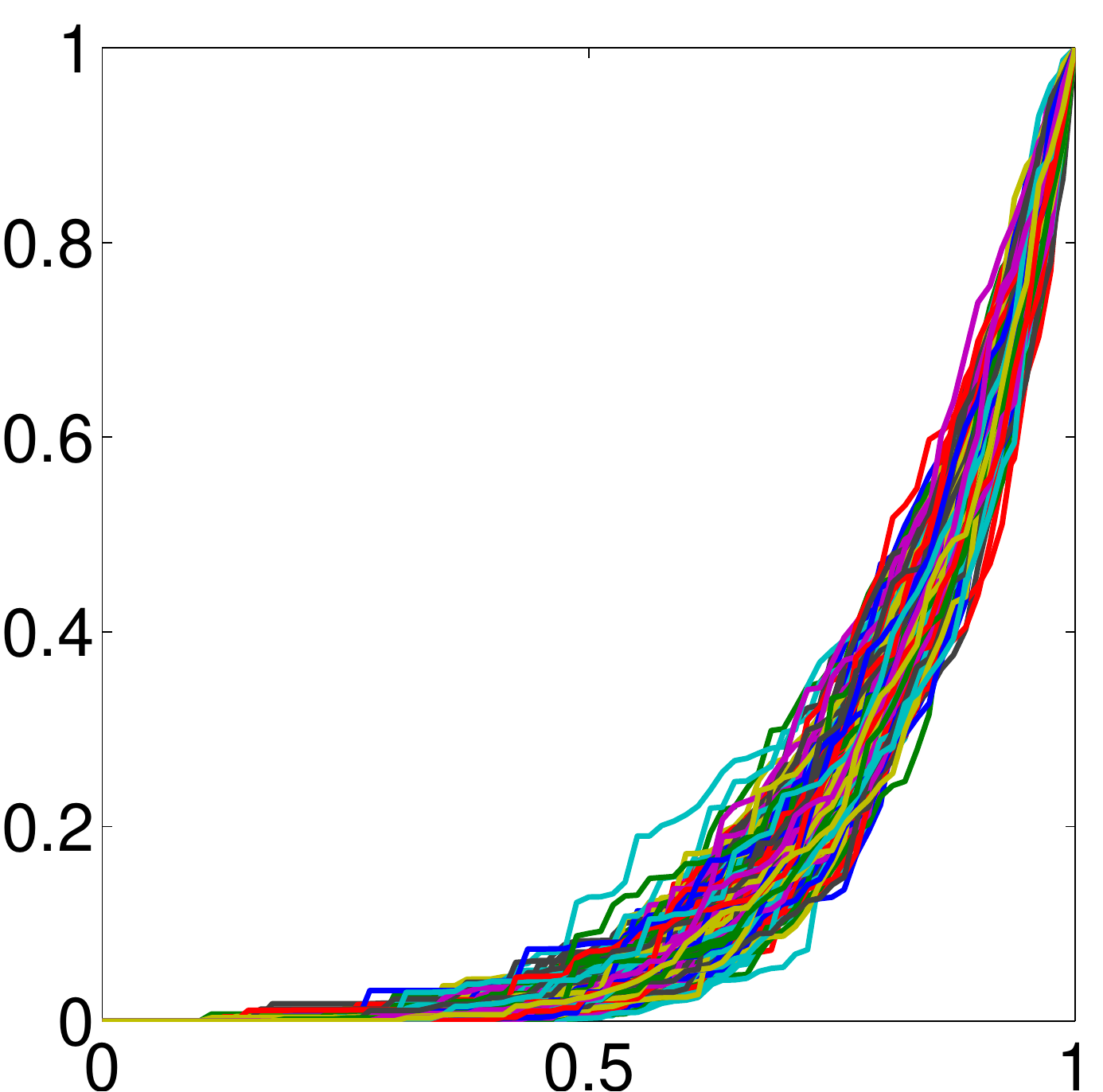}\\
				\hline
			\end{tabular}
		\end{adjustbox}
		\caption{\small {Sample warp maps from $\mathbb{D}_\theta \circ H$ based on Algorithm \ref{alg2} centred at (a) $H(t)=t$ and (b) $H$ corresponding to a Beta(5,1), under different choices of the partition size $n$ and the concentration parameter $\theta$.}}\label{fig:unif}
	\end{center}
\end{figure}

Figure \ref{degenwarp} demonstrates the degeneracy issue arising from the use of Algorithm \ref{alg1} discussed in Theorem  \ref{nowarp}.  For deterministic partitions $\Tn$, as $n \to \infty$, samples from Algorithm \ref{alg1} eventually concentrate around a deterministic warp map determined by the construction of $\Tn$. The top row of Figure \ref{degenwarp} illustrates this behavior when $\Tn$ is constructed using values of a Beta(2,1) distribution function; for increasing $n$, samples concentrate around the Beta(2,1) quantile function (see the Supplementary Material for a version of Theorem \ref{nowarp} for non-equi-spaced partitions). The bottom row illustrates the same for a uniform partition $\Tn$ wherein the samples concentrate around the identity warp map.
We now demonstrate the flexibility of $\mathbb{D}_\theta \circ H$ in modeling warp maps in $W_I$. For this purpose, we have chosen two partitions leading to two different choices of $H$: the uniform and the Beta(5,1). This results in $\mathbb{D}_\theta \circ H$ with average warp maps corresponding to uniform and Beta(5,1) distribution functions, respectively. Then, we simulated $300$ warp maps for each case under the following settings: $n=5,\ 20,\ 80$ and $\theta=0.1,\ 10,\ 100$. The simulated warp maps under the uniform partition are shown in Figure \ref{fig:unif}(a), while panel (b) shows the warp maps sampled based on the Beta(5,1) partition. In both cases, we see that the proposed distribution is very flexible, exhibiting a variety of possible shapes of warp maps under different combinations of $n$ and $\theta$. A partition created with $n=5$ results in few large jumps while a partition created with $n=80$ generates warp maps with many small jumps. As $\theta$ is increased from $0.1$ to $100$, we notice the sample tightening around the warp map corresponding to the average map induced by the partition.

\subsection{Illustrations on real data}
\label{applications}
For detailed descriptions of the datasets used throughout this section, please refer to the Supplementary Material.
We illustrate the utility of the proposed distribution and associated sampling scheme in two pairwise alignment tasks: (1) Bayesian alignment of univariate functions with and without landmark constraints, and (2) unconstrained alignment of univariate functions and higher-dimensional open and closed curves using a novel Simulated Annealing-based algorithm (see \cite{Robert:2005:MCS:1051451} for details). While the distribution $\mathbb{D}_\theta \circ H$ for sampling warp maps can be used with any alignment method, we use the framework based on the square-root velocity function (SRVF) representation of curves in $\real^d,\ d \geq 1$: $f\mapsto q:=\dot{f}(|\dot{f}|)^{-1/2}$, where $\dot{f}$ is the derivative of $f$ and $|\cdot|$ is the Euclidean norm in $\real^d$. We use this representation due to its many nice properties for the registration problem (see \cite{AS,KurtekJASA,SrivESA,LRK}). Under this representation, warping of a function $f\mapsto f\circ\gamma$ is given by $q\mapsto (q\circ\gamma)\sqrt{\dot{\gamma}}$.
\subsubsection{Bayesian alignment of curves}
\label{bayesian}
Suppose we have two functions $g_i:[0,1]\to \mathbb{R},\ i=1,2$. A Bayesian registration model can be defined using their SRVF representations. The alignment problem then centres around an $\mathbb{R}$-valued, square-integrable, separable stochastic process $X(t):=q_1(t)-q_2(t)$ for $t \in [0,1]$ with law $\mathbb{P}$ and density $p=\frac{dP}{d\mu}$ with respect to a $\sigma$-finite measure $\mu$ on $\ltwo([0,1])$. For a fixed $\gamma \in W_I$, assume that the law $\mathbb{P}_\gamma$ of the process $X_\gamma(t)=q_1(t)-q_2(\gamma(t))\sqrt{\dot{\gamma}(t)}$ is absolutely continuous with respect to $\mathbb{P}$ with density $p_\gamma$ (see Theorem 6.4.5 in \cite{bogachev} for sufficient conditions). Suppose that the SRVFs, discretized at points $[t]=\{t_{0:n},\ldots,t_{n:n}\}$, are represented as vectors $q_i([t]):=(q_i(t_{0:n}),\ldots, q_i(t_{n:n})),\ i=1,2$. Then, $X_\gamma([t]):=q_1([t])-q_2(\gamma([t]))\sqrt{\dot{\gamma}([t])} \sim p_\gamma$ prescribes a likelihood through the finite-dimensional projections of $X_\gamma$. The PL discretization of the functions, as well as the warp map $\gamma$, are theoretically supported by the work of \cite{LRK}. Thus, in this context, the finite dimensional restriction of the probability measure $\mathbb{D}_\theta \circ H$ is well-suited for defining priors on warp maps.

\begin{figure}[!t]
	\begin{center}
		\begin{adjustbox}{max width=\textwidth}
			\begin{tabular}{|ccc||ccc|}
				\hline
				\multicolumn{3}{|c||}{Unconstrained}&\multicolumn{3}{|c|}{Landmark-constrained}\\
				\hline
				(a)&(b)&(c)&(a)&(b)&(c)\\
				\hline
				\includegraphics[width=.9in]{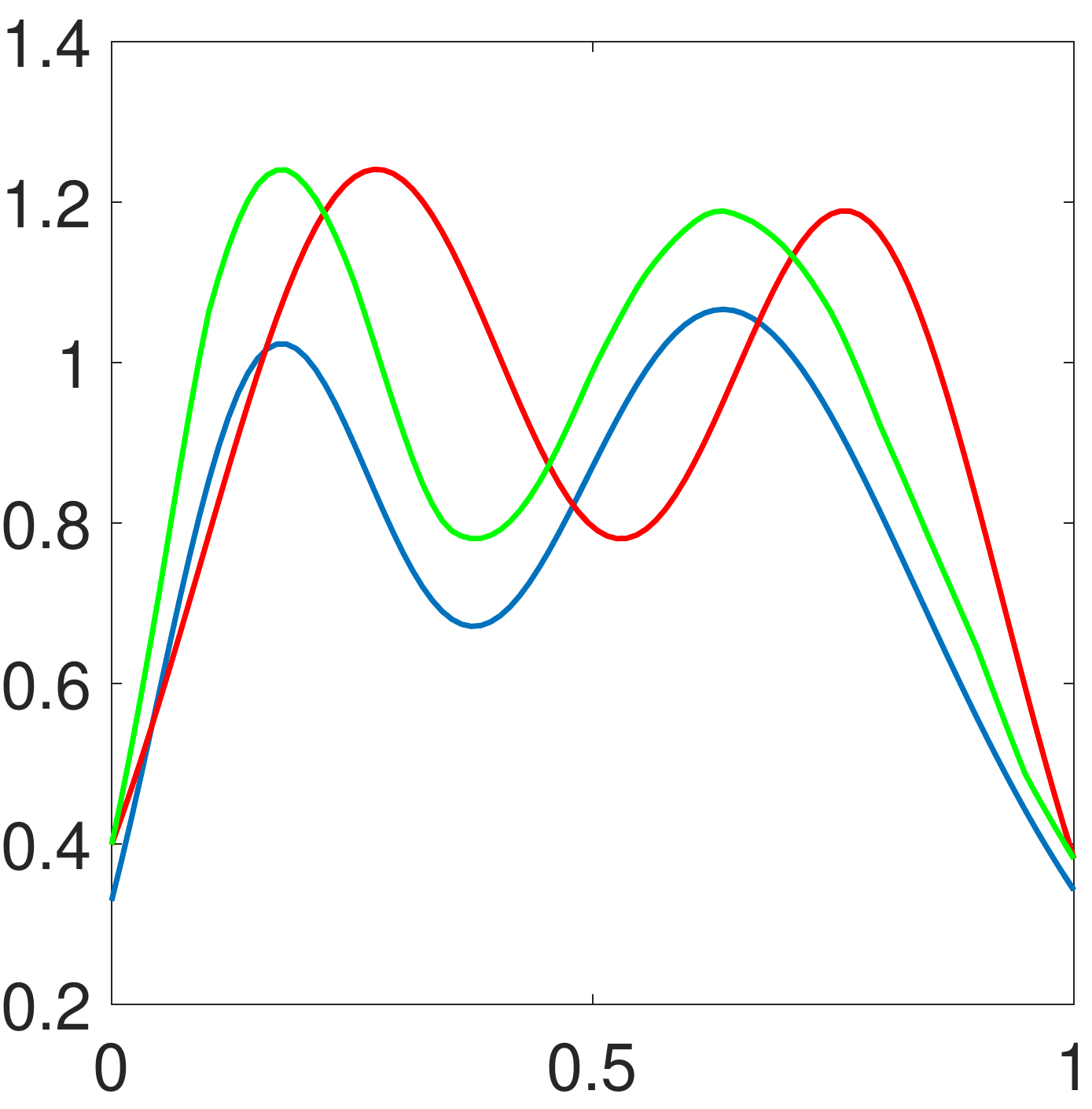}&\includegraphics[width=.9in]{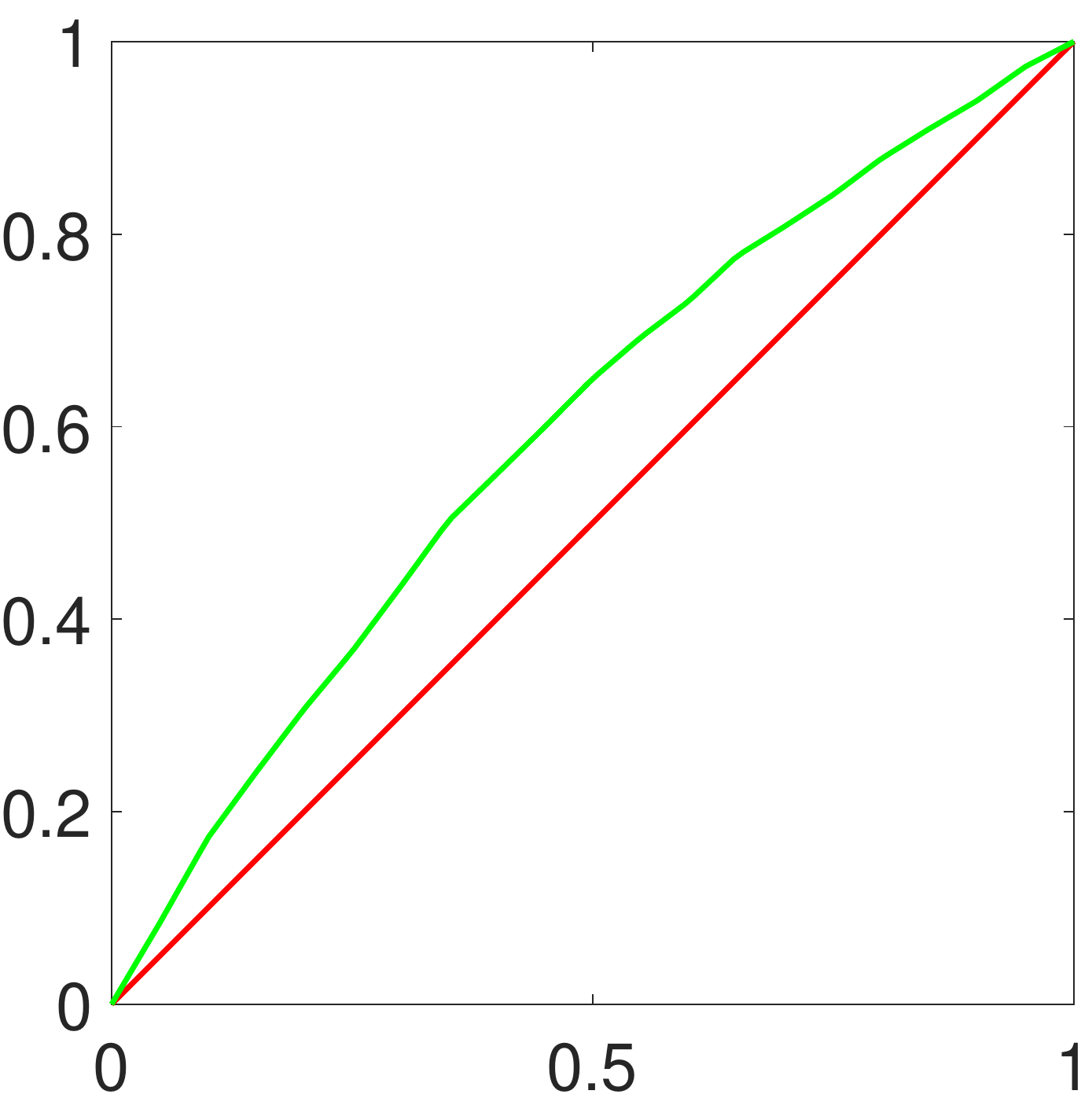}&\includegraphics[width=.9in]{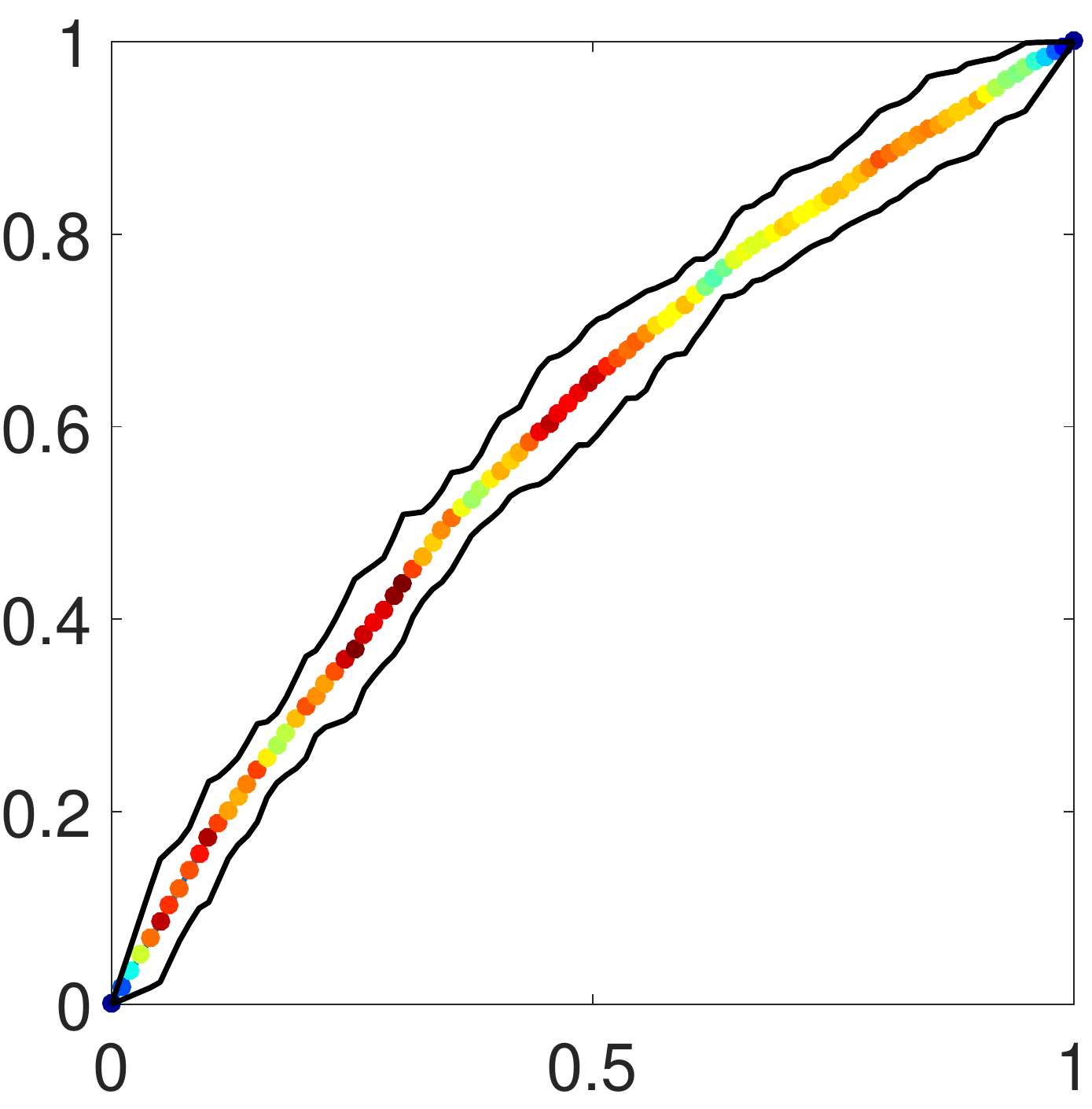}&\includegraphics[width=.9in]{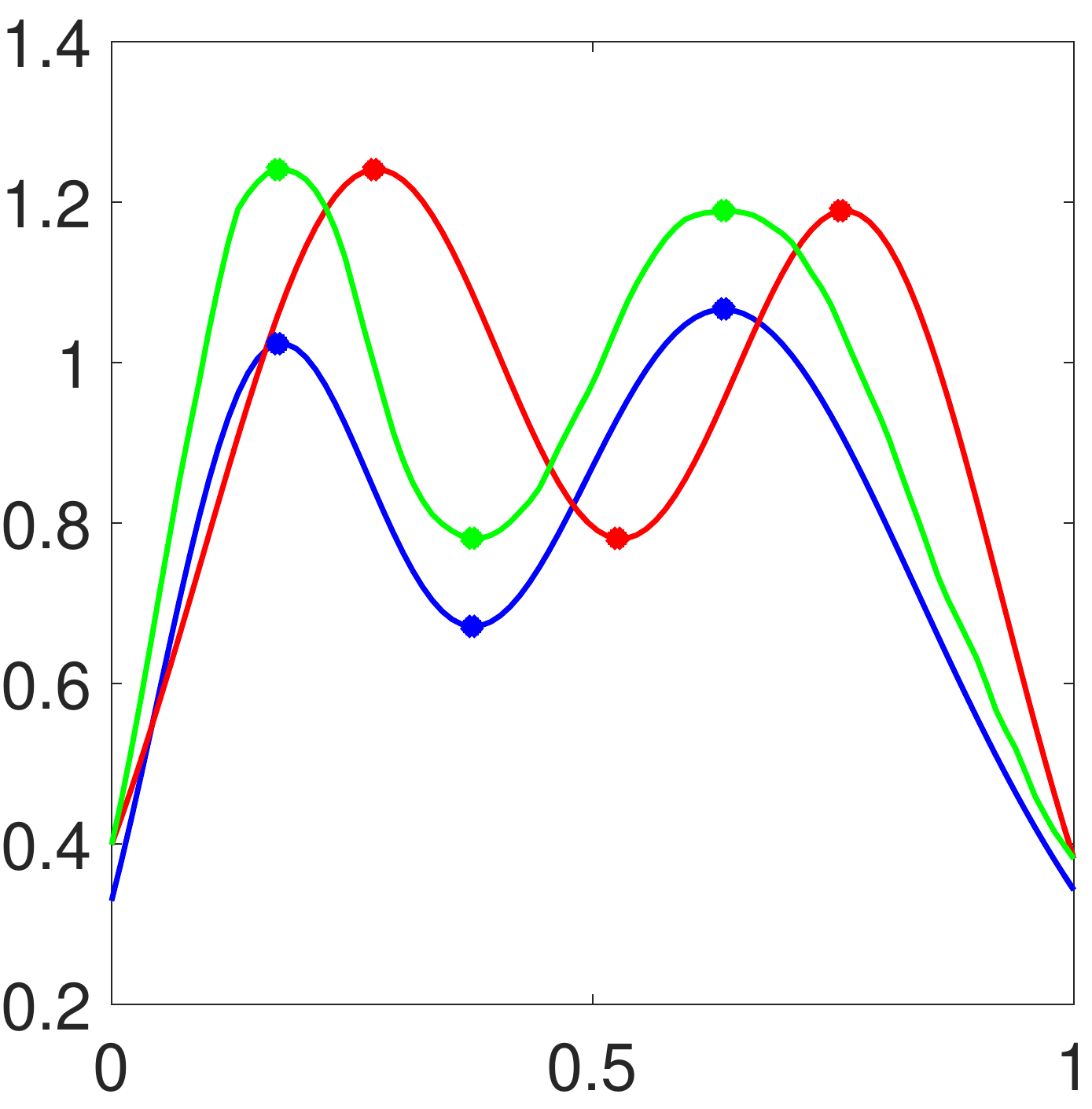}&\includegraphics[width=.9in]{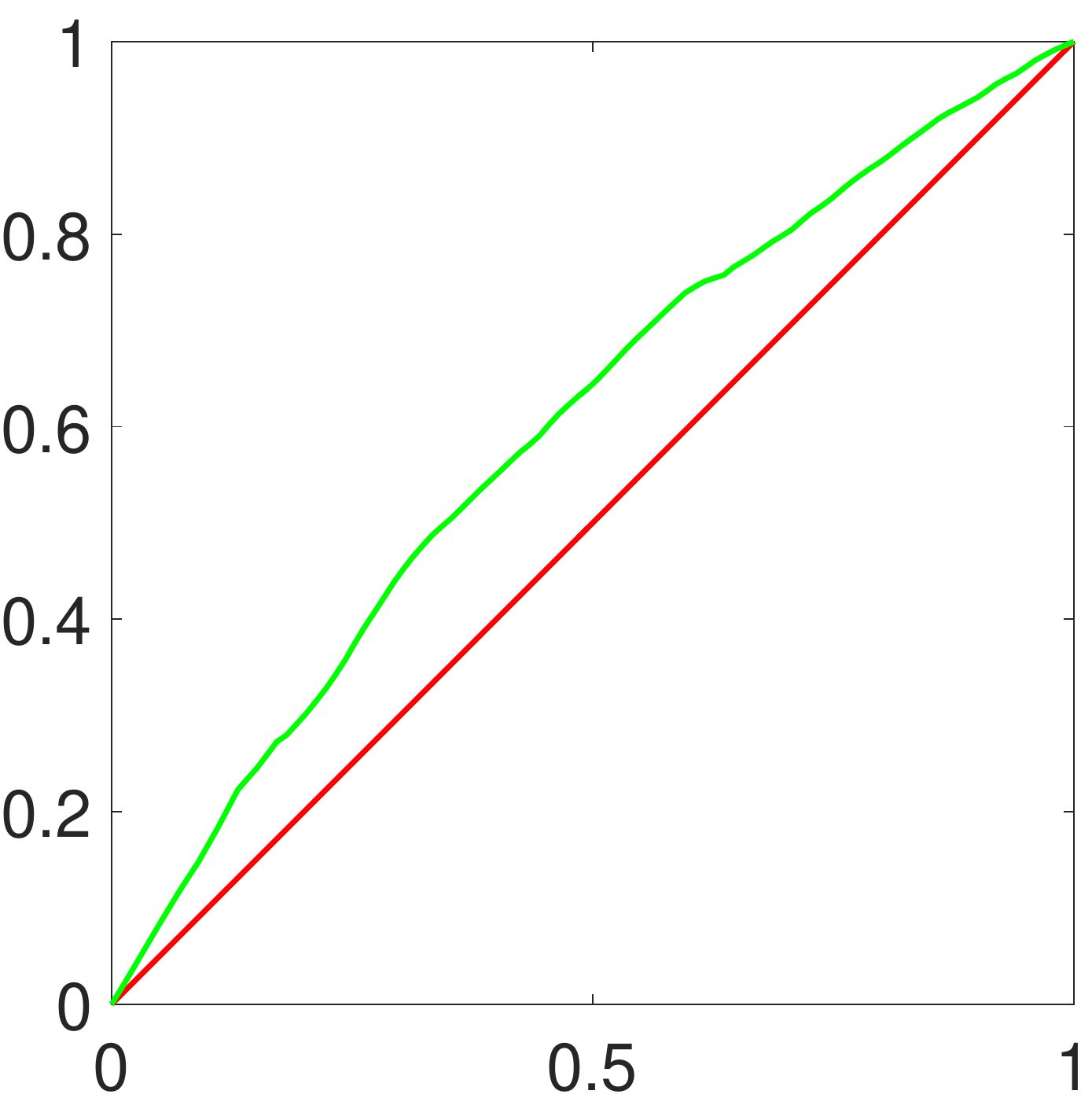}&\includegraphics[width=.9in]{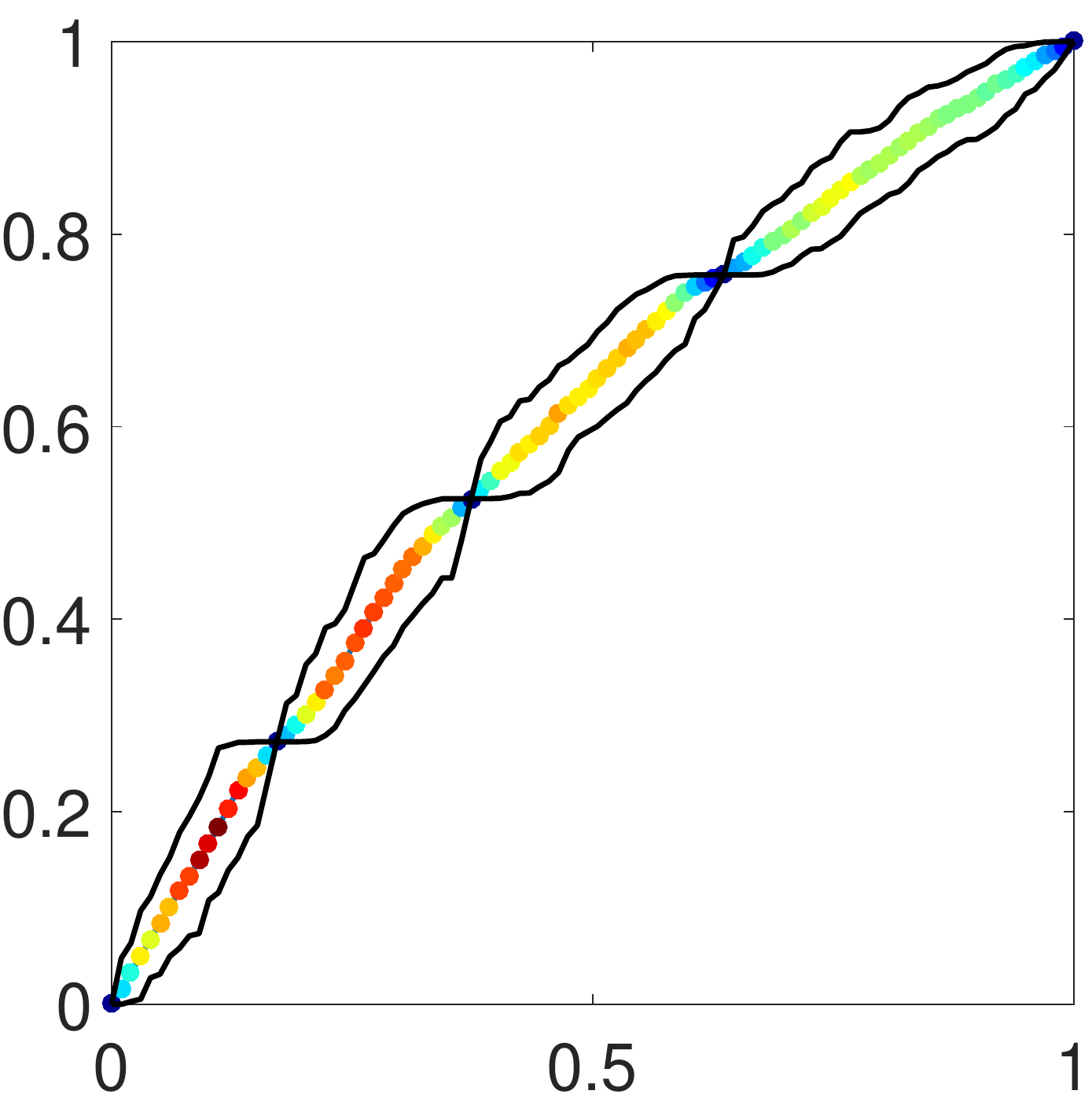}\\
				\hline
				\includegraphics[width=.9in]{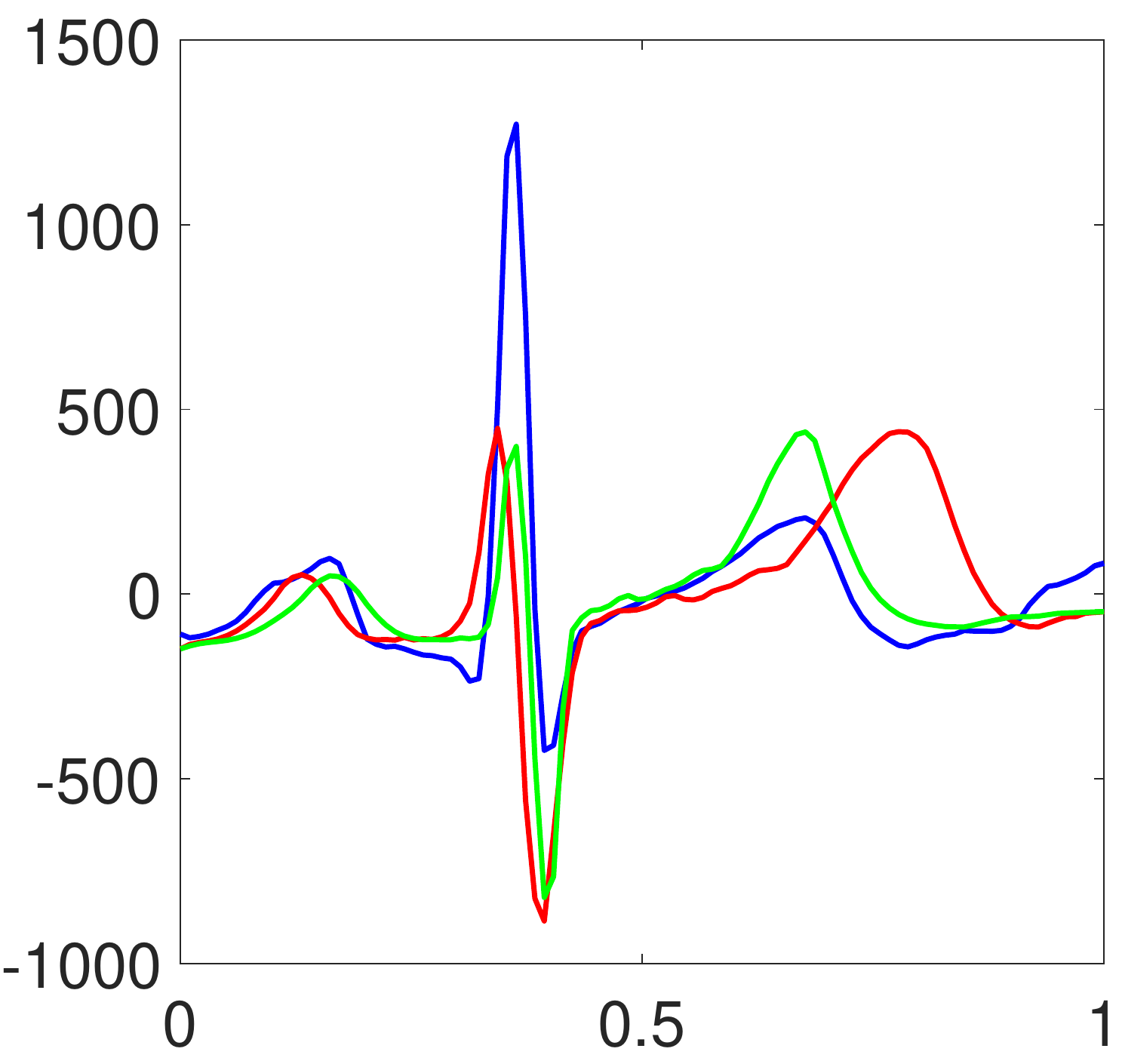}&\includegraphics[width=.9in]{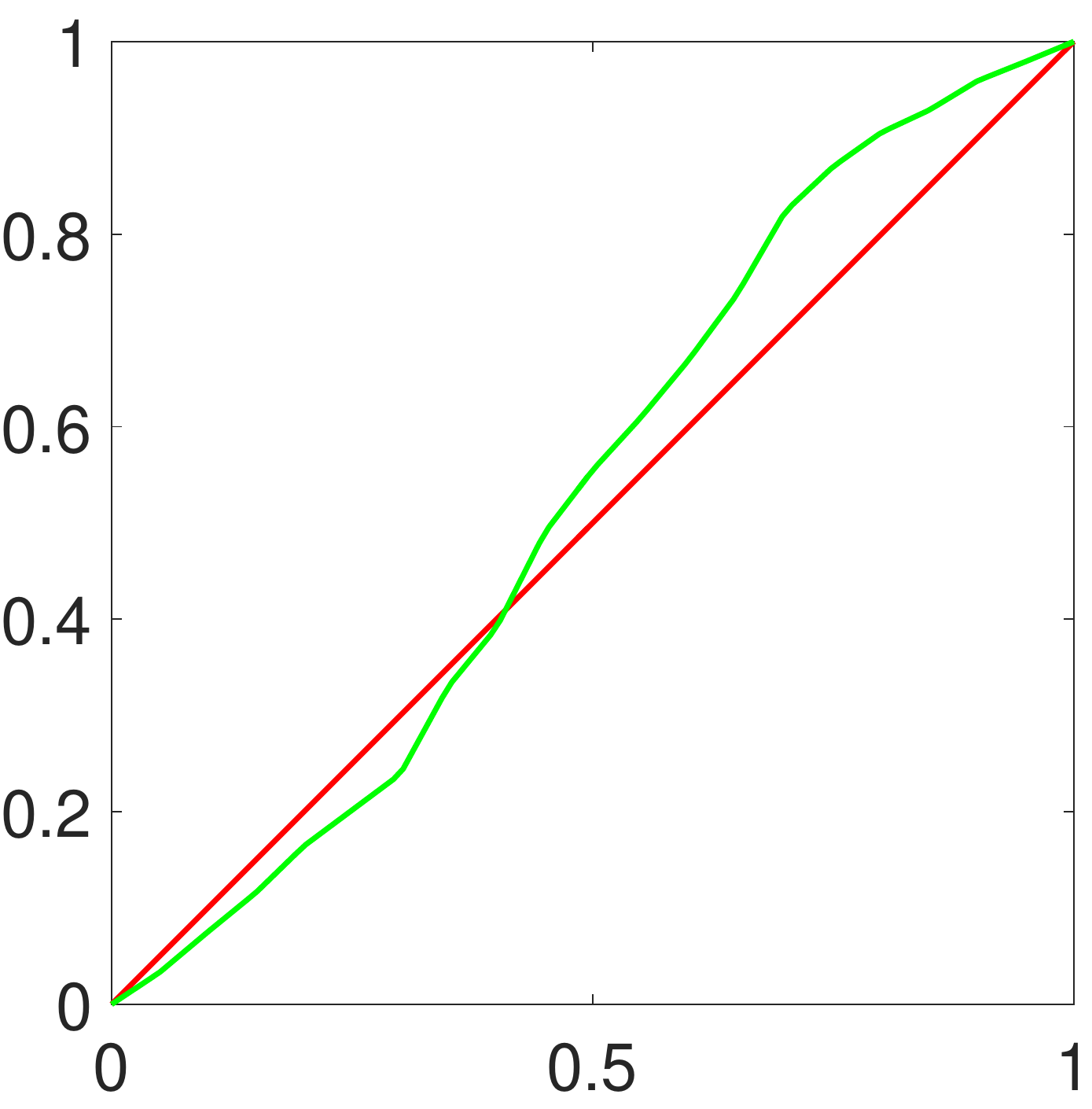}&\includegraphics[width=.9in]{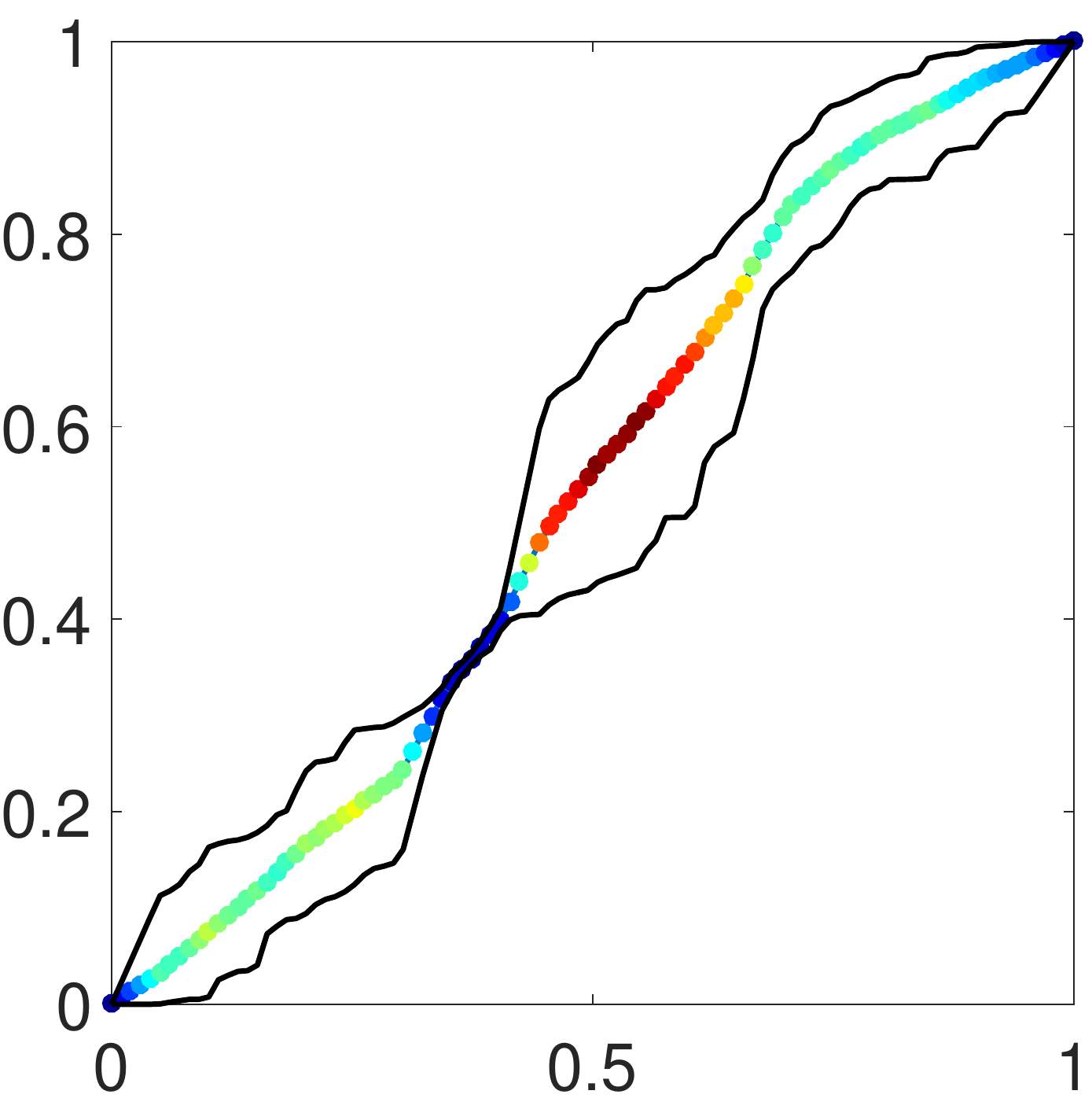}&\includegraphics[width=.9in]{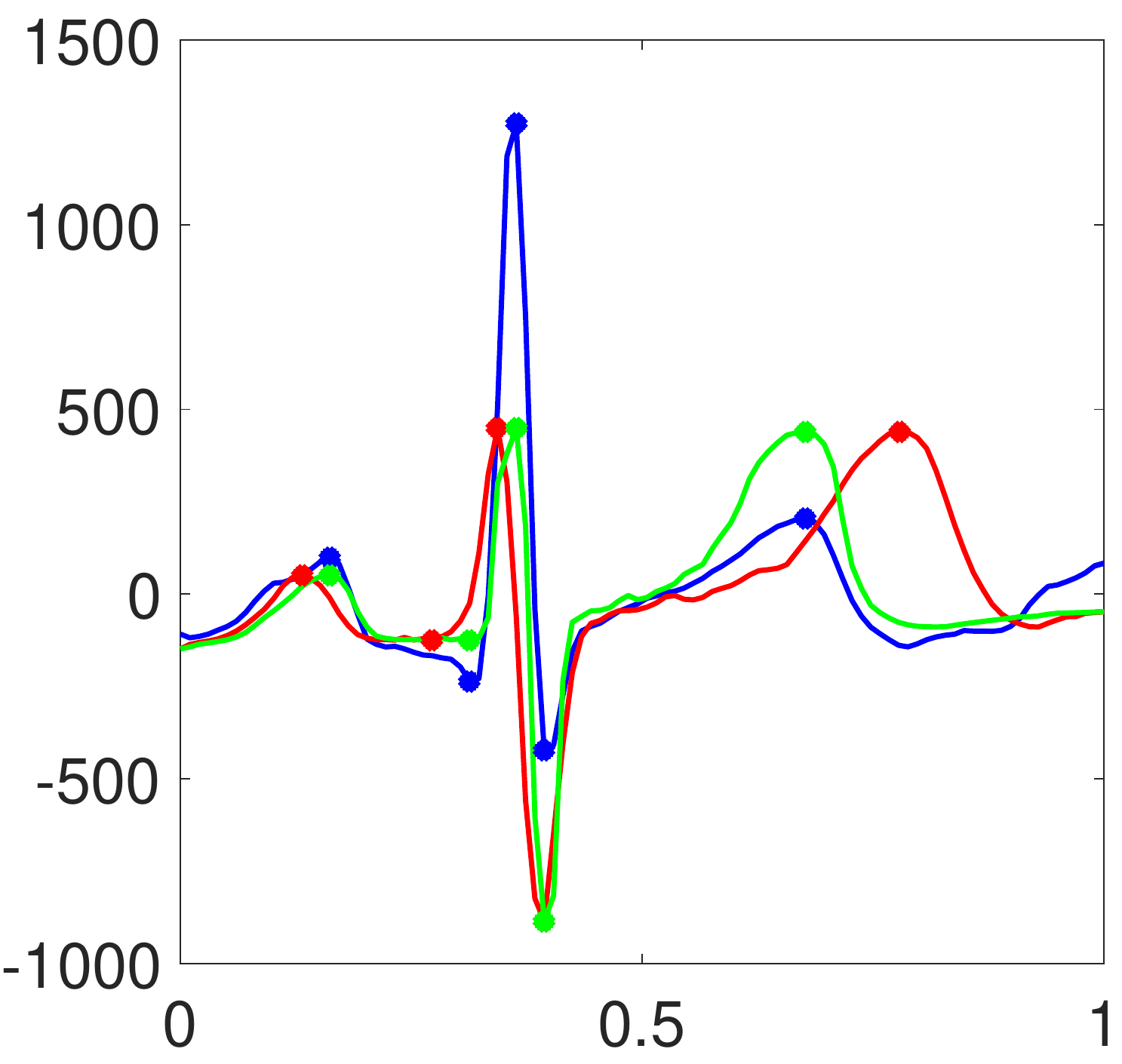}&\includegraphics[width=.9in]{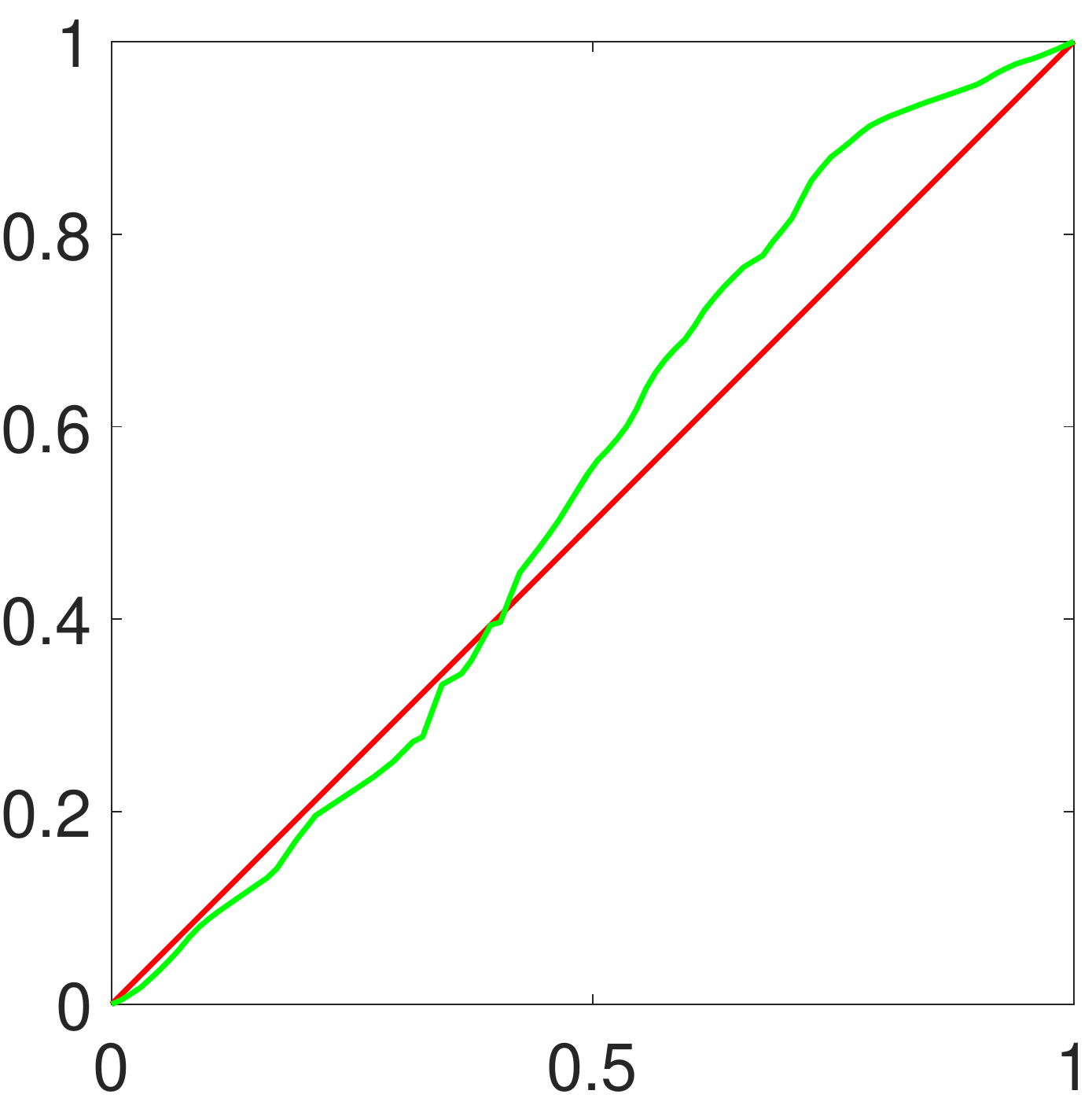}&\includegraphics[width=.9in]{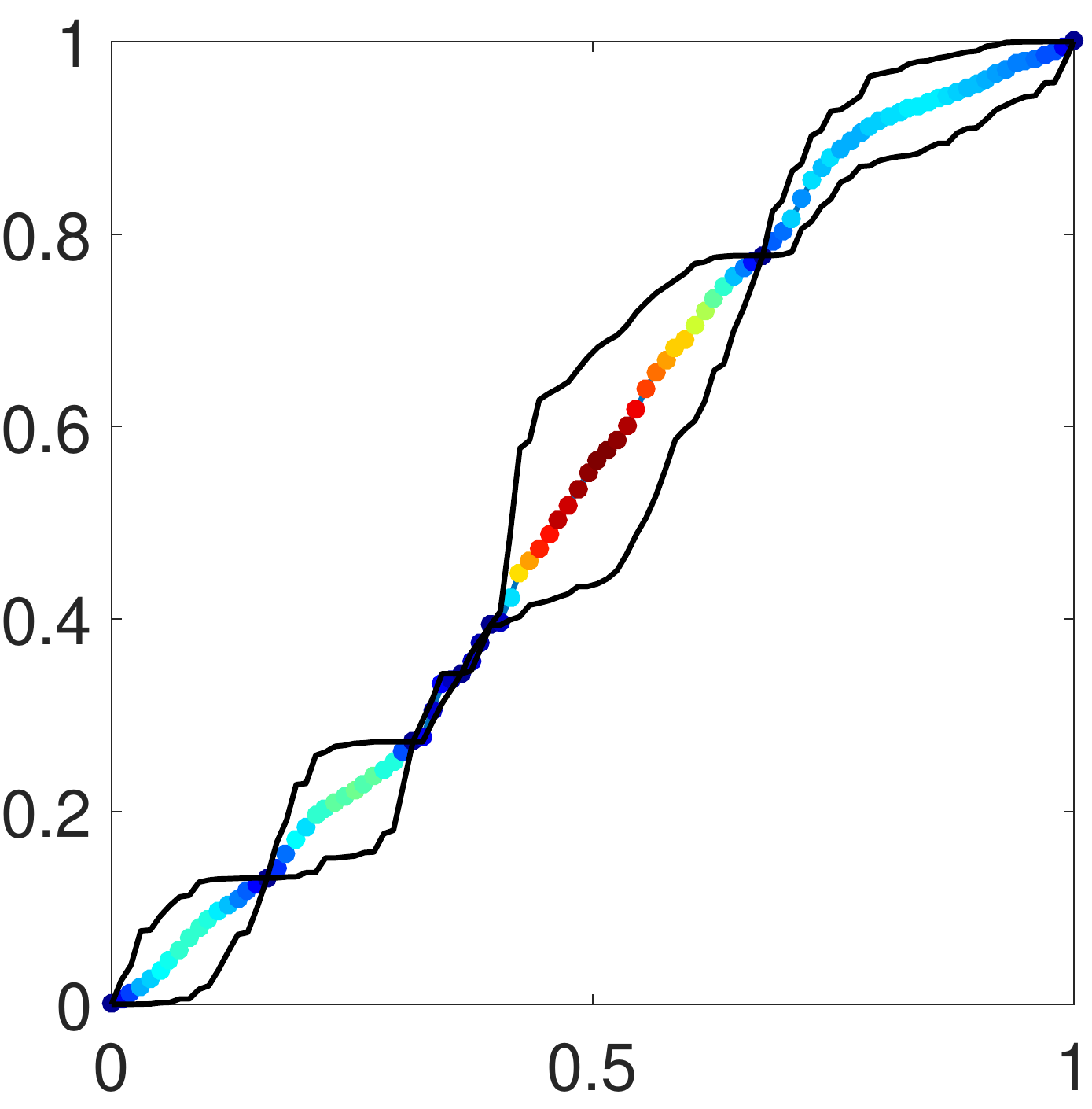}\\
				\hline
				\includegraphics[width=.9in]{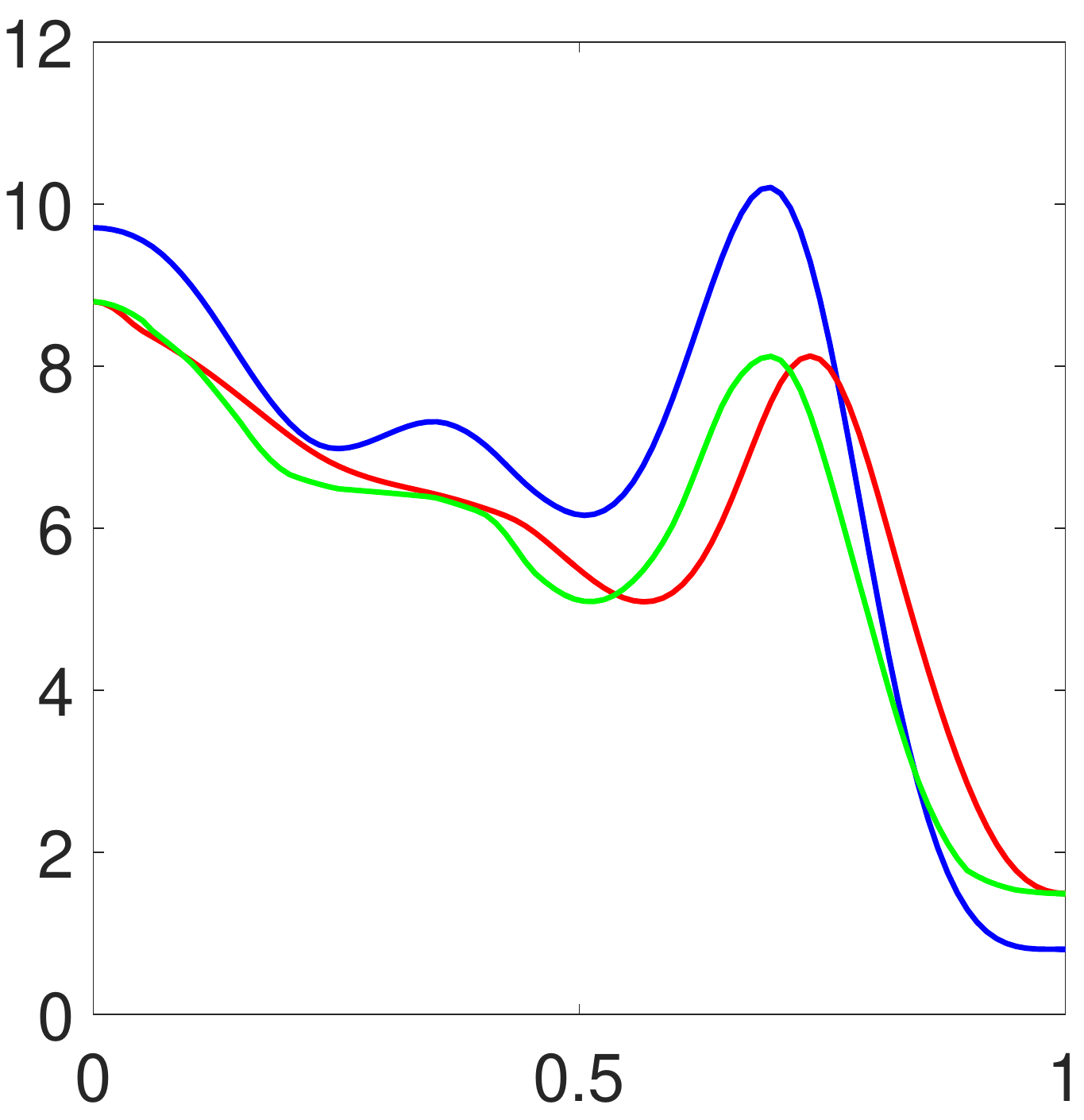}&\includegraphics[width=.9in]{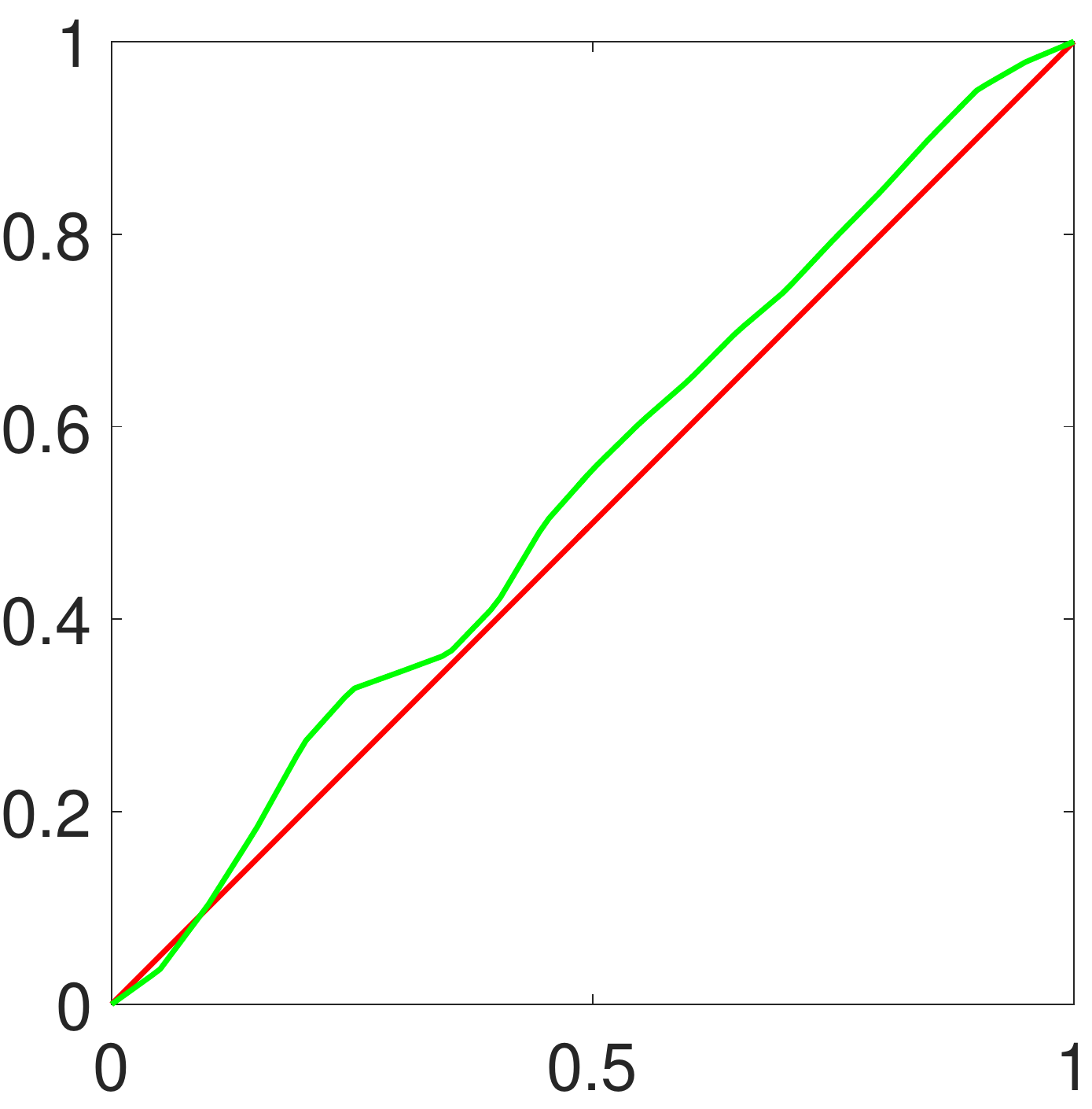}&\includegraphics[width=.9in]{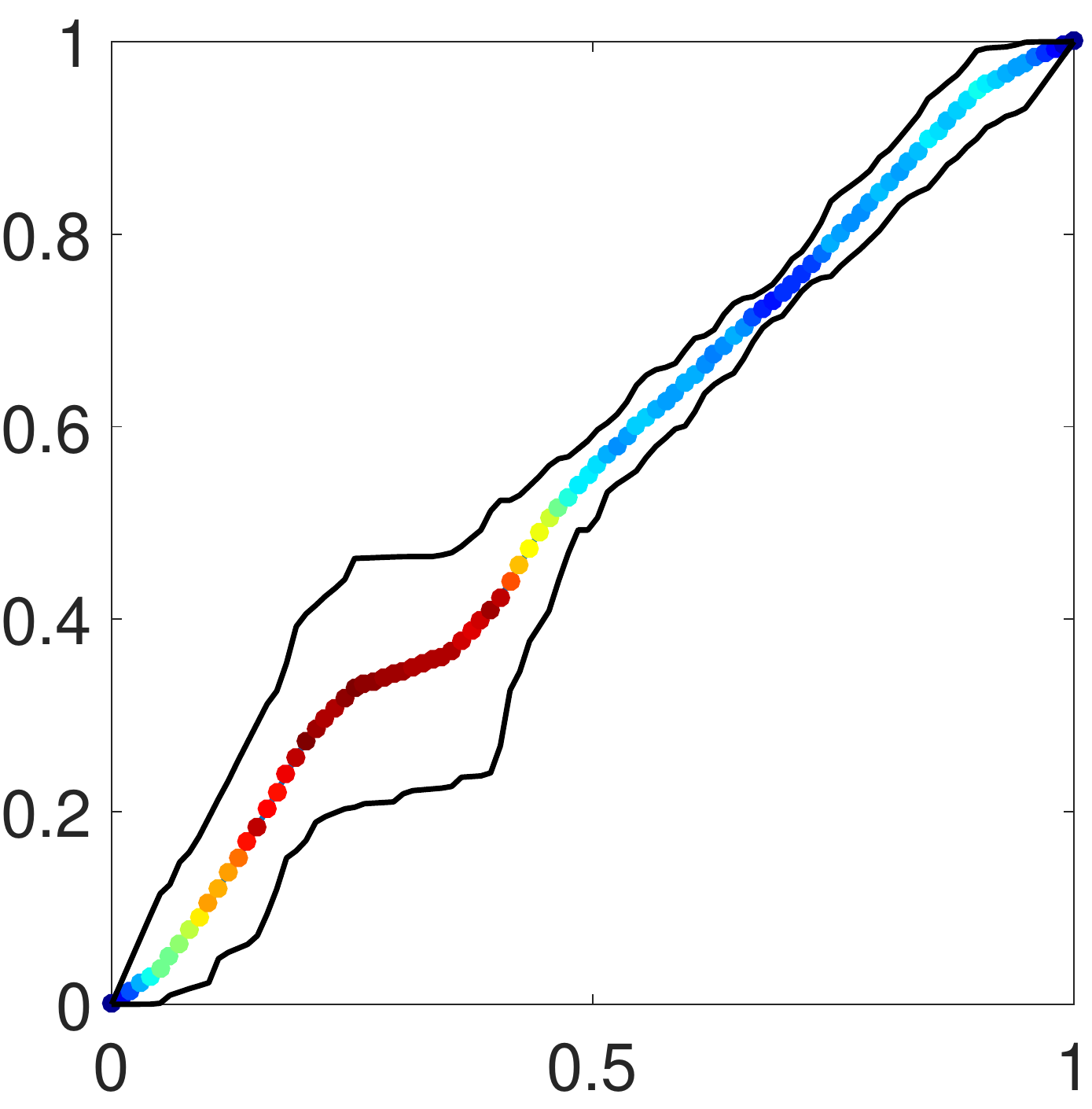}&\includegraphics[width=.9in]{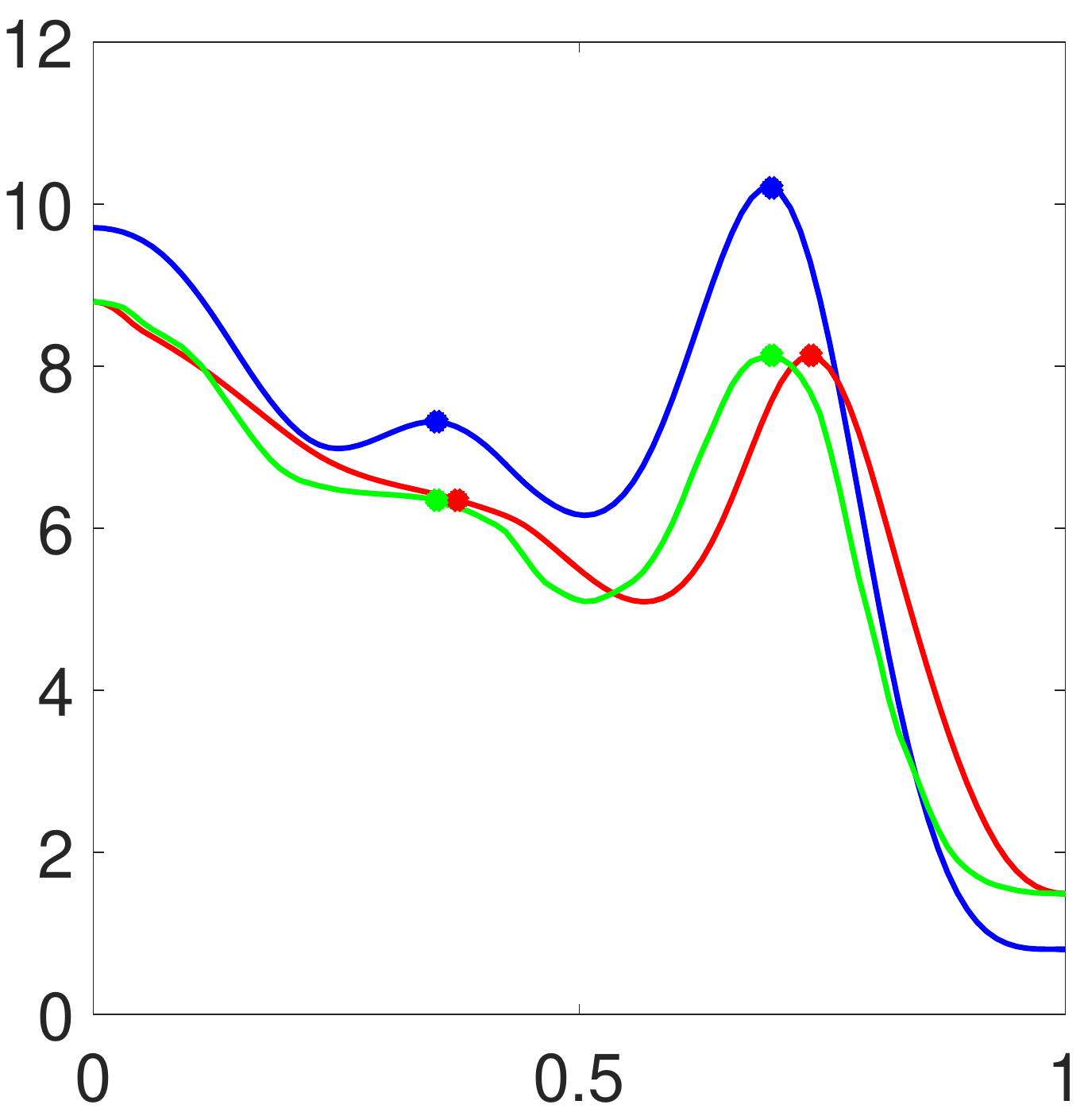}&\includegraphics[width=.9in]{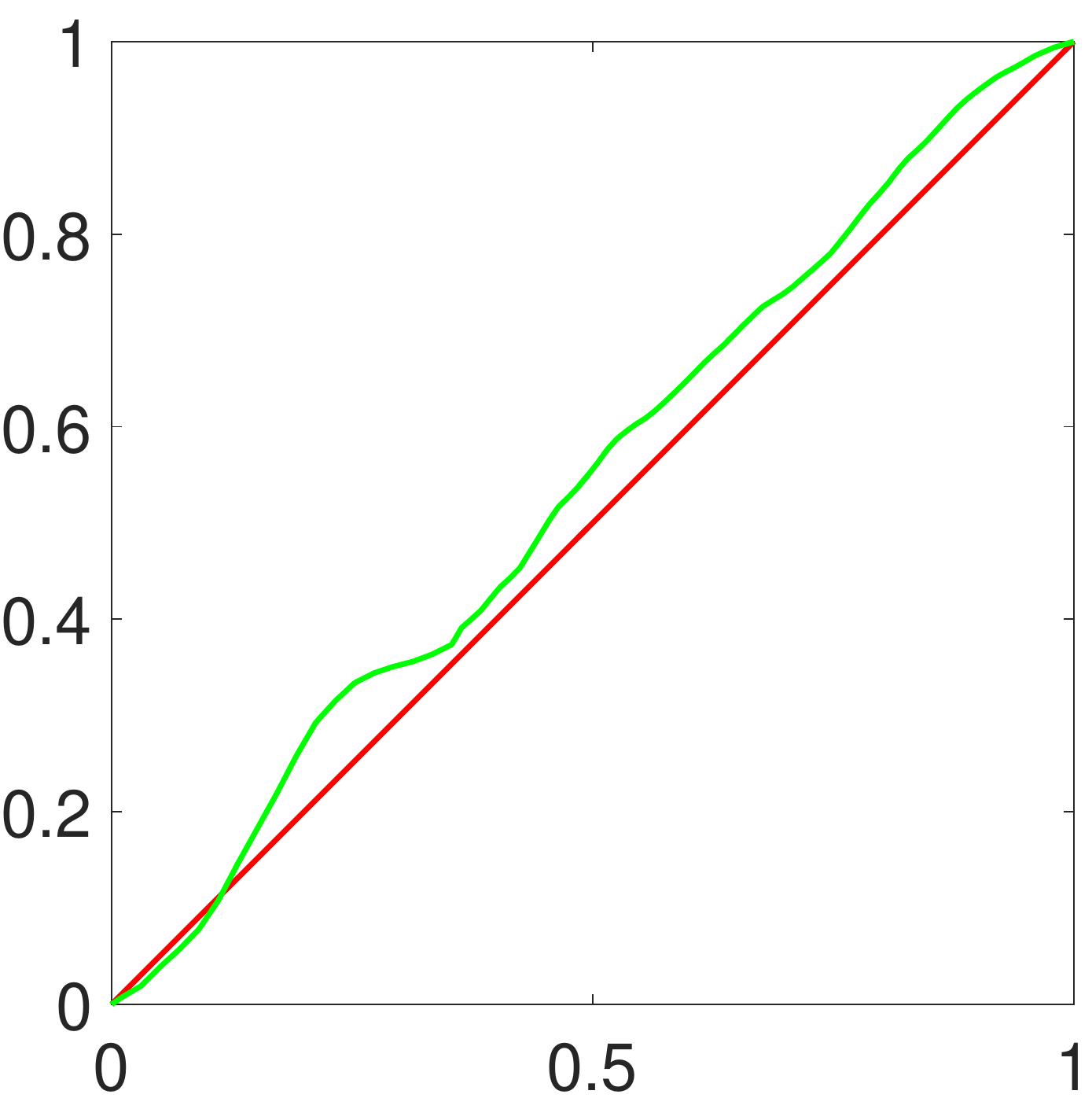}&\includegraphics[width=.9in]{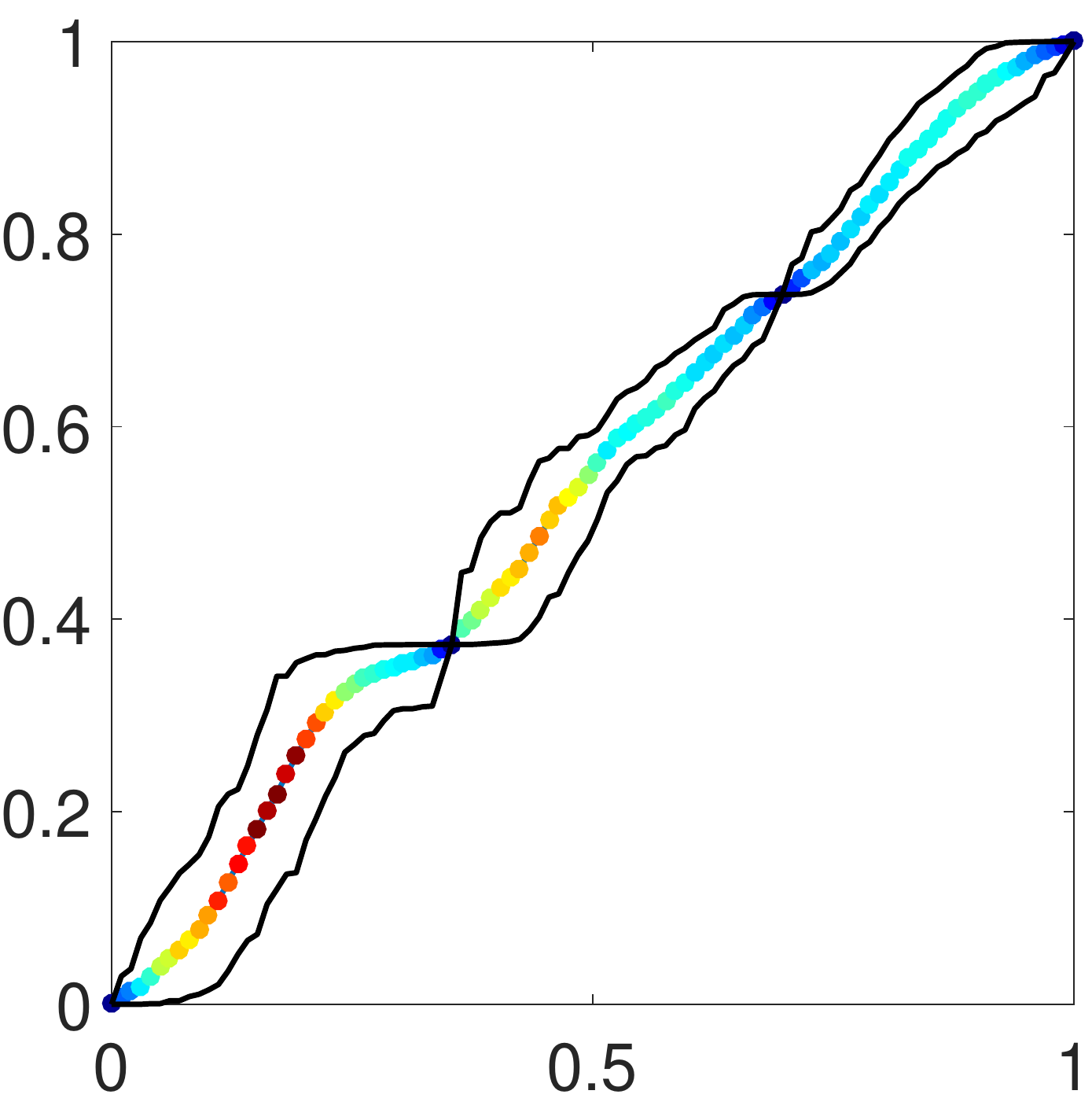}\\
				\hline
				\includegraphics[width=.9in]{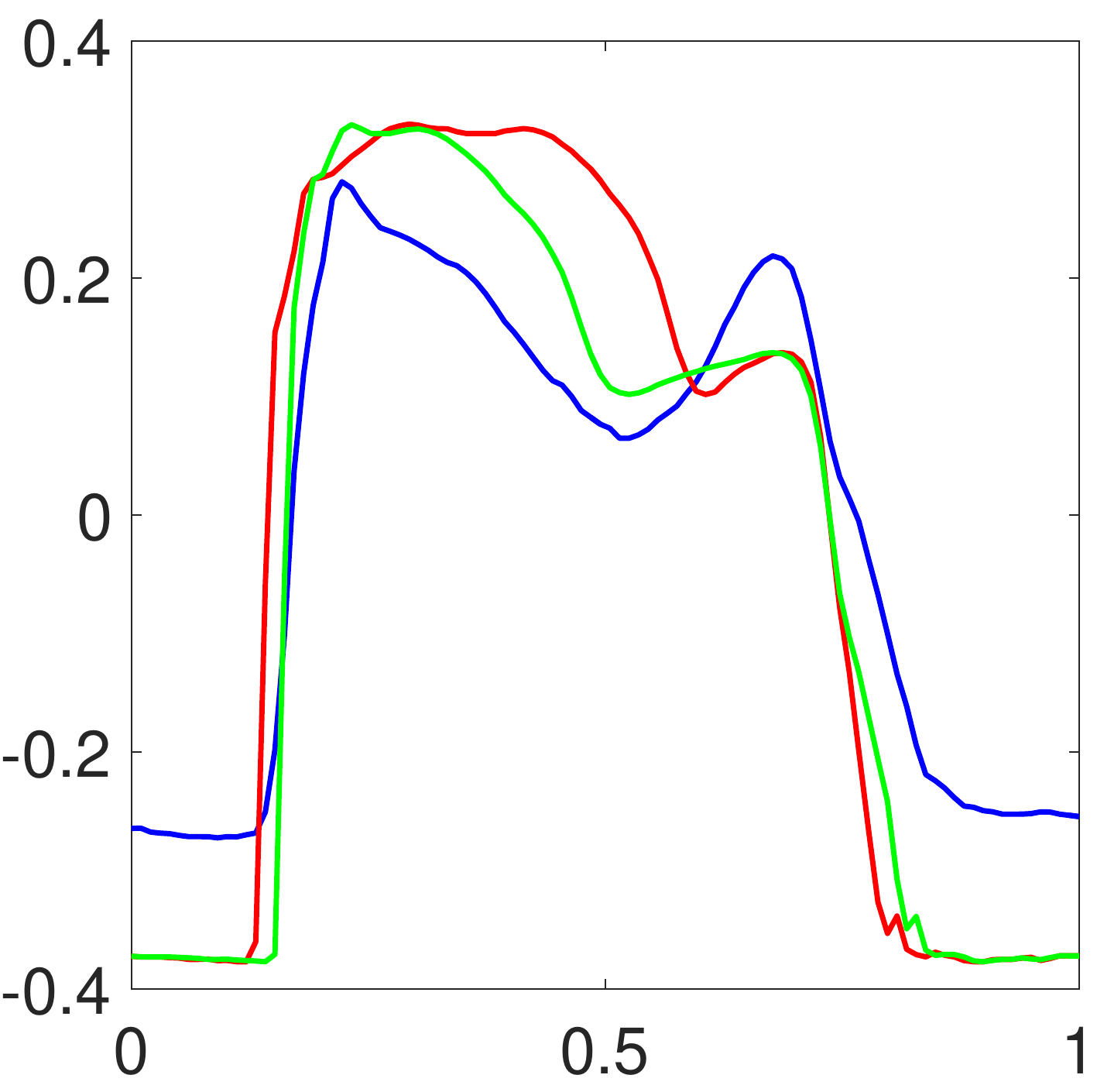}&\includegraphics[width=.9in]{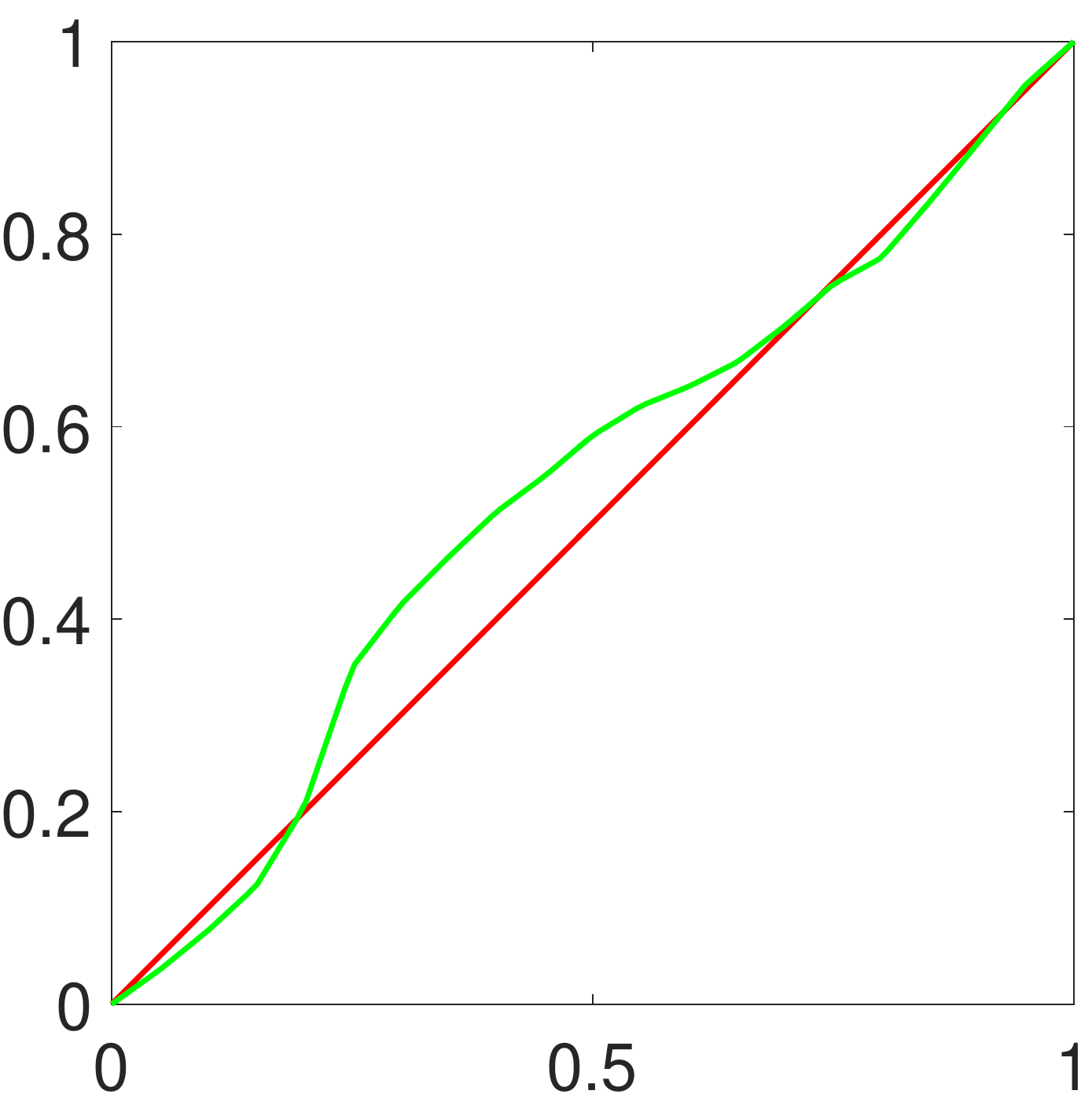}&\includegraphics[width=.9in]{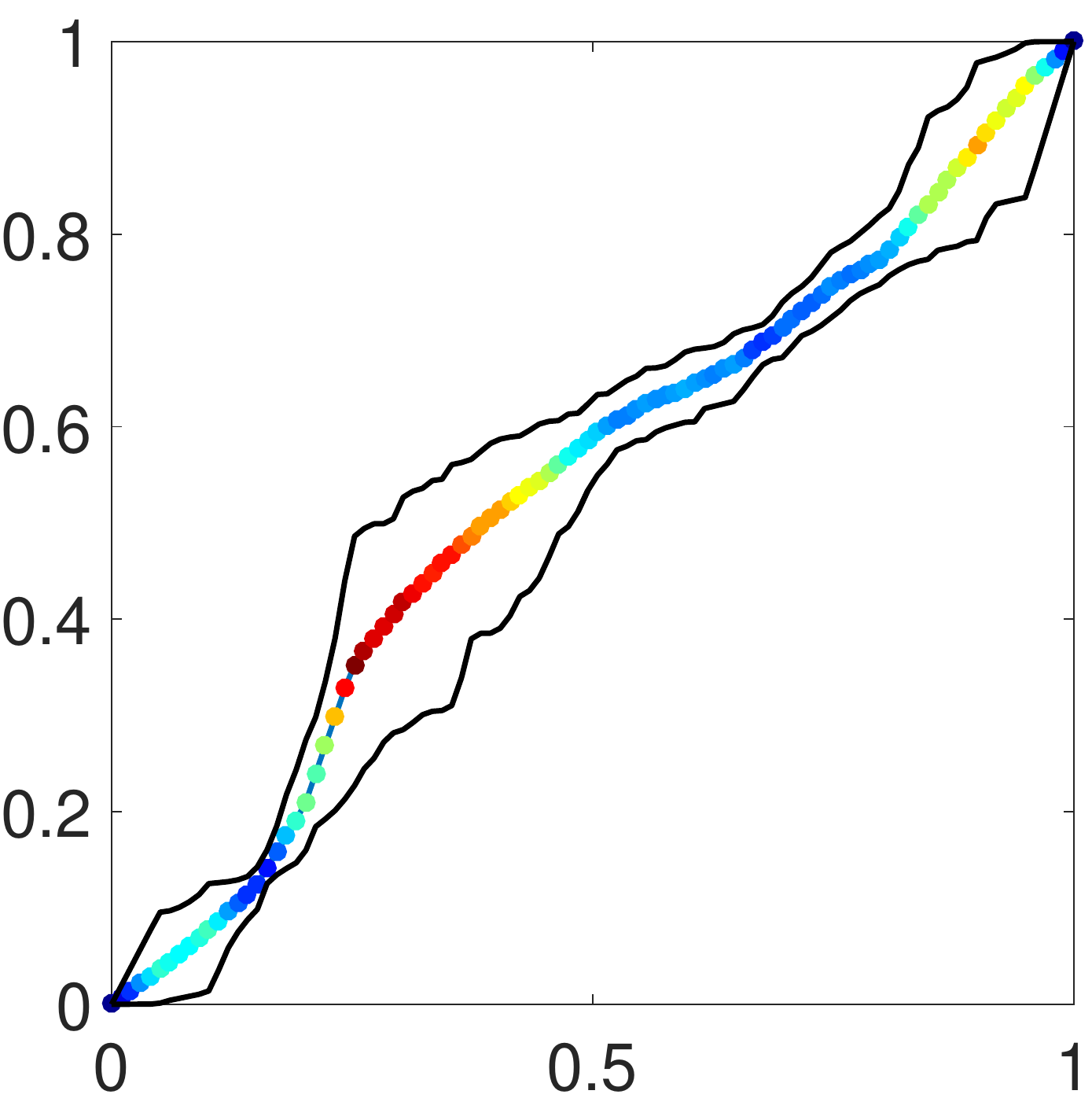}&\includegraphics[width=.9in]{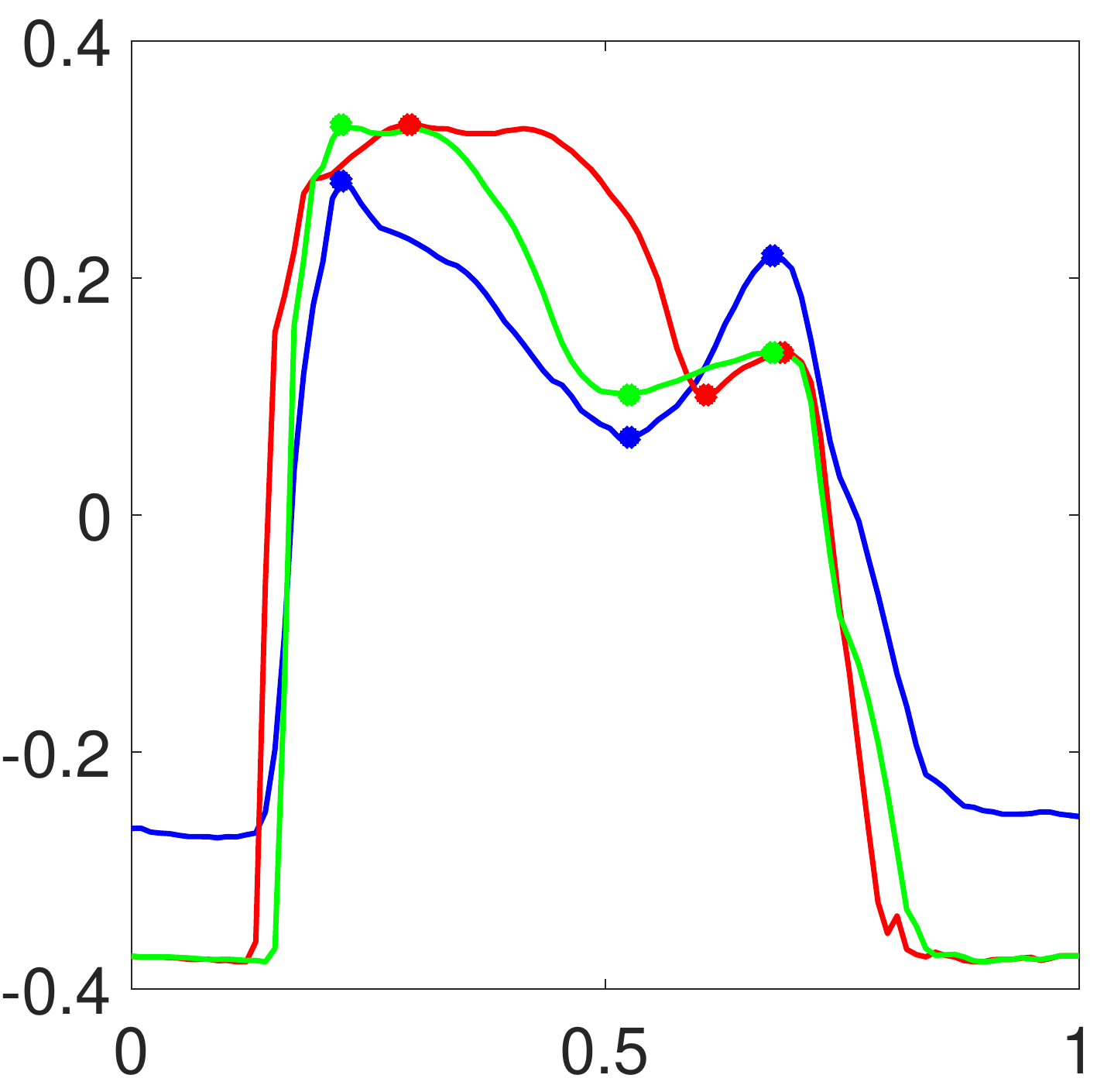}&\includegraphics[width=.9in]{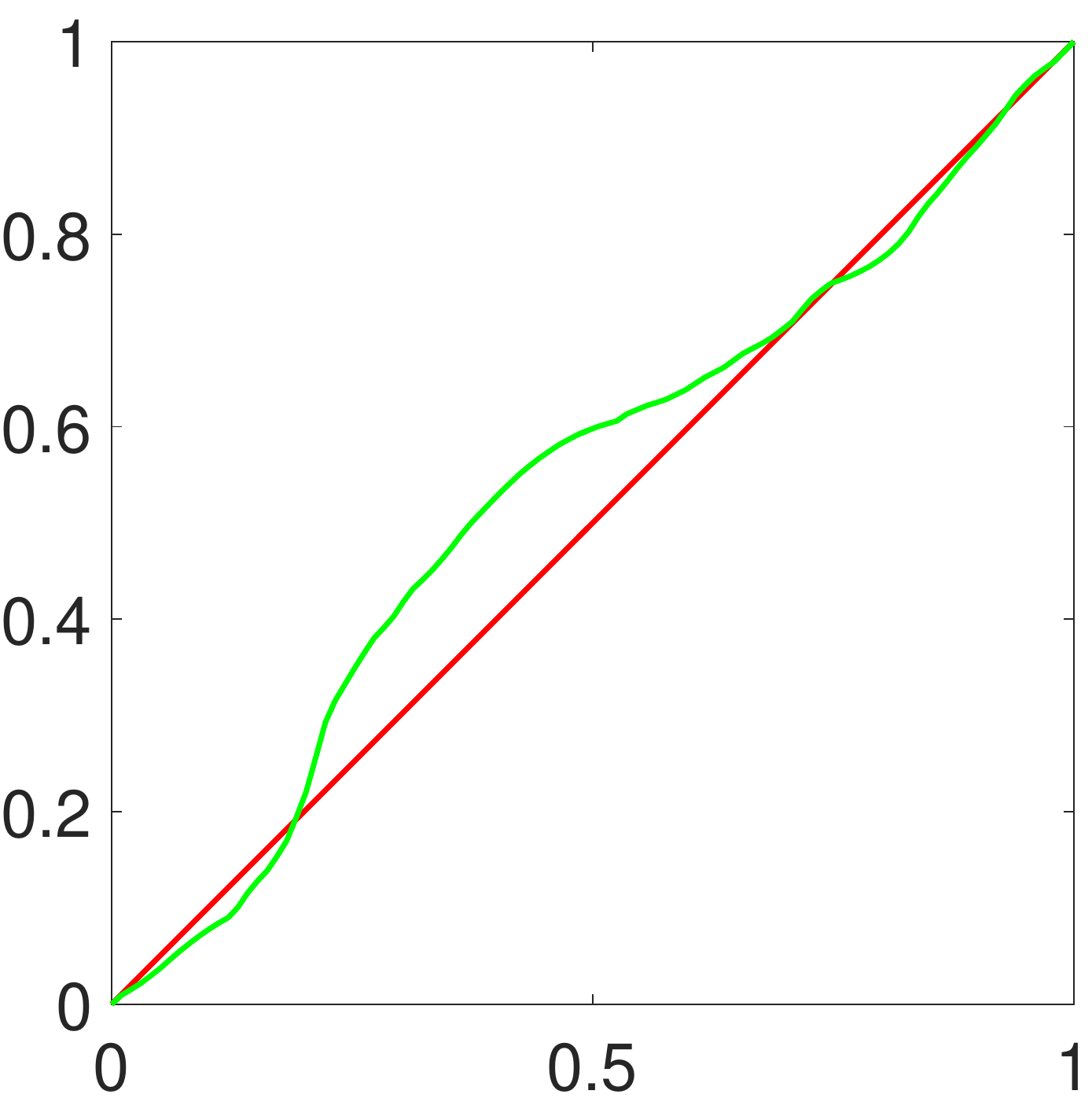}&\includegraphics[width=.9in]{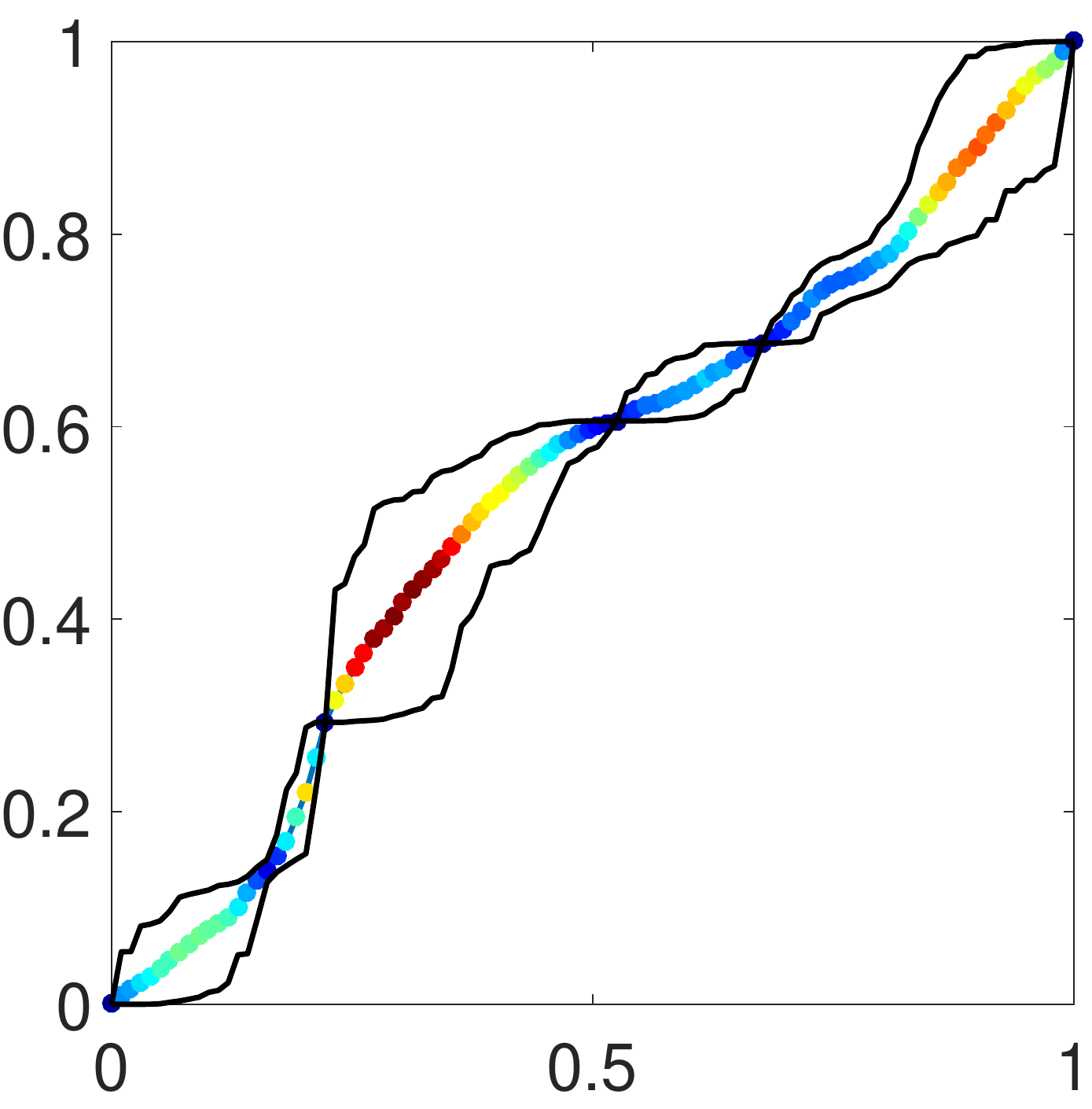}\\
				\hline
				\includegraphics[width=.9in]{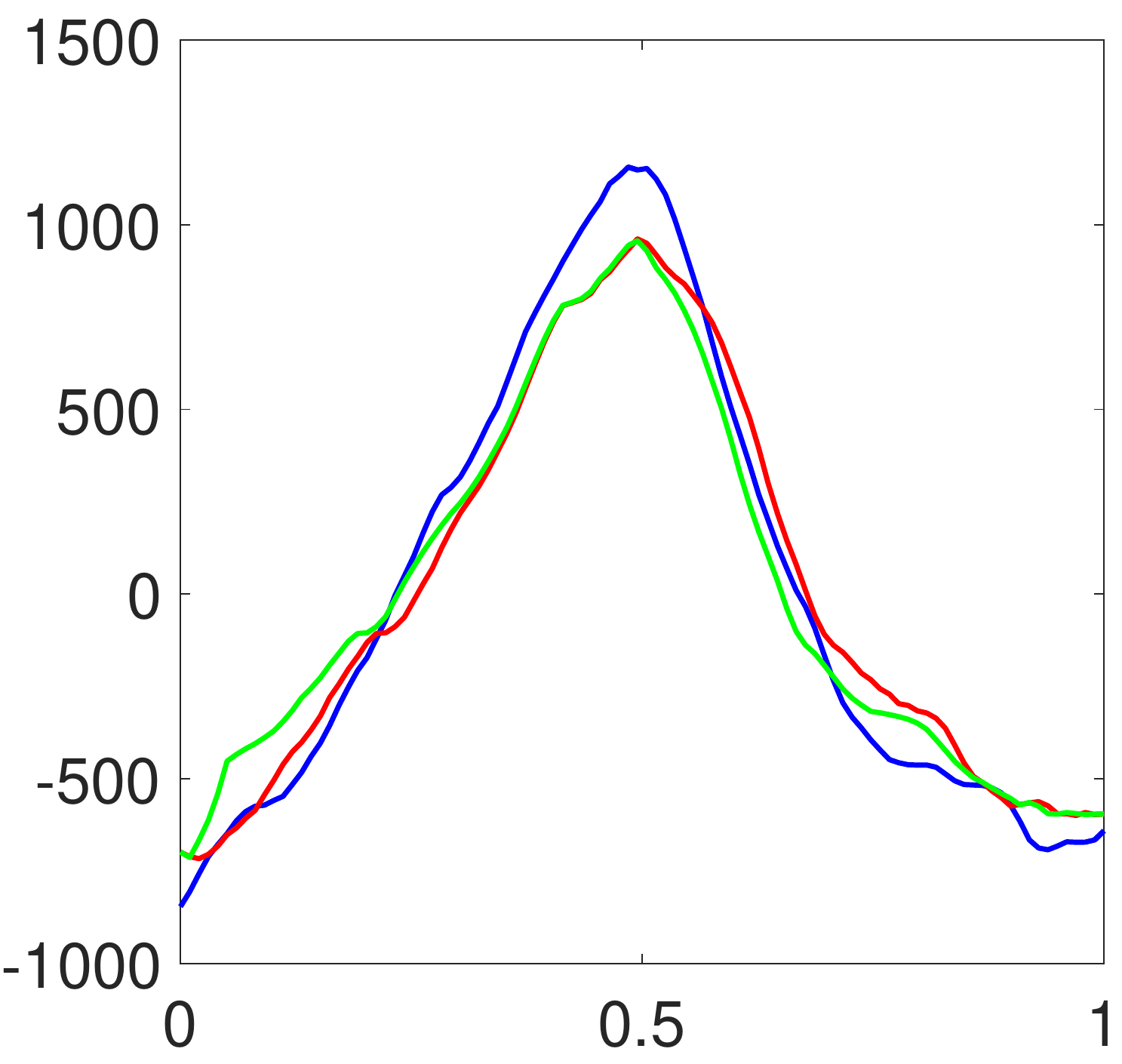}&\includegraphics[width=.9in]{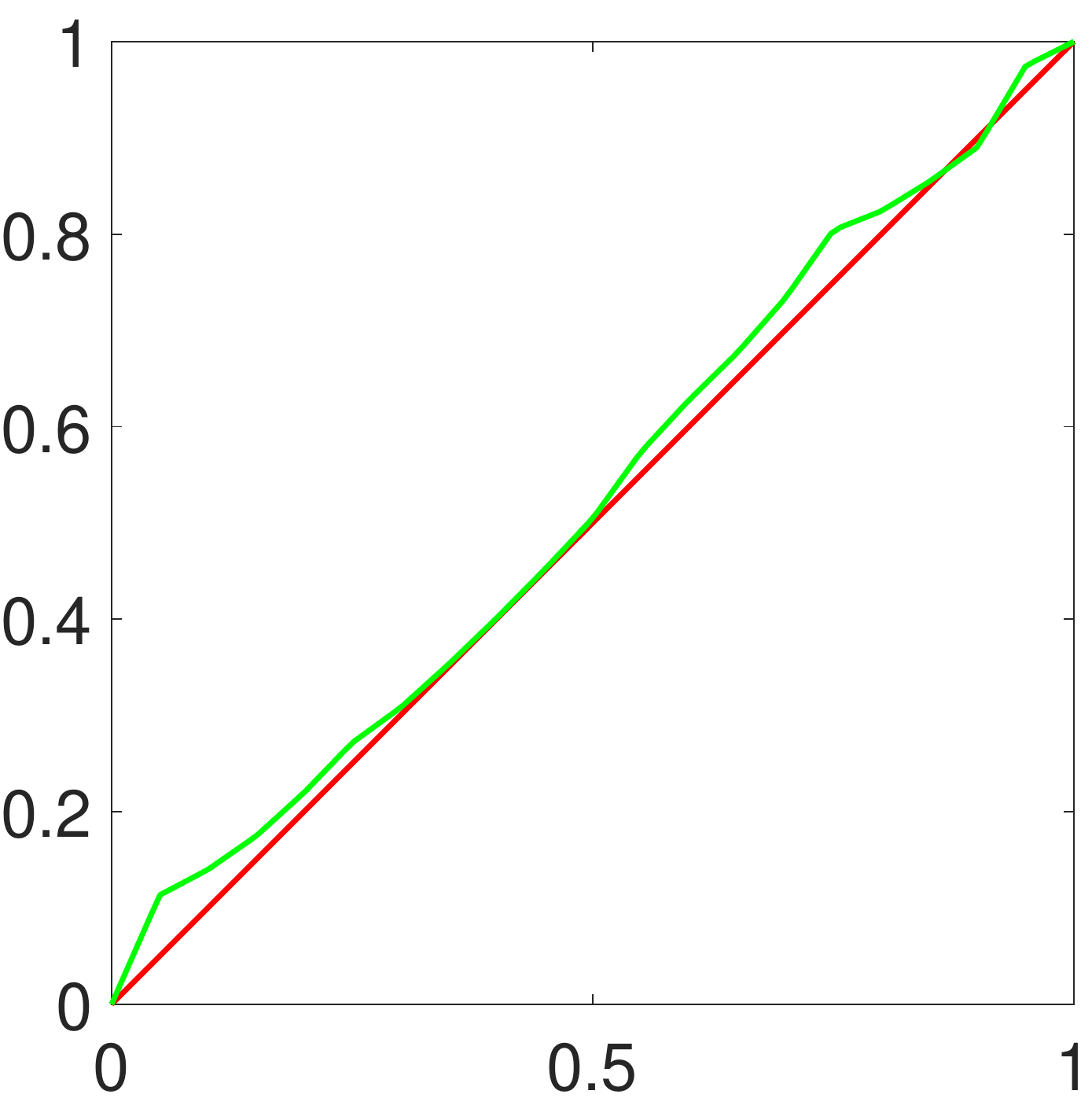}&\includegraphics[width=.9in]{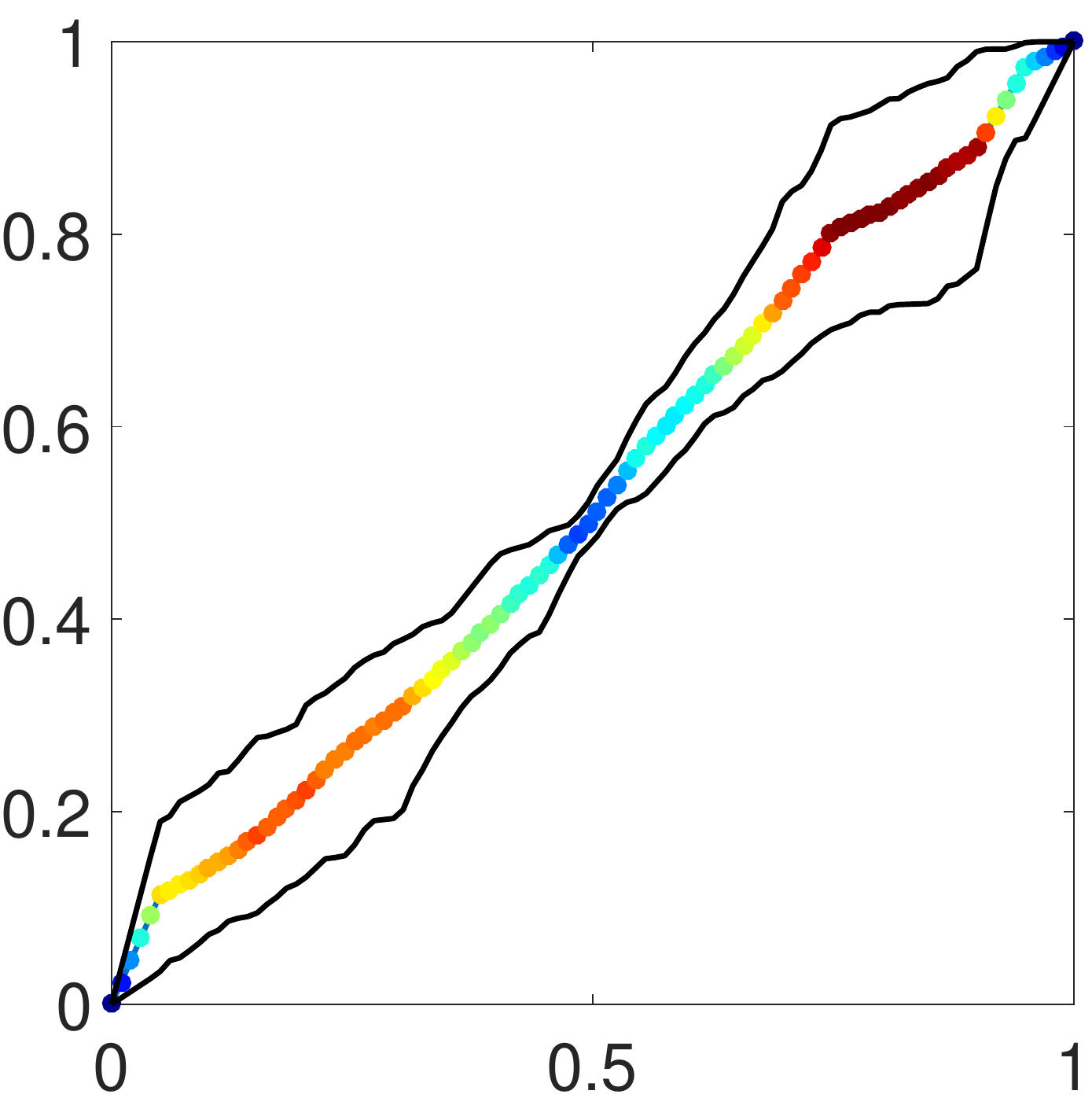}&\includegraphics[width=.9in]{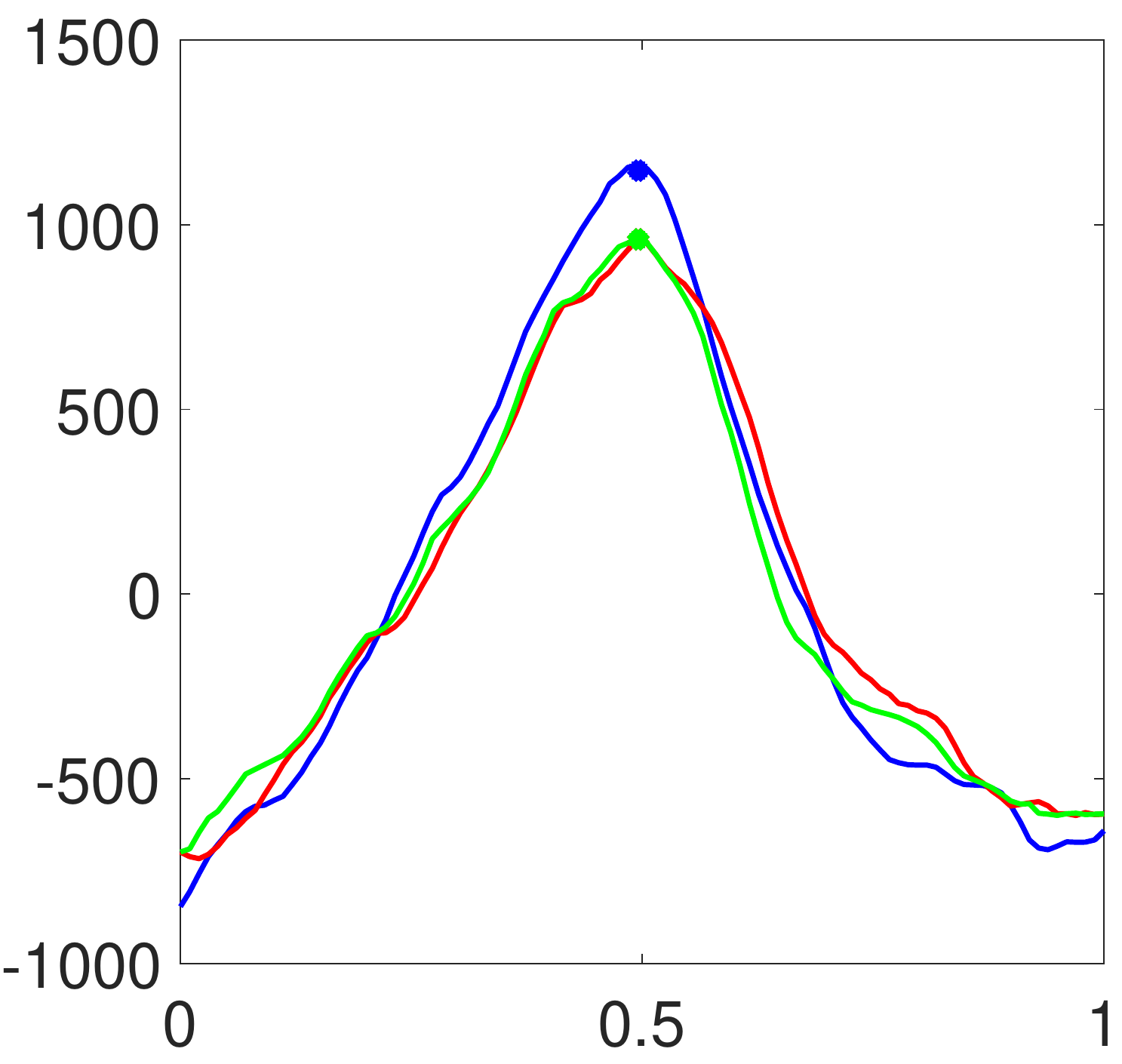}&\includegraphics[width=.9in]{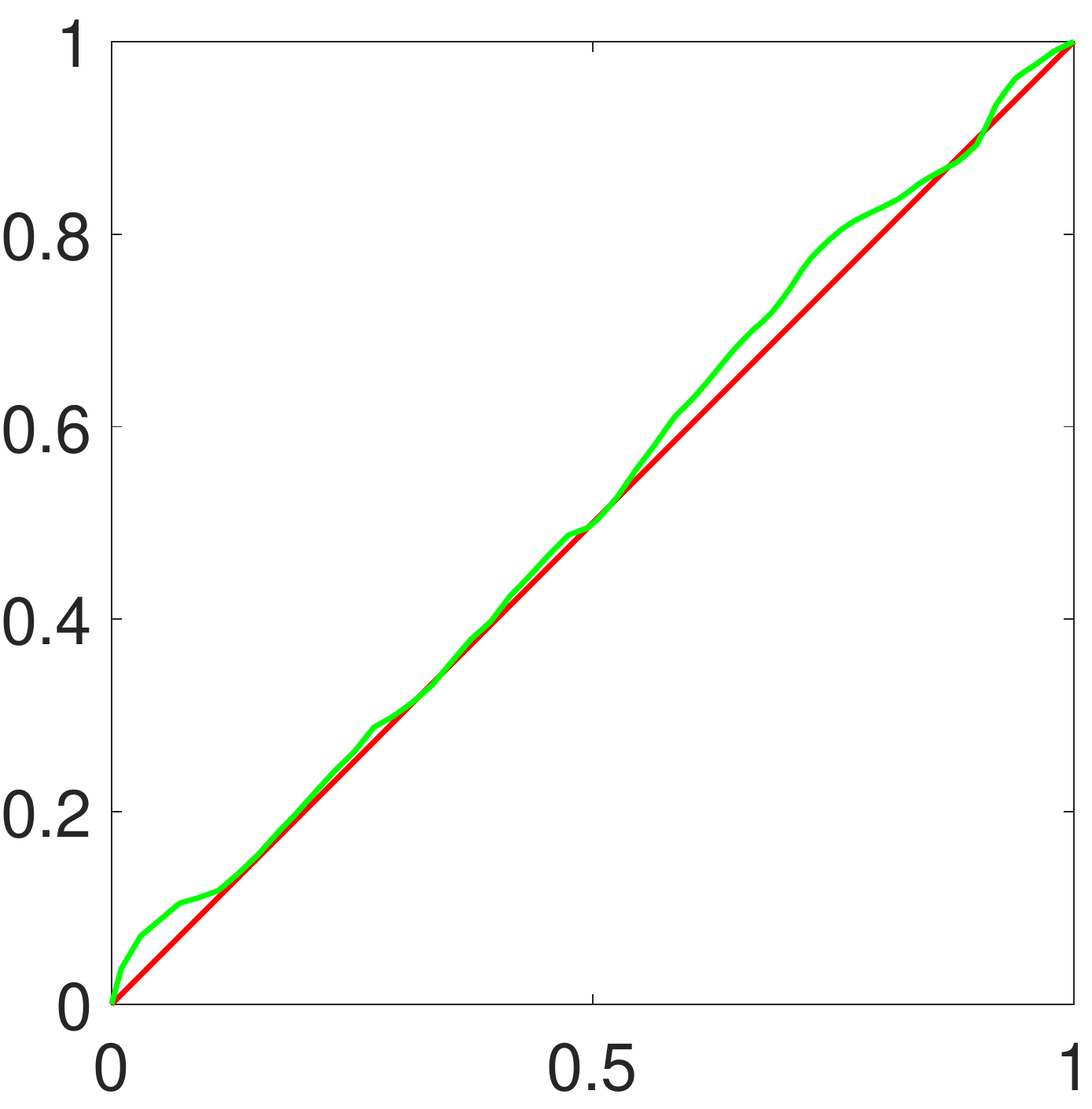}&\includegraphics[width=.9in]{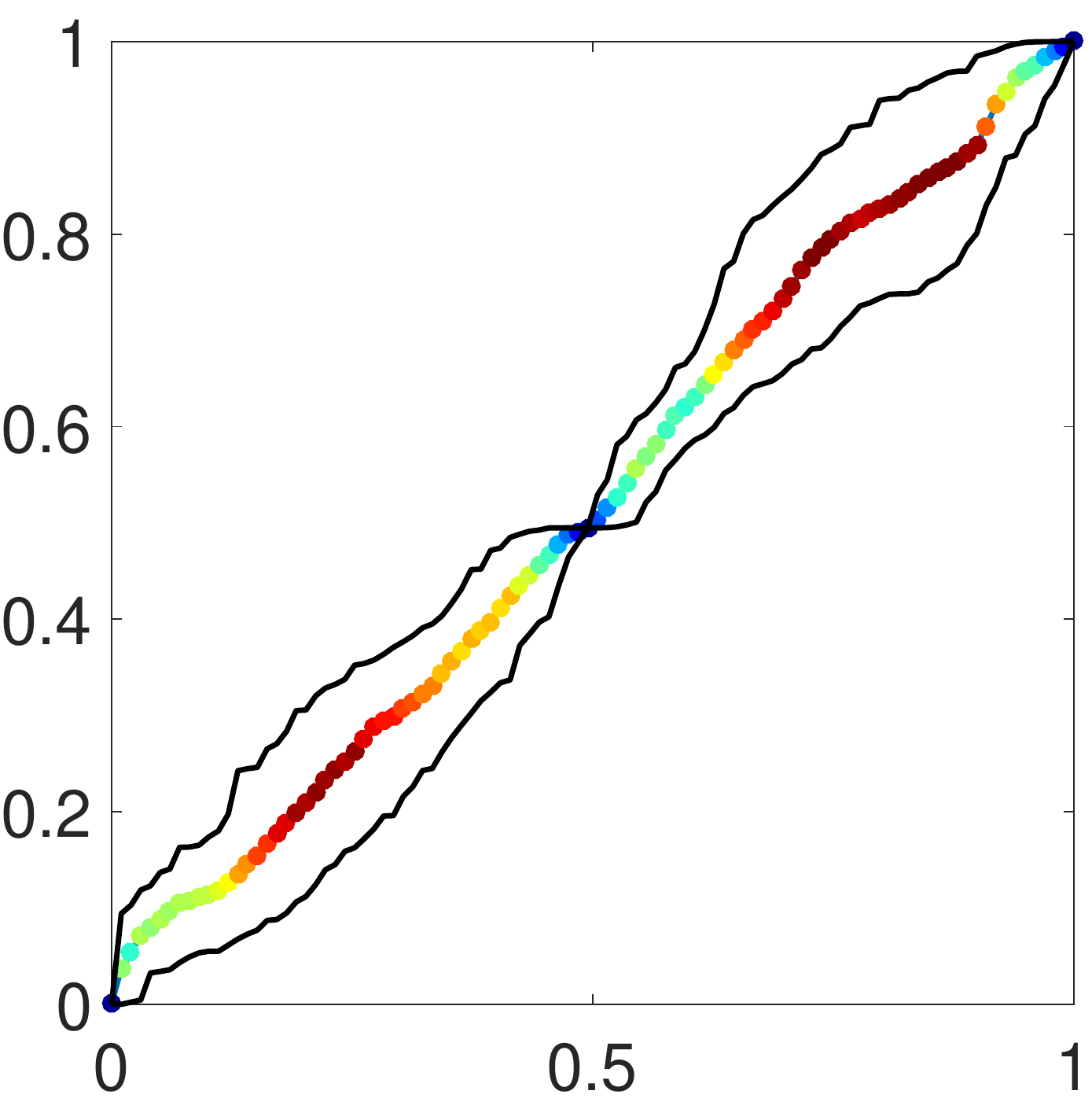}\\
				\hline
			\end{tabular}
		\end{adjustbox}
		\caption{\small Bayesian alignment of functions. From top to bottom: (1) simulated data, (2) PQRST complexes, (3) growth rate functions, (4) gait pressure cycles, (5) respiration cycles. (a) Two functions before alignment (blue and red), and red function after alignment (green); landmarks are marked in the constrained version. (b) Posterior mean (green) and identity (red) warp maps. (c) Posterior mean warp map with a pointwise $95\%$ credible interval (colors correspond to width of interval: blue=less uncertainty, red=more uncertainty).}\label{fig:BR}
	\end{center}
\end{figure}

The Bayesian model employed here is very similar to the one presented in \cite{Ian} and \cite{SK}. In short, the likelihood is a zero-mean multivariate Gaussian distribution with a diagonal covariance matrix. We select a conjugate, vague Gamma prior for the likelihood precision and analytically integrate it out of the posterior. While \cite{Ian} sampled from the posterior distribution using Markov chain Monte Carlo (MCMC), we use a simple sampling importance resampling (SIR) algorithm, with the importance function set to the prior distribution as in \cite{SK}. This may not be the best approach to sample from the posterior, but it provides a fast approximation and seems to work well in the settings we considered. We assess the utility of the proposed prior distribution on $W_I$ in two settings: (1) unconstrained, and (2) landmark-constrained function alignment. Note that the presented model can also be used for alignment of $d$-dimensional open and closed curves where $d>1$ with minor adjustments.

In order to sample from the distribution $\mathbb{D}_\theta \circ H$ on $W_I$ via the easy-to-implement Algorithm \ref{alg2}, we require three specifications: (1) choice of partition that determines the average warp map $H$, (2) $n$, which controls the size of the partition, and (3) $\theta$, which controls the spread around the average warp. We set $H(t)=t$, the uniform partition, which ensures regularization toward identity warping. We resample all functions with $100$ points and choose $n=20,\ \theta=10$. This gives flexibility in the prior to explore extreme warpings (small $\theta$) while also ensuring that the resulting warp maps are fairly regular (smaller partition prevents many small jumps). \emph{In the case of landmark-constrained alignment, owing to the properties of $\mathbb{D}_\theta \circ H$ in Proposition \ref{properties}, we re-scale the $\theta$ and $n$ proportionally to the length of each function segment, and consider each sub-problem independently.}

Figure \ref{fig:BR} presents results of unconstrained and landmark-constrained alignment for one simulated example and four real datasets. In the case of landmark-constrained alignment, the landmarks were selected either based on semantic features of the signals (e.g., PQRST points in a complex) or mathematical features (e.g., peaks and valleys). In panels (b) and (c), we show the posterior mean warp map (cross-sectional average of the posterior sample), and a pointwise $95\%$ credible interval. In all cases, the registration results are visually very good. Comparing columns (b) and (c) for the unconstrained and landmark-constrained cases, we observe intuitive differences in the posterior mean warp maps and corresponding credible intervals. The key observation is that in each of the datasets, when using additional information provided by the landmarks, the distribution $\mathbb{D}_\theta \circ H$ allows us to decompose the alignment problem into unconstrained sub-problems by enabling subset invariance; the constraints on the warp maps, and the properties of $\mathbb{D}_\theta \circ H$, ensure that there is almost zero uncertainty in regions close to the landmarks.

\subsubsection{Curve alignment with Simulated Annealing}
\label{simulated_annealing}
Based on the distributions $\mathbb{D}_\theta \circ H$ and $\mathbb{D}^s_\theta \circ H$, we present a novel stochastic algorithm for unconstrained alignment of three different types of functions: (1) univariate functions: $g:[0,1]\to\real$, (2) shapes of 3D open curves: $g:[0,1]\to\real^3$, and (3) shapes of planar closed curves: $g:\sone\to\real^2$. A crucial step of the algorithm is based on the ability to propose warp maps in the neighborhood of any other warp map. The distributions $\mathbb{D^\theta}$ and $\mathbb{D}_s^\theta$ are well-suited for this purpose: (1) we can centre the distributions at any warp map, and (2) we can control the size of the neighborhood via the parameter $\theta$. We only consider the unconstrained case; the extension to landmark-constrained alignment, as demonstrated in the previous section, is achieved through subset invariance properties of the distributions.

As with the Bayesian model, we base our stochastic alignment algorithm on the SRVF representation of curves. The energy functional that we seek to optimize is $E(\gamma)=\|q_1-(q_2\circ\gamma)\sqrt{\dot{\gamma}}\|^2$, for $\gamma$ in $W_I$ or $W_\sone$. When the domain of $q_1$ and $q_2$ is $[0,1]$, a solution to this optimization problem can be obtained using a Dynamic Programming (DP) algorithm \citep{DynPro}. The resulting solution is deterministic and depends on the fineness of the discretization and the size of the neighborhood search. In the case of closed curves, one has to either resort to a gradient descent algorithm \citep{SrivESA}, which has the obvious limitation of getting stuck in a local solution, or a DP approach with an additional seed (the point at which $\sone$ is unwrapped to $[0,1]$) search, which only gives an approximate solution. The energy functional $E$ is a natural choice for curve registration, because (1) the $\ltwo$ distance in the energy corresponds to an elastic metric on the space of curves, and (2) this elastic metric is preserved under identical warping (isometry); see \cite{AK} for details.

The Simulated Annealing alignment algorithm for functional data (i.e., $g:[0,1]\to\real$) using $\mathbb{D}_\theta \circ H$ on $W_I$ is given as Algorithm \ref{alg4}. The extension to alignment of open and closed curves for the purpose of shape analysis is commented upon below.
\begin{algorithm}
	\textbf{Alignment of functions via Simulated Annealing.}\\
	\noindent Inputs: $g_i:[0,1]\to\real,\ i=1,2$ (SRVFs $q_i:[0,1]\to\real^d,\ i=1,2$).\\
	Outputs: Optimal warp map $\gamma^*:[0,1]\to[0,1]$.\\
	Initialize: $n=20$, $\theta=100$, $T=10$, $\gamma_0=\gamma_{id}$, $E_0=\|q_1-q_2\|^2$ and $j=0$.
	\begin{enumerate}
	\itemsep 0em
		\item Generate a random $\tilde{\gamma}$ from $\mathbb{D_\theta}\circ H$ with $H$ set to the warp map $\gamma_j$. Set $\gamma_p=0.9\tilde{\gamma}+0.1\gamma_{id}$.
		\item Compute $E(\gamma_p)=\|q_1-(q_2\circ\gamma_p)\sqrt{\dot{\gamma}_p}\|^2$.
		\item Accept $\gamma_{j+1}=\gamma_p$, and set $E(\gamma_{j+1})=E(\gamma_p)$, with probability $\min\left\{1,e^{\frac{E(\gamma_j)-E(\gamma_p)}{T}}\right\}$. Otherwise, let $\gamma_{j+1}=\gamma_j$ and $E(\gamma_{j+1})=E(\gamma_j)$.
		\item Set $j=j+1$ and update the temperature to $T=T/c$ (we suggest $c=1.0001$).
	\end{enumerate}\label{alg4}
\end{algorithm}
\noindent The algorithm is fairly robust to the choices of $n, \theta, T$ and $c$. Numerical illustrations investigating robustness are presented in the Supplementary Material. In the first step of the algorithm, we propose a candidate warp map that is a linear combination of a random warping sampled from the distribution $\mathbb{D_\theta}\circ H$ centred at the previously accepted warp map and the identity warping; while not necessary, this ensures extra regularization toward identity warping. The key here is that efficient sampling from $\mathbb{D}_\theta\circ H$ centred at an arbitrary warp map is easily enabled through Algorithm \ref{alg2}. Furthermore, the relative values of the concentration parameter $\theta$ (neighborhood size) and the $T$ control the dynamics of the algorithm.

{\small
	\begin{figure}[!t]
		\begin{center}
			\begin{tabular}{|cc||cc|}
				\hline
				(a)&(b)&(a)&(b)\\
				\hline
				\multicolumn{2}{|c||}{(1)}&\multicolumn{2}{c|}{(2)}\\
				\hline
				\includegraphics[width=1in]{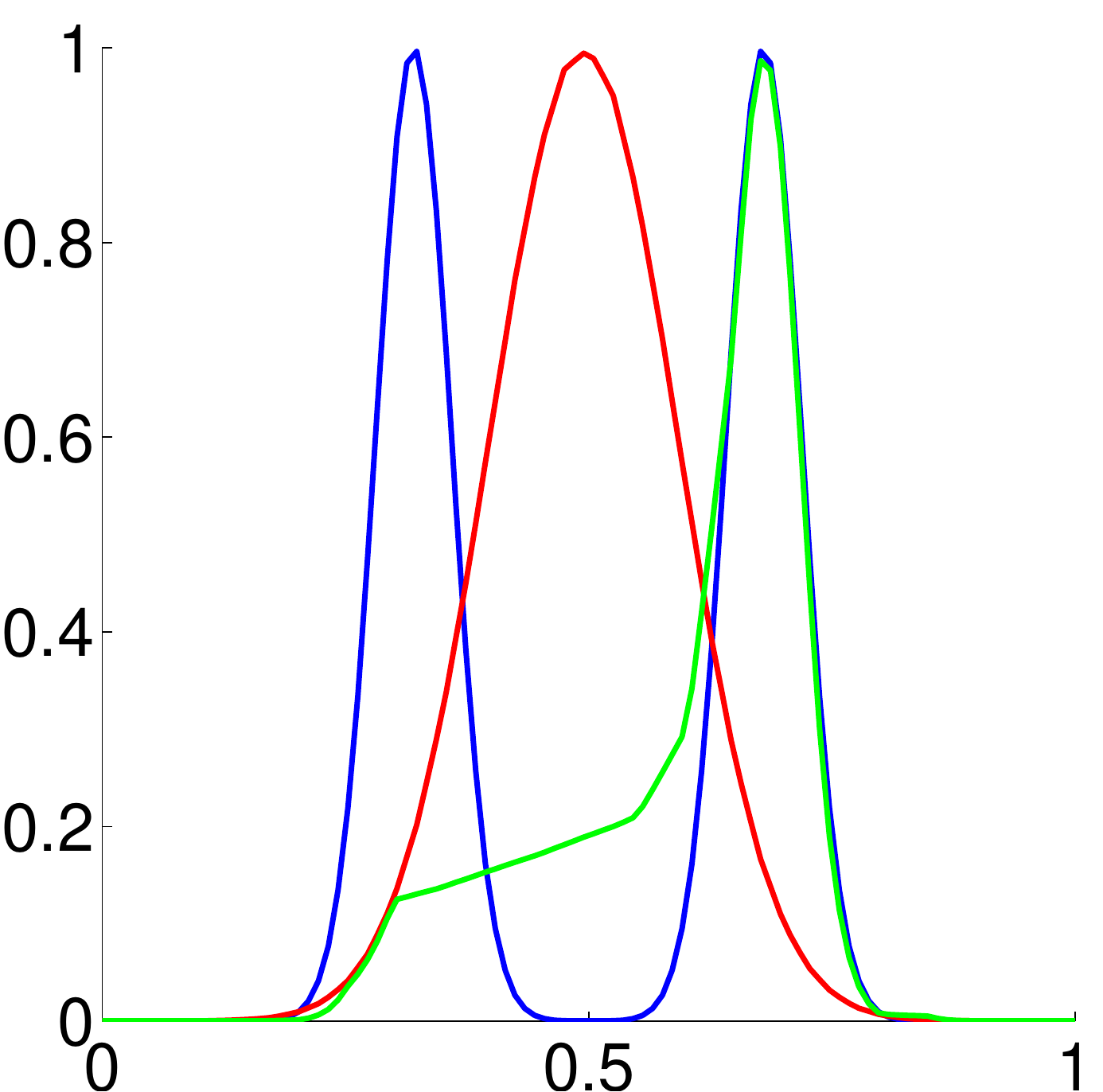}&\includegraphics[width=1in]{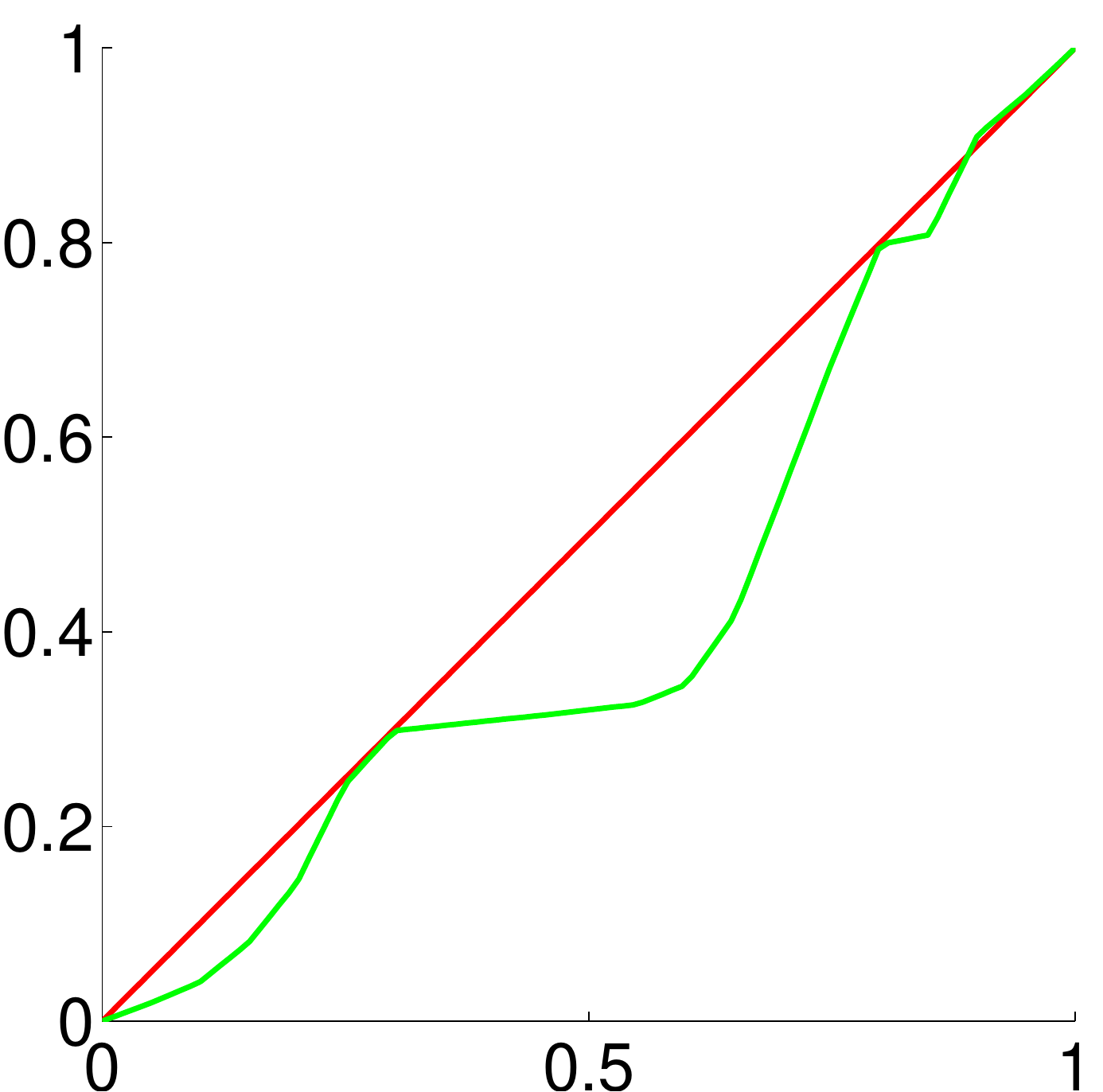}&\includegraphics[width=1in]{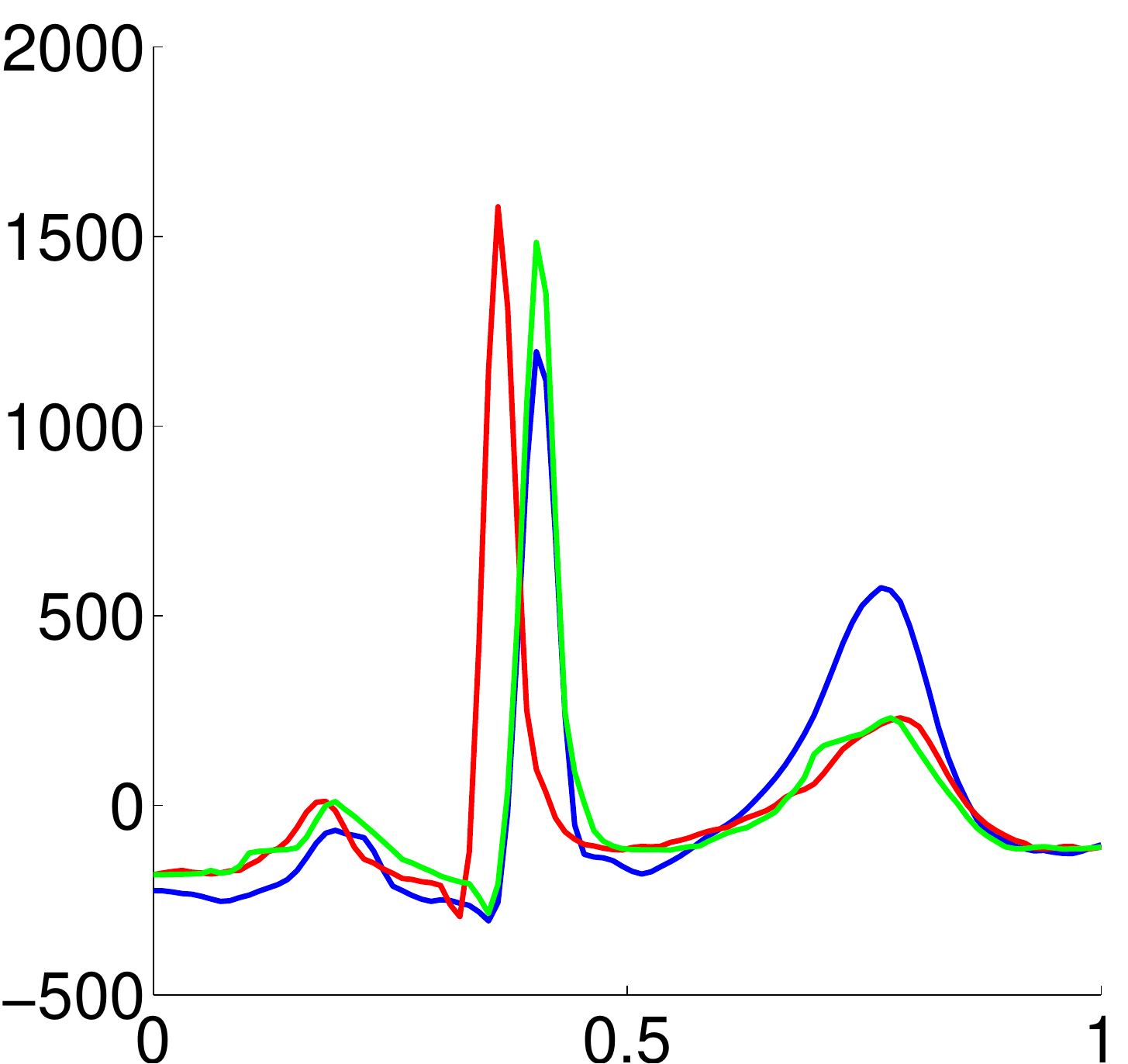}&\includegraphics[width=1in]{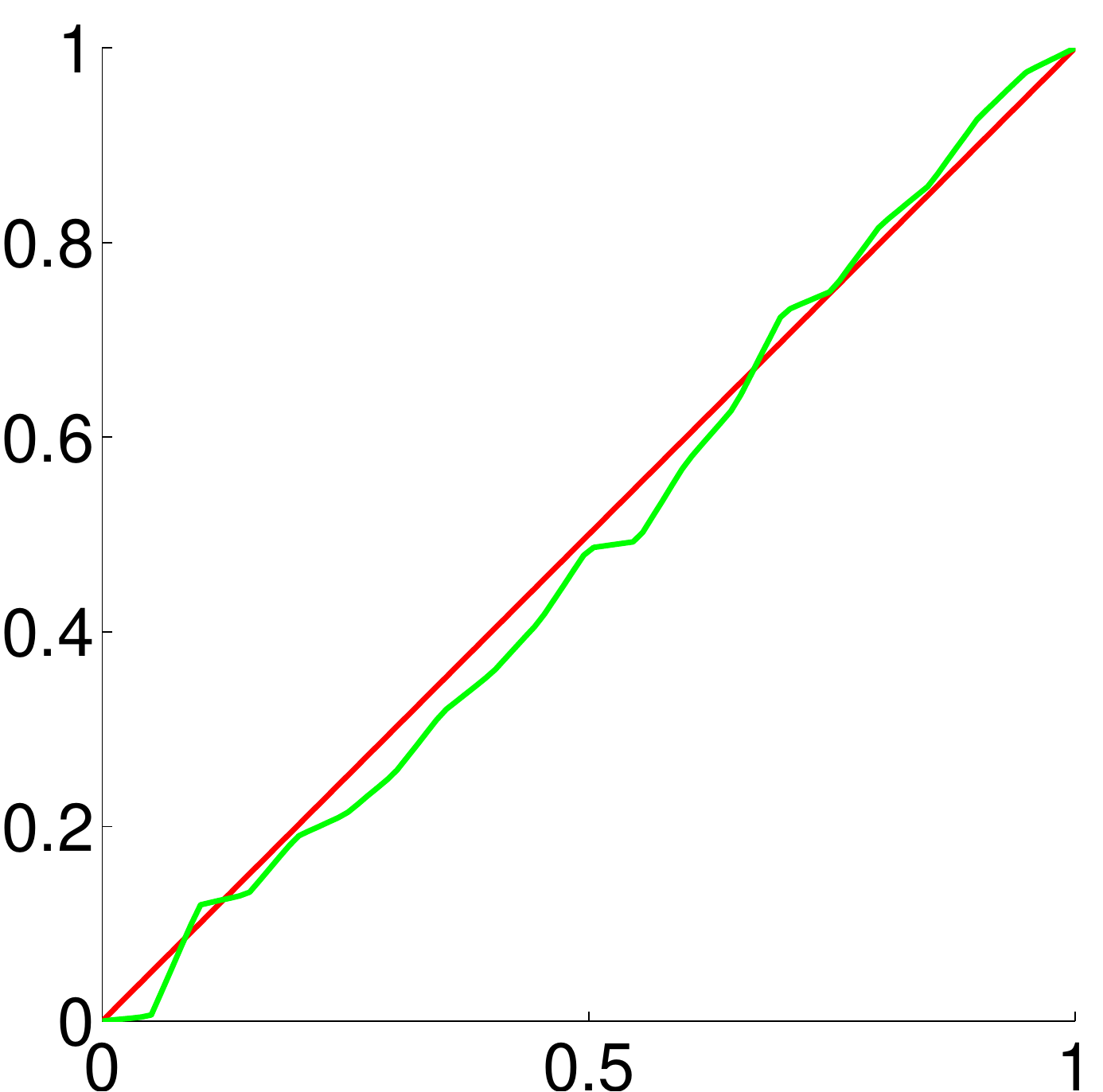}\\
				\hline
				\multicolumn{2}{|c||}{(3)}&\multicolumn{2}{c|}{(4)}\\
				\hline
				\includegraphics[width=1in]{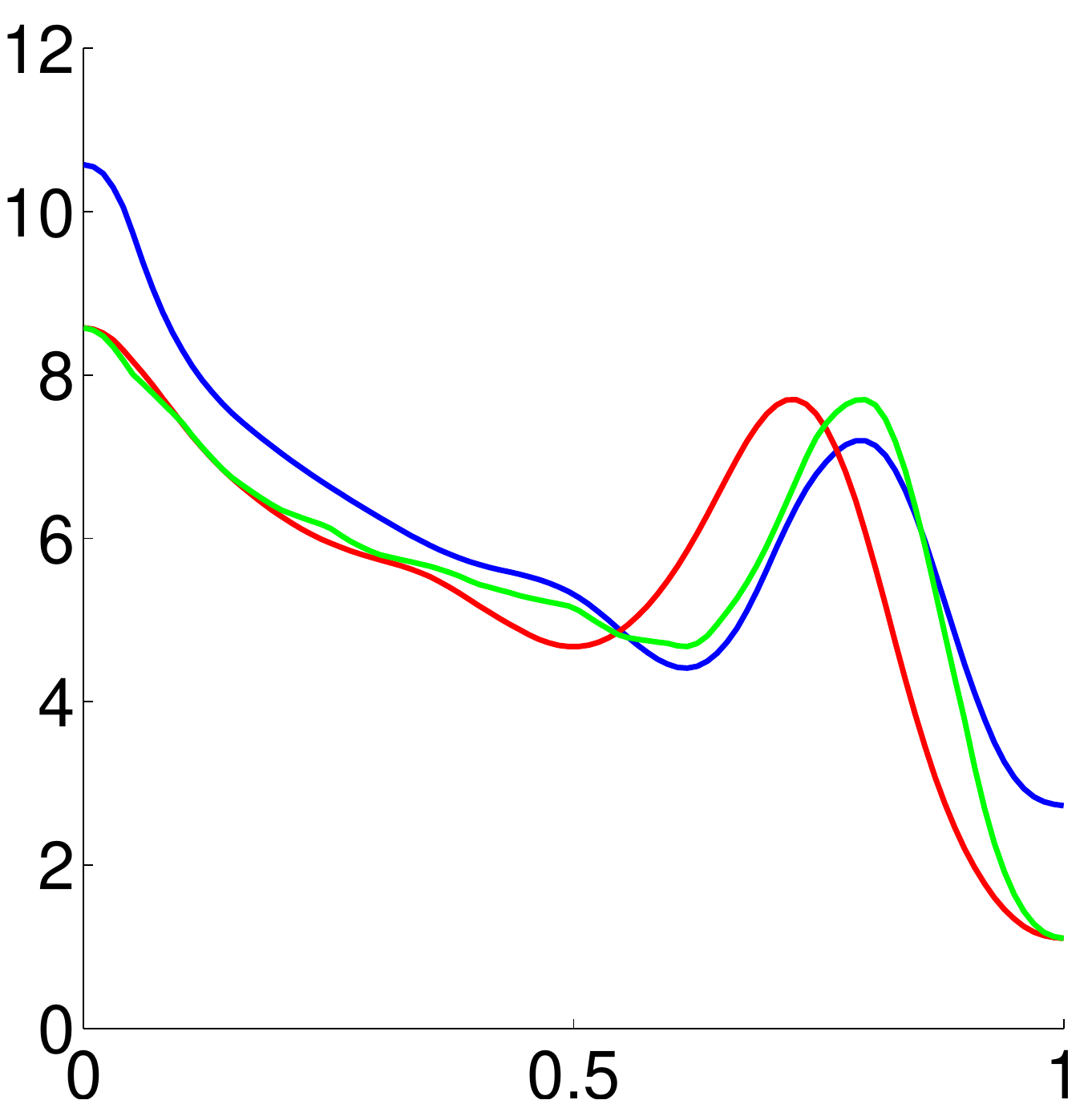}&\includegraphics[width=1in]{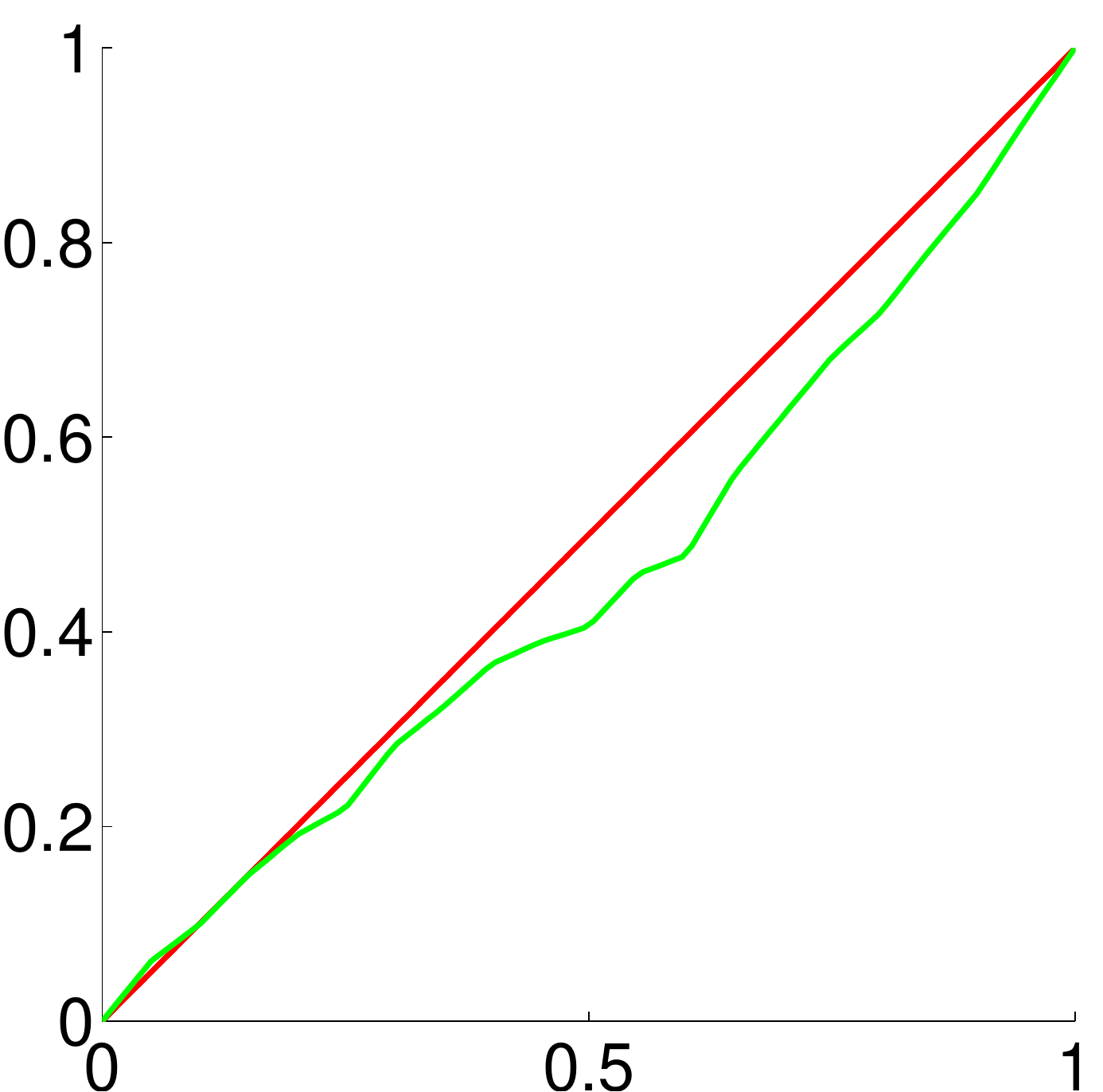}&\includegraphics[width=1in]{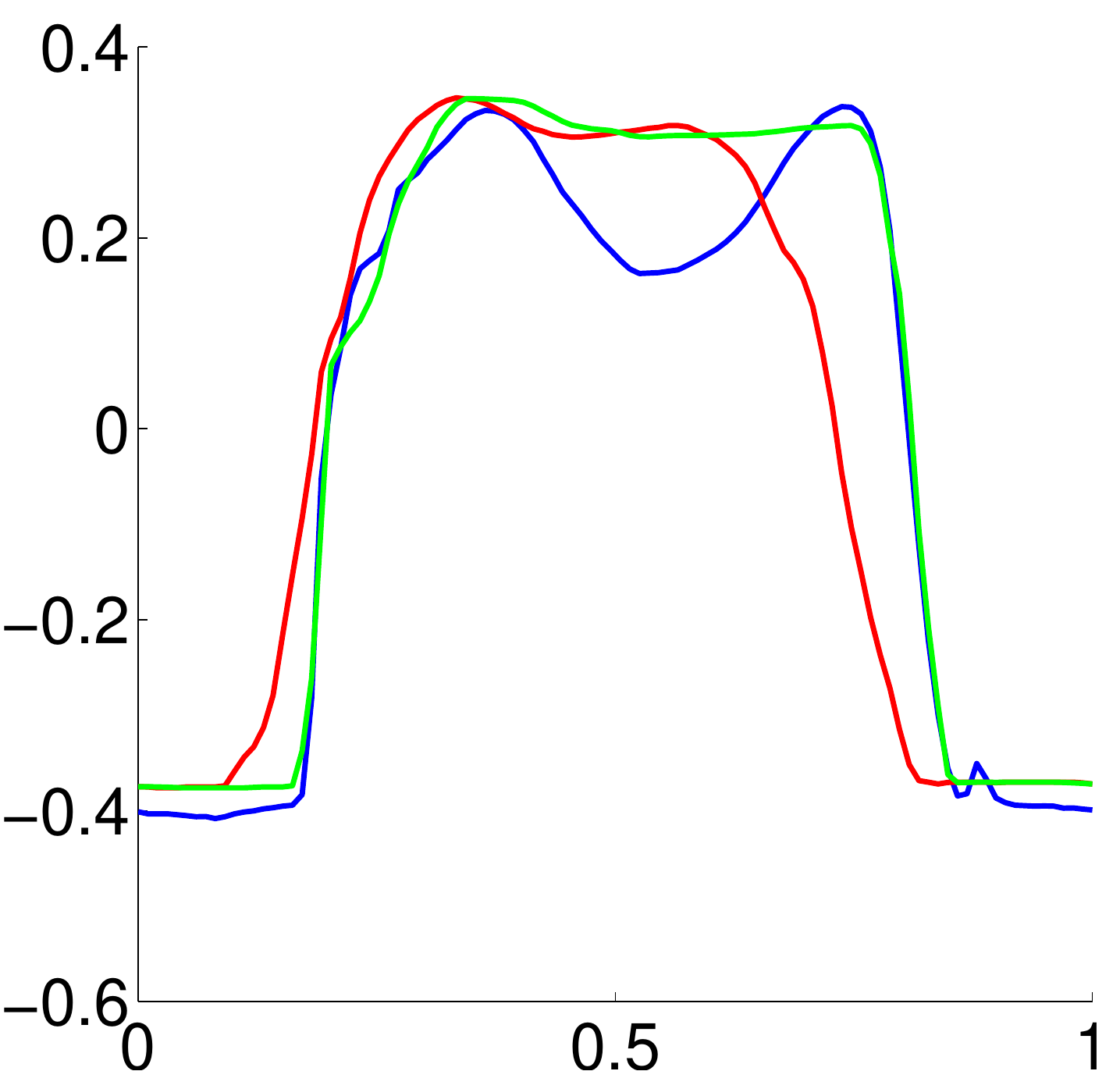}&\includegraphics[width=1in]{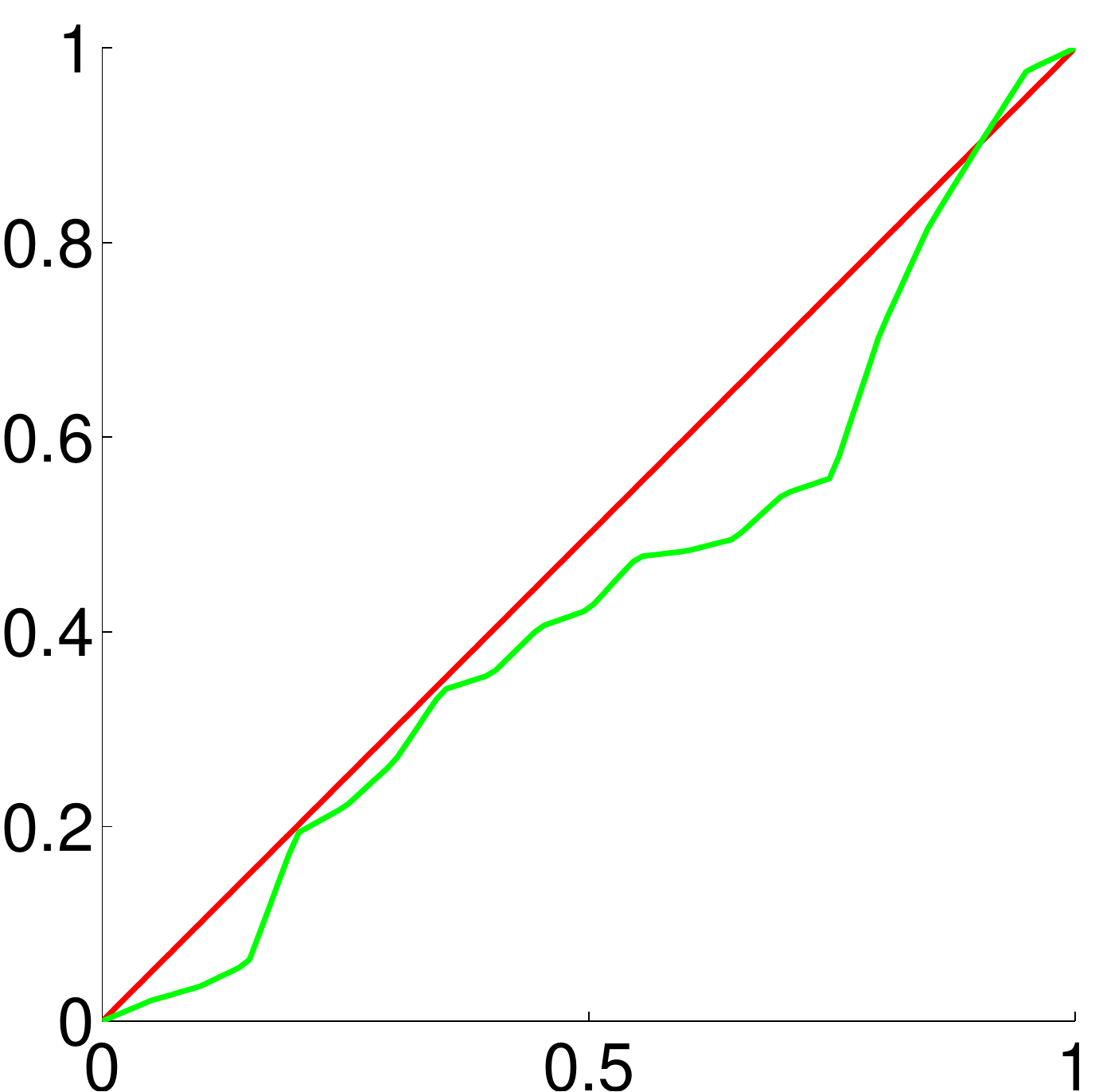}\\
				\hline
				\multicolumn{2}{|c||}{(5)}&\multicolumn{2}{c|}{(6)}\\
				\hline
				\includegraphics[width=1in]{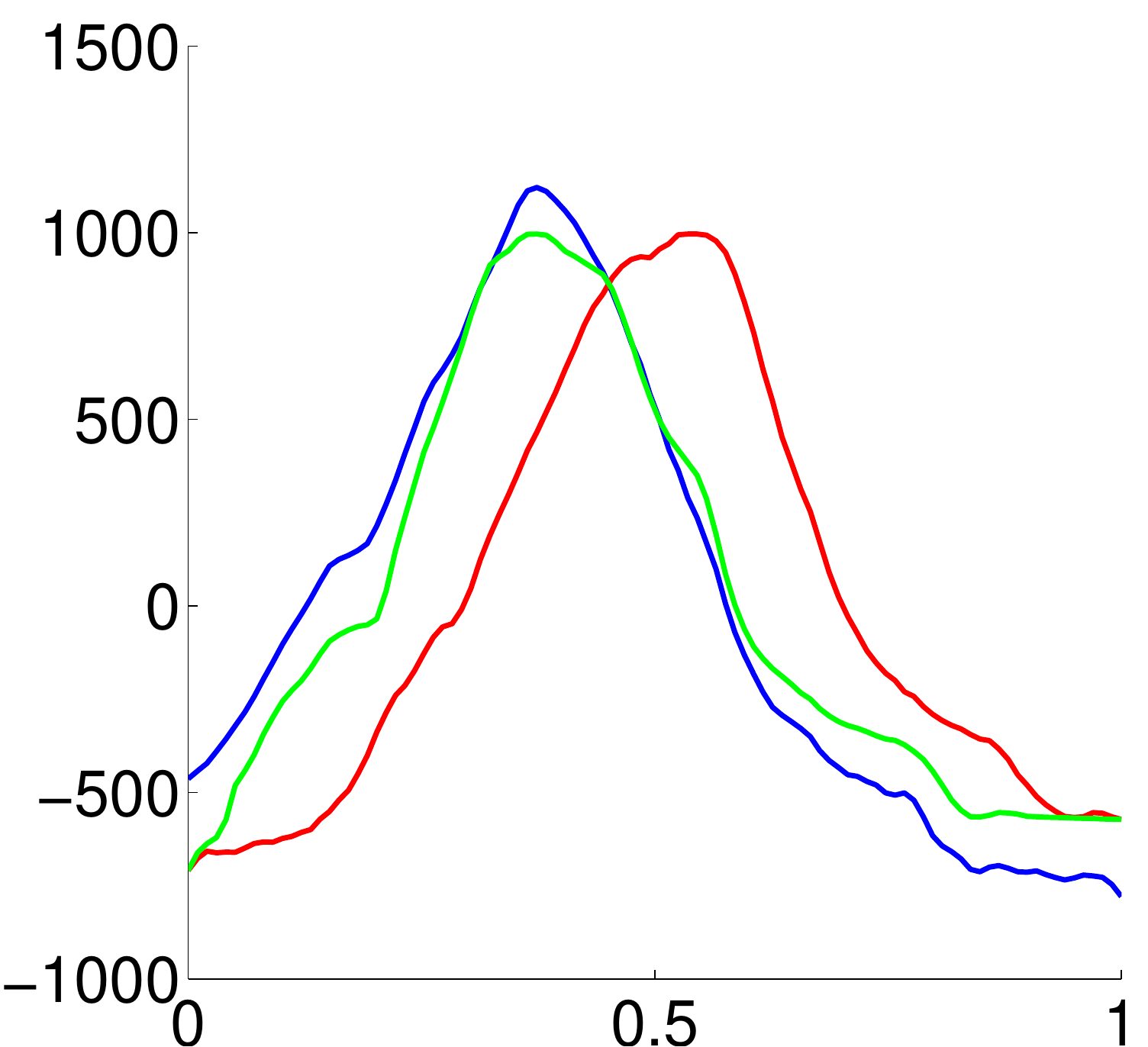}&\includegraphics[width=1in]{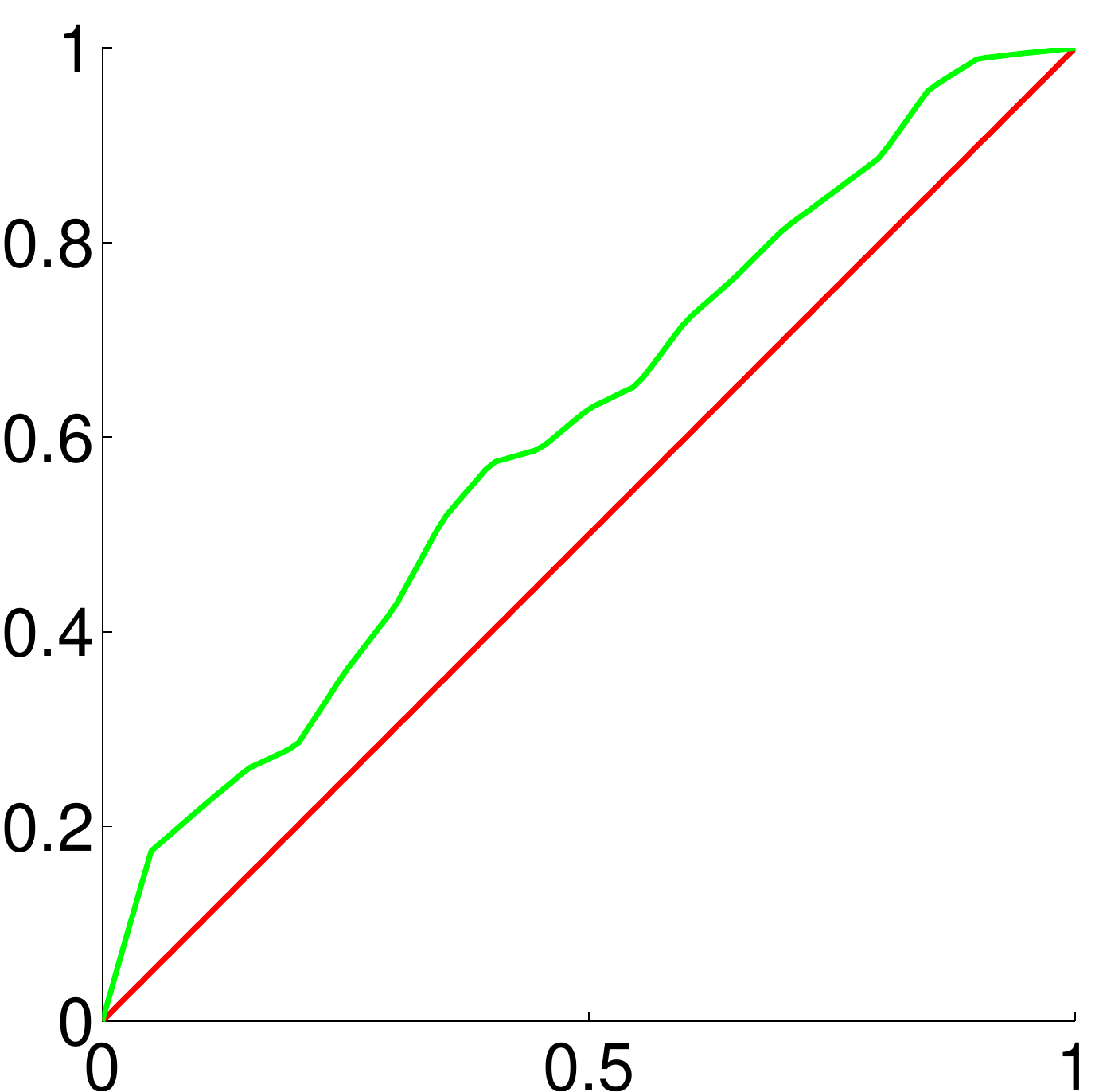}&\includegraphics[width=1in]{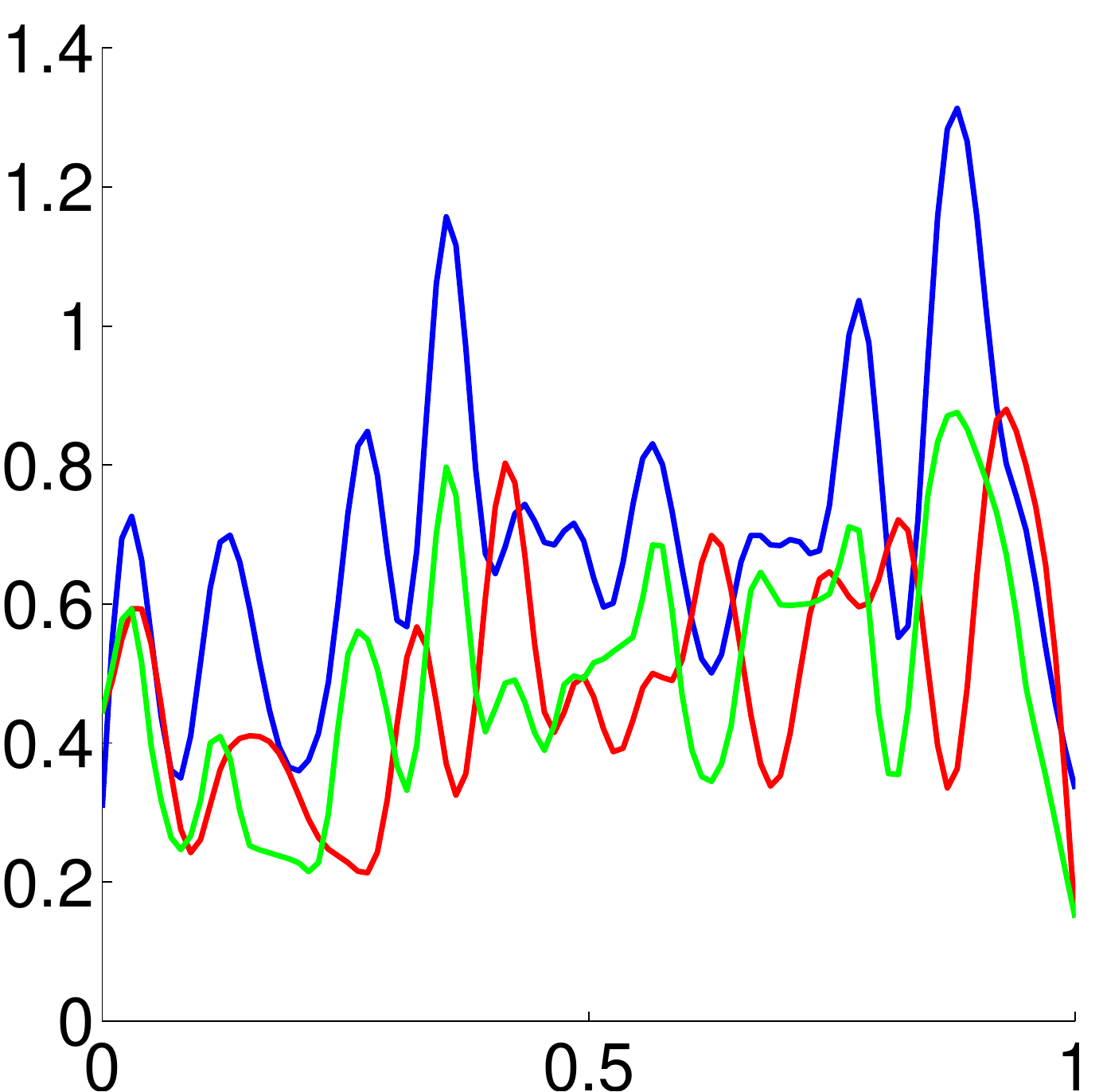}&\includegraphics[width=1in]{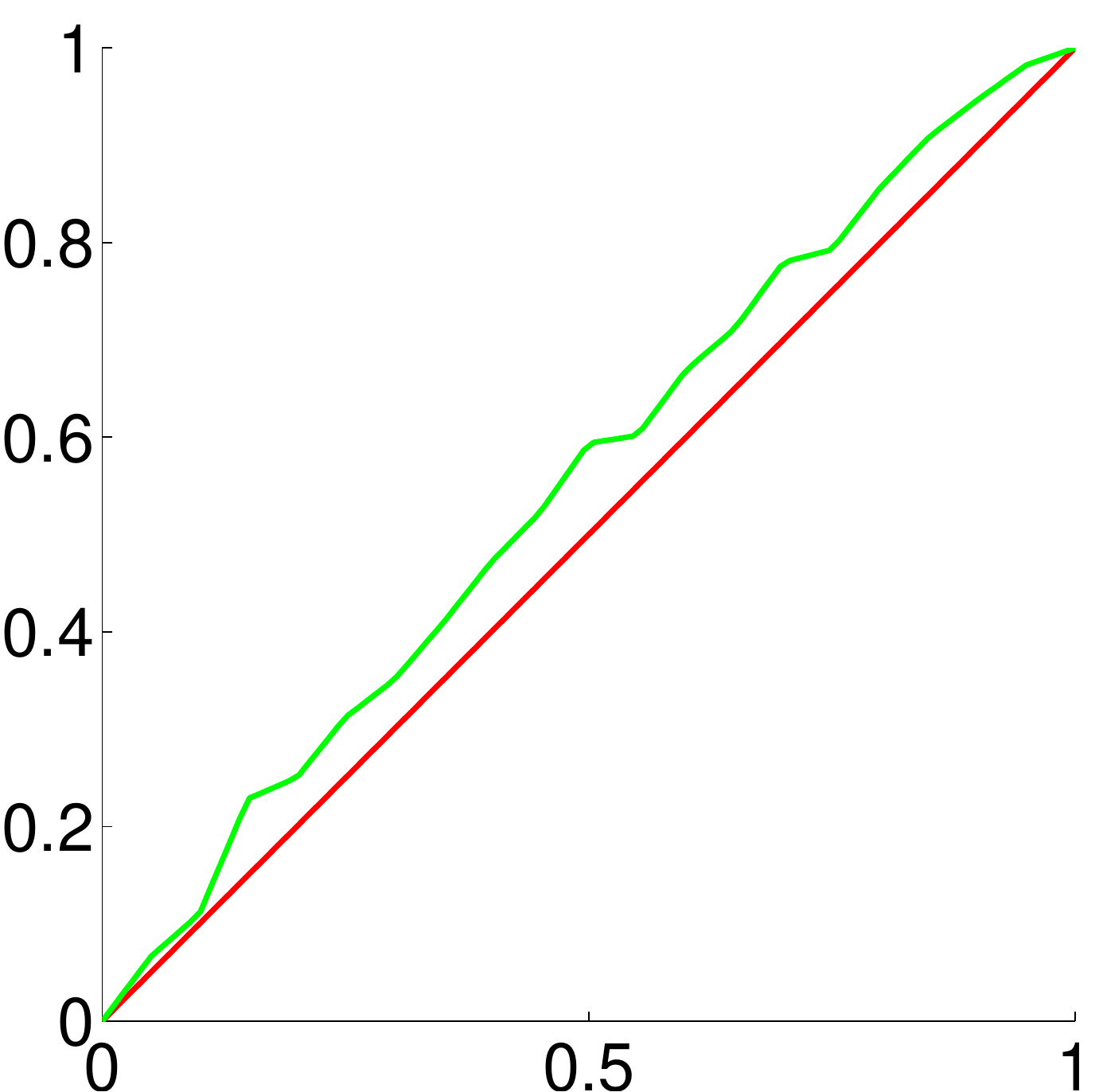}\\
				\hline
			\end{tabular}
			\caption{\small Simulated Annealing-based alignment on six examples: (1) simulated data, (2) PQRST complexes, (3) growth rate functions, (4) gait pressure cycles, (5) respiration cycles, (6) tangential acceleration of signature curves. (a) Two functions before alignment (blue and red), and red function after alignment (green). (b) Optimal (green) and identity (red) warp maps.}\label{fig:SAF}
		\end{center}
	\end{figure}
}

Next, we briefly comment on how Algorithm \ref{alg4} is extended to the case of alignment of open and closed curves for the purpose of shape analysis. In the case of shapes, translation, scale and rotation variations are nuisances and have to be removed in addition to alignment via warp maps. Translation is removed automatically through the SRVF representation; scale variation is removed by normalizing all curves to unit length. Rotation variability is accounted for by adding a Procrustes step at each iteration of the algorithm (see \cite{MD} for details on Procrustes alignment). In the case of closed curves, we must additionally propose the seed point, which is used to unwrap $\sone$ to $[0,1]$ and eventually sample from $\mathbb{D}_\theta\circ H$; this is accomplished via a random proposal, in the neighborhood of the current seed point, from the von-Mises distribution on $\sone$.

{\small
	\begin{table}[!t]
		\begin{center}
			\begin{tabular}{|c|ccc||c|ccc|}
				\hline
				&(a)&(b)&(c)&&(a)&(b)&(c)\\
				\hline
				(1)&2.74&1.43&1.37&(2)&92.06&17.58&21.77\\
				(3)&3.07&0.71&0.82&(4)&1.36&0.53&0.52\\
				(5)&57.57&11.09&15.43&(6)&3.73&0.76&1.26\\
				\hline
			\end{tabular}
		\end{center}
\vspace{-.5cm}
		\caption{\small Distances (a) before alignment, (b) after DP alignment, and (c) after Simulated Annealing alignment, for the six examples in Figure \ref{fig:SAF}.}\label{tab:SAF}
	\end{table}
}

We begin with six examples that consider univariate functional data. The alignment results are presented in Figure \ref{fig:SAF}. In all cases, the computed warp map provides a nice alignment of features across functions. For example, in panel (2), we consider alignment of two PQRST complexes without imposing landmark constraints on the warp maps. Originally, the PQRST peaks and valleys are not well aligned; this is especially evident in the case of the R peaks (highest peak in each function). The proposed method is able to align the peaks very well. This is also the case in the more complex example (6) that considers alignment of two signature tangential acceleration functions. These functions contain many peaks and valleys that are not in correspondence before alignment. The proposed method is able to effectively align all of the peaks and valleys via a suitable warp map. Table \ref{tab:SAF} provides a numerical evaluation of our approach. For each of the six examples, we compare three different distances between the functions: (a) distance before alignment, (b) distance after alignment using DP, and (c) distance after alignment using Simulated Annealing. The proposed method provides comparable performance to DP.

Next, we present several results of registering shapes of 3D open curves. In this case, we use two datasets that were previously considered in \cite{KurtekJASA}: (1) simulated spirals, and (2) fibers extracted from diffusion tensor magnetic resonance images (DT-MRIs). The results are presented in Figure \ref{fig:SAO}. For each example, we show the optimal warp map, the evolution of the energy $E(\gamma)$ as a function of the number of iterations, and the geodesic path (shortest distance deformation under the $\mathbb{L}^2$ metric on SRVF representations) between the two shapes before and after Simulated Annealing-based alignment. For the simulated spirals, the additional alignment via warp maps results in a much more natural geodesic deformation between them, where the shapes of the individual spirals are better preserved. This is also the case for the DT-MRI fibers, albeit not as clear. The top portion of Table \ref{tab:SAOC} provides a quantitative comparison of DP-based alignment and the proposed method. In three out of the four given examples, the proposed method results in a significantly shorter distance between the considered shapes (the maximum distance on this shape space is $\pi/2$).
{\small
\begin{figure}[!t]
	\begin{center}
		\begin{tabular}{|cc||cc|}
			\hline
			(a)&(b)&(a)&(b)\\
			\hline
			\multicolumn{2}{|c||}{(1)}&\multicolumn{2}{c|}{(2)}\\
			\hline
			\includegraphics[width=1in]{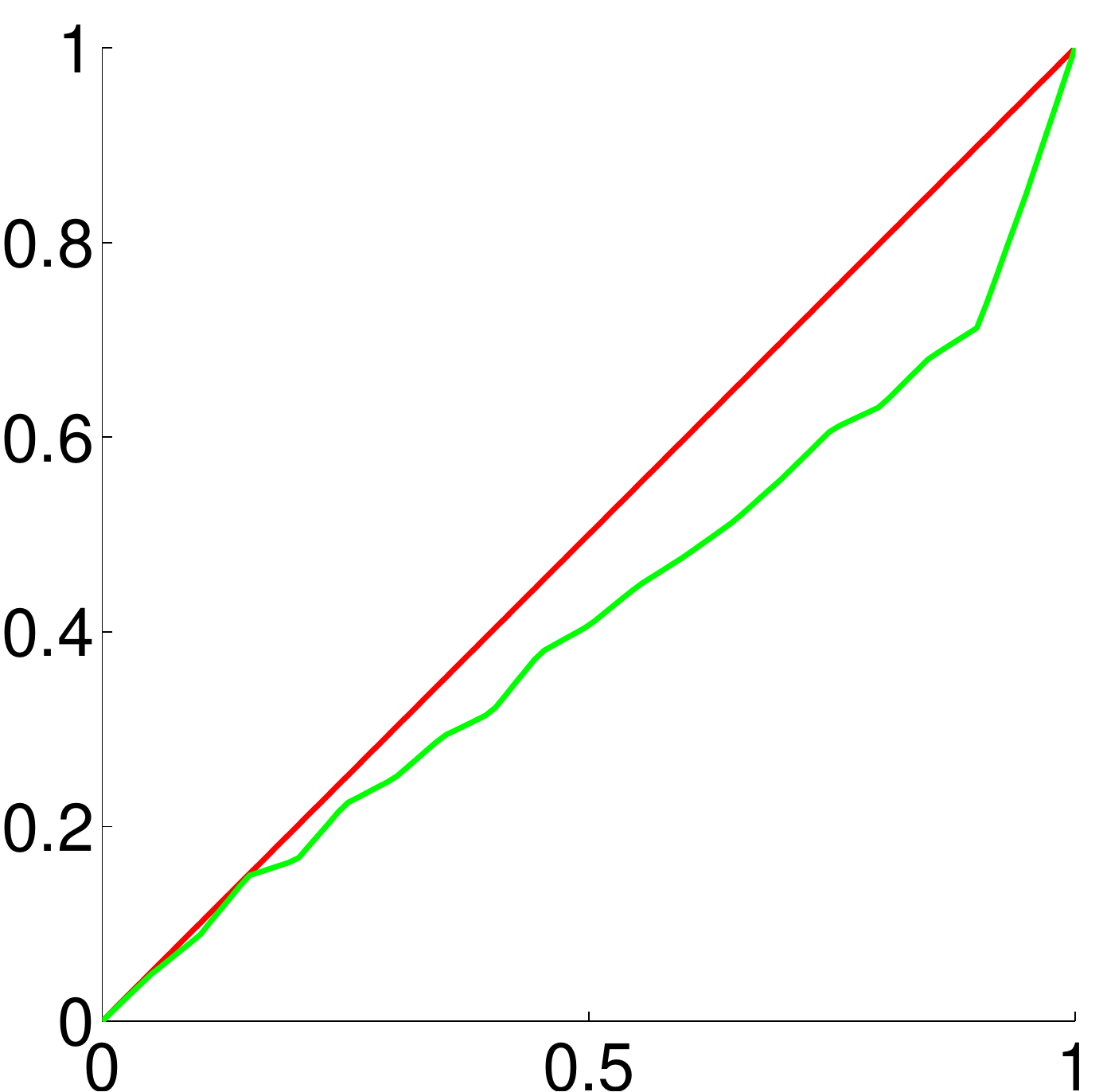}&\includegraphics[height=1in]{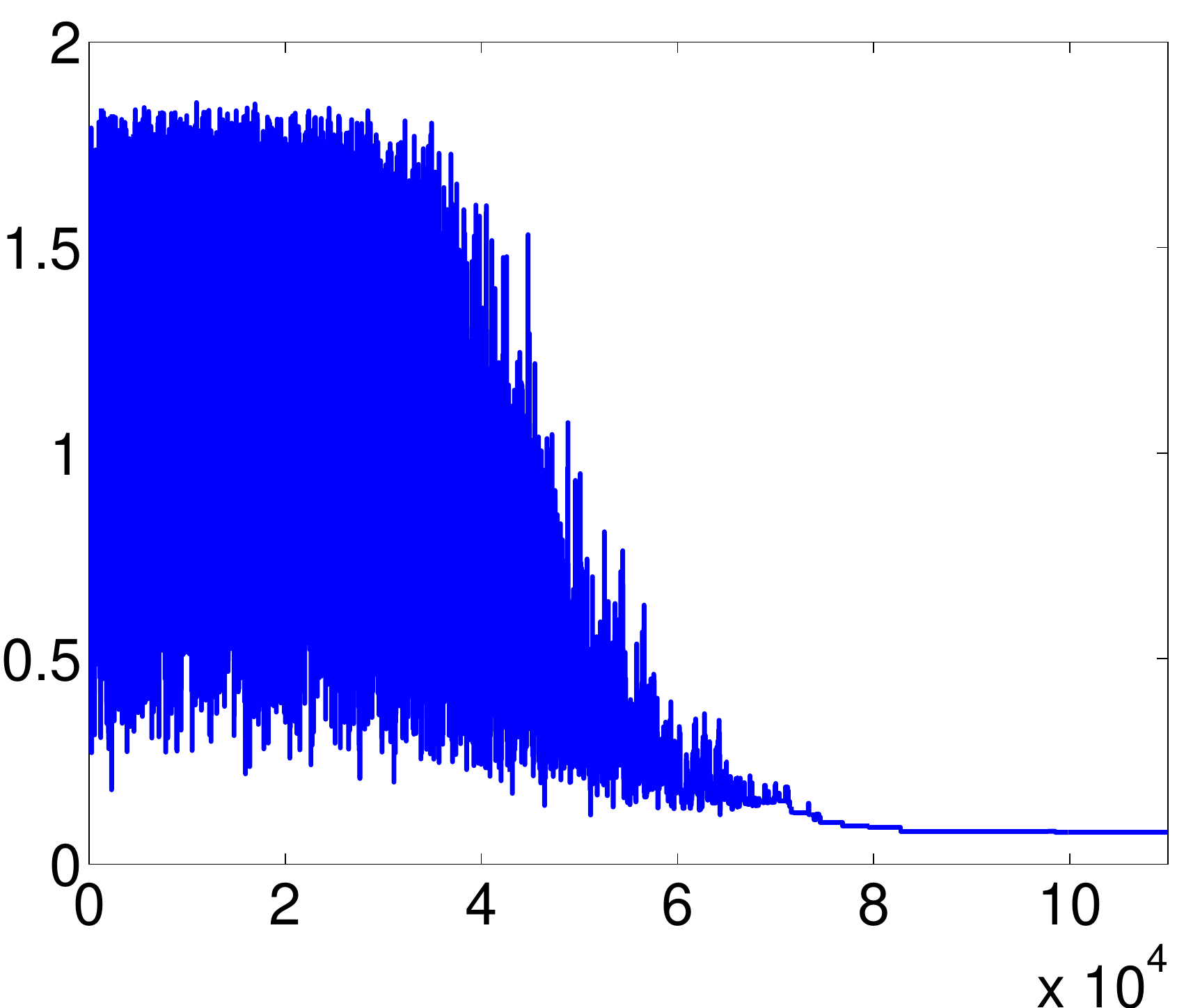}&\includegraphics[width=1in]{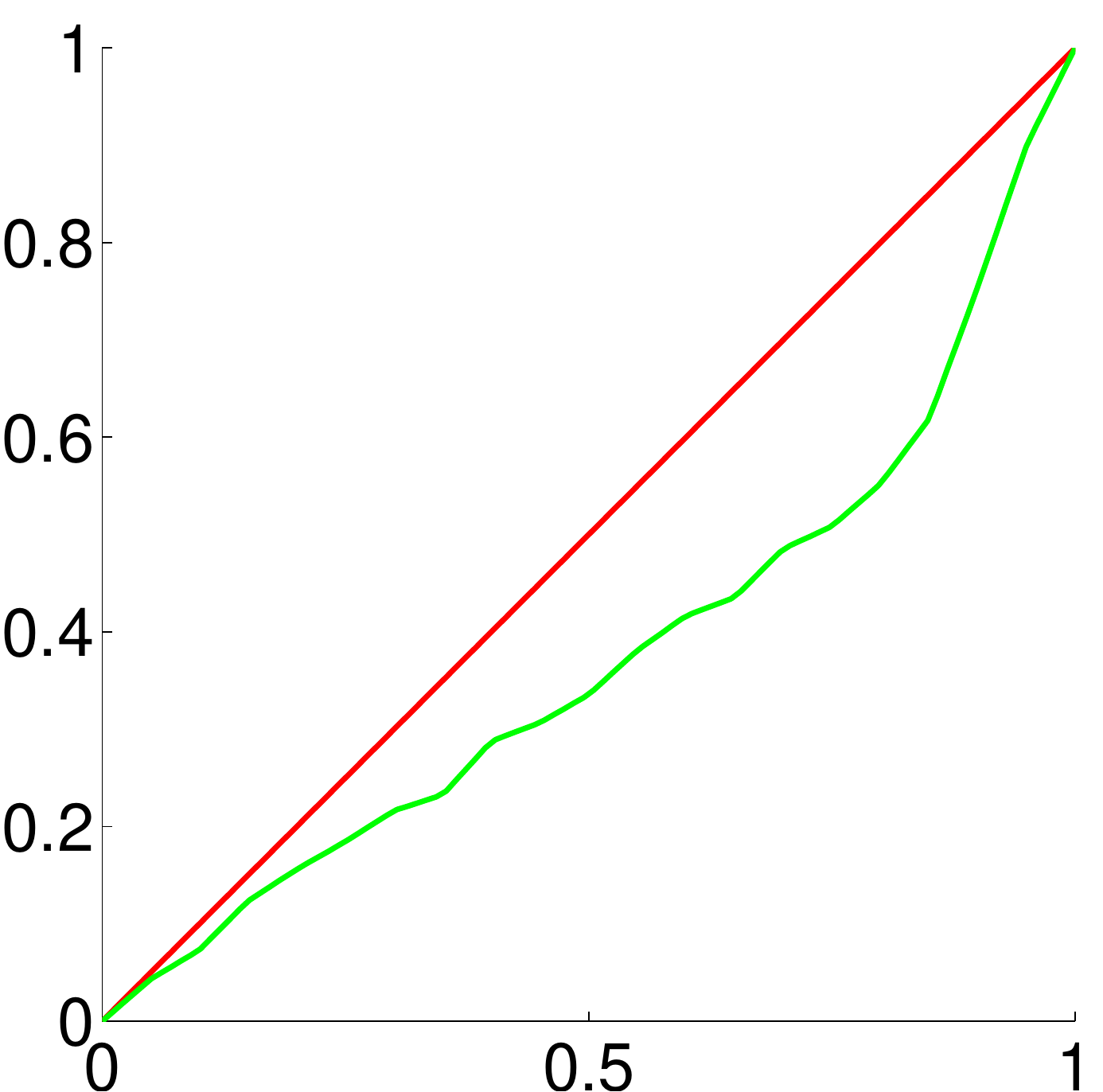}&\includegraphics[height=1in]{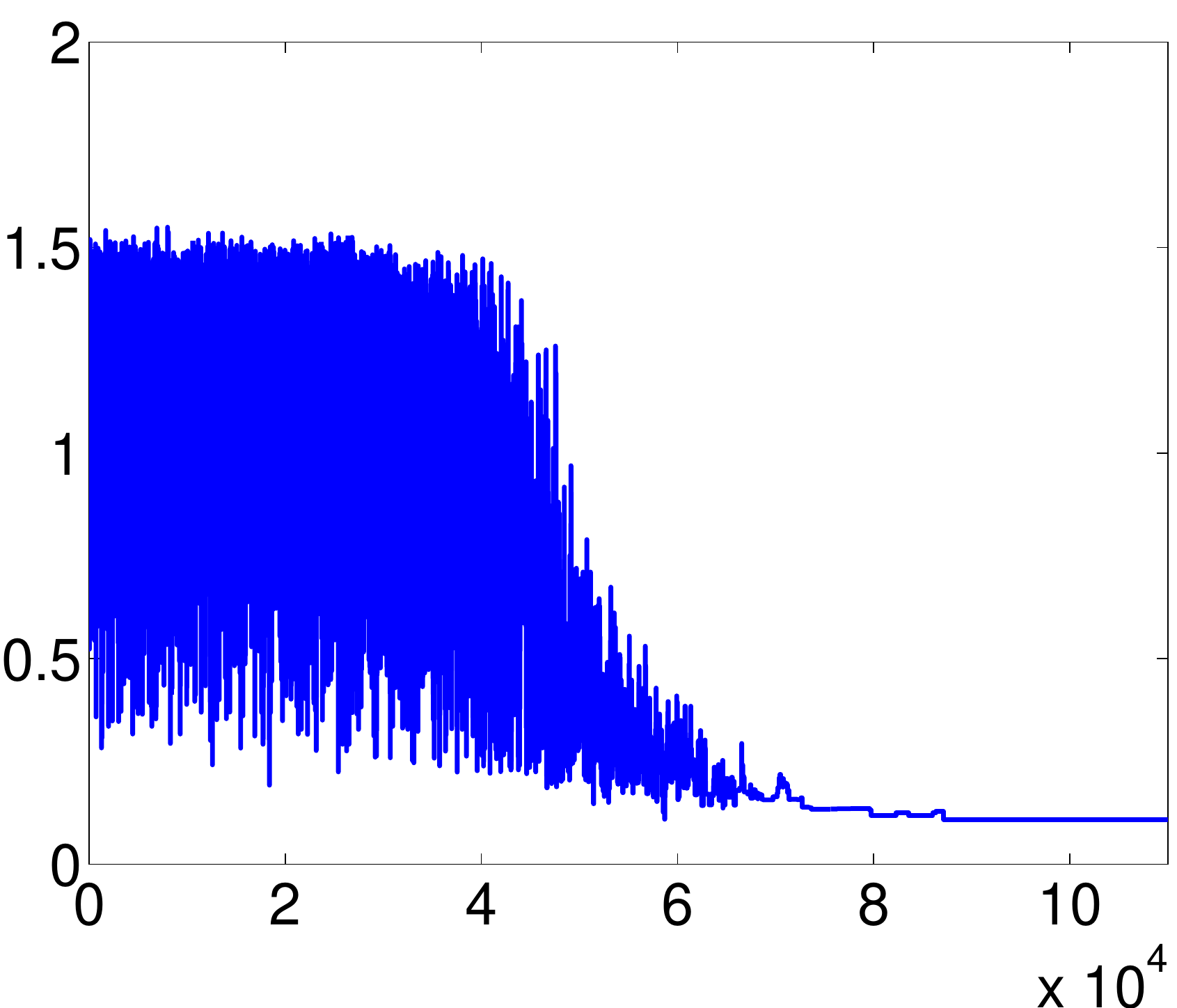}\\
			\multicolumn{4}{|c|}{Geodesic after alignment}\\
			\multicolumn{2}{|c||}{\includegraphics[width=2.1in]{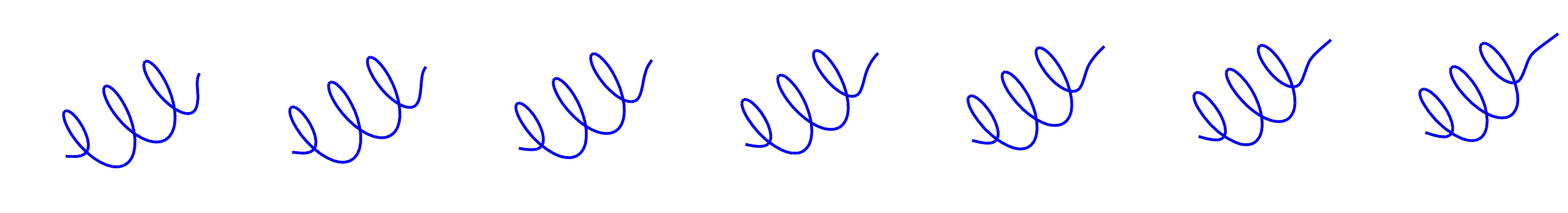}}&\multicolumn{2}{c|}{\includegraphics[width=2.1in]{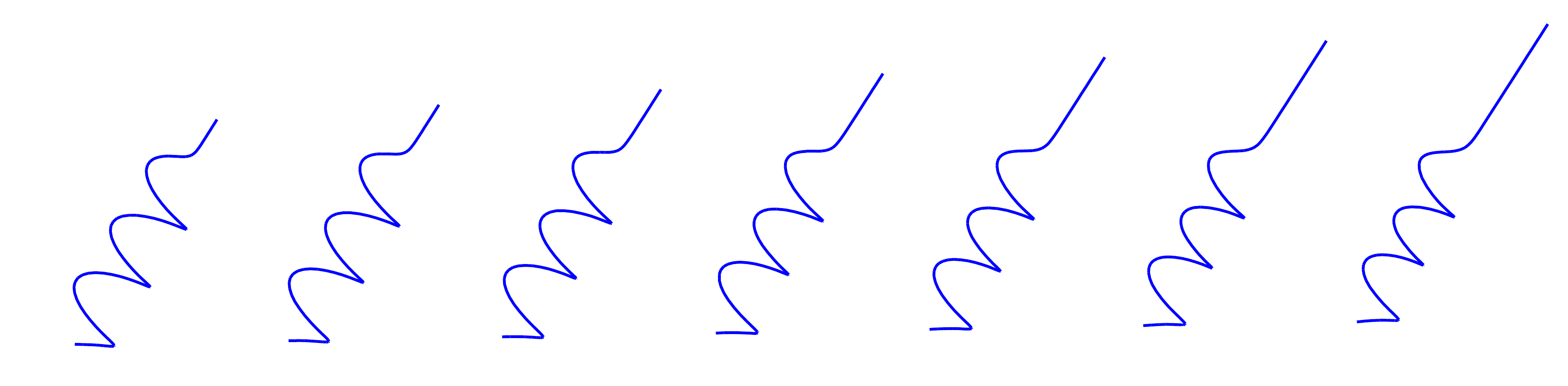}}\\
			\multicolumn{4}{|c|}{Geodesic before alignment}\\
			\multicolumn{2}{|c||}{\includegraphics[width=2.1in]{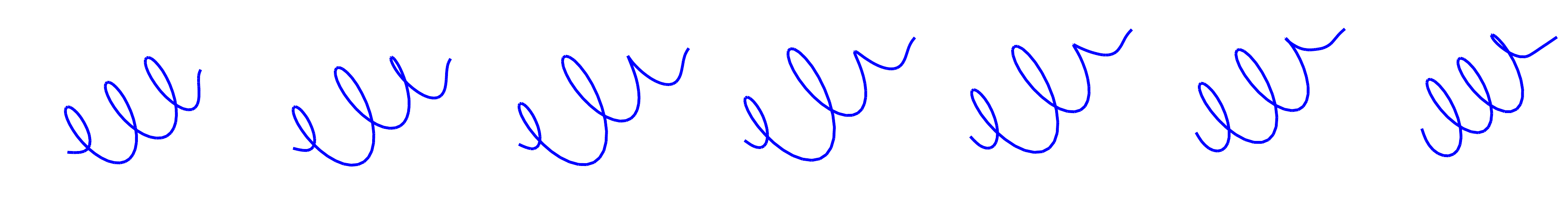}}&\multicolumn{2}{c|}{\includegraphics[width=2.1in]{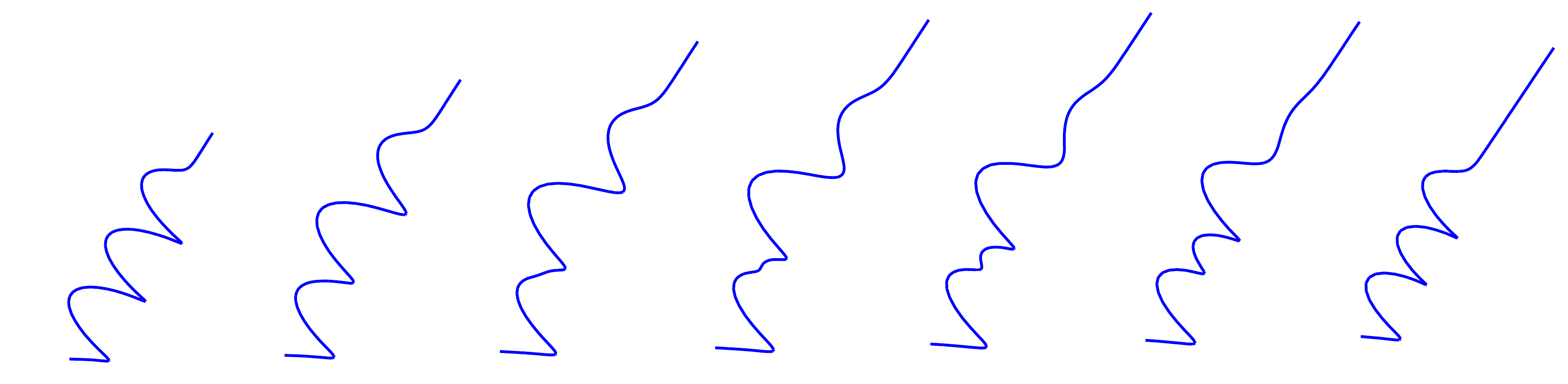}}\\
			\hline
			\multicolumn{2}{|c||}{(3)}&\multicolumn{2}{c|}{(4)}\\
			\hline
			\includegraphics[width=1in]{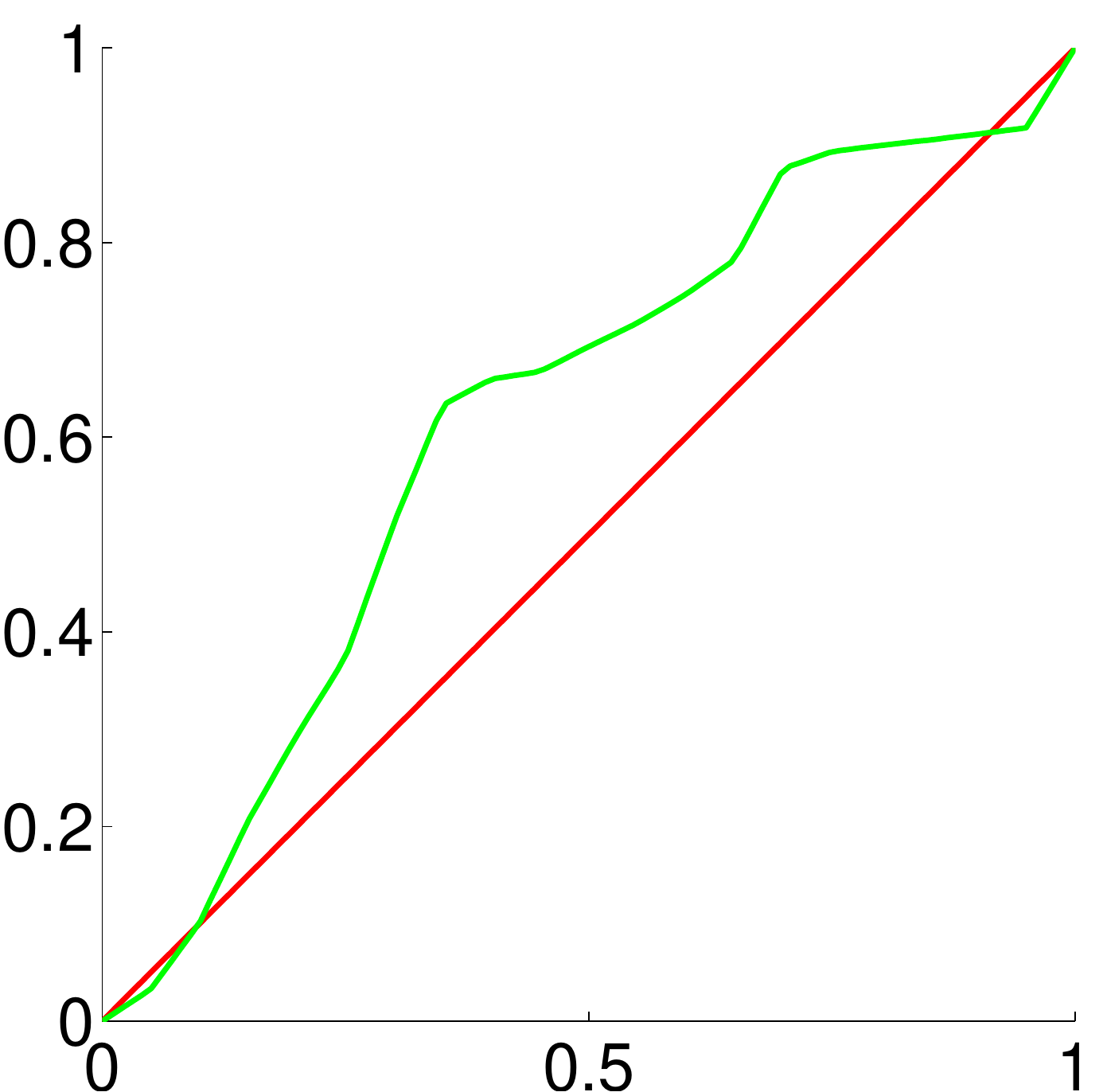}&\includegraphics[height=1in]{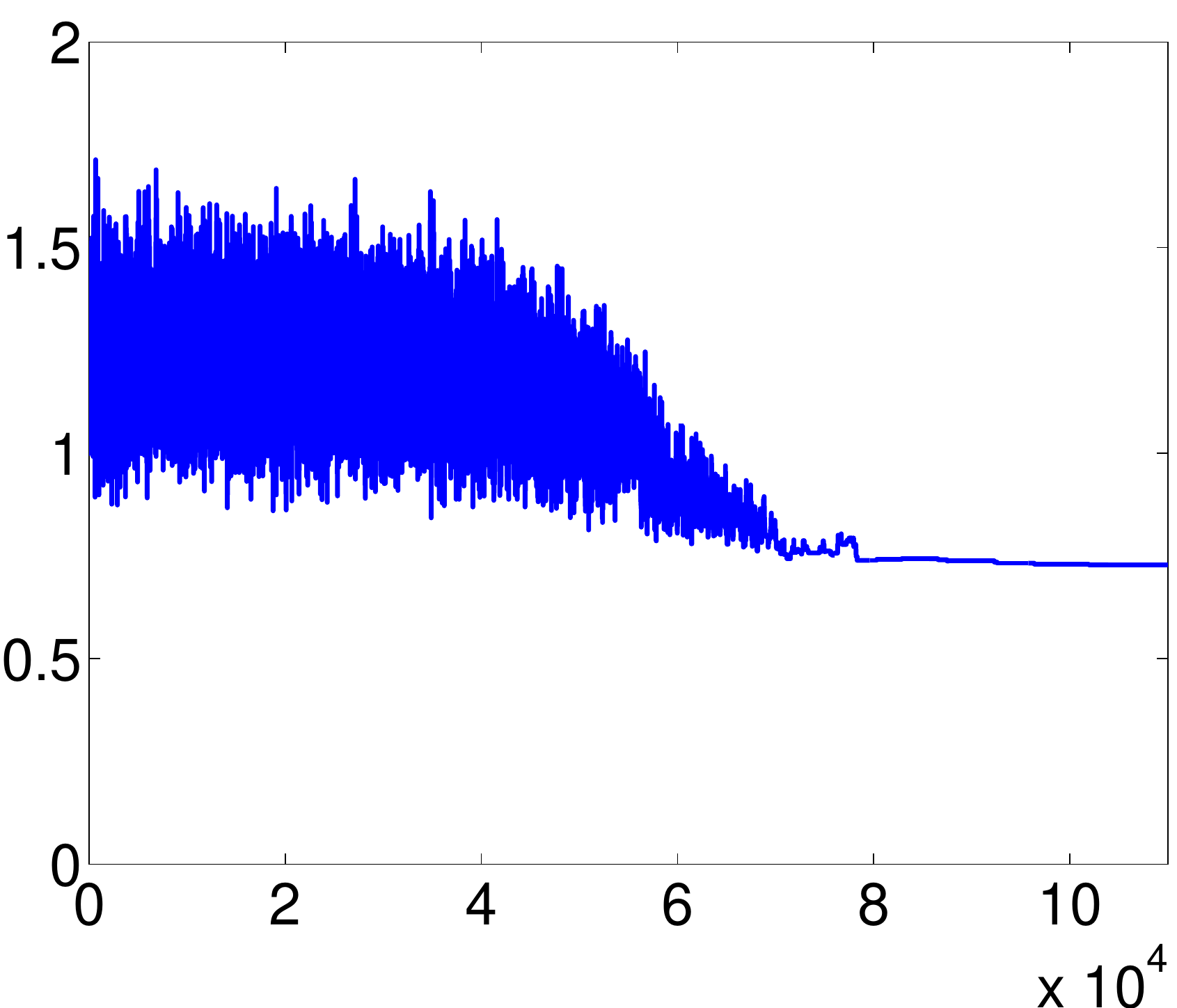}&\includegraphics[width=1in]{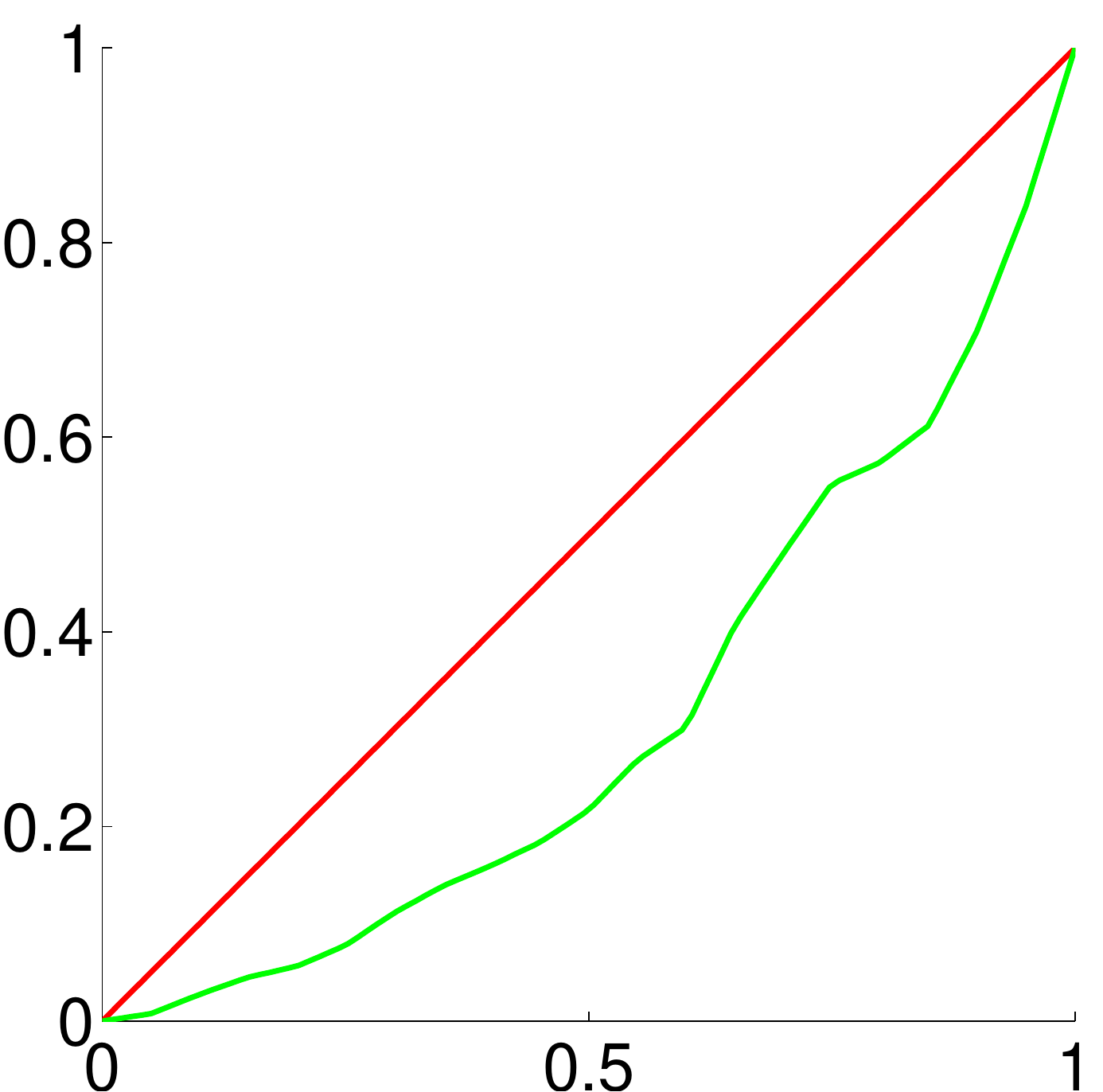}&\includegraphics[height=1in]{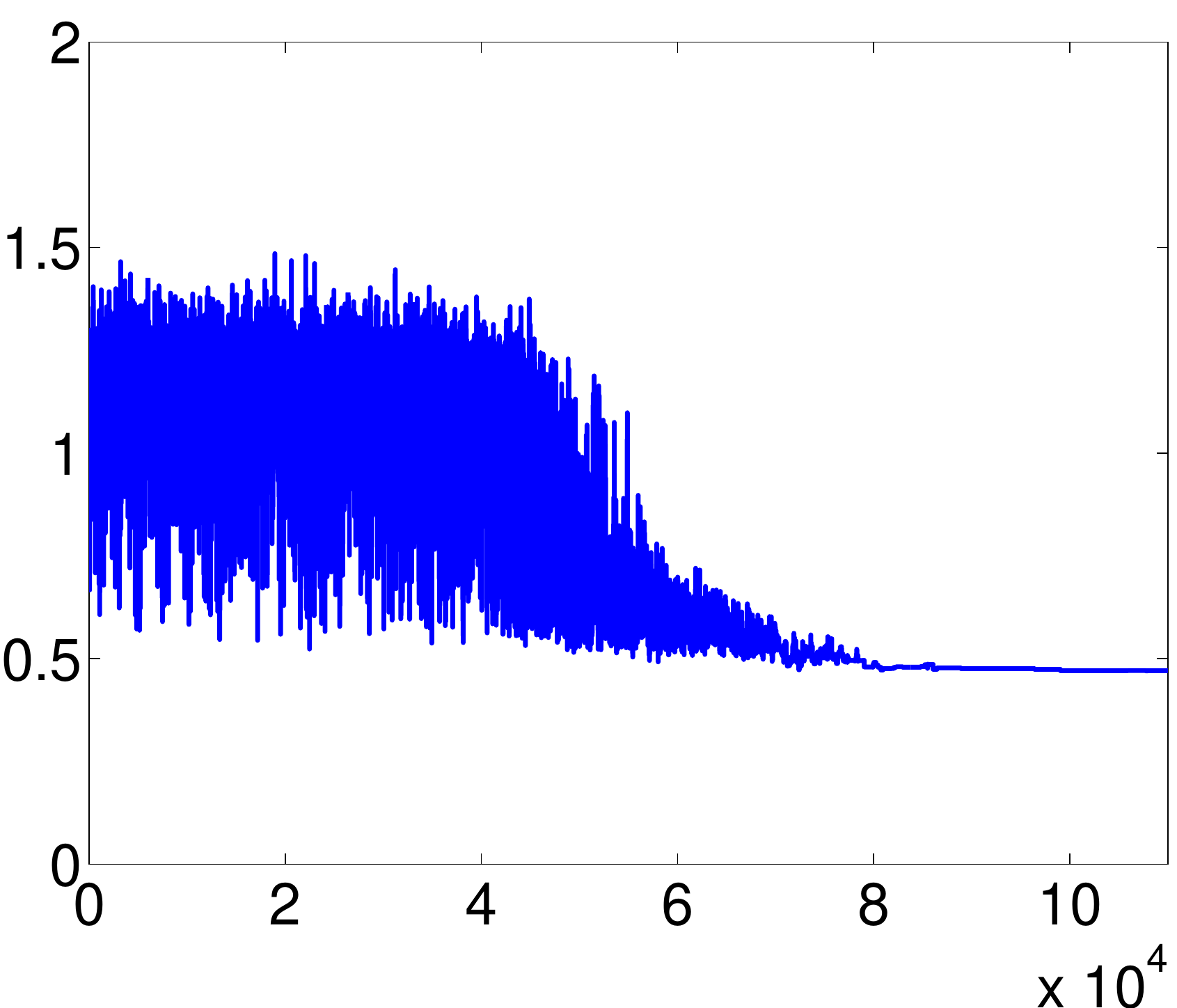}\\
			\multicolumn{4}{|c|}{Geodesic after alignment}\\
			\multicolumn{2}{|c||}{\includegraphics[width=2.1in]{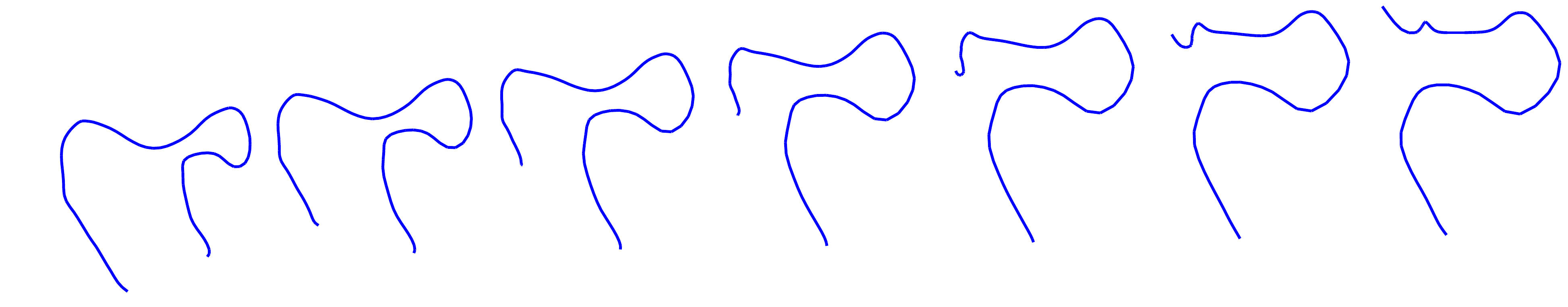}}&\multicolumn{2}{|c|}{\includegraphics[width=2.1in]{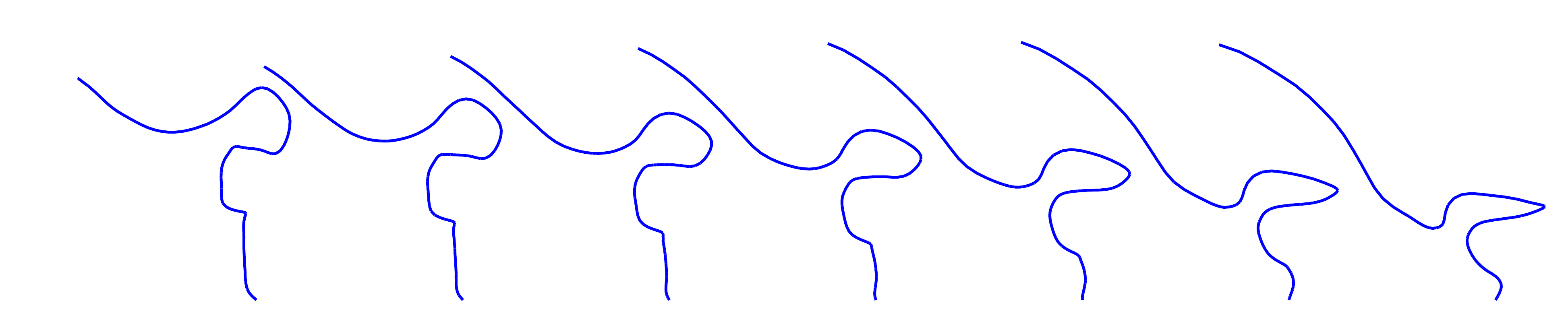}}\\
			\multicolumn{4}{|c|}{Geodesic before alignment}\\
			\multicolumn{2}{|c||}{\includegraphics[width=2.1in]{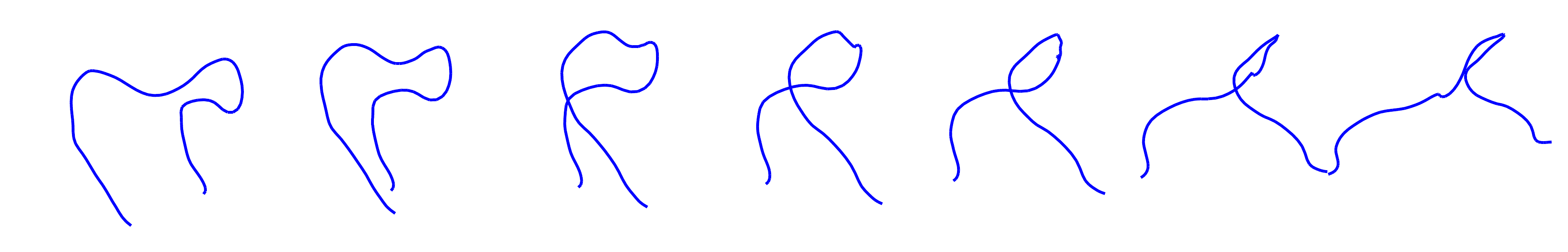}}&\multicolumn{2}{|c|}{\includegraphics[width=2.1in]{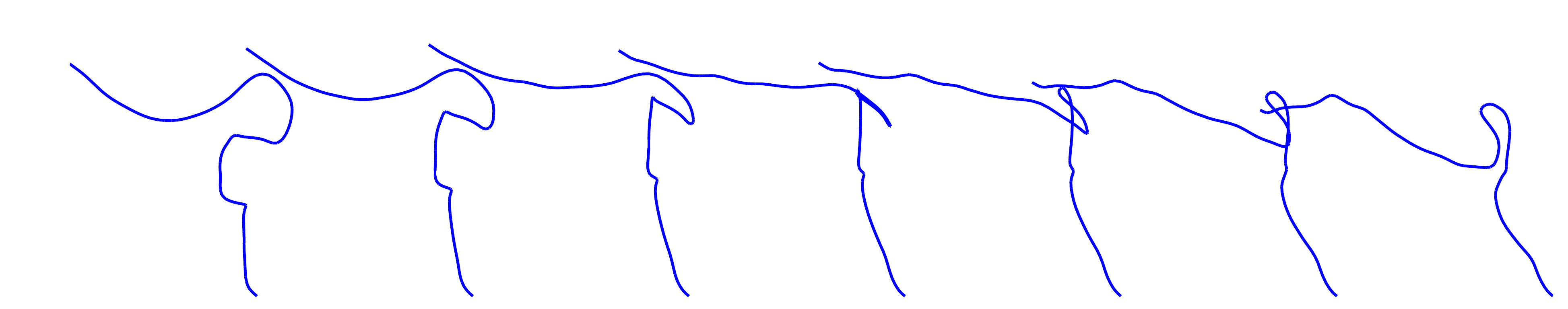}}\\
			\hline
		\end{tabular}
		\caption{\small Results of Simulated Annealing-based alignment on four different examples of open curves. (a) Optimal (red) and identity (green) warp maps. (b) Evolution of $E(\gamma)$.}\label{fig:SAO}
	\end{center}
\end{figure}
}
{\small
\begin{table}[!t]
	\begin{center}
		\begin{tabular}{|c|ccc||c|ccc|}
			\hline
			&(a)&(b)&(c)&&(a)&(b)&(c)\\
			\hline
			\multicolumn{8}{|c|}{Open Curves}\\
			\hline
			(1)&1.20&0.39&0.28&(2)&1.24&0.25&0.33\\
			(3)&1.17&0.97&0.88&(4)&1.11&0.92&0.70\\
			\hline
			\multicolumn{8}{|c|}{Closed Curves}\\
			\hline
			(1)&0.57&0.42&0.26&(2)&1.03&0.68&0.54\\
			(3)&0.63&0.42&0.45&(4)&0.77&0.50&0.52\\
			\hline
		\end{tabular}
		\caption{\small Distances (a) before alignment, (b) after DP alignment, and (c) after Simulated Annealing alignment, for the examples shown in Figures \ref{fig:SAO} and \ref{fig:SAC}.}\label{tab:SAOC}
	\end{center}
\end{table}
}

\begin{figure}[!t]
	\begin{center}
		\begin{tabular}{|cc||cc|}
			\hline
			(a)&(b)&(a)&(b)\\
			\hline
			\multicolumn{2}{|c||}{(1)}&\multicolumn{2}{c|}{(2)}\\
			\hline
			\includegraphics[width=1in]{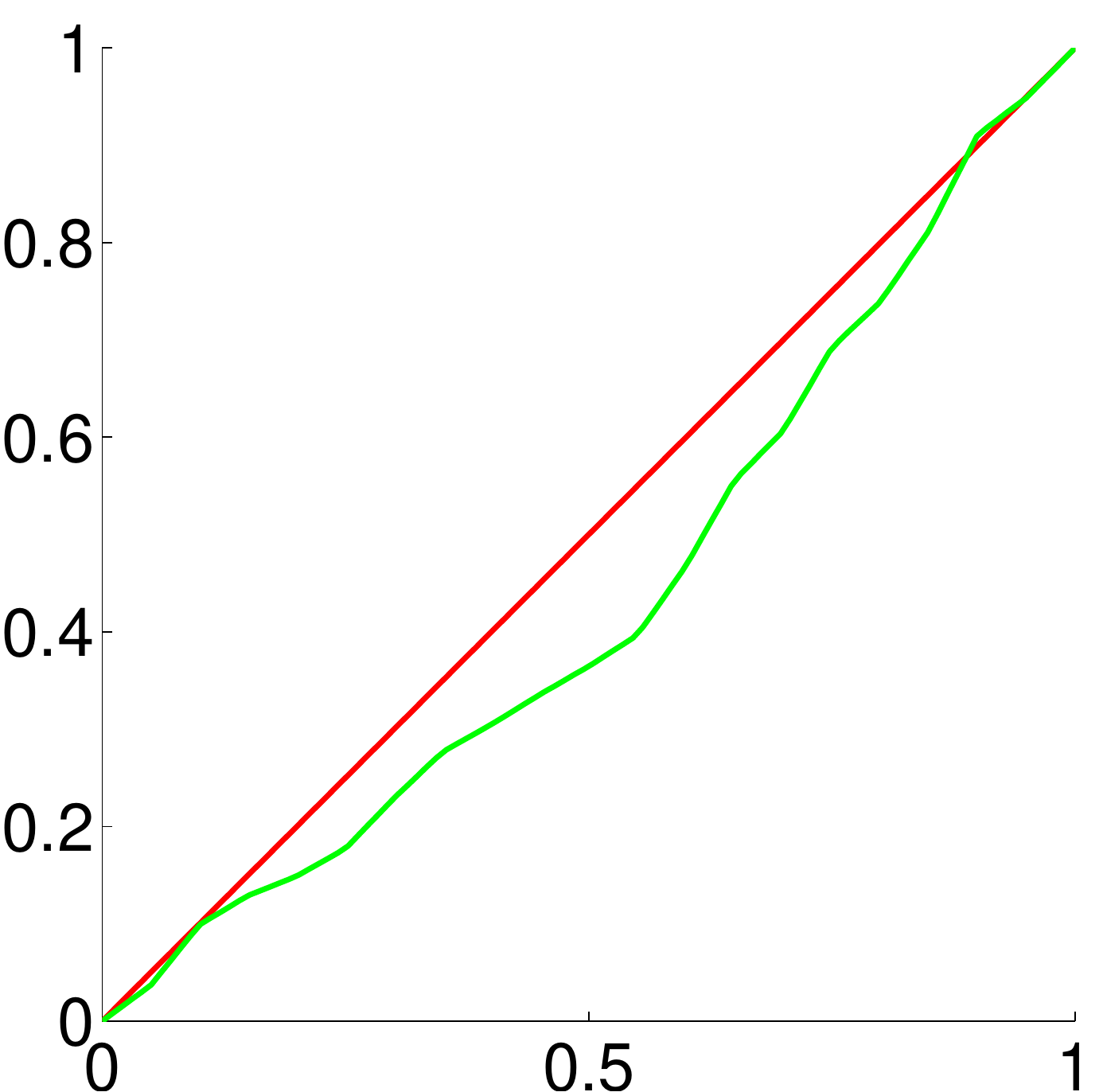}&\includegraphics[height=1in]{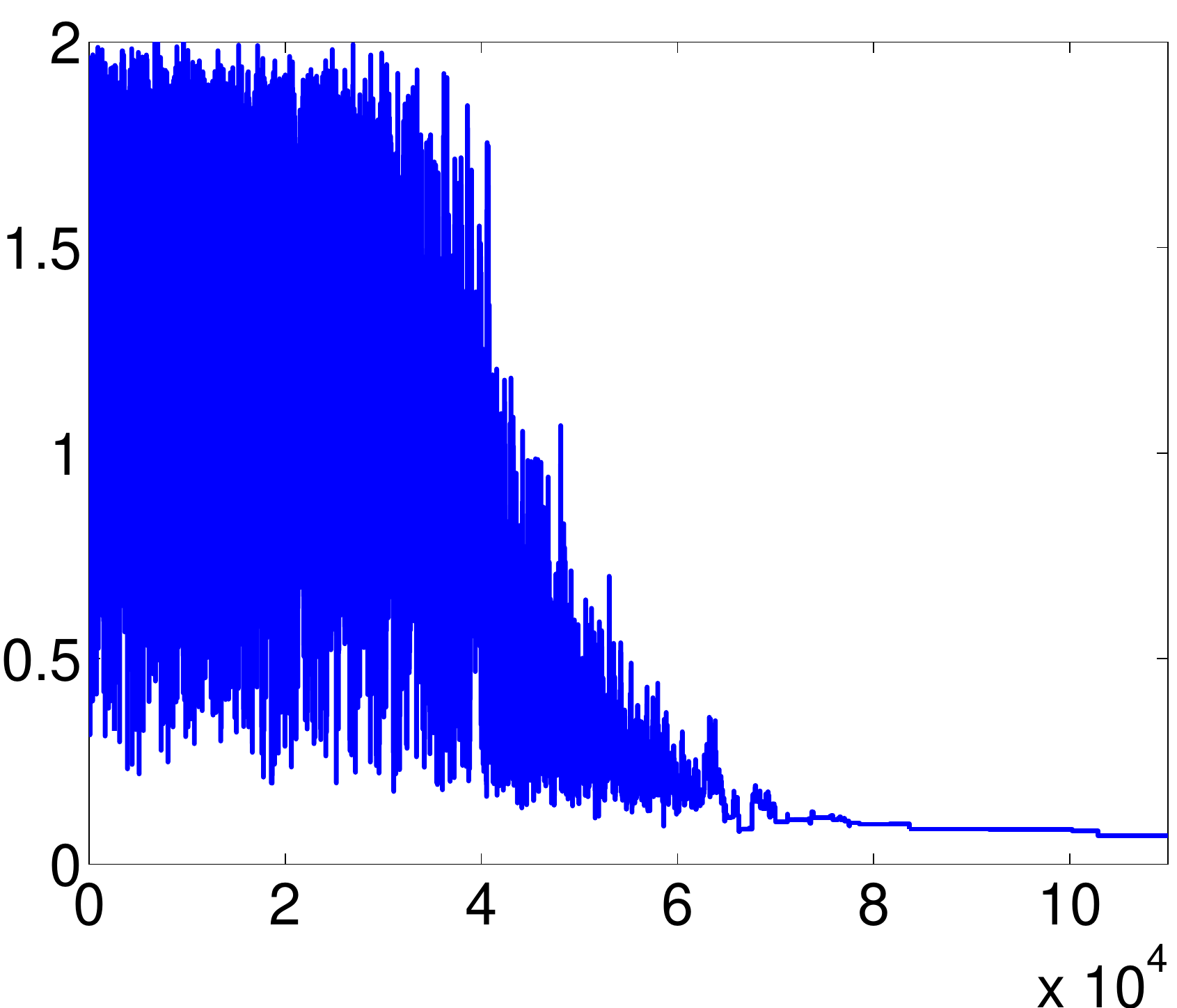}&\includegraphics[width=1in]{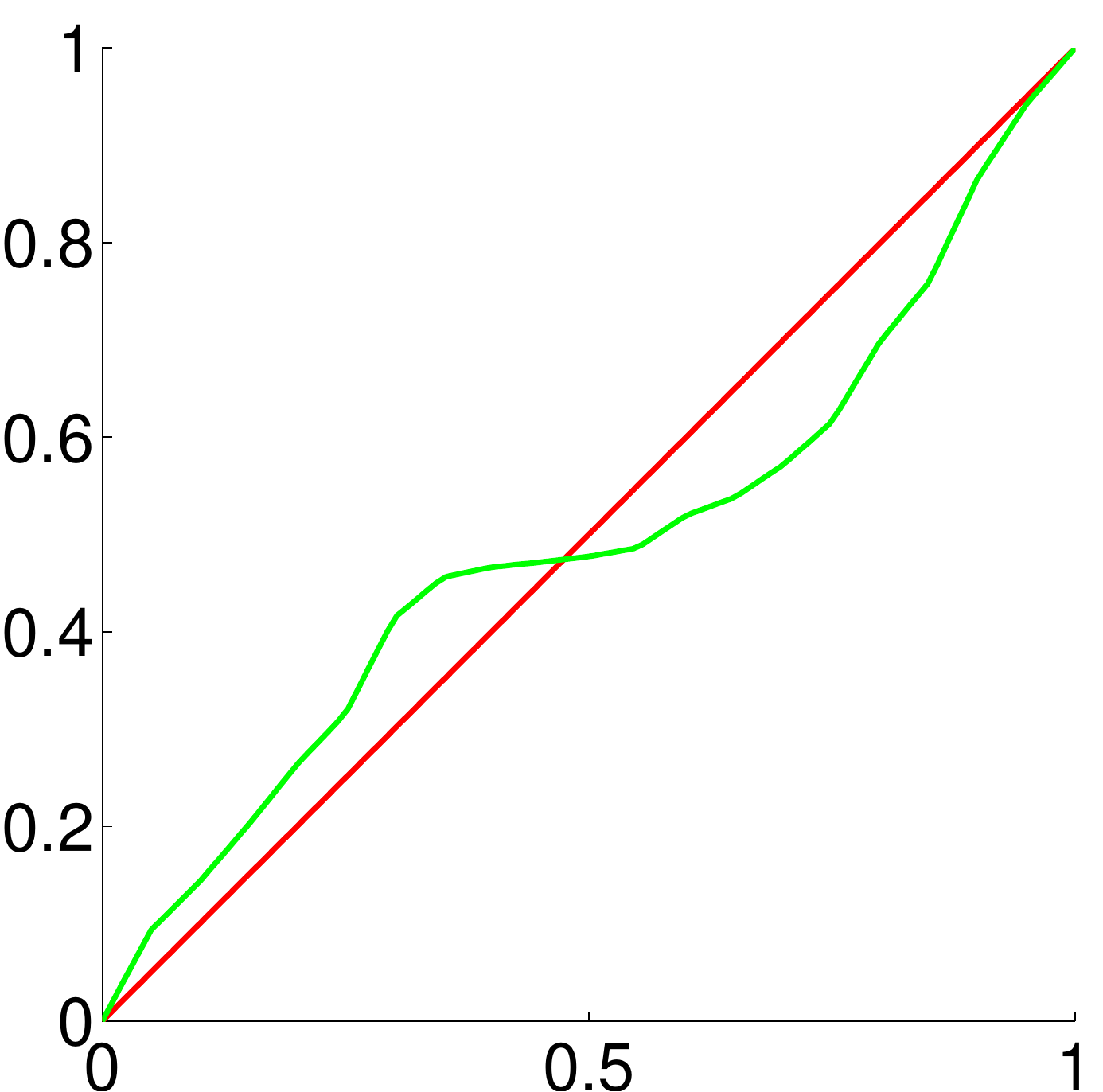}&\includegraphics[height=1in]{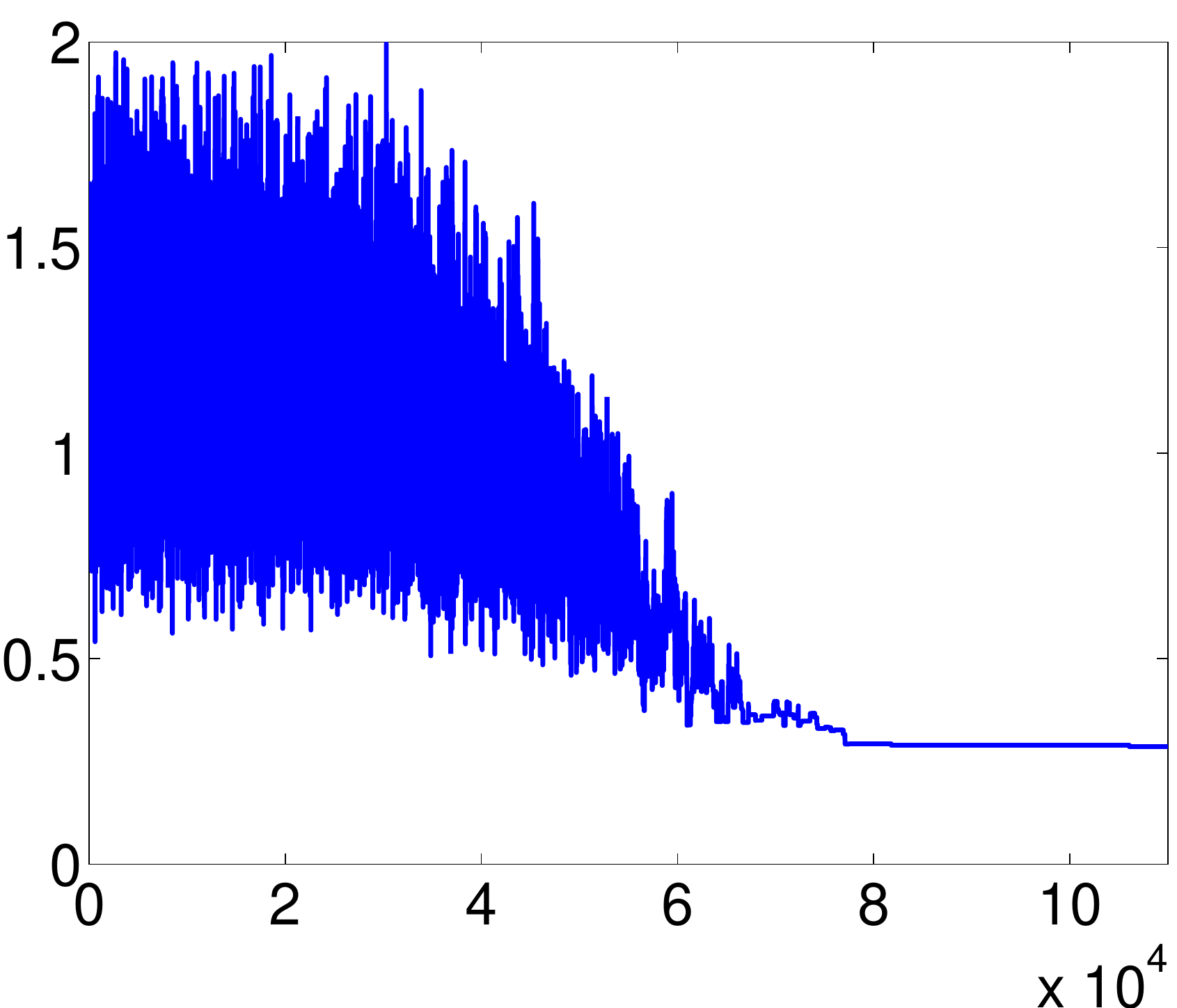}\\
			\multicolumn{4}{|c|}{Geodesic after alignment}\\
			\multicolumn{2}{|c||}{\includegraphics[width=2.1in]{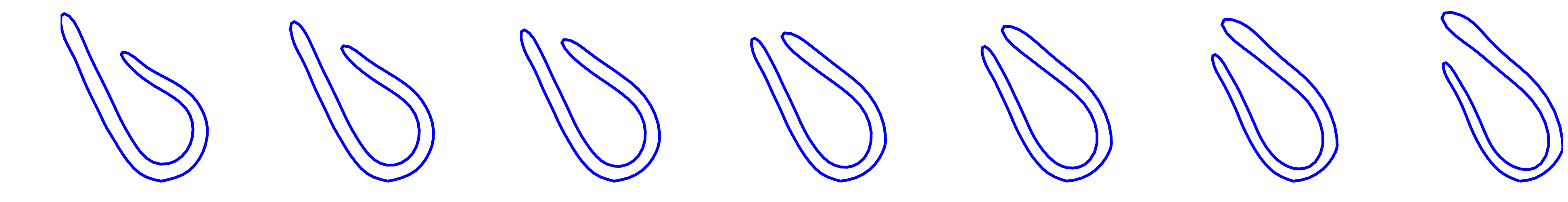}}&\multicolumn{2}{c|}{\includegraphics[width=2.1in]{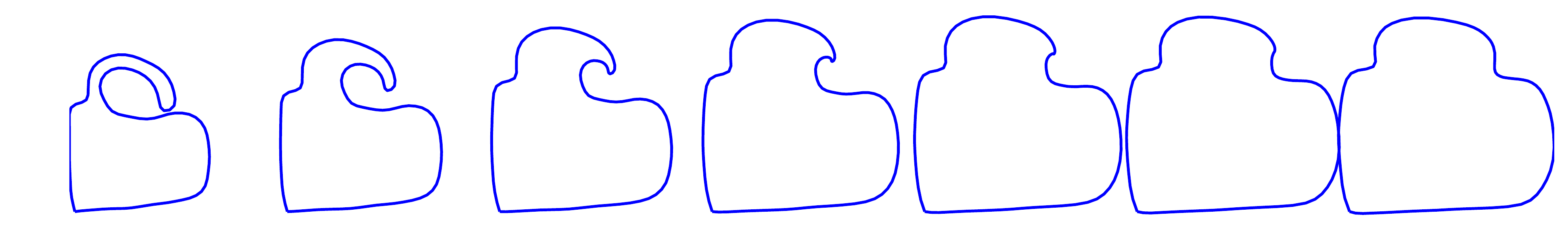}}\\
			\multicolumn{4}{|c|}{Geodesic before alignment}\\
			\multicolumn{2}{|c||}{\includegraphics[width=2.1in]{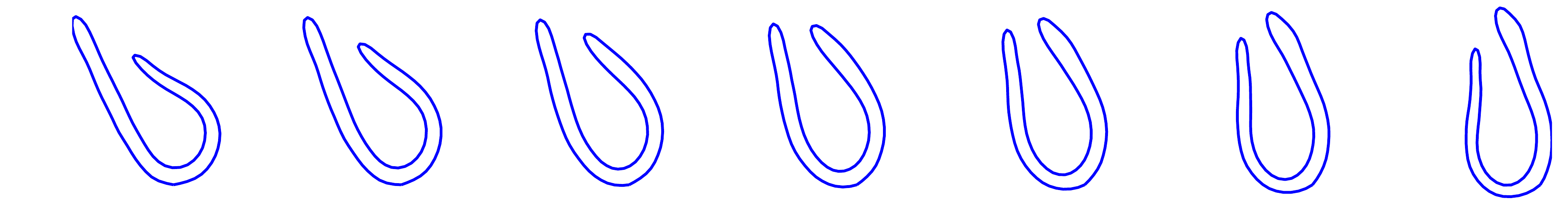}}&\multicolumn{2}{c|}{\includegraphics[width=2.1in]{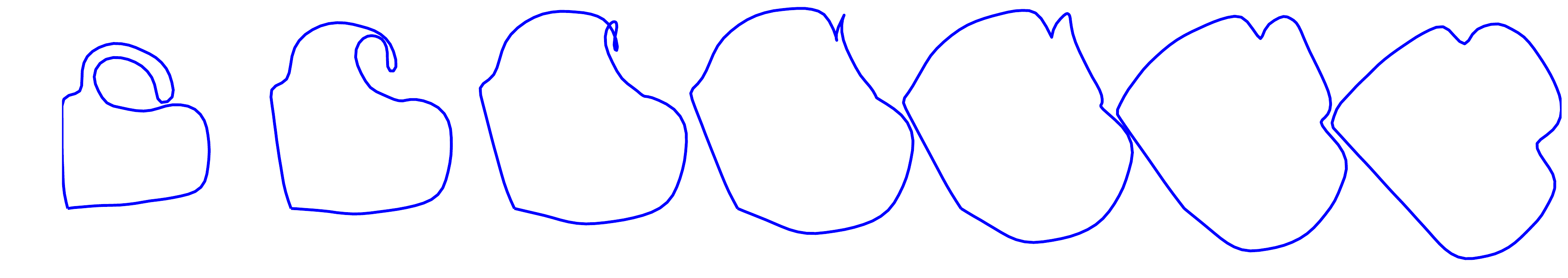}}\\
			\hline
			\multicolumn{2}{|c||}{(3)}&\multicolumn{2}{c|}{(4)}\\
			\hline
			\includegraphics[width=1in]{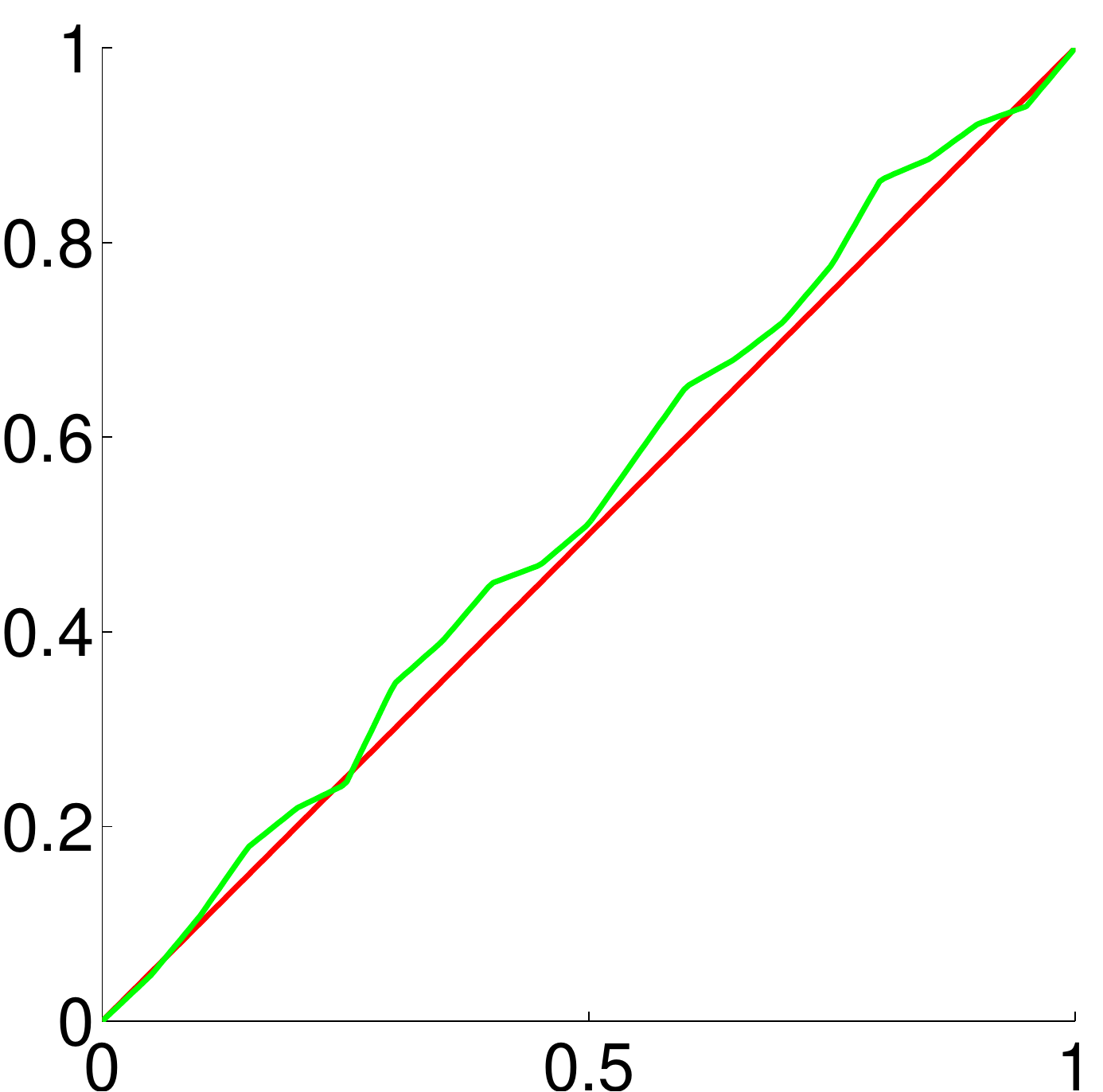}&\includegraphics[height=1in]{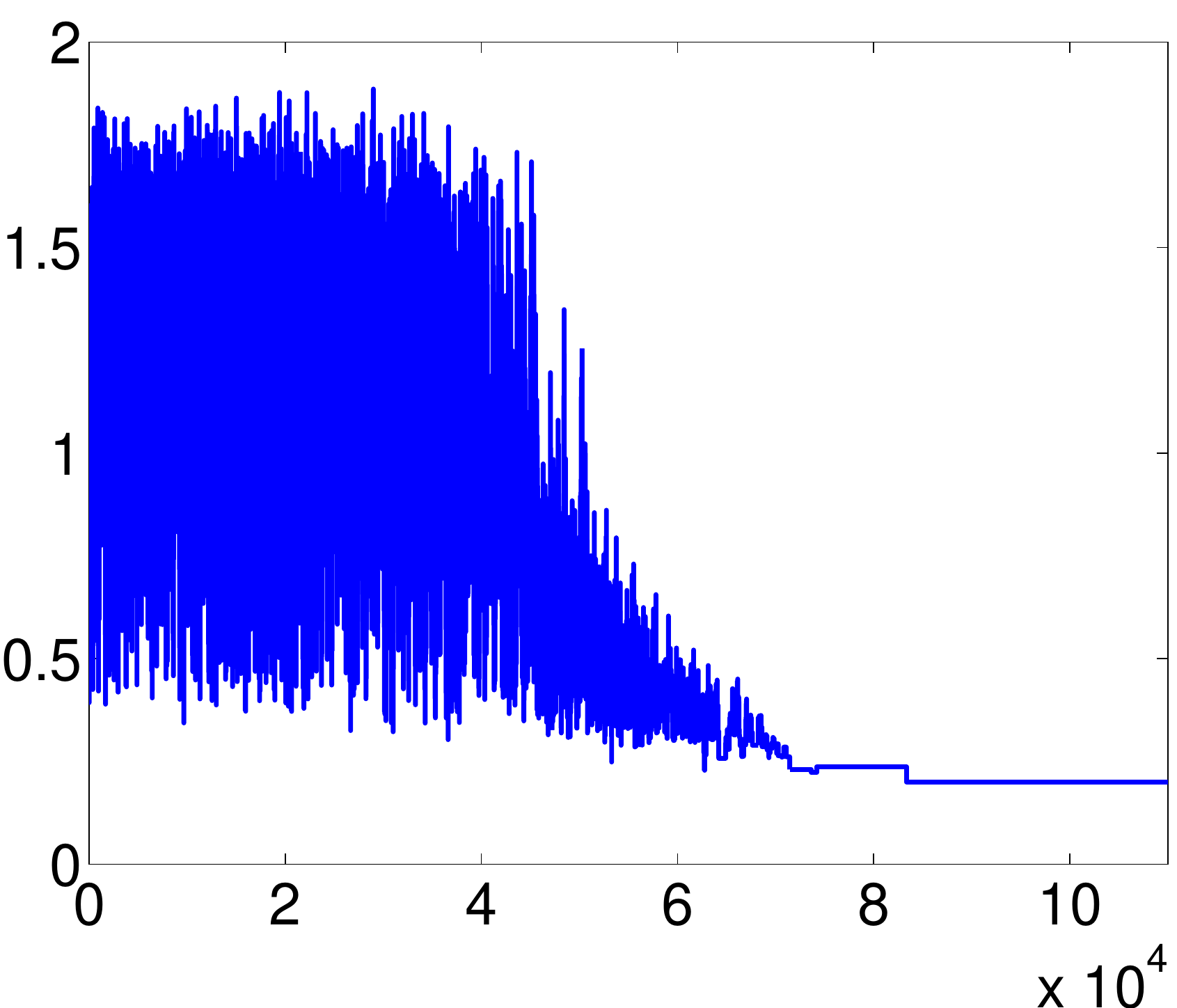}&\includegraphics[width=1in]{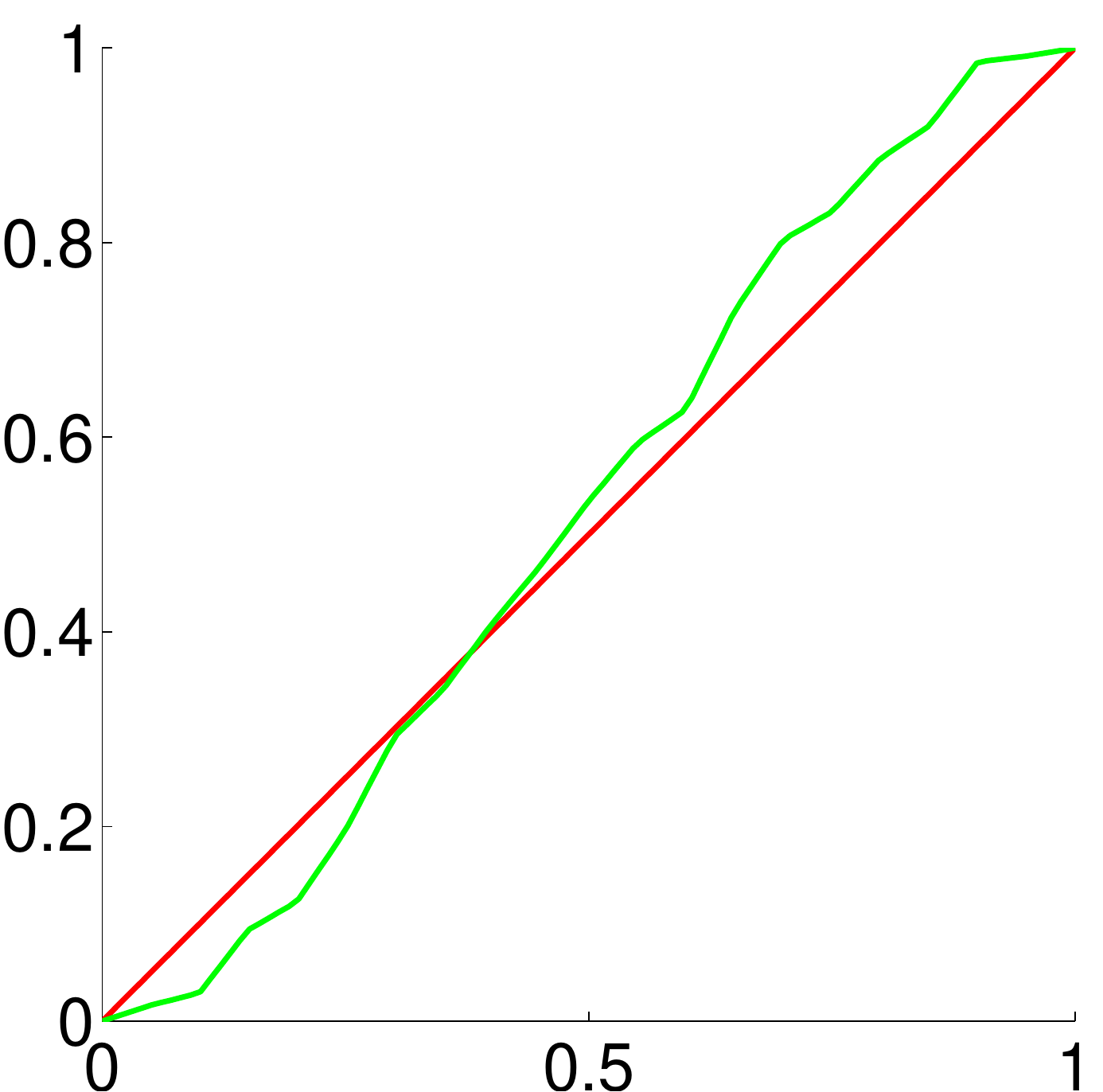}&\includegraphics[height=1in]{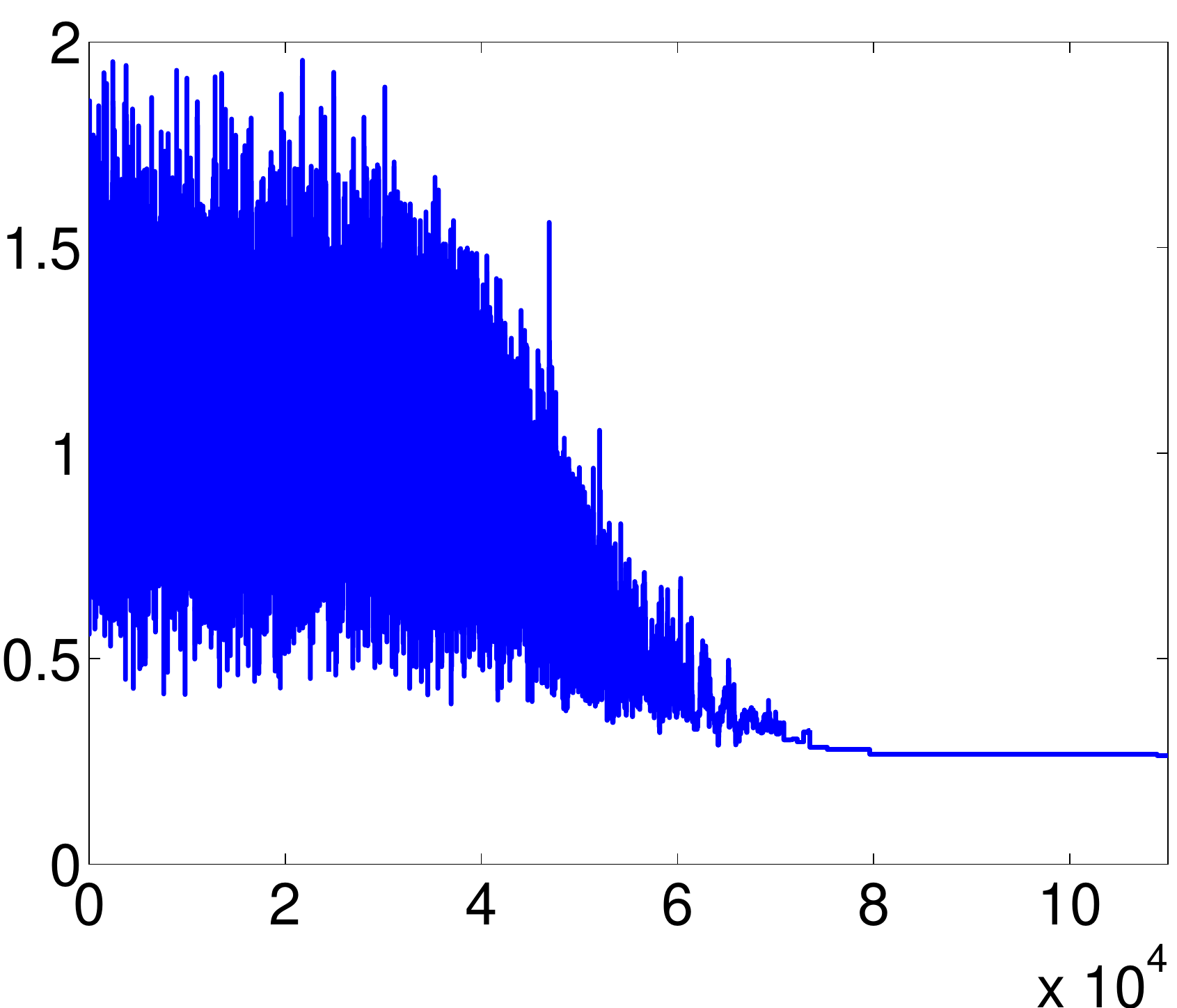}\\
			\multicolumn{4}{|c|}{Geodesic after alignment}\\
			\multicolumn{2}{|c||}{\includegraphics[width=2.1in]{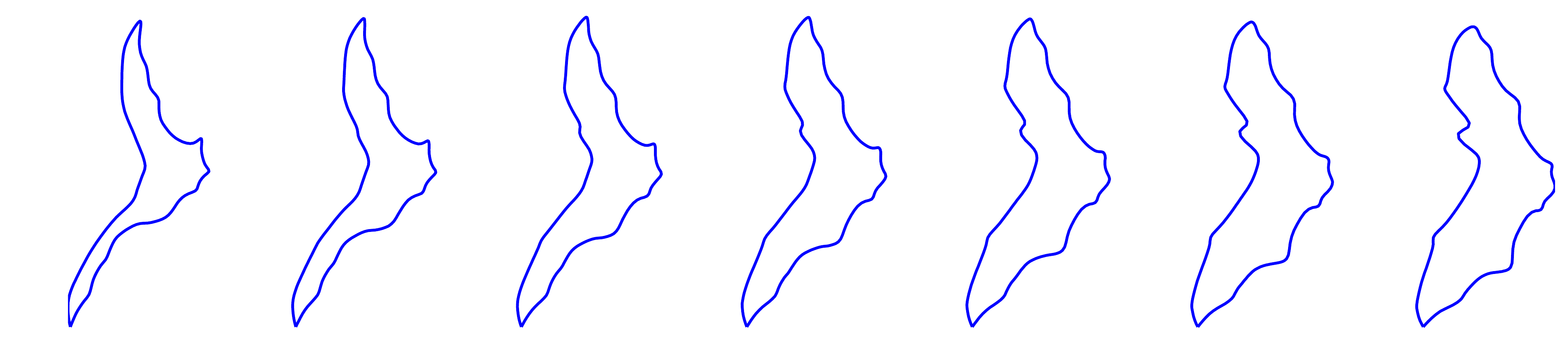}}&\multicolumn{2}{c|}{\includegraphics[width=2.1in]{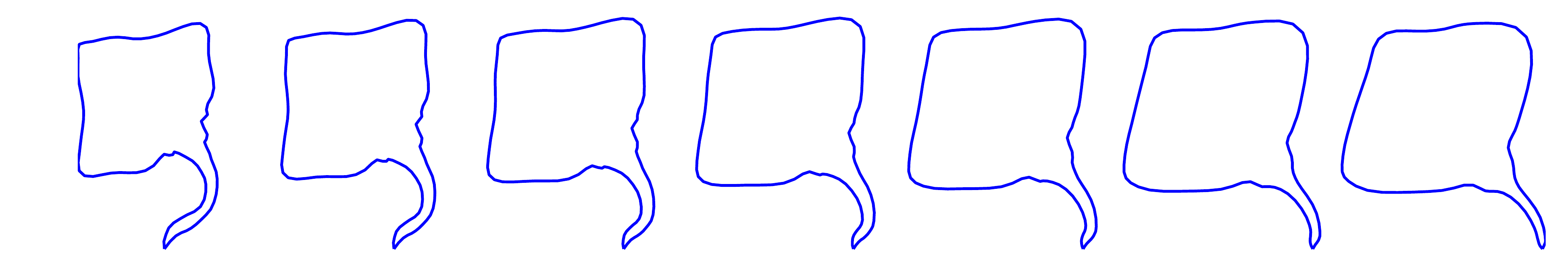}}\\
			\multicolumn{4}{|c|}{Geodesic before alignment}\\
			\multicolumn{2}{|c||}{\includegraphics[width=2.1in]{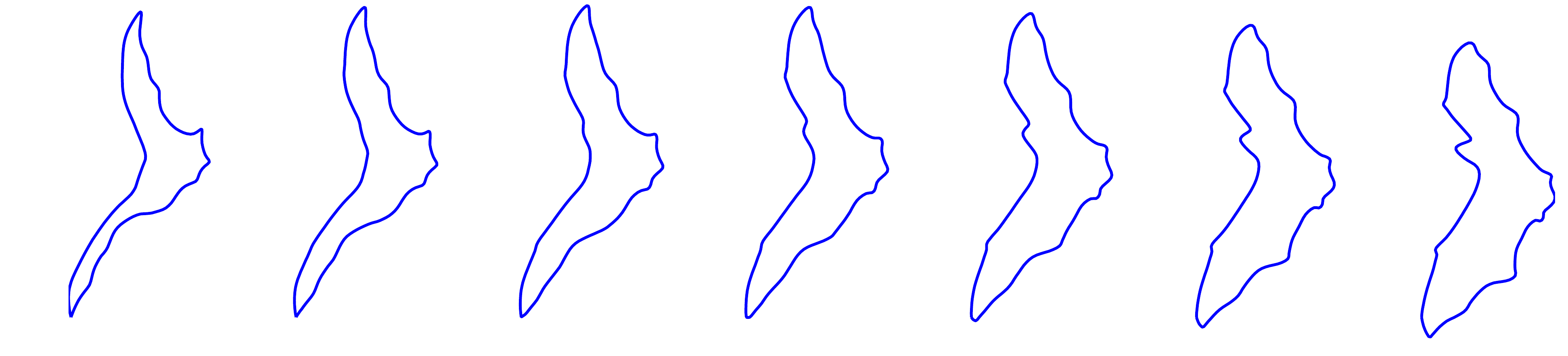}}&\multicolumn{2}{c|}{\includegraphics[width=2.1in]{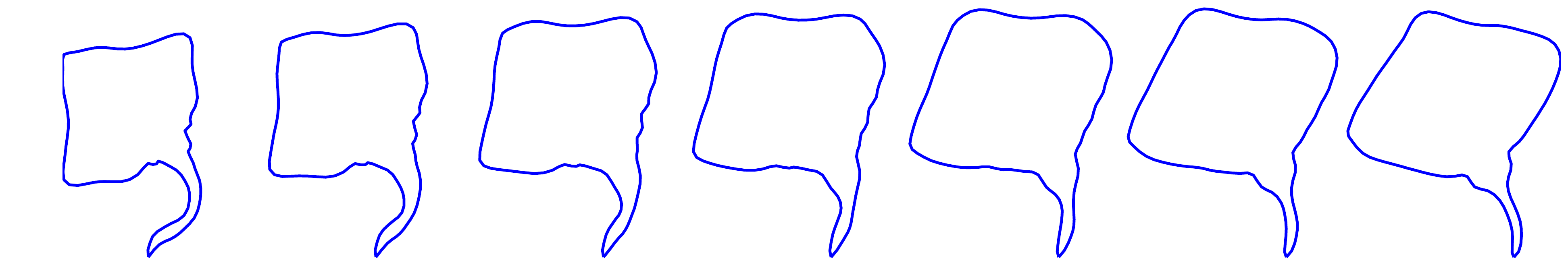}}\\
			\hline
		\end{tabular}
		\caption{\small Results of Simulated Annealing-based alignment on four different examples of closed curves. (a) Optimal (green) and identity (red) warp maps. (b) Evolution of $E(\gamma)$.}\label{fig:SAC}
	\end{center}
\end{figure}

We close this section with four examples of registering closed curves from the MPEG-7 dataset using Simulated Annealing; these curves represent fairly complex shapes including a cup with a handle and a stingray. Recall that in the case of closed curves, we must optimize over the seed placement (point at which $\sone$ is unwrapped to $[0,1]$) on the curve in addition to the warp map. The results are presented in Figure \ref{fig:SAC}. We provide the same displays as in the open curve examples. As previously, the geodesic paths after alignment represent more natural deformations between the shapes than those before alignment. This is especially evident in the cup example. We present our quantitative assessment in the bottom portion of Table \ref{tab:SAOC}. Here, we compare to a DP approach with an additional seed search. The proposed method performs significantly better on examples (1) and (2), and gives comparable performance on examples (3) and (4).

\section{Discussion}
\label{discussion}
The class of warp maps of $[0,1]$  can be identified with the set of distribution or quantile functions on $[0,1]$. Sample paths of Levy subordinators normalized to obtain normalized random measures (Dirichlet process is a special case), and variants thereof \citep{NH, RLP, BRW}, can be used to define distributions on warp maps along with numerous corresponding sampling schemes (see \cite{JG} and references therein). However, warp maps of $\sone$ cannot be identified with distribution functions, and an `intrinsic' distribution based on arc lengths (as described in the Supplementary Material) is not easily obtainable as laws of normalized random measures.

We have given two methods for {\it pairwise} matching of curves: (1) a Bayesian registration model, and (2) a stochastic search algorithm via Simulated Annealing. In many applications, it may be of interest to match multiple curves simultaneously, also termed {\it multiple} registration. This is usually accomplished by joint estimation of a template curve and an additional matching step. The proposed Simulated Annealing-based alignment can be easily incorporated into a multiple registration algorithm by replacing the commonly used DP approach. Multiple alignment via a formal Bayesian model requires a prior on the template curve in addition to the warp maps. Nonetheless, our approach can be readily built into existing Bayesian multiple registration models such as the one presented in \cite{Ian}.

Although unexplored in this paper, the proposed distribution on warp maps is well-suited for curve registration with landmarks observed with uncertainty in their placement. The desirable properties allow us to develop priors centred at piecewise linear warp maps that match the landmarks exactly. This allows for incorporating prior information into the problem by regularizing the warp maps toward a landmark induced warping.\\


\if1\blind
{
\noindent\textbf{Acknowledgements:} We thank Ian Dryden, Huiling Le and Eric Klassen for helpful discussions. We also thank two anonymous reviewers for their suggestions. This research was partially supported by NSF DMS 1613054 and NIH R01 CA214955 (to KB and SK). SK was also partially supported by NSF CCF 1740761.
}\fi


\bibliographystyle{Chicago}
\bibliography{biblio}

\end{document}


\def\spacingset#1{\renewcommand{\baselinestretch}%
{#1}\small\normalsize} \spacingset{1}


\if1\blind
{
  \title{\bf Supplementary Material for ``Distribution on Warp Maps for Alignment of Open and Closed Curves"}
  \author{Karthik Bharath$^1$ and Sebastian Kurtek$^2$\\
    $^1$\small{School of Mathematical Sciences, University of Nottingham, UK}\\
    $^2$\small{Department of Statistics, The Ohio State University, USA}}
    \date{}
  \maketitle
} \fi

\if0\blind
{
  \bigskip
  \bigskip
  \bigskip
  \begin{center}
    {\LARGE\bf Supplementary Material for ``Distribution on Warp Maps for Alignment of Open and Closed Curves"}
\end{center}
  \medskip
} \fi

\spacingset{1.45}

\section{Proofs of all results in the paper}

\noindent\textbf{Proof of Theorem 1:} Let $\tilde{Y}_n(t):=\sum_{i=1}^{\nt}p_i,\ 0 \leq t \leq 1$, an element of $D([0,1])$. Under the uniform norm $\|x \|:=\sup_{0\leq t \leq 1}|x(t)|$, consider
$\delta_n:=\|Y_n-\tilde{Y}_n\|$. Since the limit process $\id$ is deterministic and in $C([0,1])$ (in fact, it is uniformly continuous), the proof follows if we can show that $\tilde{Y}_n$ converges in probability to $\id$ in the uniform topology (Skorohod $J_1$ topology on $D([0,1])$ relativized to $C([0,1])$ coincides with the uniform topology) and $\delta_n \overset{P}\to 0$.

Uniform distribution of the vector $(p_1,\ldots,p_n)$ in $\Dn$ implies that it is equal in distribution to the vector of uniform spacings: $(U_{1:n}-U_{0:n},U_{2:n}-U_{1:n},\ldots,U_{n:n}-U_{n-1:n})$, where $0=:U_{0:n}<U_{1:n}<\ldots<U_{n-1:n}<U_{n:n}:=1$ are order statistics corresponding to $U_1,U_2,\ldots,U_n$ i.i.d. Uniform on $[0,1]$. Then $\tilde{Y}_n(t)=\sum_{i=1}^{\nt} p_i=U_{\nt:n}$. Since $\sqrt{n}(U_{\nt:n}-\id)$ converges in distribution to a standard Brownian Bridge process in $D([0,1])$ (see p. 308 of \cite{VDV}), using the continuous mapping theorem and Chebyshev's inequality,
\begin{small}$$\|U_{\nt:n}-\id\|=\sup_{0\leq t \leq 1}|U_{\nt:n}-\id|\overset{P}\to 0, \quad n\to \infty.$$\end{small}

In order to see why $\delta_n \overset{P}\to 0$, observe that $\tilde{Y}_n$ is a step function based on a partial sum of uniform sample quantiles, and hence changes values only at times $i/n,\ i=1,\ldots,n$. Since $Y_n$ is a linearly interpolated version of $\tilde{Y}_n$, their values coincide at $i/n,\ i=1,\ldots,n$, and in the intervals $\Big(\frac{i-1}{n},\frac{i}{n}\Big]$ they differ at most by the size of the uniform spacing $U_{i:n}-U_{i-1:n}$. Therefore, since $ \max_{1 \leq i \leq n}U_{i:n}-U_{i-1:n}=O_p(\log n/n)$ \citep{LD},
\begin{small}
\begin{align*}
\delta_n&=\sup_{0\leq t \leq 1}|Y_n(t)-\tilde{Y}_n(t)|\leq \max_{1 \leq i \leq n}\sup_{\frac{i-1}{n} \leq t \leq \frac{i}{n}} |Y_n(t)-\tilde{Y}_n(t)|\leq \max_{1 \leq i \leq n}U_{i:n}-U_{i-1:n} \overset{P} \to 0, \quad n \to \infty.
\end{align*}
\end{small}
\qed

\noindent\textbf{Theorem 1 when $\mathcal{T}_n$ is a non-equi-spaced partition:} We describe how a degenerate limit arises even in the case of a non-equi-spaced partition. Suppose $\mathcal{T}_n$ is constructed based on deterministic points $0=:t_0<t_1<\cdots<t_{n-1}<t_n:=1$ such that $t_i-t_{i-1}\neq 1/n$ for all $i=1,\ldots,n$. Consider an equi-spaced partition $\mathcal{T}'_n$ based on $0=:s_0<s_1<\cdots<s_{n-1}<s_n:=1$. Then, there exists a piecewise linear (PL) warp map $\gamma_n:[0,1]\to[0,1]$ establishing exact correspondence between $\mathcal{T}_n$ and $\mathcal{T}'_n$ with $\gamma_n(s_i)=t_i,\ i=1,\ldots,n$. Any such set of $\{t_i,\ i=0,\ldots,n\}$ can be associated with theoretical quantiles of a distribution or quantile function $\gamma$ of a random variable on $[0,1]$. Thus, to obtain any sensible limit under the setup of Theorem 1,  it is reasonable to assume that the continuous $\gamma_n$ satisfies $\|\gamma_n-\gamma\|\overset{P}\rightarrow 0$ as $n \to \infty$, where $\gamma$ is a homeomorphic warp map of $[0,1]$ with inverse $\gamma^{-1}$. The relevant process to consider in this setting is the (deterministic) time-changed process \begin{small}$$Y_n(\gamma_n(t)):=\sum_{i=1}^{\ngt}p_i+(n\gamma_n(t)-\ngt)p_{\ngt+1}, \quad t \in[0,1].$$\end{small} From the proof of Theorem 1 note that $(p_1,\ldots,p_n) \overset{d}=(U_{1:n}-U_{0:n},U_{2:n}-U_{1:n},\ldots,U_{n:n}-U_{n-1:n})$. The map $t \mapsto Y_n(\gamma_n(t))$ is bijective ($\gamma_n^{-1}$ exists and is unique), and we can instead consider the transformed process
\begin{small}$$Z_n(t):=\sum_{i=1}^{\nt}\left[\gamma_n^{-1}(U_{i:n})-\gamma_n^{-1}(U_{i-1:n})\right]+(nt-\nt)\left[\gamma_n^{-1}(U_{\nt+1:n})-\gamma_n^{-1}(U_{\nt:n})\right].$$\end{small}

For ease of notation define the Uniform quantile function $[0,1]\ni t \mapsto E_n^{-1}(t):=U_{\nt:n}$, with $\|E_n^{-1}-\id\| \overset{P}\to 0$, as in the proof of Theorem 1. Consider
\begin{small}$$\tilde{Z}_n:=\sum_{i=1}^{\nt}\left[\gamma_n^{-1}(U_{i:n})-\gamma_n^{-1}(U_{i-1:n})\right]=\gamma_n^{-1}(E_n^{-1}(t)),$$\end{small} an element of $D([0,1])$.
Observe that
\begin{small}$$\|\tilde{Z}_n-\gamma^{-1}\circ \id\|:=\sup_{0 \leq t \leq 1}|\gamma_n^{-1}(E_n^{-1}(t))-\gamma^{-1}(\id (t))|\overset{P}\rightarrow 0,$$\end{small}
which follows from the bound
\begin{small}
\begin{equation}
\label{bound}
 \sup_{0 \leq t \leq 1}|\gamma_n^{-1}(E_n^{-1}(t))-\gamma^{-1}(\id (t))|\leq \sup_{0 \leq t \leq 1}|\gamma_n^{-1}(t)-\gamma^{-1}(t)|+\sup_{0 \leq t \leq 1}|\gamma^{-1}(E_n^{-1}(t))-\gamma^{-1}(\id(t))|,
 \end{equation}
 \end{small}
since the first term on the right hand side converges to zero by assumption, while the second converges to zero in probability since $\|E_n^{-1}-\id\| \overset{P}\to 0$ and $\gamma^{-1}$ is uniformly continuous. We now need to show that
\begin{small} $$\delta_n:=\|Z_n-\tilde{Z}_n\|=\sup_{0 \leq t \leq 1}\left|(nt-\nt)\left[\gamma_n^{-1}(U_{\nt+1:n})-\gamma_n^{-1}(U_{\nt:n})\right]\right| \overset{P}\to 0.$$\end{small}
We have
\begin{small}
\begin{align*}
&\delta_n \leq \sup_{0 \leq t \leq 1}\left|nt-\nt\right|\Big(\sup_{0 \leq t \leq 1}|\gamma_n^{-1}(U_{\nt+1:n})-\gamma^{-1}(t)|
+ \sup_{0 \leq t \leq 1}|\gamma_n^{-1}(U_{\nt:n})-\gamma^{-1}(t)\Big)\\
&=\sup_{0 \leq t \leq 1}\left|nt-\nt\right|\Big(\sup_{0 \leq t \leq 1}|\gamma_n^{-1}(E_n^{-1}(t)+o_{up}(1/n))-\gamma^{-1}(t)|
+ \sup_{0 \leq t \leq 1}|\gamma_n^{-1}(E_n^{-1}(t))-\gamma^{-1}(t)\Big),
\end{align*}
\end{small}
since $\sup_{0 \leq t \leq 1}(U_{\nt+1:n}-U_{\nt:n}) \overset{P}\to 0$ and $(U_{\nt+1:n}-U_{\nt:n})=O_P(1/n)$; the notation $o_{up}$ denotes uniform convergence in probability over $t\in[0,1]$. Observing that $\left|nt-\nt\right|=O(1/n)$ and $\id(t)=t$ by definition, and by employing the reasoning in Equation (\ref{bound}), we conclude that $\delta_n$ converges to zero in probability as $n \to \infty$. We have thus established that $Z_n$ converges in probability to the deterministic $\gamma^{-1}$ in the uniform topology.
\qed
\vspace*{10pt}
\noindent\textbf{Proof of Corollary 1:} The result is a special case of the above result with non-equi-spaced partitions based not on a PL $\gamma_n$, but instead a smooth $Q$. As with the proof of Theorem 1, we first show that the result is true for $\tilde{Y}_n(t):=\sum_{i=1}^{\nt}p_i \in D([0,1])$ since the limit $Q$ is in $C([0,1])$. We are hence required to prove that $\sup_{0\leq t \leq 1}|\tilde{Y}_n(t)-Q(t)|\overset{P}\to 0$ as $n \to \infty$. Since $p_i=x_{i:n}-x_{i-1:n}$, $\sum_{i=1}^{\nt}p_i=x_{\nt:n}$. Then, with $U_{i:n}$ as Uniform order statistics,
\begin{small}
\begin{align*}
|\tilde{Y}_n(t)-Q(t)|&=|x_{\nt:n}-Q(t)|=|Q(U_{\nt:n})-Q(t)|=\frac{1}{f(Q(\alpha_{n,t}))}|(U_{\nt:n})-t|,
\end{align*}
\end{small}
where $\min(U_{\nt:n},t)<\alpha_{n,t}<\max(U_{\nt:n},t)$. From the assumptions on $F$, we have
\begin{small}
\begin{align*}
\sup_{0\leq t \leq 1}|\tilde{Y}_n(t)-Q(t)|& = \sup_{0 \leq t \leq 1}\frac{|U_{\nt:n}-t|}{f(Q(\alpha_{n,t}))} = \sup_{0 \leq t \leq 1}\frac{|U_{\nt:n}-t|}{f(Q(\alpha_{n,t}))-f(Q(t))+f(Q(t))}\\
&= \sup_{0 \leq t \leq 1}\frac{|U_{\nt:n}-t|}{|f(Q(t))+(\alpha_{n,t}-t)f^{'}(Q(\beta_{n,t}))/f(Q(\beta_{n,t}))|},
\end{align*}
\end{small}
where $\min(\alpha_{n,t},t)<\beta_{n,t}<\max(\alpha_{n,t},t)$. Since $\sup_{0 \leq t \leq 1}|U_{\nt:n}-t|\overset{P}\to 0$, and the denominator is bounded in probability due to assumptions on $F$, we have $\sup_{0\leq t \leq 1}|\tilde{Y}_n(t)-Q(t)|\overset{P}\to 0$.

To show that $\tilde{Y}_n$ can be uniformly approximated by $Y_n$, an identical argument as in Theorem 1 can be used once we observe that
\begin{small}
\[
\max_{1 \leq i \leq n} x_{i:n}-x_{i-1:n}=\max_{1 \leq i \leq n}\frac{U_{i:n}-U_{i-1:n}}{f(Q(\alpha_{i,n}))}\leq \max_{1 \leq i \leq n}(U_{i:n}-U_{i-1:n})\inf_{0 \leq t \leq 1}\frac{1}{f(Q(t))},
\]
\end{small}
where $U_{i-1:n}<\alpha_{i,n}<U_{i:n}$ a.s.; the numerator converges in probability to zero since it is $O_p(\log n/n)$ and the denominator is positive due to assumptions on $F$.
\qed
\\
\vspace*{10pt}

\noindent\textbf{Corollary 1 when $\mathcal{T}_n$ is a non-equi-spaced partition:} Following proof of Theorem 1 with non-equi-spaced partition $\mathcal{T}_n$, we only need to consider the time-changed process ${Y}_n\circ \gamma_n$, where $\gamma_n$ is the piece-wise linear warp map matching  $\mathcal{T}_n$ to an equi-spaced partition $\mathcal{T}'_n$. Since $\gamma_n$ is bijective with inverse $\gamma_n^{-1}$, we can equivalently consider the process $Z_n:[0,1] \to [0,1]$ where 
$$Z_n(t):=\sum_{i=1}^{\nt}\left[\gamma_n^{-1}(Q(U_{i:n}))-\gamma_n^{-1}(Q(U_{i-1:n}))\right]+(nt-\nt)\left[\gamma_n^{-1}(Q(U_{\nt+1:n}))-\gamma_n^{-1}(Q(U_{\nt:n}))\right].$$ 
The proof follows using identical arguments as with proof of Theorem 1 with non-equi-spaced partition in conjunction with the assumptions on the density $f$.
\qed
\\

\vspace*{10pt}
\noindent\textbf{Proof of Theorem 2:} Our method of proof is to first show that the result holds when $f$ is Uniform on $[0,1]$. We then show that the sequence of point processes constructed using a non-Uniform $f$ can be approximated arbitrarily well by the one with Uniform $f$. We first review some preliminary tools. Let $t_1,t_2,\ldots$ be independent Uniform random variables on $[0,1]$, $E_1,E_2,\ldots$ be independent unit-mean Exponential random variables, and $G_i=E_1+\ldots+E_i$ a Gamma random variable with shape $i$ and scale equal to one. We begin with the following observation. If $\tilde{\Po}:=\sum_{i=1}^\infty \delta_{\{t_i, G_i\}}$ is a Poisson point process on $[0,1]\times \mathbb{R}_+$ with mean measure $dt \times dx$, then $\lambda^{-1}\tilde{\Po}:=\sum_{i=1}^\infty \delta_{\{t_i, \lambda^{-1}(G_i)\}}\overset{d}{=}\Po$ \citep{FK}. Our proof uses this representation of the Gamma process.

Let $M$ denote the set of all Radon measures taking values in $\mathbb{N}$, and $\mathcal{M}$ the corresponding $\sigma$-algebra. Convergence in $\mathcal{M}$ is defined through vague convergence of measures: $\mu_n\to \mu$ if and only if $\int \phi d\mu_n \to \int \phi d\mu$ for all test functions $\phi$ which are continuous on a compact support. The vague topology is metrizable such that $\mathcal{M}$ is a complete, separable metric space, and weak convergence of probability measures is meaningful \citep{karr}. Therefore, we can view the point processes $\Po_n$ and $\Po$ as random measures on $[0,1]\times \mathbb{R}_+$, and we consider a compactified $\mathbb{R}_+$ by including $\infty$. The Prohorov distance between $\Po_n$ and $\Po$ (by an abuse of notation) is defined as
\begin{small}
\[
\rho(\Po_n,\Po):=\inf\left\{\epsilon>0: P(\Po_n \in B)\leq P(\Po \in B^\epsilon)+\epsilon, B \in \mathcal{M}\right\},
\]
\end{small}
where $B^\epsilon=\{\mu: d(\mu, B)<\epsilon\}$ based on the Prohorov distance $d$ on $M$; refer to page 28 of \cite{karr} for details. The Prohorov distance provides a metric for vague convergence: $\Po_n$ converges in distribution to $\Po$ if and only if $\rho(\Po_n,\Po)\to 0$ \citep{bill} . However, the Prohorov distance is bounded above by the total variation distance:
\begin{small}$$TV(\Po_n,\Po):=\sup\left\{|P(\Po_n \in B)-P(\Po \in B)|: B \in \mathcal{M}\right\}.$$\end{small} Consequently, our method of proof is to show that $TV(\Po_n,\Po)\to 0$ as $n \to \infty$. We next define a few quantities used in the proof.

\noindent\emph{Proof of Part (1):} Define the mapping $\lambda: t \mapsto \int_{t}^\infty \frac{e^{-u}}{u} du$, and its left continuous inverse $\lambda^{-1}$. Since $\Po_n$ is a (random) product measure on $[0,1]\times \mathbb{R}_+$, convergence of $\Po_n$ to $\Po$ in the total variation distance can be established through the convergence of the individual components of the (random) product measure (see, for example, Lemma 3.3.7 of \cite{reiss}):
\begin{small}
\begin{align}\label{TV}
&TV\Big(\sum_{k=1}^{n} \delta_{\{t_k,\lambda^{-1}(\eta_{k-1})\}}, \sum_{k=1}^\infty \delta_{\{t_k,\lambda^{-1}(G_k)\}}\Big)\leq TV\left(\sum_{k=1}^{n} \delta_{\{t_k\}}, \sum_{k=1}^\infty \delta_{\{t_k\}}\right) \nonumber \\
&+TV\left(\sum_{k=1}^{n} \delta_{\{\lambda^{-1}(\eta_{k-1})\}}, \sum_{k=1}^\infty \delta_{\{\lambda^{-1}(G_k)\}}\right).
\end{align}
\end{small}
In general, in order to establish convergence of a point process $\sum_{i=1}^{n} \delta_{z_{i,n}}$ to $\sum_i \delta_{z_i}$ in total variation distance, it suffices to show that for every $k<n$, as $n \to \infty$
\begin{small}
\begin{equation}
\sup\left\{\Big|P\left[(z_{1,n},\ldots,z_{k,n})\in B\right]-
P\left[(z_1,\ldots,z_k)\in B\right]\Big|: B \in \mathcal{B}(\mathbb{R}_+^k)\right\}\longrightarrow 0, \label{mapping}
\end{equation}
\end{small}
since the composition of the vector $(z_1,\ldots,z_k)$ with the measurable map $(z_1,\ldots,z_k)\mapsto \sum_{i=1}^k \delta_{z_i}$ provides a representation of $\sum_i \delta_{z_i}$ (p. 192 of \cite{reiss}).

Consider now the sequence $\Po_n:=\sum_{i=1}^n\delta_{t_i,\lambda^{-1}(nv_i)}$ corresponding to the case when $f$ is the Uniform density on $[0,1]$. Trivially, the first term on the right hand side of Equation (\ref{TV}) converges in total variation to zero on Borel sets on $[0,1]$.  In order to show the convergence of the second term, consider the one-dimensional point process $\mathcal{N}_n:=\sum_{k=1}^{n} \delta_{\lambda^{-1}(z_{k,n})}$. Based on Equation $(\ref{mapping})$, we need to show that $\mathcal{N}_n$ converges weakly to $\mathcal{N}:=\sum_{k=1}^\infty \delta_{\lambda^{-1}(G_k)}$. Since $(p_1,\ldots,p_n)$ is uniformly distributed on $\Dnn$, $p_i$ are distributed as the spacings of i.i.d. Uniform random variables on $[0,1]$ with $f(x)=1,\ 0\leq x \leq 1$. It is known that $(np_1,\ldots,np_k)\overset{d}{\to}(E_1,\ldots,E_{k})$ for every $k\in \{1,\ldots,n\}$. In fact (see p. 201 in \cite{reiss}),
\begin{small}
\begin{equation}\label{expo}
\sup_{B \in \mathcal{B}(\mathbb{R}_+^{k})}\Big|P\left[(np_1,\ldots,np_k)\in B\right]-
P\left[(E_1,\ldots,E_{k})\in B\right]\Big|\leq \frac{Ck}{n},
\end{equation}
\end{small}
for some constant $C>0$. This implies that $p_i$ are asymptotically independent Exponential random variables with unit mean. Therefore, for each $k$, $n(p_1+\ldots+p_k)=nv_k$ is asymptotically distributed as $G_{k}$. For our purposes
\begin{small}
\begin{equation}\label{gamma}
\sup_{B \in \mathcal{B}(\mathbb{R}_+^{k})}\Big|P\left[(nv_1,\ldots,nv_k)\in B\right]-
P\left[(G_1,\ldots,G_{k})\in B\right]\Big|\to 0,
\end{equation}
\end{small}
as $n \to \infty$. Our interest is in showing that, as $n \to \infty$,
\begin{small}
	\begin{equation}
	\label{lambda_gamma}
	\sup_{B \in \mathcal{B}(\mathbb{R}_+^{k})}\Big|P\left[(\lambda^{-1}(nv_1),\ldots,\lambda^{-1}(nv_k))\in B\right]-
	P\left[(\lambda^{-1}(G_1),\ldots,\lambda^{-1}(G_{k}))\in B\right]\Big|\to 0.
	\end{equation}
\end{small}It suffices to consider only bounded rectangles as the sets $B\in \mathcal{B}(\mathbb{R}_+^{k})$, since finite unions of bounded rectangles form a basis of relatively compact open sets closed under finite unions and intersections \citep{kallenberg}.  For a monotone function $h:\mathbb{R}\to \mathbb{R}$, consider its extension $\bar{h}: \mathbb{R}_+^{k}\to \mathbb{R}_+^{k}$ such that $\bar{h}(x_1,\ldots,x_{k})=(h(x_1),\ldots,h(x_{k}))$. For any finite union $B=\cup_{i=1}^{k}(a_i,b_i]$ of half-open rectangles, $\bar{h}(B)=\cup_{i=1}^{k}(h(a_i),h(b_i)]$. The inverse Levy measure $\lambda^{-1}$ is left continuous and monotone decreasing; hence, the preceding probability is equal to
\begin{small}
\begin{equation}
\sup_{\lambda(B) \in \mathcal{B}(\mathbb{R}_+^{k})}\Big|P\left[(nv_1,\ldots,nv_k)\in \lambda(B)\right]-
P\left[(G_1,\ldots,G_{k})\in \lambda(B)\right]\Big|\to 0,
\end{equation}
\end{small}
for all $\lambda(B)$ expressed as finite unions of half-open rectangles. Evidently $\lambda(B)\in \mathcal{B}(\mathbb{R}_+^{k})$, and by Equation (\ref{gamma}), we have shown Equation (\ref{lambda_gamma}) to be true.

We now consider the sequence of point processes when $f$ is a non-Uniform density on $[0,1]$, satisfying the assumptions in the statement of Theorem 2. Here, $(p_1,\ldots,p_n)$ is the vector of spacings of $x_i$, which are i.i.d. $f$ non-Uniform on [0,1]. Consider $(u_1,\ldots,u_n)$, uniformly distributed on $\Dnn$. Using similar arguments as before, in particular the relationship between Equation (\ref{expo}) and Equation (\ref{gamma}), it is enough to show that, for every $k<n$, the probability $P\Big[(\lambda^{-1}(f(Q(\zeta_{1,n})v_1)),\ldots,\lambda^{-1}(f(Q(\zeta_{k,n})v_k)))\in B\Big]$ converges to $P\Big[(\lambda^{-1}(E_1),\ldots,\lambda^{-1}(E_{k}))\in B\Big]$ uniformly over all subsets $B \in  \mathcal{B}(\mathbb{R}_+^{k})$. For this, it suffices to show that
\begin{small}
\[
\sup_{B \in \mathcal{B}(\mathbb{R}_+^{k})}\Big|P\left[(nu_1,\ldots,nu_k)\in B\right]-
P\left[(f(Q(\zeta_{1,n})np_1,\ldots,f(Q(\zeta_{k,n})np_k) \in B\right]\Big|\to 0,
\]
\end{small}in conjunction with the monotonicity property of $\lambda^{-1}$. In other words, we are required to show that the two sequences $(nu_1,\ldots,nu_k)$ and $(f(Q(\zeta_{1,n}))np_1,\ldots,f(Q(\zeta_{k,n}))np_k)$ are close in total variation distance over $B \in \mathcal{B}(\mathbb{R}_+^{k})$.

It is known (see Theorem 2.3 on p. 22 of \cite{anirban}) that weak convergence of a sequence of probability measures with unimodal densities to a measure with a unimodal density is equivalent to their convergence in total variation metric; this explains the unimodal assumption on densities of $p_i$ in the statement of Theorem 2. Unimodality of the density of $p_i$ is unaffected by the constant scale factor $f(Q(\zeta_{i,n}))$. We proceed in two steps. First, we show that, for each $i=1,\ldots,n$,
\begin{small}
\[
\left|f(Q(\zeta_{i,n}))np_i-nu_i\right|\overset{P}{\longrightarrow}0, \quad n \to \infty.
\]
\end{small}
Since $nu_i$ converges in distribution to a unit-mean Exponential $E_i$, we can claim that $f(Q(\zeta_{i,n})np_i$ also converges in distribution to $E_{i}$, for each $i=1,\dots,k$ (Theorem 4.1 in \cite{bill}). Second, for $i=1,\ldots,k$, if $\mathbb{P}_i$ denotes the probability measure of $f(Q(\zeta_{i,n})np_i$, $\mathbb{P}^{'}_i$  the probability measure of $nu_i$, and $\mathbb{E}_{i}$ denotes the probability measure of $E_{i}$, then, for every $k<n$, we can use an upper bound on the total variation distance between the $\otimes_{i=1}^{k} \mathbb{P}_i$ and $\otimes_{i=1}^{k} \mathbb{P}^{'}_i$. Specifically, using the following inequality for the product measures under the total variation metric (Lemma 3.3.7 of \cite{reiss})
\begin{small}
\[
TV(\otimes_{i=1}^{k} \mathbb{P}_i, \otimes_{i=1}^{k} \mathbb{P}^{'}_i)\leq \sum_{i=1}^{k} TV(\mathbb{P}_i, \mathbb{P}^{'}_i) \leq \sum_{i=1}^{k} \left[TV(\mathbb{P}_i, \mathbb{E}_{i})+TV(\mathbb{E}_{i}, \mathbb{P}^{'}_i)\right],
\]
\end{small}
we would have the required result since we have shown that each unimodal $\mathbb{P}_i$ converges to a $\mathbb{E}_i$, which is unimodal, and consequently in total variation metric over all Borel sets in $\mathbb{R}_+$; the same holds for the measure $\mathbb{P}^{'}_i$.

We now prove the first step. Bearing in mind that $u_i=U_{i:n}-U_{i-1:n}$ and $p_i=x_{i:n}-x_{i-1:n}\overset{d}{=}Q(U_{i:n})-Q(U_{i-1:n})$, note that
\begin{small}
\begin{align*}
\left|f(Q(\zeta_{i,n}))np_i-nu_i\right|&=\left|f(Q(\zeta_{i,n})n\left[Q(U_{i:n})-Q(U_{i-1:n})\right]-nu_i\right|\\
&=\left|\frac{f(Q(\zeta_{i,n}))}{f(Q(\alpha_{i,n}))}nu_i-nu_i\right|
=\left|\frac{f(Q(\zeta_{i,n}))}{f(Q(\alpha_{i,n}))}-1\right|nu_i,
\end{align*}
\end{small}
where $U_{i-1:n}<\alpha_{i,n}<U_{i:n}$ a.s., and $Q'(t)=[f(Q(t))]^{-1},\ 0\leq t \leq 1$, justified by the assumptions on $F$. Consequently,
$
\left|\frac{f(Q(\zeta_{i,n}))}{f(Q(\alpha_{i,n}))}\right|\overset{P}\longrightarrow 1,
$
since $F$, $Q$, and the density $f$ are assumed to be continuous, and by the condition that $\lim_{x \downarrow 0}f(x)<\infty$ and positive. Since $nu_i=n(U_{i:n}-U_{i-1:n})=O_P(1)$ and bounded in probability, we have the required result.
\qed\\

\vspace*{10pt}
\noindent\underline{\emph{Proof of Part (2):}} Recall that the point process or random measure is an element of $\mathcal{M}([0,1]\times \mathbb{R}_+)$, the space of Radon measures on the product space $[0,1]\times \mathbb{R}_+$ equipped with the vague topology. Consider the map $\Phi: \mathcal{M}([0,1])\times \mathbb{R}_+ \to D([0,1])$ defined by
\begin{small}
\[
\Phi\left(\sum_{i=1}^n \delta_{\{t_i,\lambda^{-1}(z_{i,n})\}}\right)(t)=\displaystyle \sum_{i} \lambda^{-1}(z_{i,n})\mathbb{I}_{t_i \leq t}, \quad t \in [0,1].
\]
\end{small}
The map $\Phi$ is a.s. continuous with respect to the distribution of the Poisson random measure (see p. 221 of \cite{resnick}). The continuous mapping theorem applied to $\Po_n$ ensures that $G_n$ converges in $D([0,1])$ to $\mathcal{G}\circ H$ equipped with the $J_1$ topology. Note that while $\lambda^{-1}(z_{i,n})\ldots\geq \lambda^{-1}(z_{1,n})>0$ are the jumps of the Gamma process ranked in decreasing order, the jumps of $G_n(t)=\sum_{i} \lambda^{-1}(z_{i,n})\mathbb{I}_{t_i \leq t}$ are re-ordered based on the locations of $t_i$ in the interval $[0,t]$. That is, if $(\lambda^{-1}(z_{i,n}),t_i)$ are bivariate random variables, and $t_{1:n}<\ldots<t_{n:n}$ are order statistics of $t_i$, then the $\lambda^{-1}(z_{i,n})$ in $G_n$ are the corresponding induced order statistics.

Now consider $G^l_n(t):=G_n(t)+(nt-\nt)\lambda^{-1}(z_{i,\nt+1})$,  the linearly interpolated version of $G_n$. Since the limit process $\mathcal{G}\circ H$ is a pure jump process, whereas $G^l_n$ is a sequence with paths in $C([0,1])$, there are unmatched jumps in the limit process: the jumps on $\mathcal{G}\circ H$ do not have any corresponding jumps in $G^l_n$, and the $J_1$ topology (and the uniform topology) is too strong and inappropriate. The appropriate topology is the weaker Skorohod's $M_1$ topology based on parametric representations of completed graphs of $G_n$. The set
\begin{small}
\[
A_{G_n}:=\Big\{(t,x)\in [0,1]\times \mathbb{R}_+: x=\alpha G_n(t-)+(1-\alpha)G_n(t) \text{ for some }\alpha \in [0,1]\Big\},
\]
\end{small}where $G_n(t-)=\lim_{s \uparrow t}G_n(s)$, is the completed graph of $G_n$ for every $n$. Therefore, the complete graph of $G_n$ besides the points of the graph $\{t,G_n(t): t \in [0,1]\}$ also contain the line segments joining $(t,G_n(t))$ and $(t,G_n(t-))$ for all points of discontinuity $t$ of $G_n$. An order on $A_{G_n}$ is then defined by saying that $(t_1,x_1) \leq (t_2,x_2)$ if either (i) $t_1<t_2$, or (ii) $t_1=t_2$ and $|G_n(t_1-)-x_1|\leq |G_n(t_2-)-x_2|$. With an order on $A_{G_n}$, a parametric representation of $A_{G_n}$ is then defined as a continuous nondecreasing function $m(u,v):[0,1] \to A_{G_n}$, and let $\Theta(G_n)$ be the set of parametric representations of $G_n$ for a fixed $n$. The metric $d_{M_1}$ on $D([0,1])$ that induces the $M_1$ topology is then defined as
\begin{small}
\[
d_{M_1}(G_n,G^l_n)=\inf_{\substack{(u_1,v_1) \in \Theta(G_n)\\(u_2,v_2) \in \Theta(G^l_n)}}\Big\{\sup_{0\leq t \leq 1}|u_1(t)-u_2(t)|\vee \sup_{0\leq t \leq 1}|v_1(t)-v_2(t)|\Big\},
\]
\end{small}
where $a \vee b =\max(a,b)$. Note that convergence in $J_1$ topology implies convergence in the $M_1$ as well. Therefore, $G_n$ converges weakly to $\mathcal{G}\circ H$ in the $M_1$ topology, and in this topology $G_n$ and $G^l_n$ are asymptotically equivalent (see p. 214 in \cite{WW}): $d_{M_1}(G_n,G^l_n)\leq 1/n$. This ensures that $G^l_n$ converges in the $M_1$ topology to $\mathcal{G}\circ H$.
\qed\\

\vspace*{10pt}
\noindent\textbf{Proof of Theorem 3:} The proof uses a point process approach in conjunction with the map $\Phi: \mathcal{M}([0,1])\times \mathbb{R}_+ \to D([0,1])$ used in Part (2) of Theorem 2. The main arguments are similar to the one used in Theorem 7.1 of \cite{resnick} to prove the weak convergence of the partial sum process to a Levy jump process based on the vague convergence of the point process of extremes. The main challenge lies in the fact that the transformed increments are no longer i.i.d. For brevity, we do not reproduce the proof, but prove two Lemmas and cite a result on a Kolmogorov-type maximal inequality for a finite sequence of exchangeable random variables, using which, minor modifications of the arguments of \cite{resnick} lead directly to the required result.
\begin{lemma}
\label{lemma1}
The sequence of point processes $N_n:=\sum_i \delta_{\{i/n,\xi_{i,n}\}}$ converges vaguely to the Poisson point process $N$ on $[0,1]\times (y_0,\infty)$ with intensity measure $dt \times e^{-y}dy$, where $y_0=\inf\{y:e^{-e^{-y}}>0\}$.
\end{lemma}
\begin{proof}
Since $p_i$ are uniform spacings, for every $n$, the joint survival function of $p_i$ is given by
\begin{small}$$P(p_1>r_1,\ldots,p_n>r_n)=\left[1-\sum_{i=1}^nr_i\right]^n.$$\end{small} This implies that the $p_i$ are exchangeable random variables, and consequently, so are $\xi_{i,n}=np_i-\log n$. Using the representation $(p_1,\ldots,p_n)\overset{d}=\left(\frac{E_1}{\sum_{i=1}^{n} E_i},\ldots,\frac{E_{n}}{\sum_{i=1}^{n} E_i}\right)$ \citep{pyke} where $E_i$ are i.i.d. unit-mean Exponential random variables, we can claim that $\max_{1 \leq i \leq n}\xi_{i,n}$ converges in distribution to the (standard) Gumbel distribution with distribution function $M(y)=e^{-e^{-y}},\ y \in \mathbb{R}$. This is because $\max_i E_i-\ln (n+1)$ converges in distribution to the Gumbel \citep{HF} using Slutsky's theorem since $1/n \sum_{i=1}^nE_i \overset{P}\to 1$. The exchangeable $\xi_{i,n}$ are strictly stationary, and  the convergence of their maximum is equivalent to $\lim_{n\to \infty}nP(\xi_{1,n}>y)\to \log M(y) = e^{-y}$.

For the i.i.d. case, the above would have sufficed to claim vague convergence of $N_n$ to $N$. Since $\xi_{i,n}$ are strictly stationary, an additional condition is required. Namely, the conditions popularly known as $D'(u_n)$ and $D(u_n)$ ensure convergence of $N_n$ to $N$ (Theorem 5.7.2 in \cite{leadbetter}). These conditions roughly require dependence between blocks of fixed sizes of the triangular array sequence $\{\xi_{i,n}\}$ to vanish with increasing $n$. From the representation of $(p_i,\ i=1,\ldots,n)\overset{d}=(E_i/\sum_{k=1}^nE_k,\ i=1,\ldots,n)$, we note that $np_i$ are asymptotically independent as $n \to \infty$. It is easy now to check that this property of the $p_i$ ensures that the sequence $\{\xi_{i,n}\}$ satisfies $D'(u_n)$ and $D(u_n)$. We leave the details to the interested reader.
\end{proof}
\begin{lemma}
\label{lemma2}
$\lim_{\epsilon \downarrow0} \lim \sup_{n \to \infty} nE(\xi_{1,n}^2\mathbb{I}_{|\xi_{1,n}|\leq \epsilon})=0.$
\end{lemma}
\begin{proof}
\begin{small}
\begin{align*}
\lim \sup_{n \to \infty} nE(\xi_{i,n}^2\mathbb{I}_{|\xi_{i,n}|\leq \epsilon})&=\lim \sup_{n \to \infty}\displaystyle \int_{\{\xi:|\xi_{1,n}|\in (-\epsilon,\epsilon)\}} \xi^2nP(\xi_{1,n}\in d\xi)
=\displaystyle \int_{\{\xi:|\xi_{1,n}|\in (-\epsilon,\epsilon)\}} \xi^2e^{-\xi},
\end{align*}
\end{small}
since $\lim_{n\to \infty}nP(\xi_{1,n}>y)$ exists and, as seen in the proof of Lemma \ref{lemma1}, equals $M(y)=e^{-y}$. The integrand on the right hand side is bounded above by $\xi^2$, and the integral hence converges to zero as $\epsilon \to 0$.
\end{proof}
The only ingredient missing in using the argument employed in Theorem 7.1 of \cite{resnick} is a Kolomogov-type maximal inequality for
the truncated partial sum $\sum_{i=1}^{\nt}[\xi_{i,n}-E(\xi_{i,n}\mathbb{I}_{\xi_{i,n}\leq \epsilon})]$ based on zero-mean exchangeable random variables $\xi_{i,n}-E(\xi_{i,n}\mathbb{I}_{\xi_{i,n}\leq \epsilon})$. This is readily available from Theorem 1 of \cite{ARP}; conditions that ensure its applicability are easily satisfied in our setup.

\section{Distribution for warp maps on $\sone$ without unwrapping}
Here, we present an alternative construction of a distribution for warp maps of $\sone$. This construction is not explicitly used in the two alignment algorithms presented in the paper. Warp maps of $\sone$ constructed without unwrapping cannot be identified with distribution or quantile functions. Our approach is based on viewing points in $\sone$ as angles expressed in radians with zero identified with $2\pi$. With an arbitrary choice $t \in \sone$ as the origin, a translation $r_s: \sone \to \sone$ of the origin by $s$, defined as the shift $r_s(t)=t-s$, amounts to the operation $t \mapsto (t-s) \mod 2\pi$. The mapping $r_s$ is the rotation operator on $\sone$. An anti-clockwise orientation with respect to the origin can be chosen as a positive orientation and the arc-length distance between two points $s,t \in \sone$ is then the arc-length of the positively oriented segment from $s$ to $t$, denoted by $|s-t|$. Therefore, the rotation $r_s$ is a distance-preserving diffeomorphism since it preserves arc-lengths.

Since the choice of the origin on $\sone$ is arbitrary, any probability measure on the set of warp maps of $\sone$ should be impervious to this choice. The partition-based approach on the unit interval described for $W_I$ can be used profitably here by constructing a point process on $[0,1)$ based on the identification of $\sone$ with $[0,1]$ along with a choice of orientation. We describe the point process-based construction for a uniform partition of $\sone$, which leads to a probability measure centered at the identity map; extension to non-uniform partitions follows along similar lines as in the case of the unit interval.

Let $0<t_{1:n}<\ldots<t_{n-1:n}<1$ be uniform order statistics that form a partition of $[0,1)$. For any $\bar{\gamma} \in W_I$ assume that $\bar{\gamma}(t_{i:n})$ are the anti-clockwise endpoints of $n$ independent arcs of equal length $a_n$, with $a_n \to 0$, that are placed randomly (according to a uniform distribution) and independently on $\sone$. Set $\gamma(t_{i:n}):=\max\{\bar{\gamma}(t_{i:n})-a_n,0\}$. Conditional on $t_{i:n}$, this ensures that $s_i:=\gamma(t_{i:n})-\gamma(t_{i-1:n})$ as $n \to \infty$ share the asymptotic properties of uniform spacings. From R\'{e}nyi's representation and Slutsky's theorem $ns_i$ converge in distribution to independent unit-mean Exponential random variables. Thus, $w_{i,n}=n(s_1+\ldots+s_i)$ are asymptotically Gamma distributed with shape $i$ and scale equal to one.

As in the previous section, let $\lambda(dy)=y^{-1}e^{-y}dy$. Note that $w_{i,n}$ can be identified with bounded subsets of $\mathbb{R}$ since $\sone=\mathbb{R}/2\pi\mathbb{Z}$. It is easy to verify that the conditions for convergence of point processes on a circle by \cite{JH} are satisfied, and we hence have the following result, proof of which is almost identical to the proof of Theorem 2.
\begin{theorem1}
\label{th3}
Let $t_i$ be i.i.d. uniform on $[0,1)$, and let $s_i$ be constructed by random and independent placement of arcs of length $a_n$ on $\sone$. Assume that $a_n \to 0$ with $a_n=o(\log n/ n)$ as $n \to \infty$. The point process $\mathcal{P}_n:=\sum_{i=1}^{n} \delta_{\left\{t_i, \lambda^{-1}(w_{i,n})\right\}}$ converges weakly to a Poisson point process $\mathcal{P}$ on $(0,1] \times \mathbb{R}_+$ with intensity measure $dx \times \lambda(dy)$, and is invariant to the choice of the origin on $\sone$.
\end{theorem1}
\noindent Since we have identified zero with one on $[0,1]$ to represent $\sone$, note that the limit Poisson process is on $(0,1]\times \mathbb{R}_+$ with 0 excluded. This ensures that $\mathcal{P}$ can be identified with a Gamma process $\mathcal{G}^s \circ H$ ($H$ is the Uniform distribution function in the above theorem), and its normalized version $\mathcal{D}^s \circ H$ is referred to as the Dirichlet process. The resulting probability measure $\mathbb{D}^s \circ H$ is invariant to the choice of the origin on $\sone$ since the distribution is completely specified by spacings that are distance-preserving.

Choosing a parameterized $\lambda_{\theta}$ leads to a parametric measure. The measure $\mathbb{D}^s_\theta \circ H$ induced by the random partition based on order statistics from $H$ is interpreted in a similar manner to that on $W_I$.
Due to the identification of zero with one in $[0,1]$, an alternate way of viewing the law of the process $\D \circ H$ on the warp maps of the unit interval $W_I$ is as the law of $\mathcal{D}^s \circ H$ conditioned on $\gamma_s(0)=0$. Consequently, the properties of the measure $\mathbb{D}^s_\theta \circ H$ are identical to the case of warp maps of the interval in Proposition 1. \emph{This ensures that $\mathbb{D}^s_\theta \circ H$ satisfies subset invariance and can be centred at any desired warp map}. Extension to the case of $s_i$ constructed using independently placed arcs according to a non-Uniform distribution (resulting in non-uniform spacings) can be achieved in a manner similar to the interval case. As with Part (2) of Theorem 2, it becomes necessary to assume some conditions on the density $f$ that generates the spacings; from the proof, it is evident that the behaviour of $f$ near zero influences the asymptotics of $s_i$. Under the assumption that $f$ is continuous, bounded, and has a minimum point $m$ with $f(m)>0$ \citep{PH}, Part (2) of Theorem 2 carries through with minor alterations (with $a_n =o(\log n/n)$). The same applies to corresponding results related to Theorem 2 and Proposition 1.

\section{Real data description}

We use small subsets of several standard datasets to illustrate the properties of the proposed distribution on warp maps in the context of function and curve alignment. (1) The PQRST complexes came from the PTB Diagnostic ECG Database \cite{PTB1} on Physionet \cite{goldberger}. Each ECG biosignal was segmented into PQRST cycles using the method of \cite{biosignals}. PQRST refers to the five waves in each cycle: the first, second and third positive peaks are the P, R and T waves, while the first and second negative peaks are the Q and S waves. The ECG is a common diagnostic and monitoring tool for various heart diseases. (2) The Berkeley growth data are height growth functions for 39 boys and 54 girls, which were measured from birth until 18 years of age \cite{tuddenham:1954}. We use growth rate functions (i.e., derivative of growth functions) for boys only. Statistical analysis of such functions is important in biometrics by providing insights into various growth patterns. (3) The foot pressure functions came from the Gait Dynamics in Neuro-Degenerative Disease database on Physionet. The raw data were obtained using force-sensitive resistors, with the output roughly proportional to the force under the foot. Each biosignal was segmented into gait cycles using the method of \cite{biosignals}. Gait signals, and their dynamics, are important in assessing Huntington's and Parkinson's disease progression, which are known to affect motor control. (4) The lung volume respiration data is described in \cite{biosignals}. Each function represents measurements of strain changes during respiration. Assessing breathing variation in imaging studies via such measurements is important in various imaging applications including 4D image reconstruction. (5) We use signature curve data from \cite{Yeung04svc2004:first}. For alignment purposes, we first compute their tangential acceleration functions as described in \cite{RS2}. (6) A DT-MRI image of a brain is composed of $3\times 3$ tensor matrices that describe the constraints on local diffusion of water molecules. Based on this imaging modality, one can extract fiber tracts, i.e., three-dimensional curves, by following principal directions of diffusion, a process called tractography. We use the data of \cite{KurtekJASA}. Such data is often used to estimate structural connectivity in the brain. (7) MPEG-7\footnote{http://www.dabi.temple.edu/\textasciitilde shape/MPEG7/dataset.html} is a large database of shapes commonly used to assess comparison, matching and retrieval algorithms in computer vision.

\section{Additional Simulated Annealing alignment results}

\begin{figure}[!t]
	\begin{center}
			\begin{tabular}{|c|}
            \hline
            \includegraphics[width=1.5in]{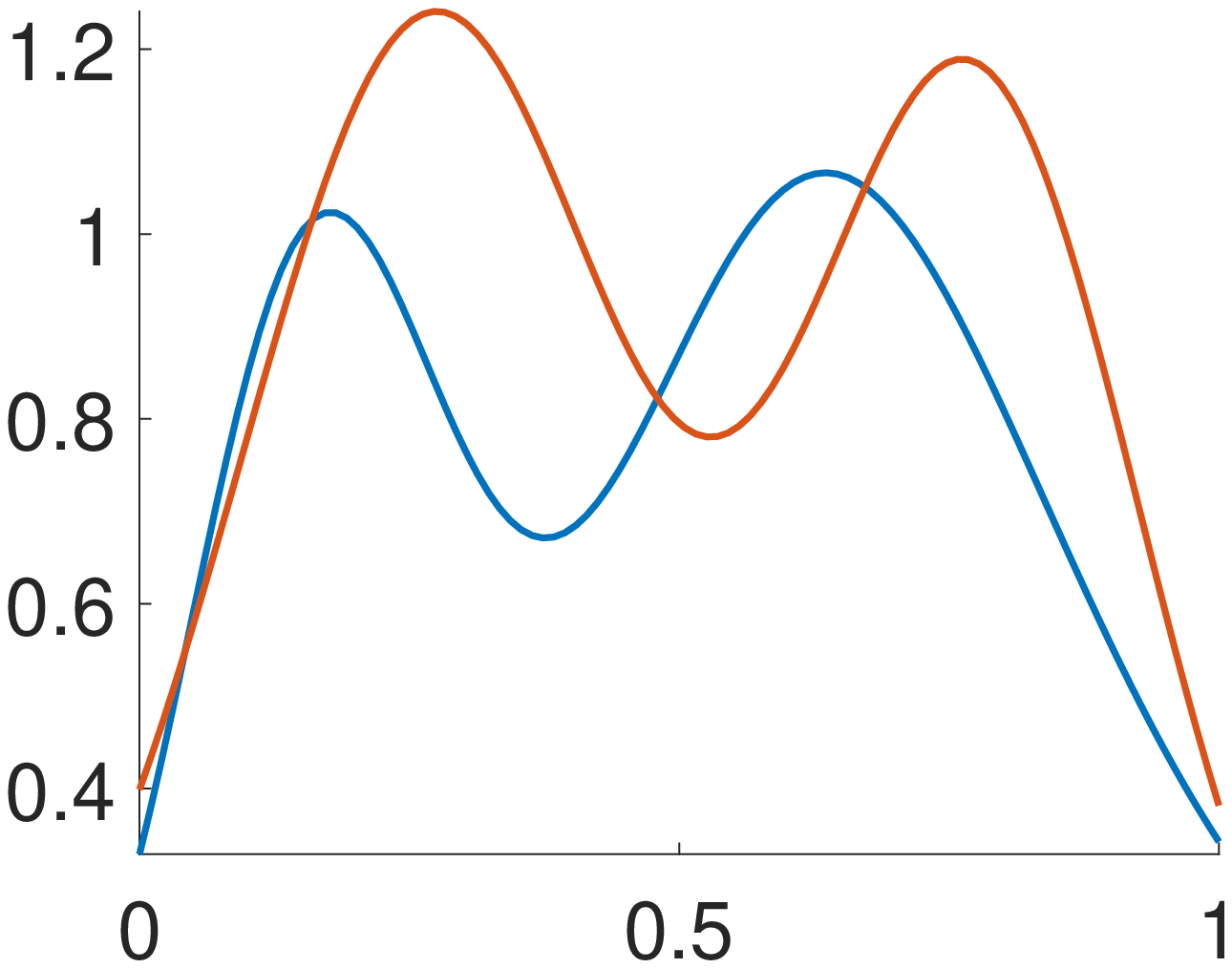}\\
            \hline
            \end{tabular}
			\begin{tabular}{|c|c|c|c|c|}
				\hline
                &$n$&$\theta$&$T$&$c$\\
                \hline
                (a)&\includegraphics[width=1.2in]{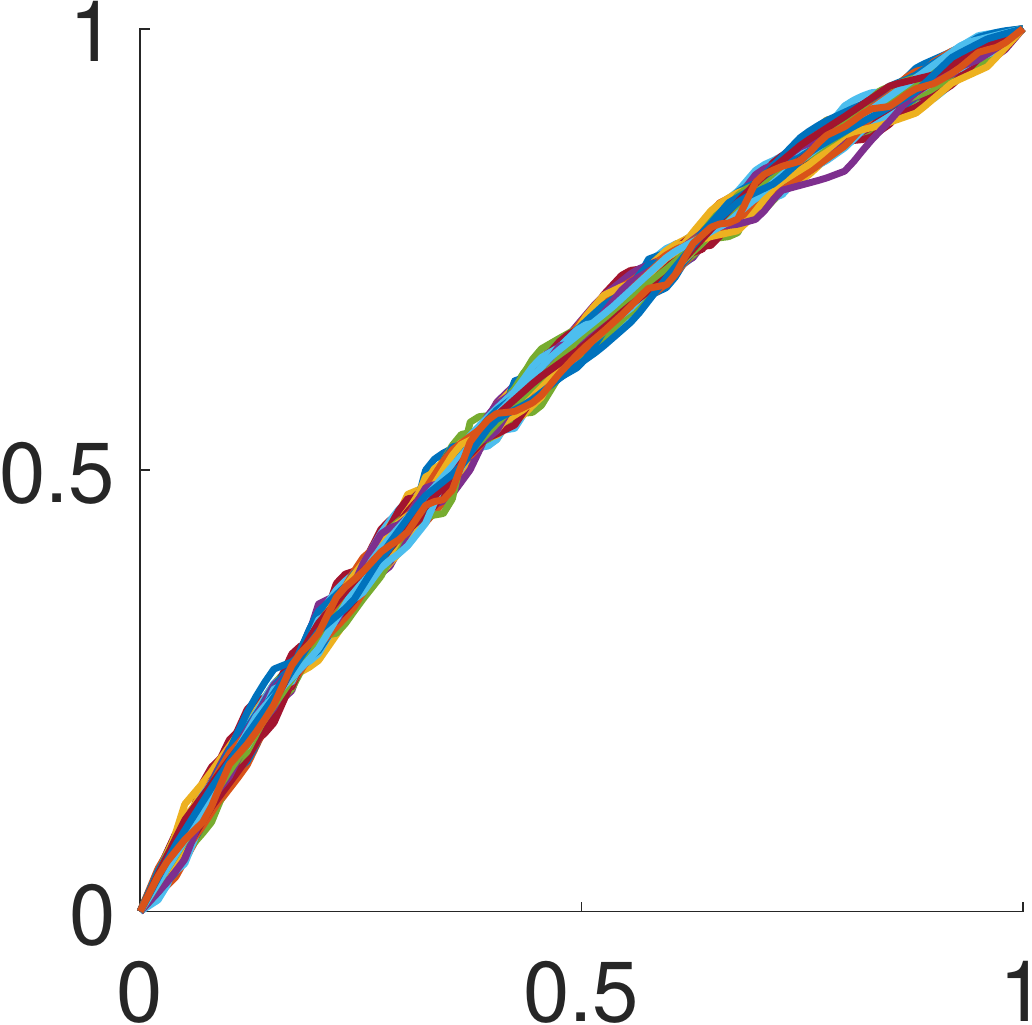}&\includegraphics[width=1.2in]{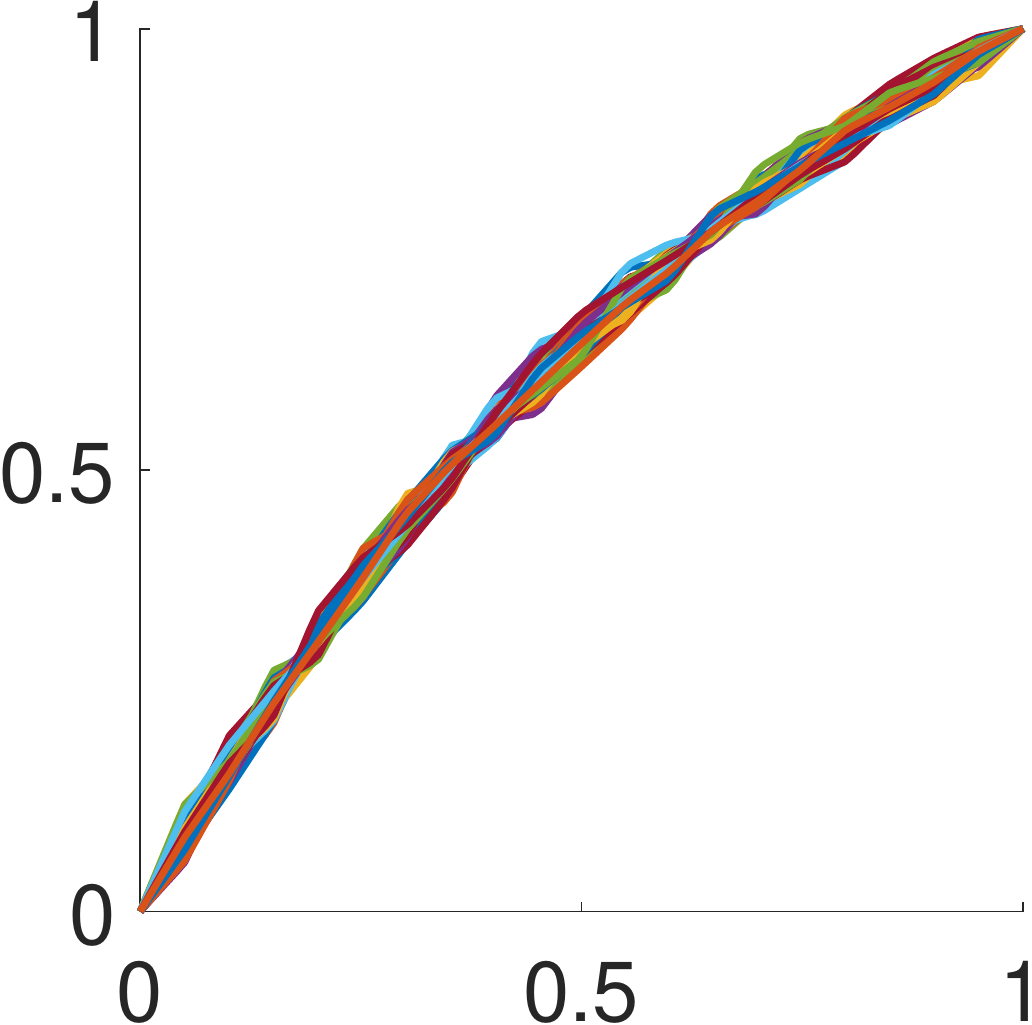}&\includegraphics[width=1.2in]{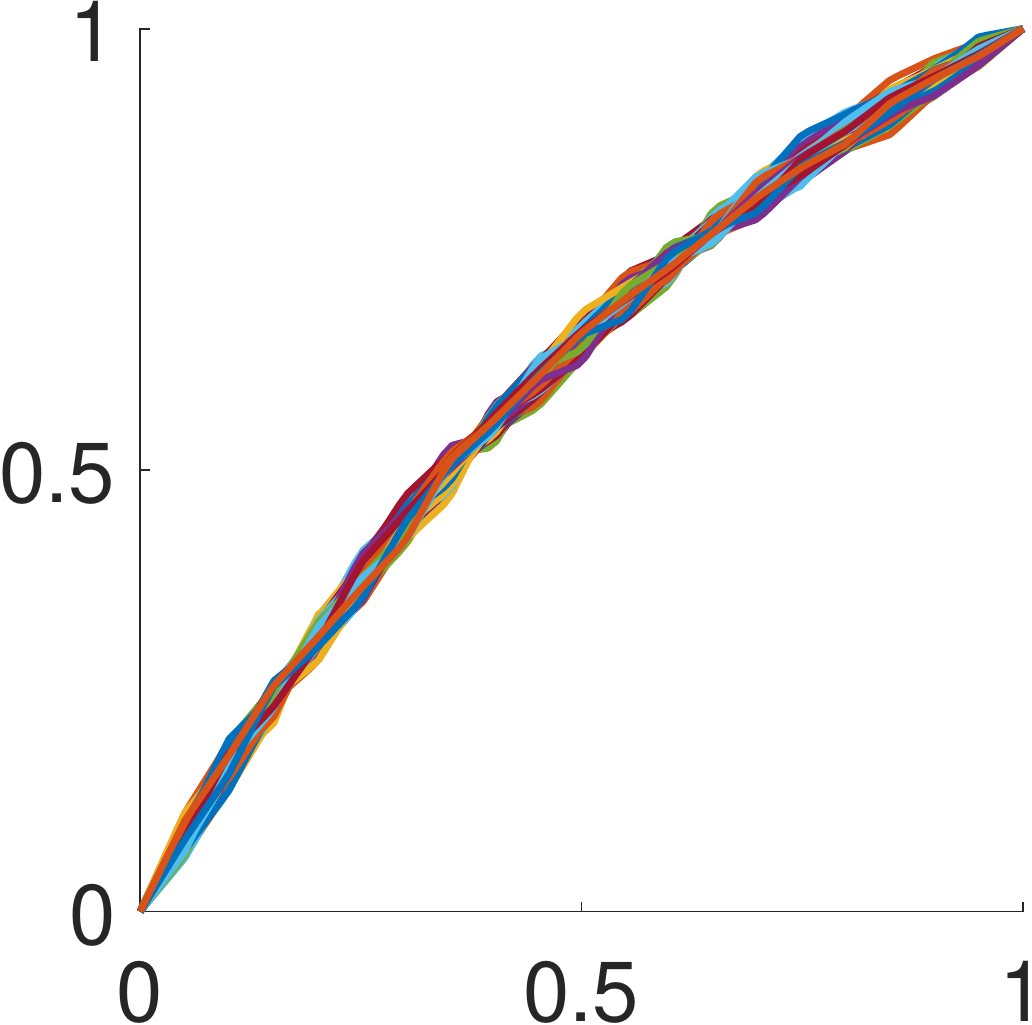}&\includegraphics[width=1.2in]{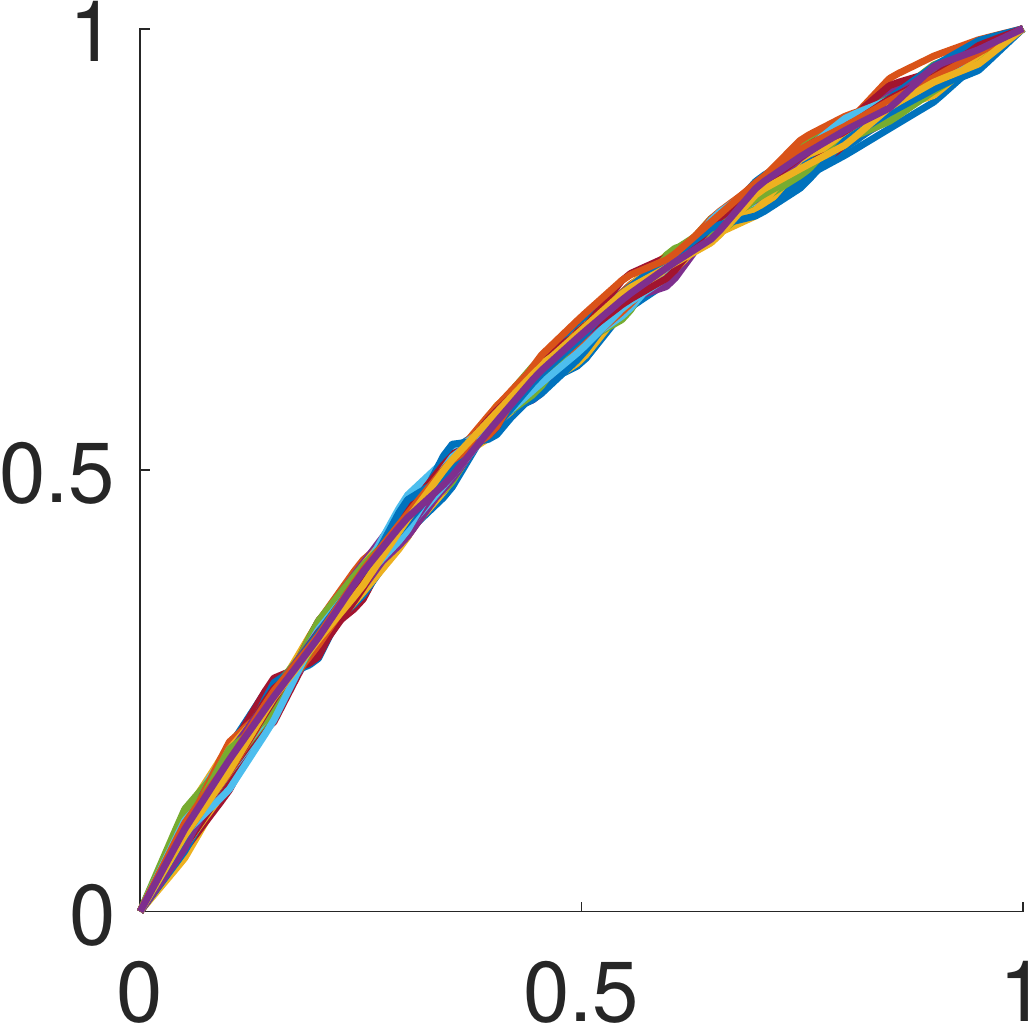}\\
                \hline
                (b)&\includegraphics[width=1.2in]{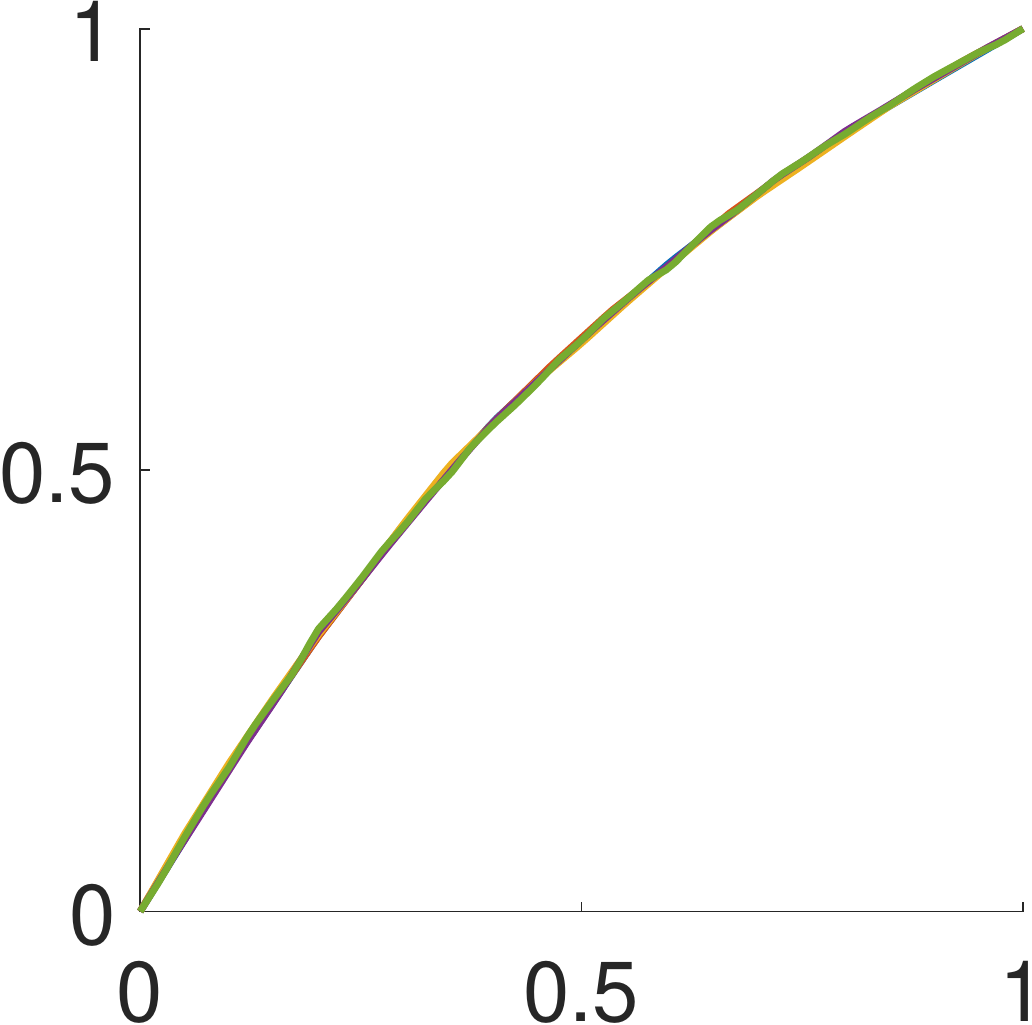}&\includegraphics[width=1.2in]{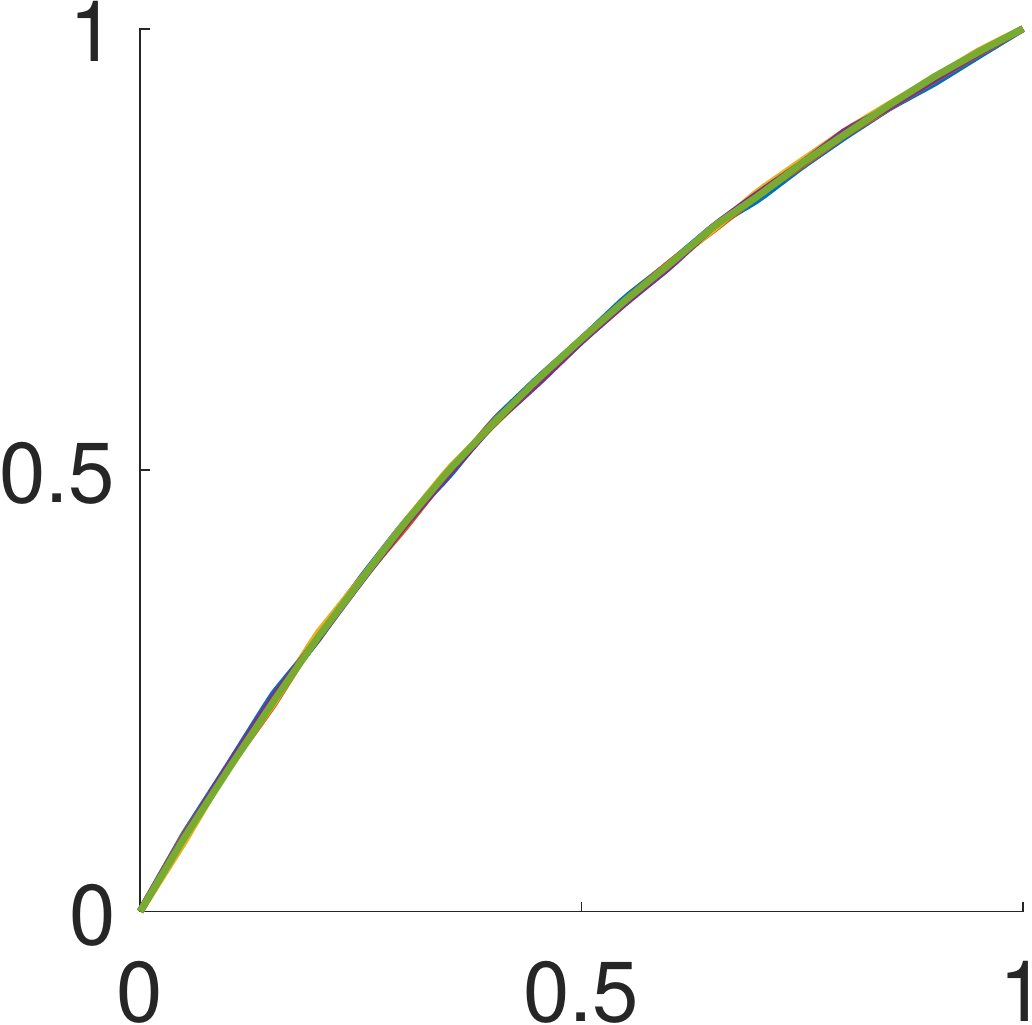}&\includegraphics[width=1.2in]{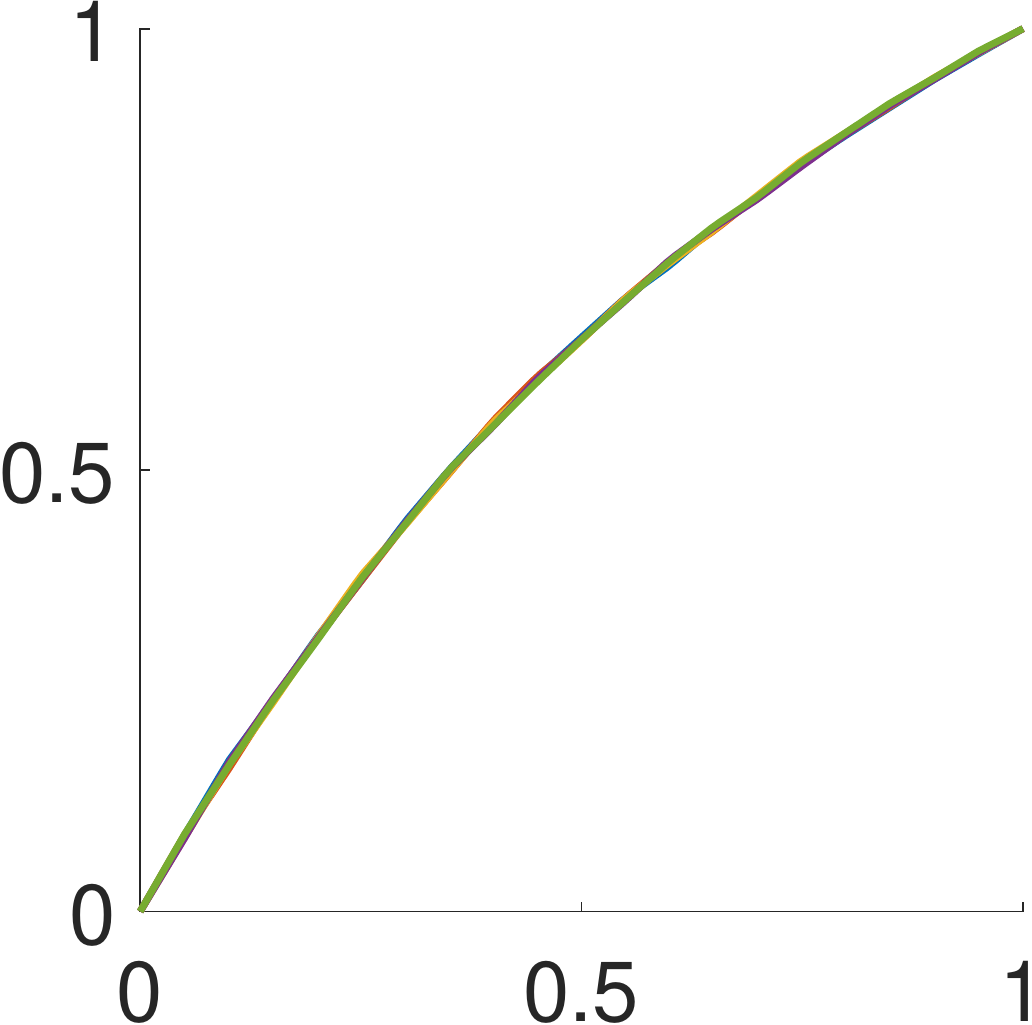}&\includegraphics[width=1.2in]{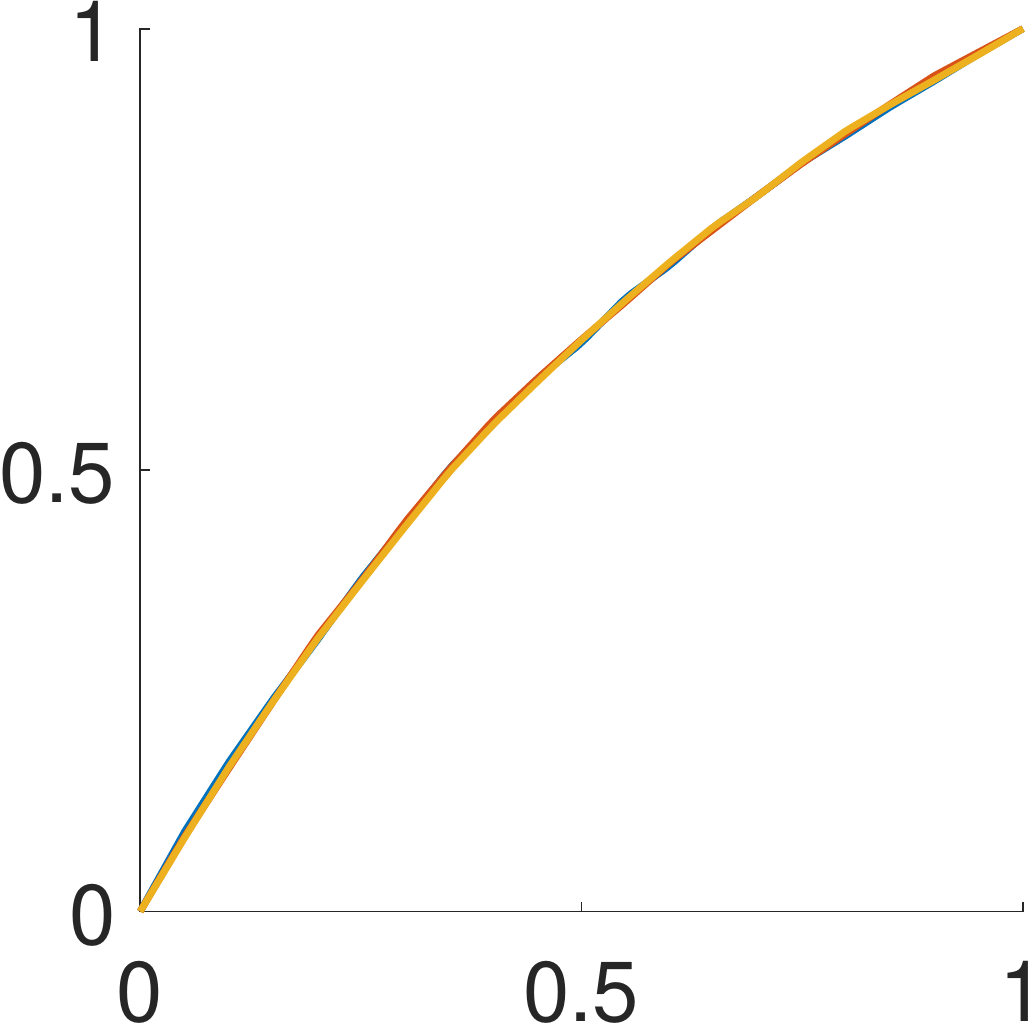}\\
                \hline
                (c)&\includegraphics[width=1.2in]{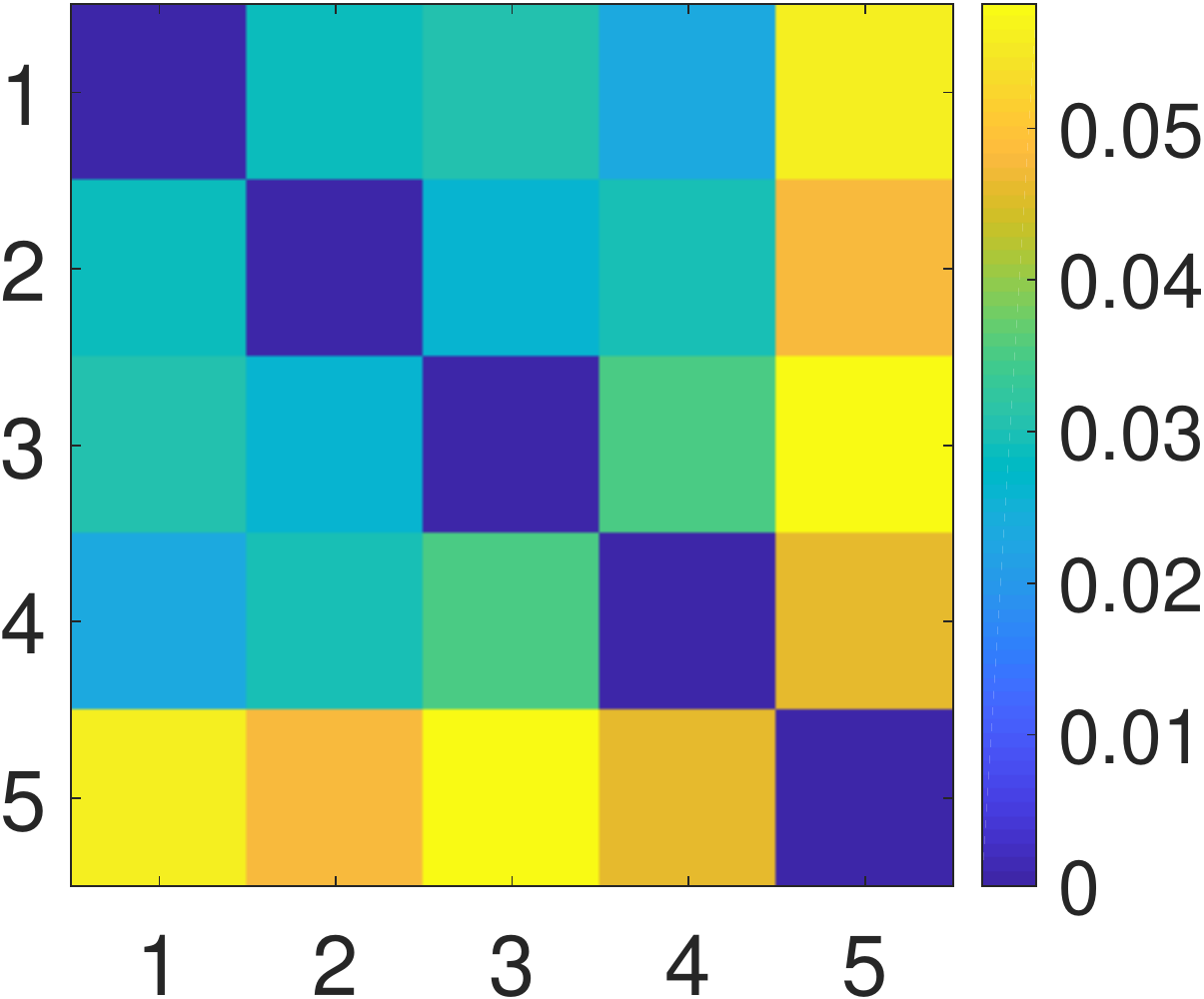}&\includegraphics[width=1.2in]{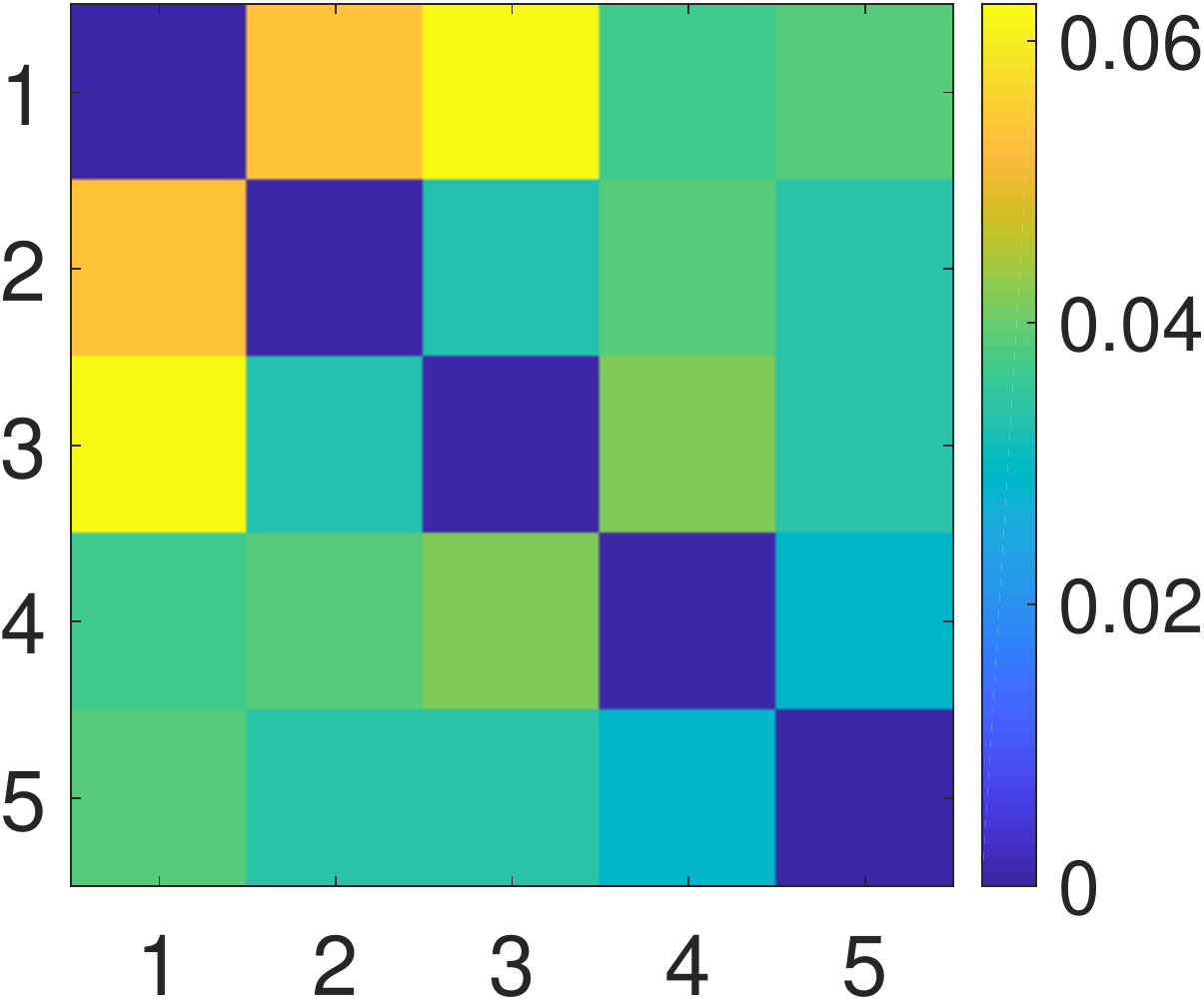}&\includegraphics[width=1.2in]{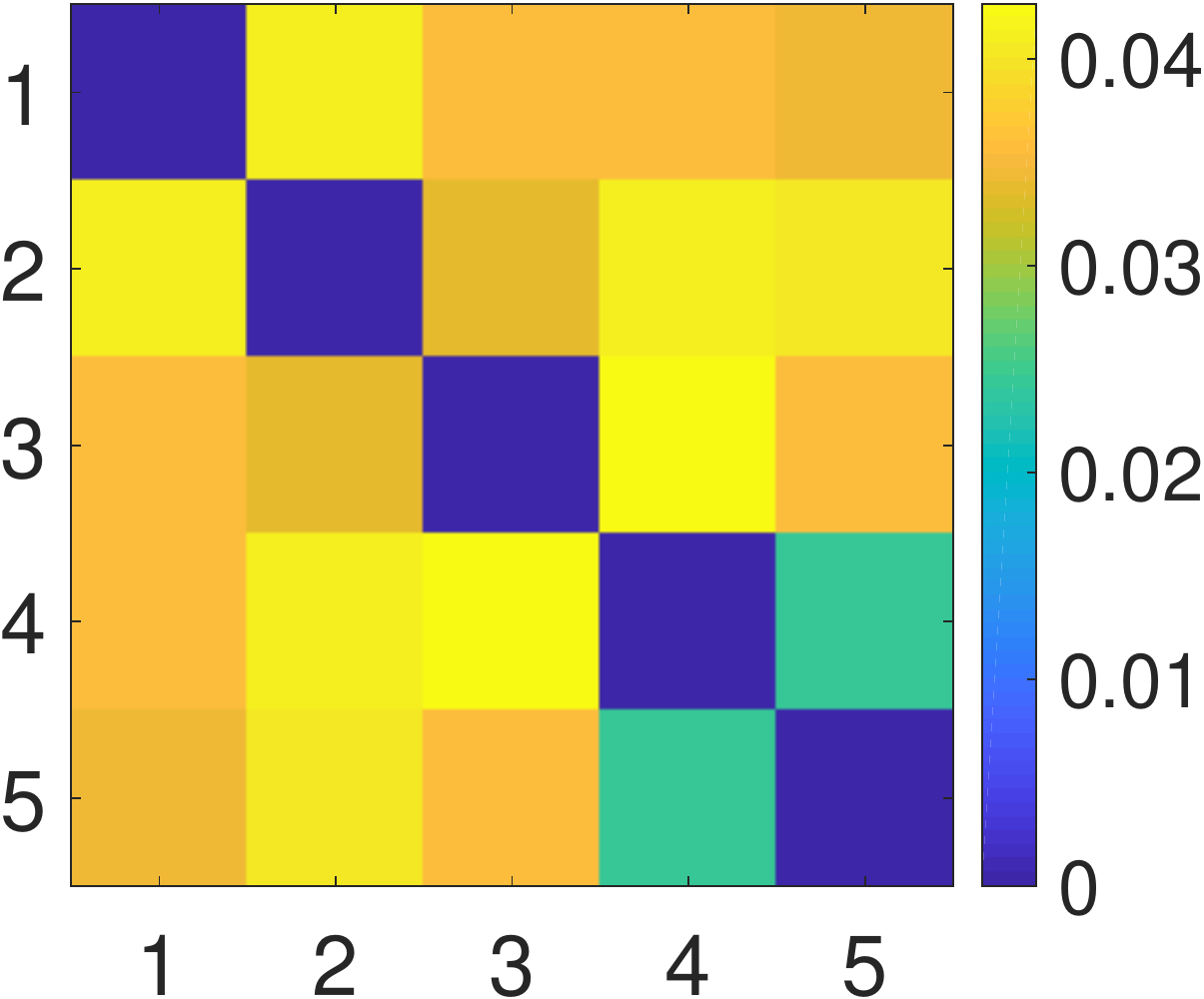}&\includegraphics[width=1.2in]{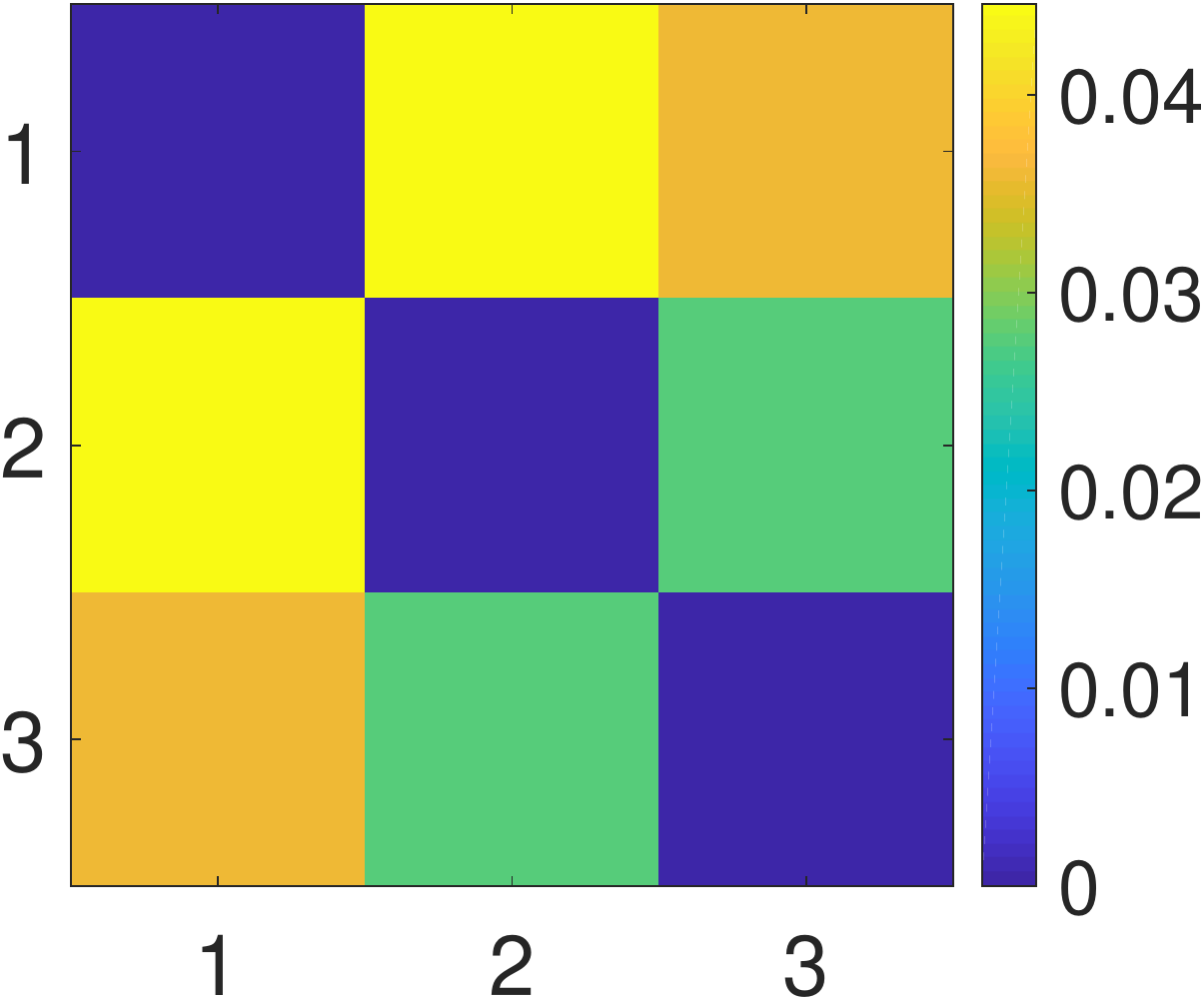}\\
			    \hline
            \end{tabular}
		\caption{\small Top: Simulated functions used to evaluate the stability and computational cost of Simulated Annealing alignment. Bottom: Alignment results across 20 replicates for different choices of parameters in Algorithm 4: $n=10,15,20,25,40$, $\theta=50,90,100,110,150$, $T=5,9,10,11,15$ and $c=1.0001,1.0005,1.001$. (a) Optimal warping functions for all parameter values and replicates ($20\times 5=100$ for $n$, $\theta$ and $T$, and $20\times 3=60$ for $c$). (b) Means of optimal warp maps across replicates (5 for $n$, $\theta$ and $T$, and 3 for $c$). (c) Pairwise distances between means in (b) (same order as parameters listed above).}\label{simdatarobres}
	\end{center}
\end{figure}

Here, we briefly study the stability and computational cost of Simulated Annealing-based alignment for two functions $g_i:[0,1]\to\mathbb{R},\ i=1,2$. The functions were selected from a simulated dataset of 21 functions, and are displayed in the top panel of Figure \ref{simdatarobres}. The results are shown in Figure \ref{simdatarobres}. The default parameter values we use in our implementation are $n=20$, $\theta=100$, $T=10$ and $c=1.0001$. In this simulation we study robustness of the alignment with respect to these parameter choices. For this purpose, we performed 20 alignments for different combinations of parameter values by changing one parameter at a time. We use the following settings: $n=10,15,20,25,40$, $\theta=50,90,100,110,150$, $T=5,9,10,11,15$ and $c=1.0001,1.0005,1.001$. Thus, for example, the left plot in Figure \ref{simdatarobres}(a) shows 100 optimal warp maps, 20 for each of the five values of $n$ considered (the other parameters are held at their default settings). This was repeated for all parameters. In Figure \ref{simdatarobres}(b), we display the cross-sectional averages of the warp maps computed across replicates. Finally, Figure \ref{simdatarobres}(c) shows the pairwise Fisher--Rao distances \citep{AS} (maximum is $\pi/2$) between the averages displayed in panel (b). In all cases, the alignment solutions are robust to the various parameter choices, which is evidenced by the very small variability in the plots in panel (a). In fact, the average warps in panel (b) are nearly indistinguishable, which is also confirmed by the very small distances reported in panel (c). A single Simulated Annealing alignment result, for two functions sampled with 100 points each, can be computed in approximately 30 seconds in Matlab 2018a on a Dell OptiPlex 7050 desktop with an Intel Core i7 processor and 16GB of RAM. This time can be further reduced by optimizing the number of iterations needed in the algorithm as well as implementing the method in C++.

\bibliographystyle{Chicago}
\bibliography{biblio}